\documentclass[11pt]{article}
\pdfoutput=1 % necessary for proper colour transparency in arXiv

\usepackage{defs/packages}
\usepackage{defs/macros}
\usepackage{defs/diagrammatics}
\usepackage{defs/environments}

\usepackage[a4paper, total={6in, 9in}]{geometry}

\counterwithin{figure}{section}
\usepackage{perpage}
\MakePerPage{footnote}

\title{Odd Khovanov homology and higher representation theory}
\newcommand{\orcid}[1]{\href{https://orcid.org/#1}{\,\orcidlink{#1}}}
\author{Léo Schelstraete\orcid{0000-0001-7167-3964
} \and Pedro Vaz\orcid{0000-0001-9422-4707}}
\date{\vspace{-5ex}}

\begin{document}
\maketitle

%%%%%%%%%%%%%%%%%%%%%%%%%%%%%%%%%%
%%%          abstract          %%%
%%%%%%%%%%%%%%%%%%%%%%%%%%%%%%%%%%
\begin{abstract}
  We introduce a super analogue of $\mathfrak{gl}_2$-foams, and use it to define an invariant of oriented tangles, shown to coincide with odd Khovanov homology when restricted to links. We then define a supercategorification of the $q$-Schur algebra of level 2 and realize our construction as a certain super-2-representation.
  This gives a representation-theoretic construction of odd Khovanov homology, where signs naturally arise as a byproduct of the super-2-categorical structure.
  In the process, we define a tensor product for chain complexes in super-2-categories, suitably compatible with homotopies. This could be of independent interest.
\end{abstract}

\tableofcontents

\section{Introduction}
\label{sec:intro}

Khovanov homology \cite{Khovanov_CategorificationJonesPolynomial_2000} is a homological invariant of links categorifying the Jones polynomial. Since its discovery, many approaches have appeared, including
singular surfaces
\cite{
  Bar-Natan_KhovanovsHomologyTangles_2005,
  CMW_FixingFunctorialityKhovanov_2009,
  Sano_FixingFunctorialityKhovanov_2021,
  Caprau_$sl2$TangleHomology_2008,
  Khovanov_Sl3LinkHomology_2004,
  Blanchet_OrientedModelKhovanov_2010,
  MSV_$mathfrakslN$linkHomology$N_2009,
},
categorical skew Howe duality
\cite{
  LQR_KhovanovHomologySkew_2015,
  QR_$mathfraksl_n$Foam2category_2016,
},
matrix factorizations
\cite{
  KR_MatrixFactorizationsLink_2008,
  % Yonezawa_QuantumMathfrakSl_2011,
  % Wu_ColoredMathfrakSl_2014,
},
categorification of tensor products of representations
\cite{Webster_KnotInvariantsHigher_2017},
category $\cO$
\cite{
  MS_CombinatorialApproachFunctorial_2009,
  % Sussan_CategorySlLink_2007
},
algebraic geometry
\cite{
  CK_KnotHomologyDerived_2008,
  % CK_KnotHomologyDerived_2008a,
},
symplectic geometry
\cite{
  SS_LinkInvariantSymplectic_2006,
},
gauge theory
\cite{
  Witten_KhovanovHomologyGauge_2012,
  % Witten_FivebranesKnots_2012
},
and mirror symmetry
\cite{
  Aganagic_KnotCategorificationMirror_2020,
  % Aganagic_KnotCategorificationMirror_2021
}.

Given the central role played by Khovanov homology, it is surprising that the Jones polynomial admits another categorification, called \emph{odd Khovanov homology} \cite{ORS_OddKhovanovHomology_2013}.\footnote{In another direction, symmetric Khovanov homology \cite{QRS_AnnularEvaluationLink_2018,RW_SymmetricKhovanovRozanskyLink_2020} is yet another categorification of the Jones polynomial.}
It agrees with Khovanov homology modulo $2$, but the two homologies are distinct in the sense that one can find pairs of knots distinguished by one but not the other, and vice versa \cite{Shumakovitch_PatternsOddKhovanov_2011}.
Odd Khovanov homology was discovered by Ozsváth, Rasmussen and Szabó in an attempt to lift to the integers the Ozsváth--Szabó spectral sequence \cite{OS_HeegaardFloerHomology_2005} from reduced Khovanov homology to the Heegaard--Floer homology of the branched double cover. In comparison, Dowlin has shown that there exists a spectral sequence from (even) Khovanov homology to knot Floer homology \cite{Dowlin_SpectralSequenceKhovanov_2024}.
Also, while (even) Khovanov homology is well-understood via the representation theory of quantum groups, odd Khovanov homology is expected to be related to the representation theory of quantum \emph{super}groups \cite{CHW_QuantumSupergroupsFoundations_2013}; see \cref{rem:super_lie_theory} for details on this last speculation.

Heuristically, an \emph{odd analogue} is an algebraic structure exhibiting anti-commutative behaviour, with the same graded rank as its even (commutative) counterpart and such that the even and odd constructions agree when reduced modulo $2$.
The appearance of odd Khovanov homology sparked interest in finding odd analogues of known categorified and geometric structures
\cite{
  EL_DGStructuresOdd_2020,
  EL_OddCategorification$U_qmathfraksl_2$_2016,%
  BE_MonoidalSupercategories_2017,BE_SuperKacMoody_2017,%
  ENW_RealSpringerFibers_2021,%
  ELV_DerivedSuperequivalencesSpin_2022,%
  EKL_OddNilHeckeAlgebra_2014,EQ_DifferentialGradedOdd_2016,%
  KKO_SupercategorificationQuantumKacMoody_2013,%
  KKO_SupercategorificationQuantumKac_2014,%
  LR_OddificationCohomologyType_2014,%
  NV_OddKhovanovsArc_2018,%
  BK_OddGrassmannianBimodules_2022,%
}, motivated by their relation with (even) Khovanov homology.
Underlying most of these constructions is the notion of a \emph{super-2-category} \cite{BE_MonoidalSupercategories_2017} (or 2-supercategory), a certain categorical structure akin to a linear 2-category, where the interchange law is twisted by the extra data of a $\bZ/2\bZ$-grading, or \emph{parity}, on the 2Hom-spaces. Diagrammatically, this is pictured as follows:
\begin{equation*}
    \xy (0,0)* {
    \begin{tikzpicture}[scale=0.8]
        \draw (0,2) node[above] {\scriptsize $g'$} to
            node[black_dot,pos=.3] {}
            node[right,pos=.3] {\scriptsize $\beta$}
            (0,0) node[below] {\scriptsize $f'$};
        \draw (1,2) node[above] {\scriptsize $g$} to
            node[black_dot,pos=.7] {}
            node[right,pos=.7] {\scriptsize $\alpha$}
            (1,0) node[below] {\scriptsize $f$};
            % \node at (2,1) {$a$};
            % \node at (-1,1) {$b$};
    \end{tikzpicture} }\endxy
    \;=\;(-1)^{\abs{\alpha}\cdot\abs{\beta}}
    \xy (0,0)* {
    \begin{tikzpicture}[scale=0.8]
        \draw (0,2) node[above] {\scriptsize $g'$} to
            node[black_dot,pos=.7] {}
            node[right,pos=.7] {\scriptsize $\beta$}
            (0,0) node[below] {\scriptsize $f'$};
        \draw (1,2) node[above] {\scriptsize $g$} to
            node[black_dot,pos=.3] {}
            node[right,pos=.3] {\scriptsize $\alpha$}
            (1,0) node[below] {\scriptsize $f$};
            % \node at (2,1) {$a$};
            % \node at (-1,1) {$b$};
    \end{tikzpicture} }\endxy
\end{equation*}
Here $\abs{\alpha}$ and $\abs{\beta}$ denote the respective parities of the 2-morphisms $\alpha$ and $\beta$.
Note that a super-2-category is \emph{not} a 2-category endowed with extra structure. By \emph{supercategorification} we mean the process of categorifying a category with a super-2-category.

\medbreak

In this article, we show that odd Khovanov homology arises from an odd analogue in categorified representation theory.
More precisely, we give a supercategorification $\sfoam_d$ of a certain integral and idempotent form $\web_d$ of the representation theory of $U_q(\glt)$, such that $\sfoam_d$ has the same graded rank as its even counterpart $\foam_d$ (\cref{sec:foams}).
We then define an invariant of oriented tangles as a certain tensor product of complexes in $\sfoam_d$, and show that it coincides with odd Khovanov homology when restricted to links (\cref{sec:topo}).
Here $\sfoam_d$ stands for \emph{super $\glt$-foams}, where $\glt$-foams are certain decorated singular surfaces used in Blanchet's approach to Khovanov homology \cite{Blanchet_OrientedModelKhovanov_2010}.
Our approach can be understood as an odd analogue of his:
\[\begin{tikzpicture}
  \node (F) at (0,1) {$\foam_d$};
  \node (SF) at (4,1) {$\sfoam_d$};
  \node (W) at (2,-.5) {$\web_d$};
  \begin{scope}[shift={(8,0)}]
    \node[align=center] (Kh) at (0,1) {
      \text{\small (even)}\\[-1ex]
      \text{\small Khovanov}\\[-1ex]
      \text{\small homology}
    };
    \node[align=center] (OKh) at (4,1) {
      \text{\small odd}\\[-1ex]
      \text{\small Khovanov}\\[-1ex]
      \text{\small homology}
    };
    \node[align=center] (J) at (2,-.5) {
      \text{\small Jones}\\[-1ex]
      \text{\small polynomial}
    };
  \end{scope}
  \node at (2+8/2,.25) {$\rightsquigarrow$};
  \draw[bend right,->] (F) to node[left]{$K_0$} (W);
  \draw[bend left,->] (SF) to node[right]{$K_0$} (W);
  \draw[bend right,->] (Kh) to node[left]{$\chi_q$} (J);
  \draw[bend left,->] (OKh) to node[right]{$\chi_q$} (J);
\end{tikzpicture}\]
In particular, $\foam_d$ and $\sfoam_d$ both categorify the same underlying category $\web_d$, just as even and odd Khovanov homologies both categorify the Jones polynomial. However, the two categorifications are structurally of different flavour, as $\foam_d$ categorifies while $\sfoam_d$ \emph{super}categorifies.

Similarly, the tensor product of complexes should be understood in a super sense. We define this super tensor product in \cref{sec:complexes}; as this is technical (although not conceptually difficult), we also give a minimal version in \cref{subsec:hypercubic_chain_complexes}, sufficient for the purpose of \cref{sec:topo}.
These results first appeared in the first author's master's thesis \cite{Schelstraete_SupercategorificationKhovanovlikeTangle_2020}.

The main result of this paper may be summarized as follows:
\begin{quote}
  \textsc{Slogan:} \emph{Odd Khovanov homology arises from the super-interplay of two categorified Kauffman brackets, one even and the other odd, respectively associated to the zip and unzip foams (see \cref{subsec:intro_super_foams}).}
\end{quote}
In particular, this gives an extension of odd Khovanov homology to oriented tangles (another extension was given in \cite{NP_OddKhovanovHomology_2020}; see \cref{rem:NP_arc_algebra}). This also gives a \emph{sign-coherent} definition, where signs naturally arise from the super-2-categorical structure. This contrasts with the original definition, where signs are fixed in a somewhat artificial way on the hypercube of resolutions.
We hope that this new definition will pave the way for further connections and applications.

Finally, we relate our construction to higher representation theory (\cref{sec:schur}). More precisely, we define a supercategorification $\cS\cS_{n,d}$ of the $q$-Schur algebra of level 2, extending the work of the second author \cite{Vaz_NotEvenKhovanov_2020}.
We then exhibit a super-2-functor
\[
  \cF_{n,d}\colon \cS\cS_{n,d} \to \sfoam_d,
\]
called the \emph{super foamation 2-functor}. In that way, $\sfoam_d$ is realized as a super-2-representation of $\cS\cS_{n,d}$.
% In particular, the invariant defined in \cite{Vaz_NotEvenKhovanov_2020} coincides with ours, and hence with odd Khovanov homology when restricted to links.
This gives a partial odd analogue of the work of Lauda, Queffelec and Rose \cite{LQR_KhovanovHomologySkew_2015}, who showed that $\foam_d$ can be viewed as a 2-representation of the categorification of $U_q(\mathfrak{gl}_n)$ \cite{KL_DiagrammaticApproachCategorification_2009,Rouquier_2KacMoodyAlgebras_2008} via categorical skew Howe duality.

\medbreak

We do not address here the problem of showing that our super-2-category of $\glt$-foams is sufficiently non-degenerate (see \cref{thm:basis_foam} for a precise statement), on which most of our results rely.
This is addressed by the first author in \cite{Schelstraete_RewritingModuloDiagrammatic_2026} (see also the first author's PhD thesis \cite{Schelstraete_OddKhovanovHomology_2024}) using novel techniques from rewriting theory.
In a nutshell, rewriting theory is the algorithmic study of presentations of algebraic structures; classical instances include Gröbner--Shirshov bases \cite{Buchberger_AlgorithmusAuffindenBasiselemente_1965,Shirshov_AlgorithmicProblemsLie_1962} and Bergman's diamond lemma \cite{Bergman_DiamondLemmaRing_1978}.
In recent years, rewriting theory has emerged as a promising tool in higher representation theory, either explicitly \cite{Alleaume_RewritingHigherDimensional_2018,Dupont_RewritingModuloIsotopies_2021,Dupont_RewritingModuloIsotopies_2022,DEL_SuperRewritingTheory_2024} or implicitly \cite{Elias_DiagrammaticTemperleyLiebCategorification_2010,Barbier_DiagramCategoriesBrauer_2024}.
In combination with \cite{Schelstraete_RewritingModuloDiagrammatic_2026}, this paper gives the first application of rewriting theory to quantum topology.

\begin{remark}[graded-2-categories and covering Khovanov homology]
  To simplify the exposition, we restricted this introduction (including the extended summary, \cref{subsec:extended_summary}) to the super, or $\bZ/2\bZ$-graded, case. In order to encompass both even and odd Khovanov homology, the rest of the article considers the more general setting of \emph{graded-2-categories}. This includes the tensor product on chain complexes in graded-2-categories.

  Hence, \cref{sec:foams} actually defines a graded-2-category of $\glt$-foams $\gfoam_d$, defined over a ring $\ringfoam$ with distinguished elements $X$, $Y$ and $Z$ such that $X^2=Y^2=1$ (\cref{defn:ring_R_bil}).
  Choosing $X=Y=Z=1$ recovers the even case, while choosing $X=Z=1$ and $Y=-1$ recovers the odd case.
  \[\begin{tikzpicture}[scale=.75]
  \node (F) at (0,1) {$\foam_d$};
  \node (SF) at (4,1) {$\sfoam_d$};
  \node (W) at (2,-.5) {$\web_d$};
  \node (GF) at (2,2.5) {$\gfoam_d$};
  \begin{scope}[shift={(10,0)}]
    \node[align=center] (Kh) at (0,1) {
      \text{\small (even)}\\[-1ex]
      \text{\small Khovanov}\\[-1ex]
      \text{\small homology}
    };
    \node[align=center] (OKh) at (4,1) {
      \text{\small odd}\\[-1ex]
      \text{\small Khovanov}\\[-1ex]
      \text{\small homology}
    };
    \node[align=center] (J) at (2,-.5) {
      \text{\small Jones}\\[-1ex]
      \text{\small polynomial}
    };
    \node[align=center] (GKh) at (2,2.5) {
      \text{\small covering}\\[-1ex]
      \text{\small Khovanov}\\[-1ex]
      \text{\small homology}
    };
  \end{scope}
  \node at (2+10/2,1) {$\rightsquigarrow$};
  \draw[bend right,->] (F) to node[left]{$K_0$} (W);
  \draw[bend left,->] (SF) to node[right]{$K_0$} (W);
  \draw[bend right,->] (GF) to node[left,align=center]{
    \footnotesize$X=Z=1$\\[-1ex]
    \footnotesize$\an Y=1$
  } (F);
  \draw[bend left,->] (GF) to node[right,align=center]{
    \footnotesize$X=Z=1$\\[-1ex]
    \footnotesize$\an Y=-1$
  } (SF);
  \draw[bend right,->] (Kh) to node[left]{$\chi_q$} (J);
  \draw[bend left,->] (OKh) to node[right]{$\chi_q$} (J);
  \draw[bend right,->] (GKh) to (Kh);
  \draw[bend left,->] (GKh) to (OKh);
\end{tikzpicture}\]
  Leaving $X$, $Y$ and $Z$ as formal parameters recovers \emph{covering Khovanov homology} \cite{Putyra_2categoryChronologicalCobordisms_2014}.
  Similarly, the supercategorification $\cS\cS_{n,d}$ of the $q$-Schur algebra of level 2 mentioned above is obtained from the graded variant defined in \cref{sec:schur} by setting $X=Z=1$ and $Y=-1$.
\end{remark}

\begin{remark}[not even Khovanov homology]
  In \cite{Vaz_NotEvenKhovanov_2020}, the second author defined a homological tangle invariant called ``not even Khovanov homology''.\footnote{This definition relied on the existence of a suitable super tensor product, which is the content of \cref{sec:complexes}.} Our construction is a foamy analogue of \cite{Vaz_NotEvenKhovanov_2020}; in particular, it follows from our result that not even Khovanov homology coincides with odd Khovanov homology when restricted to links, as was conjectured in \cite{Vaz_NotEvenKhovanov_2020}.
\end{remark}

The remainder of the introduction provides an extended summary of the paper, with extra historical and motivational notes. We suggest that the casual reader start there (\cref{subsec:extended_summary}).
The following remarks discuss related open problems.

\begin{remark}[connections with odd arc algebras]
  \label{rem:NP_arc_algebra}
  In \cite{NP_OddKhovanovHomology_2020}, Naisse and Putyra gave another extension of odd Khovanov homology to tangles based on arc algebras, building on previous work of Putyra \cite{Putyra_2categoryChronologicalCobordisms_2014} and Naisse--Vaz \cite{NV_OddKhovanovsArc_2018}. They conjectured that their tangle invariant coincides with the one in \cite{Vaz_NotEvenKhovanov_2020}.
  Following their conjecture and the previous remark, Naisse and Putyra's construction should coincide with our tangle invariant.
  This remains an open question.
  See also the introduction of \cref{sec:schur} for further connections with their work.
\end{remark}

\begin{remark}[other odd link homologies]
  At present, there are no known odd analogues of $\mathfrak{gl}_n$-Khovanov homologies \cite{KR_MatrixFactorizationsLink_2008} outside of the case $n=2$.
  In the foamy construction of these homologies \cite{Khovanov_Sl3LinkHomology_2004,MSV_$mathfrakslN$linkHomology$N_2009}, the $n$th power of the dot is zero: $(\textit{dot})^n=0$.
  On the other hand, if one requires the dot to be odd, then $(\textit{dot})^2=-(\textit{dot})^2$, and at least if $2$ is invertible in the ground ring, this implies that $(\textit{dot})^2=0$.
  This gives an obstruction to a naive construction of odd $\mathfrak{gl}_n$-Khovanov homology when $n\geq 3$.
  However, it may be that $\mathfrak{gl}_n$ is not the correct direction to look at; see the next remark.
\end{remark}

\begin{remark}[connections with super Lie theory]
  \label{rem:super_lie_theory}
  Supercategorification is known to be related to super Lie theory. For instance, consider the case of $\slt$, associated with the Cartan datum consisting of a single vertex. The Lie superalgebra $\mathfrak{osp}_{1|2}$ similarly arises from the Cartan super datum consisting of a single odd vertex. Hill and Wang introduced in \cite{HW_CategorificationQuantumKacMoody_2015} (see also the series of papers \cite{CHW_QuantumSupergroupsFoundations_2013,CH_QuantumSupergroupsBraid_2016,CHW_QuantumSupergroupsII_2014,CSW_QuantumSupergroupsVI_2019,CFL+_QuantumSupergroupsIII_2014,Clark_QuantumSupergroupsIV_2014}) a \emph{covering quantum group} $U_{q,\pi}(\slt)$ defined over $\bQ(q)[\pi]/(\pi^2-1)$, such that setting $\pi=1$ recovers $U_{q}(\slt)$ while setting $\pi=-1$ recovers $U_{q}(\mathfrak{osp}_{1|2})$.
  On the other hand, Ellis and Lauda \cite{EL_OddCategorification$U_qmathfraksl_2$_2016} constructed a supercategorification of $U_{q,\pi}(\slt)$, later reformulated (and extended to other super Cartan data) by Brundan and Ellis \cite{BE_SuperKacMoody_2017} (see also the beginning of \cref{sec:schur} for further references).
  Given those interactions, it was thought that an odd homology should correspond to a covering quantum group (resp.\ a Lie superalgebra), with odd Khovanov homology corresponding to $U_{q,\pi}(\slt)$ (resp.\ $\mathfrak{osp}_{1|2}$). However, an explicit connection between odd Khovanov homology and $\mathfrak{osp}_{1|2}$ remains an open problem (see however \cite{Clark_OddKnotInvariants_2017,Ebert_NewPresentationOsp1|2polynomial_2022}). We expect that a further careful study of our construction will lead to such a connection.

  Let us also note that under some assumptions \cite{HW_CategorificationQuantumKacMoody_2015}, the only Lie superalgebras in finite type are $\mathfrak{osp}_{1|2n}$. This suggests that an odd $\mathfrak{so}_{2n+1}$-link homology should exist.
  % This would explain why there is no $\mathfrak{sl}_n$-Khovanov homology.
  This is further corroborated by the work of Mikhaylov and Witten \cite{MW_BranesSupergroups_2015} on link homologies associated to $\mathfrak{so}_{2n+1}$, where the $\slt\cong\mathfrak{so}_3$ case is conjectured to coincide with odd Khovanov homology. See also \cite{EL_OddCategorification$U_qmathfraksl_2$_2016} for a discussion.
\end{remark}

\begin{remark}[connections with supercategorified quantum groups]
  \label{rem:intro_relation_to_supercategorified_QG}
  While \cref{sec:schur} gives a supercategorification (denoted $\cS\cS_{n,d}$) of the $q$-Schur algebra of level 2 (denoted $\qschur_{n,d}$), it does \emph{not} give a supercategorification of the algebra map $U_q(\mathfrak{gl}_n)\to\qschur_{n,d}$.
  Indeed, while there exists a conjectural (even) categorification of $U_q(\mathfrak{gl}_n)$ \cite{MSV_DiagrammaticCategorification$q$Schur_2013}, there is, at the time of writing, not even a candidate for a supercategorification of $U_q(\mathfrak{gl}_n)$.
  Relating $\cS\cS_{n,d}$ to supercategorified quantum groups \cite{BE_SuperKacMoody_2017} (also known as \emph{super Kac--Moody 2-categories}) remains an interesting open problem. See the previous two remarks for related comments.
\end{remark}

\subsection{Acknowledgments}

We thank James MacPherson, Grégoire Naisse, Kris Putyra, Louis-Hadrien Robert and Paul Wedrich for interesting discussions, and the anonymous referees for their many useful suggestions. L.S.\ was supported by the Fonds de la Recherche Scientifique -- FNRS under the Aspirant Fellowship FC 38559.
P.V.\ was supported by the Fonds de la Recherche Scientifique -- FNRS under Grant Nos. MIS-F.4536.19 and CDR-J.0189.23.

\subsection{Extended summary}
\label{subsec:extended_summary}

\subsubsection{What is odd Khovanov homology?}

The original construction of Khovanov homology consists in introducing a \emph{hypercube of resolutions} associated with every link diagram, using a suitable 2-dimensional TQFT to algebraize the hypercube, and turning the hypercube into a chain complex by assigning signs to its edges following the Koszul rule.
In that case, the commutative Frobenius algebra associated with the TQFT is the algebra $\bZ[x]/(x^2)$.

Odd Khovanov homology is constructed similarly, only with the exterior algebra $\wedge(x_1,\ldots,x_n)$ taking the role of $\bZ[x_1,\ldots,x_n]/(x_1^2,\ldots,x_n^2)$.
Of course, the exterior algebra is not a commutative Frobenius algebra, and so the associated TQFT is a \emph{projective} TQFT in the sense that it is only functorial up to signs.
Somewhat miraculously, it was shown in \cite{ORS_OddKhovanovHomology_2013} that this defect in functoriality can be balanced out when assigning signs to the hypercube. This requires a much more artificial sign assignment than the Koszul rule, based on a case-by-case analysis of possible squares in the hypercube.

As the anti-commutativity in the exterior algebra is controlled by a $\bZ/2\bZ$-grading, it is natural to wonder whether one could give a construction of odd Khovanov homology using a super-2-category.
Ideally, the superstructure would control all signs appearing in odd Khovanov homology, that is, all interchanges of saddles.
One solution, pursued by Putyra \cite{Putyra_2categoryChronologicalCobordisms_2014}, is to pull back the TQFT to linear relations on cobordisms. This provides a partial analogue for odd Khovanov homology of Bar-Natan's ``picture-world'' construction of (even) Khovanov homology.
In this context, a merge is even and a split is odd:
\begin{IEEEeqnarray*}{CCcCC}
  \mspace{50mu}&
  \xy(0,0)*{\begin{tikzpicture}[scale=.7,yscale=-1]
    \draw (0,0) to[out=-90,in=-90] (1,0) edge[out=90,in=90] (0,0);
    \draw (0,0) to[out=90,in=-90] (-1,2);
    \draw (1,0) to[out=90,in=-90] (2,2);
    \draw (0,2) to[out=-90,in=180] (.5,1) to[out=0,in=-90] (1,2);
    \draw (-1,2) edge[dashed,out=-90,in=-90] (0,2);
    \draw (0,2) to[out=90,in=90] (-1,2);
    \draw (1,2) edge[dashed,out=-90,in=-90] (2,2);
    \draw (2,2) to[out=90,in=90] (1,2);
    \fill[fill_foam1=\shop+.1] (0,0) to[out=90,in=-90] (-1,2) to[out=-90,in=-90] (0,2) to[out=-90,in=180] (.5,1) to[out=0,in=-90] (1,2) to[out=-90,in=-90] (2,2) to[out=-90,in=90] (1,0) to[out=-90,in=-90] (0,0);
    \fill[fill_foam1=\shop+.3] (0,0) to[out=90,in=-90] (-1,2) to[out=90,in=90] (0,2) to[out=-90,in=180] (.5,1) to[out=0,in=-90] (1,2) to[out=90,in=90] (2,2) to[out=-90,in=90] (1,0) to[out=90,in=90] (0,0);
  \end{tikzpicture}}\endxy
  &\mspace{150mu}&
  % UNZIP
  \xy(0,0)*{\begin{tikzpicture}[scale=.7]
    \draw (0,0) to[out=-90,in=-90] (1,0) edge[dashed,out=90,in=90] (0,0);
    \draw (0,0) to[out=90,in=-90] (-1,2);
    \draw (1,0) to[out=90,in=-90] (2,2);
    \draw (0,2) to[out=-90,in=180] (.5,1) to[out=0,in=-90] (1,2);
    \draw (-1,2) to[out=-90,in=-90] (0,2) to[out=90,in=90] (-1,2);
    \draw (1,2) to[out=-90,in=-90] (2,2) to[out=90,in=90] (1,2);
    \fill[fill_foam1=\shop+.3] (0,0) to[out=90,in=-90] (-1,2) to[out=-90,in=-90] (0,2) to[out=-90,in=180] (.5,1) to[out=0,in=-90] (1,2) to[out=-90,in=-90] (2,2) to[out=-90,in=90] (1,0) to[out=-90,in=-90] (0,0);
    \fill[fill_foam1=\shop+.1] (0,0) to[out=90,in=-90] (-1,2) to[out=90,in=90] (0,2) to[out=-90,in=180] (.5,1) to[out=0,in=-90] (1,2) to[out=90,in=90] (2,2) to[out=-90,in=90] (1,0) to[out=90,in=90] (0,0);
  \end{tikzpicture}}\endxy
  &
  \mspace{50mu}
  \xy(0,0)*{\begin{tikzpicture}[scale=.7]
    \draw[->] (0,.2) to (0,3-.2);
    \node[rotate=-90] at (.5,1.5) {\scriptsize reading direction};
  \end{tikzpicture}}\endxy
  \\*[1ex]
  &
  \textit{even}
  &&
  \textit{odd}
  &
\end{IEEEeqnarray*}
However, the superstructure only partially controls interchanges of saddles, and one still needs to use the artificial sign assignment from the original construction.
Also, and contrary to Bar-Natan's approach to Khovanov homology, it does not generalize in an obvious way to tangles---indeed, whether a saddle is a split or a merge is a global property. See however \cite{NP_OddKhovanovHomology_2020} for an answer to this question using an odd analogue of arc algebras.

This suggests that one should look further away from the original construction. Heuristically, it is plausible that different choices of TQFTs, identical up to signs, lead to the same invariant. After all, this is the take-home message of odd Khovanov homology: sign issues on the level of the TQFT can be balanced out when choosing a sign assignment on the hypercube.
As we show below, a solution is given by super $\glt$-foams.

%%%%%%%%%%%%%%%%%%%%%%%%%%%%%%%
%%%          foams          %%%
%%%%%%%%%%%%%%%%%%%%%%%%%%%%%%%
\subsubsection{Super \texorpdfstring{$\glt$}{gl2}-foams}
\label{subsec:intro_super_foams}

%   \caption{Local model for foams.}
%   \label{fig:local_model_foam}
% \end{wrapfigure}

An early and important problem in Khovanov homology also concerned signs. Namely, Khovanov's original construction is not properly functorial under link cobordisms, but only so up to signs. Solutions to this problem were provided by several authors \cite{CMW_FixingFunctorialityKhovanov_2009,Sano_FixingFunctorialityKhovanov_2021,Caprau_$sl2$TangleHomology_2008,Vogel_FunctorialityKhovanovHomology_2020}.
A solution introduced by Blanchet \cite{Blanchet_OrientedModelKhovanov_2010} used
\emph{foams}, certain decorated singular surfaces first introduced by Khovanov in his definition of $\mathfrak{sl}_3$-Khovanov homology \cite{Khovanov_Sl3LinkHomology_2004} and generalized to $\mathfrak{gl}_n$ for $n\geq 3$ in \cite{MSV_$mathfrakslN$linkHomology$N_2009}.
Adapting this construction to the $n=2$ case, Blanchet defined a functorial version of Khovanov homology using \emph{$\glt$-foams}.\footnote{Blanchet called them \emph{enhanced $\slt$-foams}, but as they were later understood to be related to $\glt$ rather than $\slt$, we call them $\glt$-foams.}
The proof of functoriality was later generalized to all $\mathfrak{gl}_n$-link homologies in \cite{ETW_FunctorialityColoredLink_2018}.
In practice, working with $\glt$ instead of $\slt$ leads to better-behaved constructions. This comes down to the fact that in the former case, the fundamental representations $\bigwedge^0(\bC(q)^2)$ and $\bigwedge^2(\bC(q)^2)$ are not isomorphic, and keeping track of this distinction leads to better control on signs (see \cite[Section~1F]{LQR_KhovanovHomologySkew_2015} for a discussion).

As it turns out, the same heuristics give control on signs in odd Khovanov homology.
One needs to work with $\glt$-foams equipped with a sort of Morse decomposition, in the sense that they decompose into a composition of the following local pictures:
\begin{IEEEeqnarray*}{CcCcCcCcC}
  \text{dot}
  &&
  \text{cup}
  &&
  \text{cap}
  &&
  \text{zip}
  &&
  \text{unzip}
  \\*[1ex]
  % DOT
  \figfoam[scale=1.3]{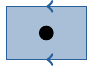}
  &\mspace{30mu}&
  % CUP
  \figfoam[scale=1.3]{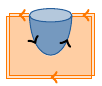}
  &\mspace{30mu}&
  % CAP
  \figfoam[scale=1.3]{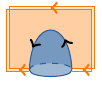}
  &\mspace{30mu}&
  % ZIP
  \figfoam[scale=1.3]{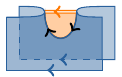}
  &\mspace{30mu}&
  % UNZIP
  \figfoam[scale=1.3]{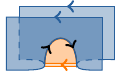}
  \\*[1ex]
  \textit{odd}
  &&
  \textit{odd}
  &&
  \textit{even}
  &&
  \textit{even}
  &&
  \textit{odd}
\end{IEEEeqnarray*}
Assigning parities to these local pictures turns a $\glt$-foam into a \emph{super} $\glt$-foam.
The super-2-category $\sfoam_d$ is generated by linear combinations of super $\glt$-foams modulo relations, some of which encode \emph{super} isotopies (see \cref{fig:rel_foam} in the text, setting $X=Z=1$ and $Y=-1$). Note that (super) $\glt$-foams can be decorated with dots: in the super case, this dot is odd, just as the variable $x_i$ is in the exterior algebra $\wedge(x_1,\ldots,x_n)$.

Boundaries of (super) $\glt$-foams are certain trivalent graphs called \emph{webs} (or \emph{spider webs}). Their edges have a \emph{thickness}
$\xy(0,0)*{\begin{tikzpicture}[scale=.6,rotate=90]
  \draw[web1] (0,0) to (0,1);
\end{tikzpicture}}\endxy$
and
$\xy(0,0)*{\begin{tikzpicture}[scale=.6,rotate=90]
  \draw[web2] (0,0) to (0,1);
\end{tikzpicture}}\endxy$. They provide an integral and idempotent form $\web_d$ for intertwiners of fundamental representations of $U_q(\glt)$ \cite{CKM_WebsQuantumSkew_2014}: in this correspondence, $\xy(0,0)*{\begin{tikzpicture}[scale=.6,rotate=90]
  \draw[web1] (0,0) to (0,1);
\end{tikzpicture}}\endxy$ corresponds to the standard representation $\bC(q)^2$ and $\xy(0,0)*{\begin{tikzpicture}[scale=.6,rotate=90]
  \draw[web2] (0,0) to (0,1);
\end{tikzpicture}}\endxy$ to the determinant representation $\bigwedge^2(\bC(q)^2)$.

\begin{bigtheorem}[see \cref{sec:foams}]
  There exists a $\bZ$-graded super-2-category $\sfoam_d$ which supercategorifies the $\bZ[q,q^{-1}]$-linear category $\web_d$.
\end{bigtheorem}

%%%%%%%%%%%%%%%%%%%%%%%%%%%%%%%%%%
%%%          homology          %%%
%%%%%%%%%%%%%%%%%%%%%%%%%%%%%%%%%%
\subsubsection{Odd Khovanov homology is a super tensor product of chain complexes}
Once given the super-2-category of super $\glt$-foams, defining a homological invariant of oriented tangles is straightforward in most aspects.
It follows the usual scheme of a categorified Kauffman bracket \cite{Bar-Natan_KhovanovsHomologyTangles_2005}: we assign length-two complexes to positive and negative crossings and take an appropriate tensor product (more precisely, a horizontal composition). The differentials are respectively given by the even saddle (zip) and the odd saddle (unzip):
\begin{IEEEeqnarray*}{rCrl}
  {}\xy(0,0)*{\begin{tikzpicture}[scale=-1*.7,xscale=1.2,rotate=90,transform shape]
  \draw[black,line width=1.5pt] (1,0) to[out=90,in=-90] (0,1);
  \draw[black,line width=1.5pt,overdraw=8pt] (0,0) to[out=90,in=-90] (1,1);
  \end{tikzpicture}}\endxy
  \qquad&\mapsto&\qquad
  % \left[\mspace{10mu}
  qt^{-1}\;
  \xy(0,0)*{\begin{tikzpicture}[scale=.7,rotate=90,transform shape]
    \draw[web1] (0,0) to (0,2);
    \draw[web1] (1,0) to (1,2);
    % \draw[snakecd] (-.5,0) to (-.5,2);
  \end{tikzpicture}}\endxy
  \quad&\xrightarrow{
    \figfoam{zip}
  }\quad
  \xy(0,0)*{\begin{tikzpicture}[scale=1*.7,transform shape,rotate=90]
    \pic at (0,0) {web+};
    \pic at (0,1) {web-};
  \end{tikzpicture}}\endxy
  % \mspace{10mu}\right]
  \\[2ex]
  {}\xy(0,0)*{\begin{tikzpicture}[scale=-1*.7,xscale=1.2,rotate=90,transform shape]
  \draw[black,line width=1.5pt] (0,0) to[out=90,in=-90] (1,1);
  \draw[black,line width=1.5pt,overdraw=8pt] (1,0) to[out=90,in=-90] (0,1);
  \end{tikzpicture}}\endxy
  \qquad&\mapsto&\qquad
  % \left[\mspace{10mu}
  qt^{-1}\;
  \xy(0,0)*{\begin{tikzpicture}[scale=1*.7,transform shape,rotate=90]
    \pic at (0,0) {web+};
    \pic at (0,1) {web-};
  \end{tikzpicture}}\endxy
  \quad&\xrightarrow{
    \figfoam{unzip}
  }\quad
  \;\xy(0,0)*{\begin{tikzpicture}[scale=1*.7,rotate=90,transform shape]
    \draw[web1] (0,0) to (0,2);
    \draw[web1] (1,0) to (1,2);
  \end{tikzpicture}}\endxy
  % \mspace{10mu}\right]
\end{IEEEeqnarray*}
Here $q$ denotes a shift in quantum grading, and $t$ indicates the homological grading. A similar assignment can be defined for cups and caps; see \cref{subsec:topo_defn_invariant} for details. Taking the super tensor product of these complexes and renormalizing assigns a complex to every sliced tangle diagram. Its homotopy type is shown to be an invariant of the associated oriented tangle (\cref{thm:invariance}).

The only step requiring extra work is the last step: taking an appropriate tensor product.
This is done in \cref{sec:complexes} for a subclass of chain complexes called \emph{homogeneous polycomplexes}. A \emph{homogeneous complex} is a chain complex in a super-2-category such that each differential is homogeneous (although the parity can differ at distinct homological degrees), and a homogeneous polycomplex is a tensor product of homogeneous complexes.
This tensor product is coherent with homotopies in the following sense:

\begin{bigtheorem}[\cref{thm:invariance_homotopy_classes}]
  \label{bigthm:intro_tensor_product}
  In any super-2-category, there exists a well-defined tensor product on homogeneous polycomplexes such that if $A_1^\bullet$ and $B_1^\bullet$ (resp.\ $A_2^\bullet$ and $B_2^\bullet$) are homotopic homogeneous polycomplexes, then so are $A_1^\bullet\otimes A_2^\bullet$ and $B_1^\bullet\otimes B_2^\bullet$.
\end{bigtheorem}

If all differentials are even, this recovers the usual Koszul rule for chain complexes in linear 2-categories.
Here the notion of homotopy equivalence is the usual notion of homotopy equivalence. Indeed, a homogeneous polycomplex is a genuine chain complex; it is only when taking the tensor product that its ``super nature'' plays a role.

In \cref{sec:topo}, we detail the construction of a tangle invariant using super $\glt$-foams and show the following theorem, which is the main result of this paper:

\begin{bigtheorem}[\cref{thm:equivalence_with_odd}]
  \label{bigthm:equivalence_with_odd}
  Our construction coincides with odd Khovanov homology when restricted to links.
\end{bigtheorem}

%%%%%%%%%%%%%%%%%%%%%%%%%%%%%%%
%%%          schur          %%%
%%%%%%%%%%%%%%%%%%%%%%%%%%%%%%%
\subsubsection{Skew Howe duality}

Let $\bC_q\coloneqq\bC(q)$, and consider the space ${\textstyle \bigwedge^d}(\bC_q^n\otimes\bC_q^2)$ equipped with the actions of $U_q(\mathfrak{gl}_n)$ and $U_q(\glt)$. These actions commute, and in fact each describes the other's intertwiners. This fact is known as \emph{skew Howe duality}, and was used (generalizing from $\glt$ to $\mathfrak{gl}_k$) by Cautis, Kamnitzer and Morrison to describe $U_q(\mathfrak{gl}_k)$-intertwiners on fundamental representations by generators and relations \cite{CKM_WebsQuantumSkew_2014}.
In particular, the map
\[
\phi\colon U_q(\mathfrak{gl}_n)\to\End_{U_q(\glt)}\left({\textstyle \bigwedge^d}(\bC_q^n\otimes\bC_q^2)\right)
\]
is surjective.
Taking the quotient by the kernel of this map defines the \emph{$q$-Schur algebra of level 2}:
\[
S_{n,d}\coloneqq U_q(\mathfrak{gl}_n)/\ker \phi\cong \End_{U_q(\glt)}\left({\textstyle \bigwedge^d}(\bC_q^n\otimes\bC_q^2)\right).
\]
While the above isomorphism holds over $\bC(q)$, the two sides have distinct interpretations, and hence different natural choices of integral idempotent form.
On the one hand, $U_q(\glt)$-inter\-twiners of fundamental representations are generated by the following intertwiners:
\begin{gather*}
  \begin{tikzpicture}
    % \draw[web2,directed] (0,0) to (-.8,0);
    % \draw[web1,directed=.7] (.8,-.5) to[out=180,in=-90] (0,0);
    % \draw[web1,directed=.7] (.8,.5) to[out=180,in=90] (0,0);
    \draw[web2] (0,0) to (-.8,0);
    \draw[web1] (.8,-.5) to[out=180,in=-90] (0,0);
    \draw[web1] (.8,.5) to[out=180,in=90] (0,0);
    \node[left] at (-.8,0) {$\bigwedge^2(\bC_q^2)$};
    \node[right] at (.8,-.5) {$\bigwedge^1(\bC_q^2)$};
    \node[right] at (.8,.5) {$\bigwedge^1(\bC_q^2)$};
    \node[right=10pt] at (.8,0) {$\otimes$};
    \begin{scope}[yshift=-1.5cm]
      \node[right] at (.8,0) {$v_1\otimes v_2$};
      \node[left] at (-.8,0) {$v_1\wedge v_2$};
      \node[rotate=180] at (0,0) {$\mapsto$};
    \end{scope}
  \end{tikzpicture}
  \mspace{80mu}
  \begin{tikzpicture}
    % \draw[web2,rdirected] (0,0) to (.8,0);
    % \draw[web1,rdirected=.7] (-.8,-.5) to[out=0,in=-90] (0,0);
    % \draw[web1,rdirected=.7] (-.8,.5) to[out=0,in=90] (0,0);
    \draw[web2] (0,0) to (.8,0);
    \draw[web1] (-.8,-.5) to[out=0,in=-90] (0,0);
    \draw[web1] (-.8,.5) to[out=0,in=90] (0,0);
    \node[right] at (.8,0) {$\bigwedge^2(\bC_q^2)$};
    \node[left] at (-.8,-.5) {$\bigwedge^1(\bC_q^2)$};
    \node[left] at (-.8,.5) {$\bigwedge^1(\bC_q^2)$};
    \node[left=10pt] at (-.8,0) {$\otimes$};
    \begin{scope}[yshift=-1.5cm]
      \node[right] at (.8,0) {$v_1\wedge v_2$};
      \node[left] at (-.8,0) {$qv_1\otimes v_2-v_2\otimes v_1$};
      \node[rotate=180] at (0,0) {$\mapsto$};
    \end{scope}
  \end{tikzpicture}
\end{gather*}
These are the (spider) webs discussed before, leading to the integral and idempotent form $\web_d$ (recall that a small linear category is the same data as an algebra with a choice of orthogonal idempotents).
On the other hand, $U_q(\mathfrak{gl}_n)$ admits an integral and idempotent form denoted $\dot{U}_q(\mathfrak{gl}_n)$ (\emph{Lusztig's idempotent form}), which also provides an integral and idempotent form $\qschur_{n,d}$ for the $q$-Schur algebra of level 2.
Switching back to the categorical terminology, these $\bZ[q,q^{-1}]$-linear categories are related by $\bZ[q,q^{-1}]$-linear functors:
\[
  \dot{U}_q(\mathfrak{gl}_n)\to \qschur_{n,d}\to \web_d.
\]
The relationship between $\qschur_{n,d}$ and $\web_d$ is diagrammatically described by so-called \emph{ladder webs}; see \cite{CKM_WebsQuantumSkew_2014} for details.

%%%%%%%%%%%%%%%%%%%%%%%%%%%%%%%
%%%          schur          %%%
%%%%%%%%%%%%%%%%%%%%%%%%%%%%%%%
\subsubsection{A supercategorification of the \texorpdfstring{$q$}{q}-Schur algebra of level 2}

As shown by Lauda, Queffelec and Rose, skew Howe duality can be categorified:
% https://q.uiver.app/#q=WzAsNixbMCwxLCJcXGRvdHtVfV9xKFxcbWF0aGZyYWt7Z2x9X24pIl0sWzEsMSwiXFxkb3R7U31fe24sZH0iXSxbMiwxLCJcXHdlYl9kIl0sWzIsMCwiXFxtYXRoYmZ7Rm9hbX1fZCJdLFsxLDAsIlxcY1Nfe24sZH0iXSxbMCwwLCJcXGNVKFxcbWF0aGZyYWt7Z2x9X24pIl0sWzAsMV0sWzEsMl0sWzMsMiwiS18wIiwwLHsic3R5bGUiOnsiYm9keSI6eyJuYW1lIjoic3F1aWdnbHkifX19XSxbNSwwLCJLXzAiLDAseyJzdHlsZSI6eyJib2R5Ijp7Im5hbWUiOiJzcXVpZ2dseSJ9fX1dLFs0LDEsIktfMCIsMCx7InN0eWxlIjp7ImJvZHkiOnsibmFtZSI6InNxdWlnZ2x5In19fV0sWzQsM10sWzUsNF1d
\[\begin{tikzcd}[ampersand replacement=\&,cramped]
	{\cU(\mathfrak{gl}_n)} \& {\cS_{n,d}} \& {\foam_d} \\
	{\dot{U}_q(\mathfrak{gl}_n)} \& {\qschur_{n,d}} \& {\web_d}
	\arrow[from=1-1, to=1-2]
	\arrow["{K_0}", squiggly, from=1-1, to=2-1]
	\arrow[from=1-2, to=1-3]
	\arrow["{K_0}", squiggly, from=1-2, to=2-2]
	\arrow["{K_0}", squiggly, from=1-3, to=2-3]
	\arrow[from=2-1, to=2-2]
	\arrow[from=2-2, to=2-3]
\end{tikzcd}\]
The top arrows are 2-functors, $\cU(\mathfrak{gl}_n)$ is the \emph{categorified quantum $\mathfrak{gl}_n$} as introduced by Khovanov--Lauda and Rouquier \cite{KL_CategorificationQuantum$sln$_2010,Rouquier_2KacMoodyAlgebras_2008}, and $\cS_{n,d}$ is the \emph{2-Schur algebra} as introduced by Mackaay--Stošić--Vaz \cite{MSV_DiagrammaticCategorification$q$Schur_2013}.
The composition of the two top arrows is called the \emph{foamation 2-functor} in \cite{LQR_KhovanovHomologySkew_2015}, exhibiting $\foam_d$ as a 2-representation of $\cU(\mathfrak{gl}_n)$. In that sense, it provides an understanding of Khovanov homology in terms of higher representation theory.
As this 2-representation factors through $\cS_{n,d}$, we may equivalently view $\foam_d$ as a 2-representation of $\cS_{n,d}$, and define the \emph{foamation 2-functor} as the 2-functor ${\cS_{n,d}} \to {\foam_d}$.

In \cref{sec:schur}, we give an odd analogue of this latter result.
In \cite{Vaz_NotEvenKhovanov_2020}, the second author defined a superalgebra which can be seen as an odd analogue of the KLR algebra \cite{KL_DiagrammaticApproachCategorification_2009,Rouquier_2KacMoodyAlgebras_2008} of level 2 for the $A_n$ quiver. Taking a cyclotomic quotient of this construction leads to a supercategorification of the negative half of the $q$-Schur algebra of level 2.
In \cref{sec:schur}, we extend it to a supercategorification of the $q$-Schur algebra of level 2, giving an odd analogue of \cite{MSV_DiagrammaticCategorification$q$Schur_2013} in the level-2 case.
We then define a \emph{super foamation 2-functor}, fitting into the following commutative square:
% https://q.uiver.app/#q=WzAsNCxbMCwxLCJcXGRvdHtTfV97bixkfSJdLFsxLDEsIlxcd2ViX2QiXSxbMSwwLCJcXG1hdGhiZntGb2FtfV9kIl0sWzAsMCwiXFxjU197bixkfSJdLFswLDFdLFsyLDEsIktfMCIsMCx7InN0eWxlIjp7ImJvZHkiOnsibmFtZSI6InNxdWlnZ2x5In19fV0sWzMsMCwiS18wIiwwLHsic3R5bGUiOnsiYm9keSI6eyJuYW1lIjoic3F1aWdnbHkifX19XSxbMywyXV0=
\[\begin{tikzcd}[ampersand replacement=\&,cramped]
	{\cS\cS_{n,d}} \& {\sfoam_d} \\
	{\qschur_{n,d}} \& {\web_d}
	\arrow[from=1-1, to=1-2]
	\arrow["{K_0}", squiggly, from=1-1, to=2-1]
	\arrow["{K_0}", squiggly, from=1-2, to=2-2]
	\arrow[from=2-1, to=2-2]
\end{tikzcd}\]

In a nutshell:

\begin{bigtheorem}
  \label{bigthm:schur_and_OKh}
  The supercategorification of the $q$-Schur algebra of level 2 together with the super foamation 2-functor provides a representation-theoretic construction of odd Khovanov homology.
\end{bigtheorem}

% \subsubsection{relationship with Heegaard--Floer homology}

\section{Graded \texorpdfstring{$\glt$}{gl2}-foams}
\label{sec:foams}

This section defines the graded-2\nbd-cate\-gory $\gfoam_d$ of graded $\glt$-foams that categorifies the category $\web_d$ of $\glt$-webs. This can be viewed as a graded analogue of $\glt$-foams \emph{à la} Blanchet \cite{Blanchet_OrientedModelKhovanov_2010}. See also \cite{BHP+_Functoriality$mathfraksl_2$Tangle_2023,EST_Generic$mathfrakgl_2$foamsWeb_2020,EST_BlanchetKhovanovAlgebras_2017} for further studies of $\glt$-foams.
More precisely, the graded-2\nbd-cate\-gory $\gfoam_d$ gives a graded analogue of the category of \emph{directed} $\glt$-foams as defined in \cite{QR_$mathfraksl_n$Foam2category_2016} (see \cref{rem:chfoam_recover_even_foam}).

\medbreak

\Cref{subsec:graded_two_cat} defines the notion of graded-2-categories. The category $\web_d$ and the graded-2\nbd-cate\-gory $\gfoam_d$ are then respectively defined in \cref{subsec:webs,subsec:graded_foams}.
In \cref{subsec:string_diag_foam}, we define a string diagrammatics for graded $\glt$-foams; see also \cite{Schelstraete_RewritingModuloDiagrammatic_2026,Schelstraete_OddKhovanovHomology_2024} for another diagrammatics using shadings.
This gives a computation-friendly counterpart of the topological definition given in \cref{subsec:graded_foams}.
Finally, \cref{sec:categorification_foam} shows that $\gfoam_d$ categorifies $\web_d$, relying on the main result of \cite{Schelstraete_RewritingModuloDiagrammatic_2026}.

%%%%%%%%%%%%%%%%%%%%%%%%%%%%%%%%%%%%%%%%%%%
%%%          2-supercategories          %%%
%%%%%%%%%%%%%%%%%%%%%%%%%%%%%%%%%%%%%%%%%%%
\subsection{Graded-2-categories}
\label{subsec:graded_two_cat}

We define the notion of graded-2-cate\-go\-ries to allow gradings with generic abelian\footnote{We take groups to be abelian for simplicity; this article does not investigate the case of more general groups.} groups. Taking this abelian group to be $\bZ/2\bZ$ recovers the notion of super-2-cate\-gories~\cite{BE_SuperKacMoody_2017}.
A related definition appeared in~\cite{Putyra_2categoryChronologicalCobordisms_2014} in the context of algebras, under the name of \emph{chronological algebras}.
See also~\cite{BE_MonoidalSupercategories_2017} for an in-depth study of superstructures, including non-strict versions.
Our exposition is a direct generalization to the graded case of the exposition given in~\cite{BE_SuperKacMoody_2017} for super-2-categories (which they call \emph{2-supercategories}).

Throughout the section we fix an abelian group $G$ and a unital commutative ring $\ring$. We denote $\ring^\times$ the abelian group of invertible elements in $\ring$ equipped with the multiplicative structure.
We also fix a $\ring^\times$-valued pairing on $G$, that is, a $\bZ$-bilinear map
\[\bil\colon G\times G\to \ring^\times.\]

\begin{definition}
  \label{defn:bil_symmetric}
  A pairing $\bil$ is \emph{symmetric} if $\bil(g,h)\bil(h,g) = 1$ for all $g,h\in G$.  
\end{definition}

For completeness, this section does not assume $\bil$ to be symmetric. However, every example considered in this article will be for $\bil$ symmetric; furthermore, \cref{sec:complexes} is stated for $\bil$ symmetric, for simplicity.
Note that $\bil$ is automatically symmetric in the super case $G=\bZ/2\bZ$.

\medbreak

Recall that a $\ring$-module $V$ is said to be \emph{$G$\nbd-graded} if it is equipped with the data of a direct sum decomposition $V=\bigoplus_{g\in G}V_g$. If $v\in V_g$ for some $g\in G$, the vector $v$ is said to be \emph{homogeneous}; if furthermore $v\neq 0$, \emph{it has degree $g$}, which we write as $\deg(v)=g$.
Note that while the zero vector is homogeneous, it does not have a well-defined degree.

The hom-space $\Hom_\ring(V,W)$ between two $G$\nbd-graded $\ring$-modules inherits a structure of $G$\nbd-graded $\ring$-module, stating that a non-zero $\ring$-linear map
\[f\colon \bigoplus_{g\in G}V_g\to \bigoplus_{g\in G}W_g\]
is of degree $\deg f=h$ if $f(V_g)\subset W_{g+h}$ for all $g\in G$. If $f=0$ or if $\deg f=0$, we say that $f$ is \emph{degree-preserving}.

Denote by $\GMod$ the category of $G$\nbd-graded $\ring$-modules, and by $\GModev$ the wide subcategory\footnote{Recall that a subcategory is said to be \emph{wide} if it contains all objects.} of $\GMod$ consisting only of degree-preserving $\ring$-linear maps.
The category $\GMod$ admits a standard monoidal structure, defined on objects $V=\bigoplus_{g\in G}V_g$ and $W=\bigoplus_{g\in G}W_g$ by the formula
\[(V\otimes W)_g = \bigoplus_{\substack{(g_1,g_2)\in G\times G\\ g_1+g_2=g}}V_{g_1}\otimes_\ring W_{g_2},\]
and on morphisms $f$ and $g$ by the formula $(f\otimes g)(v\otimes w)=f(v)\otimes f(w)$.
However, one can give an alternative definition using the data of the grading and the bilinear form $\bil$:

\begin{definition}
\label{defn:graded_tensor_product}
  The \emph{$(G,\bil)$-graded tensor product} is defined on objects as $V\otimes_{G,\bil} W\coloneqq V\otimes W$ and on morphisms following the Koszul rule associated to $\bil$:
  \[ 
    (f\otimes_{G,\bil} g)(v\otimes_{G,\bil} w) = \bil(\deg v,\deg g) f(v)\otimes_{G,\bil} g(w).
  \]
\end{definition}

Equipped with this graded tensor product, $\GMod$ is \emph{not} in general a monoidal category. Indeed, morphisms respect the following \emph{graded interchange law}:
\begin{equation}
\label{eq:exchange_relation}
(f\otimes_{G,\bil} g)\circ (h\otimes_{G,\bil} k) = \bil(\deg g,\deg h)(f\circ h)\otimes_{G,\bil} (g\circ k).
\end{equation}
We now define the proper categorical structures that encompass this behaviour.

\begin{definition}
\label{defn:graded_cat}
  A \emph{$G$\nbd-graded $\ring$-linear category} is a $(\GModev,\otimes)$-enriched category.
  A \emph{$G$\nbd-graded $\ring$-linear functor} is a $(\GModev,\otimes)$-enriched functor.
\end{definition}
In other words, a $G$\nbd-graded $\ring$-linear category is a category such that each Hom is a $G$\nbd-graded $\ring$-module, and such that composition is $\ring$-bilinear and preserves the grading in the sense that $\deg(f\circ g)=\deg f+\deg g$. A $G$\nbd-graded $\ring$-linear functor is a functor between two $G$\nbd-graded $\ring$-linear categories that restricts to a degree-preserving $G$\nbd-graded $\ring$-linear map on Homspaces.

Denote by $\GCat$ the category of small $G$\nbd-graded $\ring$-linear categories and $G$\nbd-graded $\ring$-linear functors. For $\cA$ and $\cB$ two $G$\nbd-graded $\ring$-linear categories, their \emph{$(G,\bil)$-graded-cartesian tensor product} is the $G$\nbd-graded $\ring$-linear category $\cA \boxtimes_{G,\bil}\cB$ such that $\mathrm{ob}(\cA \boxtimes_{G,\bil}\cB)\coloneqq\mathrm{ob}(\cA)\times \mathrm{ob}(\cB)$ and
\[
  \Hom_{\cA \boxtimes_{G,\bil}\cB}((a,b),(c,d))\coloneqq\Hom_{\cA}(a,c)\otimes_{G,\bil}\Hom_{\cB}(b,d).
\]
Most importantly, the composition in $\cA \boxtimes_{G,\bil}\cB$ is given by the graded interchange relation \eqref{eq:exchange_relation}. This makes $\cA \boxtimes_{G,\bil}\cB$ a $G$\nbd-graded $\ring$-linear category. Moreover, one can check that $\boxtimes_{G,\bil}$ is associative and unital, giving $\GCat$ the structure of a monoidal category.

\begin{definition}
\label{defn:graded_2_cat}
  A \emph{$(G,\bil)$-graded-2\nbd-cate\-gory} is a $\GCat$-enriched category. A \emph{$(G,\bil)$-graded-2-functor} is a $\GCat$-enriched functor.
\end{definition}

In particular:

\begin{definition}
\label{defn:graded_monoidal_cat}
  A \emph{$(G,\bil)$-graded-monoidal category} is a $G$\nbd-graded $\ring$-linear category $\cC$ together with the data of a unit object and a unital and associative $G$\nbd-graded $\ring$-linear functor $\otimes_{G,\bil}\colon{\cC\boxtimes_{G,\bil}\cC}\to \cC$.
\end{definition}

For instance, the $(G,\bil)$-graded tensor product from \cref{defn:graded_tensor_product} assembles into a $G$\nbd-graded $\ring$-linear functor
\[\otimes_{G,\bil}\colon \GMod\boxtimes_{G,\bil}\GMod\to  \GMod,\]
making $(\GMod,\otimes_{G,\bil})$ a $(G,\bil)$-graded-monoidal category. (More precisely, we should consider the strictification of $\GMod$, where $\otimes_{G,\bil}$ has strict associativity and unitality.)
Note that $(\GMod,\otimes_{G,\bil})$ is \emph{not} a linear monoidal category with an extra $G$-grading; indeed, the interchange law \eqref{eq:exchange_relation} holds only ``up to scalars''.
In general, a $(G,\bil)$-graded-2-category is \emph{not} a linear 2-category with extra structure.
Indeed, the $G$-grading interacts with the 2-categorical structure, twisting the interchange law.
However, any 2\nbd-cate\-gory that is both $\ring$-linear and $G$\nbd-graded as a $\ring$-linear 2\nbd-cate\-gory can be seen as a $(G,\bil)$-graded-2\nbd-cate\-gory, with $\bil$ the trivial map.

\begin{remark}
  \label{rem:gray_categories}
  Graded-2-categories are closely related to Gray categories, certain 3-dimen\-sio\-nal categorical structures where the interchange law for 2\nbd-morphisms holds weakly.
  We use this point of view in \cite{Schelstraete_RewritingModuloDiagrammatic_2026} (see also \cite{Schelstraete_OddKhovanovHomology_2024}) to develop a rewriting theory suitable for graded-2-categories.
\end{remark}

\begin{remark}
  If $\mu$ is symmetric, the monoidal category $(\GModev,\otimes)$ has a symmetric structure given by $v\otimes w\mapsto \mu(\deg(v),\deg(w))w\otimes v$.
  In general, if $\cV$ is a symmetric monoidal category, then $\cV\md\cC at$, the category of small $\cV$\nbd-enriched categories and $\cV$\nbd-enriched functors, is itself symmetric monoidal (see e.g.\ \cite[Proposition~6.2.9]{Borceux_HandbookCategoricalAlgebra_1994a}).
  This allows an inductive definition of $\cV$\nbd-enriched $n$-categories.
  In that general framework, $(G,\mu)$-graded-2-categories are precisely
  $(\GModev)$-enriched 2-categories.
\end{remark}

\subsubsection{Diagrammatics}

Let $\cC$ be a $(G,\bil)$-graded-2\nbd-cate\-gory. Unpacking \cref{defn:graded_2_cat}, $\cC$ consists of objects $\cC_0$ together with a $G$\nbd-graded $\ring$-linear category $\cC(a,b)$ for each pair of objects $(a,b)$. We denote $\cC_1(a,b)$ its objects and for $f,g\in\cC_1(a,b)$, we denote $\cC_2(f,g)$ the Hom-space of morphisms from $f$ to $g$. Elements of $\cC_1(a,b)$ and $\cC_2(f,g)$ are respectively called 1\nbd-morphisms and 2\nbd-morphisms.
A 2\nbd-morphism $\alpha\in\cC_2(f,g)$ can be pictured using string diagrams, akin to the string diagrammatics of 2-categories:
\begin{equation*}
  \xy (0,0)* {
  \begin{tikzpicture}[scale=0.5]
      \draw   (0,2) node[above] {\footnotesize $g$} to
                  node[black_dot] {}
                  node[right=1pt] {\footnotesize $\alpha$ }
              (0,0) node[below] {\footnotesize $f$};
      \node at (1.5,1) {$a$};
      \node at (-1,1) {$b$};
  \end{tikzpicture} }\endxy
\end{equation*}
The vertical composition~$\starop_1$ denotes all the compositions in the $G$\nbd-graded $\ring$-linear categories $\cC(a,b)$. It is pictured by stacking the 2\nbd-morphisms atop each other:
\begin{equation*}
  % vertical composition
  {}\xy (0,0)* {
  \begin{tikzpicture}[scale=0.5]
      \draw (0,2) node[above] {\footnotesize $h$} to
          node[black_dot,near start] {}
          node[right,near start] {\footnotesize $\beta$}
          node[right=-2pt] {\footnotesize $g$}
          node[black_dot,near end] {}
          node[right,near end] {\footnotesize $\alpha$}
          (0,0) node[below] {\footnotesize $f$};
          \node at (1.5,1) {$a$};
          \node at (-1,1) {$b$};
  \end{tikzpicture} }\endxy
  \;\coloneqq\;
  \left(
  \xy (0,0)* {
  \begin{tikzpicture}[scale=0.5]
      \draw   (0,2) node[above] {\footnotesize $h$} to
                  node[black_dot] {}
                  node[right=1pt] {\footnotesize $\beta$ }
              (0,0) node[below] {\footnotesize $g$};
      \node at (1.5,1) {$a$};
      \node at (-1,1) {$b$};
  \end{tikzpicture} }\endxy
  \right)
  \starop_1
  \left(
  \xy (0,0)* {
  \begin{tikzpicture}[scale=0.5]
      \draw   (0,2) node[above] {\footnotesize $g$} to
                  node[black_dot] {}
                  node[right=1pt] {\footnotesize $\alpha$ }
              (0,0) node[below] {\footnotesize $f$};
      \node at (1.5,1) {$a$};
      \node at (-1,1) {$b$};
  \end{tikzpicture} }\endxy
  \right)
\end{equation*}
The horizontal composition~$\starop_0$ comes from the enrichment. It is pictured by putting the two 2\nbd-mor\-phisms side-by-side:
\begin{equation*}
  % horizontal composition
  \xy (0,0)* {
    \begin{tikzpicture}[scale=0.65]
        \draw (0,2) node[above] {\footnotesize $g'$} to
            node[black_dot] {}
            node[right] {\footnotesize $\beta$}
            (0,0) node[below] {\footnotesize $f'$};
        \draw (1,2) node[above] {\footnotesize $g$} to
            node[black_dot] {}
            node[right] {\footnotesize $\alpha$}
            (1,0) node[below] {\footnotesize $f$};
            \node at (2,1) {$a$};
            \node at (-1,1) {$c$};
    \end{tikzpicture} }\endxy
    \;\coloneqq\;
    \left(
  \xy (0,0)* {
  \begin{tikzpicture}[scale=0.5]
      \draw   (0,2) node[above] {\footnotesize $g'$} to
                  node[black_dot] {}
                  node[right=1pt] {\footnotesize $\beta$ }
              (0,0) node[below] {\footnotesize $f'$};
      \node at (1.5,1) {$b$};
      \node at (-1,1) {$c$};
  \end{tikzpicture} }\endxy
  \right)
  \starop_0
  \left(
  \xy (0,0)* {
  \begin{tikzpicture}[scale=0.5]
      \draw   (0,2) node[above] {\footnotesize $g$} to
                  node[black_dot] {}
                  node[right=1pt] {\footnotesize $\alpha$ }
              (0,0) node[below] {\footnotesize $f$};
      \node at (1.5,1) {$a$};
      \node at (-1,1) {$b$};
  \end{tikzpicture} }\endxy
  \right)
\end{equation*}
To avoid clutter, we leave regions and lines unlabelled for now.
The rule ``compose first horizontally, then vertically'' resolves the ambiguity for the order of compositions:
\begin{equation*}
  \xy (0,0)* {
    \begin{tikzpicture}[scale=0.5]
        \draw (0,2.5) to
            node[black_dot,near end] {}
            node[right,near end] {\footnotesize $\gamma$}
            node[black_dot,near start] {}
            node[right,near start] {\footnotesize $\delta$}
            (0,0);
        \draw (1,2.5) to
            node[black_dot,near end] {}
            node[right,near end] {\footnotesize $\alpha$}
            node[black_dot,near start] {}
            node[right,near start] {\footnotesize $\beta$}
            (1,0);
    \end{tikzpicture} }\endxy
    \;\coloneqq\;
   (\delta\starop_0\beta)\starop_1(\gamma\starop_0\alpha)
\end{equation*} 
Contrary to 2-categories, it is in general \emph{not} the same as ``compose first vertically, then horizontally''. Indeed, in a graded-2\nbd-cate\-gory sliding a 2\nbd-morphism past another vertically comes at the price of an invertible scalar. This is the graded interchange law:
\begin{equation}
  \label{eq:gradedcat_interchange}
  {}\xy (0,0)* {
    \begin{tikzpicture}[scale=0.5]
        \draw (0,2.5) to
            node[black_dot] {}
            node[right] {\footnotesize $\beta$}
            (0,0);
        \draw (1,2.5) to
            node[black_dot] {}
            node[right] {\footnotesize $\alpha$}
            (1,0);
            % \node at (2,1) {$a$};
            % \node at (-1,1) {$b$};
    \end{tikzpicture} }\endxy
    \;=\;
    \xy (0,0)* {
    \begin{tikzpicture}[scale=0.5]
        \draw (0,2.5) to
            node[black_dot,pos=.3] {}
            node[right,pos=.3] {\footnotesize $\beta$}
            (0,0);
        \draw (1,2.5) to
            node[black_dot,pos=.7] {}
            node[right,pos=.7] {\footnotesize $\alpha$}
            (1,0);
            % \node at (2,1) {$a$};
            % \node at (-1,1) {$b$};
    \end{tikzpicture} }\endxy
    \;=\;\bil(\deg\beta,\deg\alpha)\;
    \xy (0,0)* {
    \begin{tikzpicture}[scale=0.5]
        \draw (0,2.5) to
            node[black_dot,pos=.7] {}
            node[right,pos=.7] {\footnotesize $\beta$}
            (0,0);
        \draw (1,2.5) to
            node[black_dot,pos=.3] {}
            node[right,pos=.3] {\footnotesize $\alpha$}
            (1,0);
            % \node at (2,1) {$a$};
            % \node at (-1,1) {$b$};
    \end{tikzpicture} }\endxy
\end{equation}
In particular, one must be careful with the relative vertical positions of 2\nbd-mor\-phisms. To avoid confusion, in this article we always draw string diagrams such that no two non-trivially graded 2\nbd-mor\-phisms lie on the same vertical level. Trivially graded 2-morphisms can safely be drawn at the same vertical level, as those can slide vertically without adding scalars.
Finally, as customary already in the string diagrammatics of 2\nbd-categories, we usually do not picture identities.

If $\bil$ is symmetric (\cref{defn:bil_symmetric}), the scalar appearing in the graded interchange law depends only on the relative vertical positions of the 2\nbd-morphisms:
\begin{equation*}
  \xy (0,0)* {
  \begin{tikzpicture}[scale=0.5]
      \draw (0,2.5) to
          node[black_dot,pos=.7] {}
          node[right,pos=.7] {\footnotesize $\alpha$}
          (0,0);
      \draw (1,2.5) to
          node[black_dot,pos=.3] {}
          node[right,pos=.3] {\footnotesize $\beta$}
          (1,0);
          % \node at (2,1) {$a$};
          % \node at (-1,1) {$b$};
  \end{tikzpicture} }\endxy
  \;=\;\bil(\deg\beta,\deg\alpha)\;
  \xy (0,0)* {
  \begin{tikzpicture}[scale=0.5]
      \draw (0,2.5) to
          node[black_dot,pos=.3] {}
          node[right,pos=.3] {\footnotesize $\alpha$}
          (0,0);
      \draw (1,2.5) to
          node[black_dot,pos=.7] {}
          node[right,pos=.7] {\footnotesize $\beta$}
          (1,0);
          % \node at (2,1) {$a$};
          % \node at (-1,1) {$b$};
  \end{tikzpicture} }\endxy
  \qquad
  \text{ if $\mu$ is symmetric.}
\end{equation*}
Compare with \cref{eq:gradedcat_interchange}. In other words, passing $\beta$ down $\alpha$ always adds $\bil(\deg\beta,\deg\alpha)$, regardless of their respective horizontal positions.

\subsubsection{Grothendieck ring}
\label{subsec:more_cat_facts}

Let $\cC$ be a $(G,\bil)$-graded-2\nbd-cate\-gory equipped with an extra $H$-grading, where $H$ is an abelian group.
In this subsection, we define the Grothendieck ring of $\cC$ with respect to $H$, denoted $K_0(\cC)\vert_H$.
This is analogous to the usual notion of Grothendieck ring for an $H$-graded linear 2-category.
Crucially, we consider the $H$-grading as independent of the $G$-grading, at least on a formal level; this implies that taking $H$-envelopes does not affect the $(G,\bil)$-graded interchange law.

In the remainder of the article, we will deal with $H$-graded $(G,\mu)$-graded-2-categories where $G=\bZ\times\bZ$ and $H=\bZ$.
In this context, the $H$-grading is called the ``quantum grading'', and we write
\[
K_0(\cC)\vert_q \coloneqq K_0(\cC)\vert_H.
\]

\bigbreak
\noindent \emph{1-dimensional case: $H$-graded categories}
\smallbreak

We review some basic notions; see also \cite[chap.~11]{EMT+_IntroductionSoergelBimodules_2020}.

\noindent Let $H$ be an abelian group.
An \emph{$H$-graded category} is a category enriched over $H$-graded abelian groups. If $H=\{*\}$ is the trivial group, the category is said to be \emph{pre-additive}.
A pre-additive category $C$ is a \emph{category with $H$-shifts} if it is equipped with a group morphism $H\to\cA ut(C)$ where $\cA ut(C)$ is the group of invertible endofunctors; that is, $C$ is equipped with an action of $H$ by autofunctors.
We denote $f\{x\}$ the action of $x\in H$ on an object $f\in\ob(C)$.
When $H=\bZ$, we also write $q^n f\coloneqq f\{n\}$.

\medbreak

Let $C$ be a pre-additive category. The category $C$ naturally embeds in an additive category $C^{\oplus}$, its \emph{additive closure}. Objects in $C^{\oplus}$ are formal direct sums of objects in $C$, and morphisms are matrices whose entries are morphisms in $C$.
If $C$ is a category with $H$-shifts, then so is $C^{\oplus}$ by extending the $H$-action additively.

Let $C$ be an additive category.
The \emph{(split) Grothendieck ring of $C$} is the abelian group $K_0(C)$ generated by elements $[f]$ for each object $f\in\ob(C)$, and subject to the relations $[f\oplus g]=[f]+[g]$; in particular, if $f\cong g$ then $[f]=[g]$.
If $C$ is with $H$-shifts, then $K_0(C)$ has the structure of a $\bZ[H]$-module, where the $H$-action is given by $x\cdot_H [f]\coloneqq [f\{x\}]$.

Let $C$ be a pre-additive category. The \emph{Grothendieck ring of $C$} is the Grothendieck ring of its additive closure: $K_0(C)\coloneqq K_0(C^\oplus)$.

\medbreak

Let $C$ be an $H$-graded category.
The \emph{$H$-envelope of $C$} is the category $C_H$ whose objects are formal \mbox{$H$-shifts} $f\{x\}$, where ${x\in H}$ and ${f\in\ob(C)}$, and there is a morphism
\begin{gather*}
  \alpha\colon f\{x\}\to g\{y\}
\end{gather*}
in $C_H$ for each morphism ${\alpha\colon f\to g}$ in $C$ and each pair of elements $x,y\in H$.
If $\deg_H \alpha$ denotes the $H$-degree of $\alpha$, we set
$
  \deg_H\left(\alpha\colon f\{x\}\to g\{y\}\right) =\deg_H \alpha-x+y.
$
The $H$-envelope is both $H$-graded and with $H$-shifts.
  
Let $C$ be an $H$-graded category.
The \emph{underlying category of $C$} is the category $\ul{C}$ consisting of degree-preserving morphisms.
The \emph{$H$-shifted closure of $C$}, denoted by $\ul{C}_H$, is the underlying category of its $H$-envelope; note that it has $H$-shifts.
Finally, the \emph{$H$-Grothendieck ring of $C$} is the Grothendieck ring of its $H$-shifted closure $\ul{C}_H$:
\begin{gather*}
  K_0(C)\vert_H\coloneqq K_0(\ul{C}_H)=K_0((\ul{C}_H)^\oplus).
\end{gather*}
As before, the $H$-Grothendieck ring is a $\bZ[H]$-module.

\medbreak

All of the above extend to the case where $C$ is an $(G\times H)$-graded category, ignoring the $G$\nbd-grading throughout. Note that the $H$-envelope is canonically $G$\nbd-graded, setting
\begin{gather*}
  \deg_G\left(\alpha\colon f\{x\}\to g\{y\}\right)=\deg_G\alpha.
\end{gather*}

\bigbreak
\noindent\textit{2-dimensional case: $H$-graded $(G,\bil)$-graded-2-categories}
\smallbreak

\noindent 
Let $\cC$ be a $(G,\bil)$-graded-2\nbd-cate\-gory and $H$ an abelian group.
We say that $\cC$ is \emph{additive} if each hom-category $\cC(a,b)$ is additive and horizontal composition is bilinear.
We say that $\cC$ is \emph{with $H$\nbd-shifts} if it is equipped with an action of $H$ by automorphisms which is the identity on objects. In particular, for each pair of objects $a,b\in\ob(\cC)$ the $G$\nbd-graded category $\cC(a,b)$ is with $H$\nbd-shifts. Its Grothendieck ring is the $\bZ[H]$-linear category obtained by taking the Grothendieck ring of each hom-categories $\cC(a,b)$, and inducing a bilinear composition by setting $[f]\circ [g]\coloneqq [f\starop_0 g]$.

Let $\cC$ be an $H$\nbd-graded $(G,\bil)$-graded-2\nbd-cate\-gory.
The \emph{$H$\nbd-envelope of $\cC$} is the $(G,\bil)$-graded-2-cate\-gory with \mbox{$H$-shifts} defined by taking the $H$-envelope of each hom-category $\cC(a,b)$.
The \emph{underlying linear 2\nbd-cate\-gory of  $\cC$} is the sub-2\nbd-cate\-gory $\ul{\cC}$ obtained by taking the underlying category of each hom-category $\cC(a,b)$.
Notions of $H$\nbd-shifted closure and additive closure are defined analogously to the 1-dimensional case. Finally:

\begin{definition}
  Let $\cC$ be an $H$-graded $(G,\bil)$-graded-2\nbd-cate\-gory.
  Its \emph{$H$\nbd-Gro\-then\-dieck ring} is the Grothendieck ring of the additive closure of its $H$\nbd-shifted closure:
  \begin{gather*}
    K_0(\cC)\vert_H\coloneqq K_0((\ul{\cC}_H)^\oplus).
  \end{gather*}
  The $H$\nbd-Grothendieck ring is a $\bZ[H]$-linear category.
\end{definition}

\begin{remark}
  \label{rem:G_envelope}
  One could give a more general definition, choosing to decategorify with respect to the $G$-grading as well, and not only some (formally) distinct $H$-grading.
  This generalizes the notion of Grothendieck ring for super-2-categories given in \cite[Definition~1.16]{BE_MonoidalSupercategories_2017}.
  The main difficulty lies in properly defining the horizontal composition in the $G$-envelope, ensuring that the graded interchange law holds.
  For completeness and future reference, we give below the definition of the horizontal composition in the $G$-envelope of a $(G,\mu)$-graded-2\nbd-cate\-gory, following the conventions of \cite{BE_MonoidalSupercategories_2017}:
  \begin{gather}
    \nonumber
    {}\xy(0,0)*{\begin{tikzpicture}[scale=.5]
        \draw (0,-1) to (0,1);
        \node[black_dot] at (0,0){};
        \node[right,thick] at (0,0) {$\alpha$};
        \draw[draw=red,thick] (-1,-1) to (1,-1);
        \draw[draw=red,thick] (-1,1) to (1,1);
        \node[right] at (1,1) {$y$};
        \node[right] at (1,-1) {$x$};
    \end{tikzpicture}}\endxy
    \starop_0
    \xy(0,0)*{\begin{tikzpicture}[scale=.5]
      \draw (0,-1) to (0,1);
      \node[black_dot] at (0,0){};
      \node[right,thick] at (0,0) {$\beta$};
      \draw[draw=red,thick] (-1,-1) to (1,-1);
      \draw[draw=red,thick] (-1,1) to (1,1);
      \node[right] at (1,1) {$w$};
      \node[right] at (1,-1) {$z$};
  \end{tikzpicture}}\endxy
  % \mspace{300mu}
  % \\
  % \mspace{60mu}
  =\mu(-x,\deg_G \beta-z+w)\mu(\deg_G \alpha,w)
  \xy(0,0)*{\begin{tikzpicture}[scale=.5]
    \begin{scope}[xshift=-1.5cm]
      \draw (0,-1) to (0,1);
      \node[black_dot] at (0,0){};
      \node[right,thick] at (0,0) {$\alpha$};
      \draw[draw=red,thick] (-1,-1) to (1,-1);
      \draw[draw=red,thick] (-1,1) to (1,1);
    \end{scope}
    \draw (0,-1) to (0,1);
    \node[black_dot] at (0,0){};
    \node[right,thick] at (0,0) {$\beta$};
    \draw[draw=red,thick] (-1,-1) to (1,-1);
    \draw[draw=red,thick] (-1,1) to (1,1);
    \node[right] at (1,1) {$y+w$};
    \node[right] at (1,-1) {$x+z$};
  \end{tikzpicture}}\endxy
  % \label{eq:graded_closure_hor_comp}
  \end{gather}
\end{remark}

%%%%%%%%%%%%%%%%%%%%%%%%%%%%%%
%%%          webs          %%%
%%%%%%%%%%%%%%%%%%%%%%%%%%%%%%
\subsection{\texorpdfstring{$\mathfrak{gl}_2$}{gl2}-webs}
\label{subsec:webs}

In our context, a $\glt$-web is a certain trivalent graph, smoothly embedded in $\bR\times[0,1]$, that we view as a morphism in a certain category defined below.
Objects, called \emph{weights}, are elements of the following set:
\begin{equation}
  \label{eq:defn_uLam}
  \uLam_d\coloneqq\bigsqcup_{k\in\bN}\{\lambda\in\{1,2\}^k\mid \lambda_1+\ldots+\lambda_k=d\}.
\end{equation}
For each $\lambda\in\uLam_d$ with $k$ coordinates, we can define a label on its coordinates
\[l_\lambda\colon\{1,\ldots,k\}\to\{1,\ldots,d\}\]
by setting $l_\lambda(i) = \sum_{j<i}\lambda_i+1$. For instance, $l_{(1,1,2,1)}=(1,2,3,5)$. Foreseeing the string diagrammatics, we call this label the \emph{colour} of the coordinate.
The identity web of a weight is pictured as a juxtaposition of straight vertical lines in $\bR\times[0,1]$, decorated as \emph{single} (black) or \emph{double} (orange) lines:
\begin{equation*}
  \id_{(1,1,2,1)} \quad=\quad
  \xy(0,-4)*{\begin{tikzpicture}[yscale=.5]
    \draw[web1,directed] (0,0)node[below]{\tiny $1$} to ++(-2,0);
    \draw[web1,directed] (0,1)node[below]{\tiny $2$} to ++(-2,0);
    \draw[web2,directed] (0,2)node[below]{\tiny $3$} to ++(-2,0);
    \draw[web1,directed] (0,3)node[below]{\tiny $5$} to ++(-2,0);
    \node at (1,0){\small $1$};
    \node at (1,1){\small $1$};
    \node at (1,2){\small $2$};
    \node at (1,3){\small $1$};
    \node at (-3,0){\small $1$};
    \node at (-3,1){\small $1$};
    \node at (-3,2){\small $2$};
    \node at (-3,3){\small $1$};
    \node (A) at (1,-1) {\small $\bR\times\{0\}$};
    \node (B) at (-3,-1) {\small $\bR\times\{1\}$};
    \draw[->] (2,1) to (2,2) node[above right]{$\bR$};
  \end{tikzpicture}}\endxy
\end{equation*}
Note that we read webs from right to left. Note also that each line has a colour, given by the colour of the corresponding coordinate (this colour is unrelated to whether the line is pictured ``black'' or ``orange''). Effectively, the colour of a line counts the number of preceding lines (plus one), counting twice double lines.
A generic \emph{$\glt$-web} is either an identity web or a composition of the following generating webs:
\begin{equation*}
\label{eq:generating_webs}
  W_{i,-}=
  \xy(0,0)*{\begin{tikzpicture}[scale=.7]
    \draw[web2,directed] (.7,0) node[below]{\scriptsize $i$} to (0,0);
    \draw[web1,directed=.7] (0,0) to[out=-90,in=0] (-.7,-.5)
      node[below]{\scriptsize $i$};
    \draw[web1,directed=.7] (0,0) to[out=90,in=0] (-.7,.5)
      node[above]{\scriptsize $i+1$};
    \node at (0,-1) {$\vdots$};
    \node at (0,1) {$\vdots$};
  \end{tikzpicture}}\endxy
  \qquad\an\qquad
  W_{i,+}=
  \xy(0,0)*{\begin{tikzpicture}[scale=.7]
    \draw[web2,rdirected] (-.7,0) node[below]{\scriptsize $i$} to (0,0);
    \draw[web1,rdirected=.7] (0,0) to[out=-90,in=180] (.7,-.5)
      node[below]{\scriptsize $i$};
    \draw[web1,rdirected=.7] (0,0) to[out=90,in=180] (.7,.5)
      node[above]{\scriptsize $i+1$};
    \node at (0,-1) {$\vdots$};
    \node at (0,1) {$\vdots$};
  \end{tikzpicture}}\endxy
\end{equation*}
Hereabove the dots denote possibly additional vertical single or double lines, and composition is given by stacking webs atop each other (reading from right to left).
Note that a web $W$ has an underlying unoriented flat tangle diagram, denoted $\undcomp(W)$, given by forgetting the double lines and the orientations.

\begin{remark}[Orientation]
  \label{rem:directed_webs}
  Our webs are given an orientation ``flowing'' from $\bR\times \{0\}$ to $\bR\times \{1\}$, which is more restrictive than some definitions in the literature.
  Such webs are sometimes called \emph{acyclic} or \emph{left-directed} (e.g.\ in \cite{QR_$mathfraksl_n$Foam2category_2016}).
  As this orientation is canonical, we often omit it.
\end{remark}

\begin{definition}
\label{defn:webs}
  The category $\web_d$ has objects $\uLam_d$, and morphisms are $\bZ[q,q^{-1}]$-linear combinations of $\glt$-webs, up to the following \emph{web relations}:
  \begin{IEEEeqnarray*}{CcC}
    \stackanchor{W_{i,s_1}W_{j,s_2}=W_{j,s_2}W_{i,s_1}}{(\text{for all }s_1,s_2\in\{-,+\}\an \abs{i-j}>1)}
    &\mspace{40mu}&
    \stackanchor{\text{\small interchange}}{\text{\small relations}}
    \\[1ex]
    {}\xy(0,0)*{\begin{tikzpicture}[scale=.6]
      \draw[web2] (0,0) to (.5,0);
      \draw[web2] (1.5,0) to (2,0);
      \draw[web1,rdirected] (.5,0) to[out=90,in=180] ++(.5,.5) to[out=0,in=90] ++(.5,-.5);
      \draw[web1,rdirected] (.5,0) to[out=-90,in=180] ++(.5,-.5) to[out=0,in=-90] ++(.5,.5);
    \end{tikzpicture}}\endxy
    \mspace{5mu}=\mspace{5mu}
    (q+q^{-1})\mspace{5mu}
    \xy(0,0)*{\begin{tikzpicture}[scale=.6,transform shape,rotate=90]
      \draw[web2,directed] (0,0) to (0,2);
    \end{tikzpicture}}\endxy
    &&
    \stackanchor{\text{\small circle}}{\text{\small evaluation}}
    \\[1ex]
    {}\xy(0,0)*{\begin{tikzpicture}[scale=.35,transform shape,rotate=90]
      \pic at (0,0) {web-};
      \pic at (1,1) {web+};
      \pic at (1,2) {web-};
      \pic at (0,3) {web+};
      \draw[web1,directed] (1.5,-.5) to (1.5,.5);
      \draw[web1,directed] (-.5,.5) to (-.5,2.5);
      \draw[web1,directed] (1.5,2.5) to (1.5,3.5);
    \end{tikzpicture}}\endxy
    % \mspace{5mu}=\mspace{5mu}
    =
    \xy(0,0)*{\begin{tikzpicture}[scale=.35,transform shape,rotate=90]
      \draw[web2,directed] (0,-.5) to (0,3.5);
      \draw[web1,directed] (1,-.5) to (1,3.5);
    \end{tikzpicture}}\endxy
    \quad\an\quad
    \xy(0,0)*{\begin{tikzpicture}[scale=.35,xscale=-1,transform shape,rotate=-90]
      \pic at (0,0) {web-};
      \pic at (1,1) {web+};
      \pic at (1,2) {web-};
      \pic at (0,3) {web+};
      \draw[web1,directed] (1.5,-.5) to (1.5,.5);
      \draw[web1,directed] (-.5,.5) to (-.5,2.5);
      \draw[web1,directed] (1.5,2.5) to (1.5,3.5);
    \end{tikzpicture}}\endxy
    =
    % \mspace{5mu}=\mspace{5mu}
    \xy(0,0)*{\begin{tikzpicture}[scale=.35,xscale=-1,transform shape,rotate=-90]
      \draw[web2,directed] (0,-.5) to (0,3.5);
      \draw[web1,directed] (1,-.5) to (1,3.5);
    \end{tikzpicture}}\endxy
    &&
    \stackanchor{\text{\small isotopies of flat}}{\text{\small tangle diagrams}}
  \end{IEEEeqnarray*}
  % This ends the definition.
\end{definition}

We shall need the following notion:

\begin{definition}
  \label{defn:spatial_isotopy}
  A \emph{spatial isotopy}\footnote{The terminology is taken from \cite{Selinger_SurveyGraphicalLanguages_2011}.} of flat tangle diagrams is a usual isotopy with the addition of the following \emph{spatial move}:
  \begin{equation*}
    \xy(0,0)*{\begin{tikzpicture}[scale=.25,rotate=90]
      \draw[web1] (0,0) circle (1cm);
      \draw[web1] (2,-2) to (2,2);
    \end{tikzpicture}}\endxy
    \quad\sim\quad
    \xy(0,0)*{\begin{tikzpicture}[scale=.25,rotate=90]
      \draw[web1] (0,0) circle (1cm);
      \draw[web1] (-2,-2) to (-2,2);
    \end{tikzpicture}}\endxy
  \end{equation*}
\end{definition}

The relations in $\web_d$ fully capture the spatial isotopy classes of the underlying flat tangle diagrams, in the sense of the following lemma.

\begin{lemma}
  \label{lem:web_spatial_isotopy}
  Let $W$ and $W'$ be two $\glt$-webs with the same domain and codomain. Then $W$ and $W'$ are equal in $\web_d$ if and only if there exists a spatial isotopy between their underlying flat tangle diagrams $\undcomp(W)$ and $\undcomp(W')$.
\end{lemma}

Using relations in \cref{defn:webs}, any closed strand in $\undcomp(W)$ evaluates to $q+q^{-1}$ in $W$. Moreover, if $\undcomp(W)$ does not have any closed strand, then the last two relations in \cref{defn:webs} are enough to capture all isotopies of $\undcomp(W)$.
These two facts essentially constitute the proof of \cref{lem:web_spatial_isotopy}; a formal proof can be found in the first author's PhD thesis \cite[subsubsection~6.6.3]{Schelstraete_OddKhovanovHomology_2024}.

%%%%%%%%%%%%%%%%%%%%%%%%%%%%%%%
%%%          foams          %%%
%%%%%%%%%%%%%%%%%%%%%%%%%%%%%%%
\newcommand{\interval}{[0,1]}

\subsection{Graded \texorpdfstring{$\mathfrak{gl}_2$}{gl2}-foams}
\label{subsec:graded_foams}
% %
\begin{wrapfigure}[13]{r}{.38\textwidth}
  % \vspace{-8pt}
  \centering
  \includegraphics[scale=1]{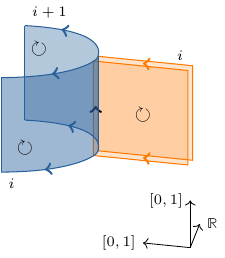}

  \caption{Local model for foams.}
  \label{fig:local_model_foam}
\end{wrapfigure}
Foams provide a suitable notion of cobordisms between webs.
They are certain singular surfaces locally modelled on the product of the interval $[0,1]$ with the letter ``Y'' (see \cref{fig:local_model_foam}), embedded in $(\bR\times \interval)\times \interval$. In \cref{fig:local_model_foam}, $\bR$ is pictured from front to back, while the last interval is pictured from bottom to top.
In this context, we refer to the singular curves as \emph{seams}, and call \emph{facets} the components of the complement of the set of seams. Facets have a thickness, either single or double, and we refer to them respectively as \emph{1\nbd-facets} (or \emph{single facets}, shaded blue) and \emph{2\nbd-facets} (or \emph{double facets}, shaded orange).
We refer to \emph{cross-sections} as cross-sections of the projection $\pi\colon(\bR\times \interval)\times \interval\to \interval$ onto the last coordinate (pictured vertically).

\begin{definition}
\label{defn:graded_foam}
  A \emph{graded $\glt$-foam}, or simply \emph{foam}, is a topological singular surface embedded in $(\bR\times \interval)\times \interval$, such that generic cross-sections are webs and all non-generic cross-sections respect one of the following local behaviours:
  \begingroup
  \def\scl{.3}
    \begin{IEEEeqnarray*}{CcCcCcCcC}
      \text{dot}
      &&
      \text{cup}
      &&
      \text{cap}
      &&
      \text{zip}
      &&
      \text{unzip}
      \\*[1ex]
      % DOT
      \figfoam[scale=1.3]{dot}
      &\mspace{30mu}&
      % CUP
      \figfoam[scale=1.3]{cup}
      &\mspace{30mu}&
      % CAP
      \figfoam[scale=1.3]{cap}
      &\mspace{30mu}&
      % ZIP
      \figfoam[scale=1.3]{zip}
      &\mspace{30mu}&
      % UNZIP
      \figfoam[scale=1.3]{unzip}
      \\*[1ex]
      \degd
      &&
      \degrcu
      &&
      \deglca
      &&
      \deglcu
      &&
      \degrca
    \end{IEEEeqnarray*}
  \endgroup
  where the axes are oriented as in \cref{fig:local_model_foam}.
  Each such local behaviour $F$ is endowed with a $\bZ^2$-degree $\deg_{\bZ^2}(F)$, denoted below the associated picture.
  Extending additively, this induces a $\bZ^2$-grading on graded $\glt$-foams.
\end{definition}

We call the above local behaviours \emph{generating foams}, or simply \emph{generators}.
A graded $\glt$-foam $F\subset (\bR\times \interval)\times\interval$ is viewed as a morphism from the web $F\cap\pi^{-1}(\{0\})$ to the web $F\cap\pi^{-1}(\{1\})$, reading from bottom to top (recall that $\pi$ is the projection onto the third coordinate):
\begin{IEEEeqnarray*}{CcCcCcCcC}
  \text{dot}
  &\mspace{30mu}&
  \text{cup}
  &\mspace{30mu}&
  \text{cap}
  &\mspace{30mu}&
  \text{zip}
  &\mspace{30mu}&
  \text{unzip}
  \\*[1ex]
  \id_{(1)}\to\id_{(1)}
  &&
  \id_{(2)}\to W_+W_-
  &&
  W_+W_-\to\id_{(2)}
  &&
  \id_{(1,1)}\to W_-W_+
  &&
  W_-W_+\to\id_{(1,1)}
\end{IEEEeqnarray*}
The \emph{dot} in the first picture of \cref{defn:graded_foam} is a formal decoration that any 1-facet may carry. By removing the 2-facets, a foam $F$ has an underlying surface, denoted $\undcomp(F)$. We assume that $\undcomp(F)$ is smooth and that the vertical projection $\pi$ defines a separative Morse function for $\undcomp(F)$, considering dots as critical points. Note that the local behaviours in \cref{defn:graded_foam} dictate the local behaviours around critical points of~$\pi$.
Because $\pi$ is separative, each such local behaviour lies on a distinct vertical position.

The facets of a foam $F$ admit a canonical orientation, induced by the canonical orientation on the web $F\cap\pi^{-1}(\{0\})$. That is, we endow the facets incident to $F\cap\pi^{-1}(\{0\})$ with an orientation compatible with the canonical orientation of $F\cap\pi^{-1}(\{0\})$, and extend globally with the condition that at a given seam, the orientation induced by the 2-facet is opposite to the orientation induced by the two 1-facets (see \cref{fig:local_model_foam}).
In particular, the facets incident to $F\cap\pi^{-1}(\{1\})$ induce an orientation on $F\cap\pi^{-1}(\{1\})$ opposite to its canonical orientation. 
Similarly, facets are endowed with a colour induced by the colour on webs; see also the string diagrammatics in \cref{subsec:string_diag_foam}.
This also defines an orientation and colour on each seam, induced from the orientation and colour of its incident 2-facet.

\medbreak

As in the non-graded case, we will consider graded $\glt$-foams up to isotopies. In our context, an \emph{isotopy} is a boundary-preserving isotopy which generically preserves the cross-section condition in \cref{defn:graded_foam}. Sliding a dot along its 1-facet is also considered as an isotopy. Moreover, we assume that the restriction of an isotopy to the underlying surface is a diffeotopy.

\begin{figure}
  \def\scl{.7}
  \newcommand{\ins}{3pt}% space between diagram and text
  \newcommand{\outs}{12pt}% space between two (set of) relation(s)
  \begin{gather*}
    %%%%% DOT-TYPE RELATIONS %%%%%
    \begin{IEEEeqnarraybox}{cCc}
      %%%%% DOT ANNIHILATION %%%%%
      \figfoam[scale=1.2]{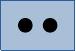}
      = 0
      &\mspace{80mu}&
      %%%%% DOT SLIDE %%%%%
      \figfoam[scale=1.2]{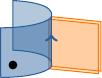}
      \;=\;\norma{}{-\;}{}\;
      \figfoam[scale=1.2]{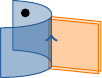}
      \\[\ins]
      \text{\footnotesize dot annihilation}
      &&
      \text{\footnotesize dot migration}
    \end{IEEEeqnarraybox}
    \\[\outs]
    %%%%% ZIGZAGS %%%%%
    \begin{IEEEeqnarraybox}{rClCrCl}
      \figfoam[scale=1.2]{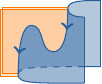}
      \;&=&\;
      \figfoam[scale=1.2]{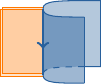}
      &\mspace{80mu}&
      \figfoam[scale=1.2]{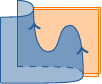}
      \;&=&\;X\;\;
      \figfoam[scale=1.2]{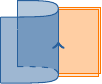}
      \\[2ex]
      \figfoam[scale=1.2]{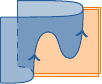}
      \;&=&\;\norma{Z\;\;}{Z\;\;}{Z^2\;\;}
      \figfoam[scale=1.2]{id_W_minus}
      &&
      \figfoam[scale=1.2]{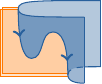}
      \;&=&\;\norma{YZ}{YZ}{YZ^2}\;\;
      \figfoam[scale=1.2]{id_W_plus}
    \end{IEEEeqnarraybox}
    \\[\ins]
    \text{\footnotesize zigzag relations}
    \\[\outs]
    %%%%% BUBBLES %%%%%
    \figfoam[scale=1.2]{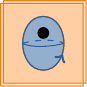}
    \;=\;
    \figfoam[scale=1.2]{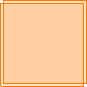}
    \mspace{80mu}
    \figfoam[scale=1.2]{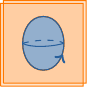}
    \;=\; 0
    \\[\ins]
    \text{\footnotesize bubble evaluations}
    \\[\outs]
    %%%%% HORIZONTAL NECK-CUTTING RELATION %%%%%
    \figfoam[scale=1.2]{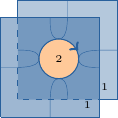}
    \;=\;
    \;\norma{}{}{Z}\;\;
    \figfoam[scale=1.2]{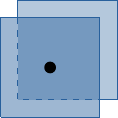}
    \;\norma{+XY}{-XY}{+XYZ}\;\;
    \figfoam[scale=1.2]{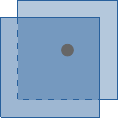}
    \\[\ins]
    \text{\footnotesize horizontal neck-cutting relation}
    \\[\outs]
    %%%%% NECK-CUTTING AND SQUEEZING %%%%%
    \begin{IEEEeqnarraybox}{cCc}
      %%%%% VERTICAL NECK-CUTTING RELATION %%%%%
      \figfoam[scale=1.2]{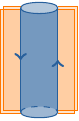}
      \;=\;
      \figfoam[scale=1.2]{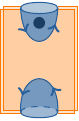}
      \;+\;
      \figfoam[scale=1.2]{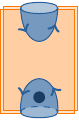}
      & \mspace{80mu} &
      %
      %%%%% SQUEEZING %%%%%
      \figfoam[scale=1.2]{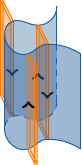}
      \mspace{10mu}=\norma{}{-}{Z^{-1}}\mspace{10mu}
      \figfoam[scale=1.2]{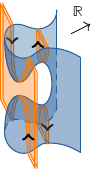}
      \\[\ins]
      \text{\footnotesize vertical neck-cutting relation}
      &&
      \text{\footnotesize squeezing relation}
    \end{IEEEeqnarraybox}
    \end{gather*}

  \caption{Relations in $\gfoam_d$. One of the two dots in the horizontal neck-cutting relation is grey, emphasizing that it sits on the back 1-facet. The $\bR$ axis is pictured from front to back, except for the squeezing relation for which it is pictured from left to right for better readability.}
  \label{fig:rel_foam}
\end{figure}

We distinguish different kinds of isotopies, in analogy with the property of the corresponding diffeotopy for the underlying surface:
\begin{itemize}
  \item If the isotopy preserves the relative vertical positions of the generating foams, we say that it is a \emph{Morse-preserving isotopy}.
  \item If the isotopy only interchanges the vertical positions of generating foams, we say that it is a \emph{Morse-singular isotopy}.
  \item If the isotopy is such that the underlying diffeotopy is a birth-death diffeotopy,\footnote{That is, collapsing the two singularities associated to a saddle and a cup or cap by ``smoothing out'' the surface, or the converse; see the zigzag relations in \cref{fig:rel_foam} for the underlying surfaces.} we say that the isotopy is a \emph{birth-death isotopy}.
\end{itemize}
We refer to \cite{Putyra_2categoryChronologicalCobordisms_2014} for relevant details on Morse theory.

In the graded case, isotopies only hold up to invertible scalars. Some are controlled by a graded-2-categorical structure:

\begin{definition}
  \label{defn:ring_R_bil}
  Let $\ringfoam$ be a commutative ring together with three invertible elements $X$, $Y$ and $Z$ in $\ringfoam^\times$ such that $X^2=Y^2=1$.
  Given this data, let $\bilfoam$ be the following bilinear form for the abelian group $G\coloneqq\bZ^2$:
  \begin{align*}
    \bilfoam\colon\bZ^2\times \bZ^2&\to \ringfoam^\times,\\
    ((a,b),(c,d)) &\mapsto X^{ac}Y^{bd}Z^{ad-bc}.
  \end{align*}
\end{definition}

Note that $\bilfoam$ is symmetric in the sense of \cref{defn:bil_symmetric}.
There are three standard choices for $\ringfoam$ and elements $X$, $Y$ and~$Z$:
\begin{itemize}
  \item \emph{even case:} leave $\ringfoam$ arbitrary and choose $X=Y=Z=1$. This recovers (even) $\glt$-foams with coefficients in $\ringfoam$ (see \cref{rem:chfoam_recover_even_foam}).
  \item \emph{odd case:} leave $\ringfoam$ arbitrary and choose $X=Z=1$ and $Y=-1$. This defines so-called \emph{super $\glt$-foams} (with coefficients in $\ringfoam$), as described in the introduction. One could also choose $Y=Z=1$ and $X=-1$, leading to an essentially identical theory.
  \item \emph{graded case:} $X$, $Y$ and $Z$ are formal parameters. More precisely, choose another commutative ring $\underline{\ringfoam}$ and set
  $\ringfoam\coloneqq\underline{\ringfoam}[X,Y,Z]/(X^2=Y^2=1).$
\end{itemize}
In what follows, we will often say that we ``choose $X=Y=Z=1$'' or ``choose $X=Z=1$ and $Y=-1$'' to mean that we consider the even or the odd case, respectively.

% Note that $\bilfoam$ is symmetric in the sense of \cref{defn:bil_symmetric}.

\begin{definition}
\label{defn:cat_foam}
   $\gfoam_d$ is the $(\bZ^2,\bilfoam)$-graded-2\nbd-cate\-gory whose objects are elements of $\uLam_{d}$, whose 1-mor\-phisms are $\glt$-webs and whose 2\nbd-morphisms are $\ringfoam$\nbd-linear combinations of graded $\glt$-foams, regarded up to the following relations:
  \begin{enumerate}[label=(\roman*)]
    \item If $\varphi\colon F_1\to F_2$ is a Morse-preserving isotopy, then $F_1=F_2$ in ${\gfoam}_d$.

    \item If $\varphi\colon F_1\to F_2$ is a Morse-singular isotopy interchanging the vertical positions of exactly two critical points $p$ and $q$, with $p$ vertically above $q$ in $F_1$, then
    \[F_1=\bilfoam (\deg p,\deg q) F_2\]
    in ${\gfoam}_d$.

    \item All the local relations in \cref{fig:rel_foam} below.
  \end{enumerate}
\end{definition}

% Setting $X=Y=Z=1$ recovers (even) $\glt$-foams (see \cref{rem:chfoam_recover_even_foam}), while setting $X=Z=1$ and $Y=-1$ gives so-called \emph{super $\glt$-foams} as discussed in the introduction. One could also set $Y=Z=1$ and $X=-1$, leading to an essentially identical theory.

%   \caption{Relations in $\gfoam_d$. A grey dot denotes a dot behind a facet. The $\bR$ axis is pictured from front to back, except for the squeezing relation for which it is pictured from left to right for better readability.}
%   \label{fig:rel_foam}
% \end{figure}

\begin{remark}
  \label{rem:foam_grading}
  The $\bZ^2$-grading induces a $\bZ$\nbd-grading on graded $\glt$-foams, the \emph{quantum grading} (or \emph{$q$-grading}). If $q\colon\bZ^2\to\bZ$ denotes the map $q(a,b)=a+b$, we set for each foam $F$:
  \begin{equation*}
    \qdeg(F)\coloneqq q(\deg_{\bZ^2}(F)).
  \end{equation*}
  Note that while symmetry with respect to the horizontal plane does not preserve the $\bZ^2$-grading, it does preserve the quantum grading.
  Topologically, the quantum grading is given by
  \[\qdeg(F) = 2\#\{\text{dots}\} - e(\undcomp(F)),\]
  where $\#\{\text{dots}\}$ is the number of dots and $e$ is the Euler measure, viewing $\undcomp(F)$ as a surface with acute right-angled corners. (Recall that the Euler measure of a surface $S$ with Euler characteristic $\chi(S)$ and with $k$ acute right-angled corners is $e(S)\coloneqq\chi(S)- k/4$; see e.g.\ \cite{Lipshitz_CylindricalReformulationHeegaard_2006}.) The Euler measure is additive under disjoint union and gluing of surfaces with corners, in accordance with the additivity of the quantum grading.

  Following \cref{subsec:more_cat_facts}, we write $\underline{(\gfoam_d)}^{\oplus}_q$ for the additive $q$\nbd-shifted closure of $\gfoam_d$: it allows formal shifts of webs in the quantum grading, restricts graded foams to those preserving the quantum grading, and allows formal direct sums.
  In this construction, we view the quantum grading as independent from the $\bZ^2$-grading; hence $\underline{(\gfoam)}^{\oplus}_q$ is still a $(\bZ^2,\bilfoam)$-graded-2\nbd-cate\-gory, and shifts in quantum degree do not affect the $\bZ^2$-degree.
\end{remark}

Similar to webs (\cref{lem:web_spatial_isotopy}), the following lemma shows that relations on foams capture the diffeotopy classes of the underlying surfaces; or rather, the underlying dotted surfaces, where sliding a dot along a connected component is considered to be a diffeotopy.
Below we write $F\projrel F'$ whenever there exists an invertible scalar $r\in \ringfoam^\times$ such that $F=r F'$.
Recall that $\undcomp(F)$ denotes the underlying surface of a foam $F$.

\begin{lemma}
  \label{lem:foam_all_isotopy_are_rels}
  Let $F$ and $F'$ be two foams in $\gfoam_d$ with the same domain and codomain. If $\undcomp(F)$ and $\undcomp(F')$ are isotopic, then $F\projrel F'$ in $\gfoam_d$.
\end{lemma}

\begin{proof}
  By Cerf theory (\cite{Cerf_DiffeomorphismesSphereDimension_1968}; see \cite[appendix~A]{Putyra_2categoryChronologicalCobordisms_2014} for a review), diffeotopic surfaces are related by diffeotopies preserving the relative vertical positions of critical points, Morse-singular diffeotopies interchanging two critical points, and birth-death diffeotopies. These correspond to the relations (i), (ii) and the zigzag relations in $\gfoam_d$, respectively. Finally, if $\undcomp(F)$ and $\undcomp(F')$ are diffeotopic by sliding a dot along a connected component, then $F$ and $F'$ are related by sliding the dot along 1-facets, possibly crossing seams using dot migration.
\end{proof}

When considering only isotopies of dots sliding along 1-facets, the above lemma can be made more precise:

\begin{lemma}
  \label{lem:moving_dots}
  Let $F\colon W\to W'$ and $F'\colon W'\to W''$ be two foams and let $D_1$ and $D_2$ be two foams identical to $\id_{W'}$ except for a dot sitting on a 1-facet. If the two dots belong to the same closed component of $\undcomp(F'\circ F)$, then
  \[ F'\circ D_1\circ F=F'\circ D_2\circ F\]
  in $\gfoam_d$.
\end{lemma}

\begin{proof}
  This is a consequence of the fact that a dot can slide across 2-facets at no cost of scalar, and along 1-facets past generators depending on $\bilfoam$. Because $\bilfoam$ is symmetric (see \cref{defn:bil_symmetric}), the latter scalar only depends on the relative vertical position of the dot.
\end{proof}

\begin{remark}
  \label{rem:chfoam_recover_even_foam}
  Recall the quantum grading from \cref{rem:foam_grading}, and $\underline{(\gfoam_d)}^{\oplus}_q$ the additive $q$\nbd-shifted closure of $\gfoam_d$.
  One can check that the $\ringfoam$-linear 2\nbd-cate\-gory
  \begin{equation*}
    \underline{(\gfoam_d)}^{\oplus}_q\vert_{X=Y=Z=1}
  \end{equation*}
  is precisely the category $n\mathbf{Foam}(N)^\bullet$ from \cite[p.~1322]{QR_$mathfraksl_n$Foam2category_2016} with $n=2$ and $N=d$, with the same quantum grading.
  This follows from \cref{lem:foam_all_isotopy_are_rels} and renormalizing the $i$-labelled dot, the $i$-labelled cap and the $i$-labelled zip by $(-1)^i$.
  In particular, we warn the reader familiar with the even setting that, while the dot migration relation in \cref{fig:rel_foam} has no sign, this is \emph{not} an essential feature of our construction, but rather a choice of normalization, allowed by the restriction to directed $\glt$-foams (see \cref{rem:directed_webs}).
  % orientation on seams same as \cite{BHP+_FunctorialityMathfrakSl_2019}, but orientation on facets is opposite
\end{remark}

%%%%%%%%%%%%%%%%%%%%%%%%%%%%%%%%%%%%%%
%%%          diagrammatics          %%%
%%%%%%%%%%%%%%%%%%%%%%%%%%%%%%%%%%%%%%
\subsection{String diagrammatics}
\label{subsec:string_diag_foam}

We introduce string diagrammatics for graded $\glt$-foams. This is based on the observation that \emph{a foam is fully described by its seams}; more precisely, by its domain object, its seams, and their orientations and colours.
For identity foams, this holds since webs are generated by $W_{i,\pm}$ and $\id_{W_{i,\pm}}$ is determined by its domain and the seam, together with its orientation and colour.
Thus, we can represent a web with an \emph{identity foam diagram} (or simply \emph{identity diagram}), a horizontal juxtaposition of oriented vertical strands coloured with elements in $\{1,2,\ldots,d-1\}$. For instance:
\begin{IEEEeqnarray*}{clCc}
  {}\xy(0,0)*{\begin{tikzpicture}[xscale=.8,scale=.8,yscale=.8]
    \draw[diag1,->] (0,0) node[below]{\tiny $1$} to (0,1.5);
    \draw[diag3,<-] (1,0) node[below]{\tiny $3$} to (1,1.5);
    \draw[diag1,<-] (2,0) node[below]{\tiny $1$} to (2,1.5);
    \draw[diag2,->] (3,0) node[below]{\tiny $2$} to (3,1.5);
  \end{tikzpicture}}\endxy
  &\mspace{-15mu}(1,2,1)&\mspace{100mu}&
  \xy(0,0)*{\begin{tikzpicture}[scale=.7*.8,yscale=.8*.8]
    \draw[web2] (.5,.5) to node[below,pos=.1] {\tiny $1$} node[below,pos=.9] {\tiny $1$} (2.5,.5);
    \draw[web2] (0,2.5) to node[below,pos=.9] {\tiny $3$} ++(1.5,0);
    \draw[web2] (3.5,1.5) to node[below,pos=.1] {\tiny $2$} ++(.5,0);
    \draw[web1] (0,0) to[out=0,in=-90] ++(.5,.5);
    \draw[web1] (0,1) to[out=0,in=90] ++(.5,-.5);
    \draw[web1] (1.5,2.5) to[out=90,in=180] ++(.5,.5) to ++(2,0);
    \draw[web1] (1.5,2.5) to[out=-90,in=180] ++(.5,-.5) to ++(1,0);
    \draw[web1] (2.5,.5) to[out=90,in=180] ++(.5,.5);
    \draw[web1] (2.5,.5) to[out=-90,in=180] ++(.5,-.5) to ++(1,0);
    \draw[web1] (3,1) to[out=0,in=-90] ++(.5,.5);
    \draw[web1] (3,2) to[out=0,in=90] ++(.5,-.5);
    \node at (5,0) {\small $1$};
    \node at (5,1.5) {\small $2$};
    \node at (5,3) {\small $1$};
  \end{tikzpicture}}\endxy
  \\*[-1pt]
  \text{\footnotesize a labelled identity diagram}
  &&&
  \text{\footnotesize the web it represents}
\end{IEEEeqnarray*}
where we remind the reader that webs are read and oriented from left to right. In the example above, we have labelled its rightmost region to specify the domain of the web.
Only the bigger numbers $(1,2,1)$ (``label'') convey a notion of thickness; the smaller numbers {\scriptsize $1$}, {\scriptsize $2$} and {\scriptsize $3$} (``colours'') capture vertical position.

Labelling one region of an identity foam diagram induces a label on all of its regions thanks to the following rule (here the value $l_\lambda$ (see \cref{subsec:webs}) is given below the corresponding coordinate):
\begin{equation}
  \label{eq:foam_diag_label}
  (\ldots,\underset{i}{1},\underset{i+1}{1},\ldots)
  % \mspace{10mu}
  \xy(0,-2)*{\begin{tikzpicture}[scale=.8]
    \draw[diag1,->] (0,0) node[below=-2pt]{\small $i$} to (0,1.5);
  \end{tikzpicture}}\endxy
  % \mspace{10mu}
  (\ldots,\underset{i}{2},\ldots)
  \mspace{10mu}\an\mspace{10mu}
  (\ldots,\underset{i}{2},\ldots)
  % \mspace{10mu}
  \xy(0,-2)*{\begin{tikzpicture}[scale=.8]
    \draw[diag1,<-] (0,0) node[below=-2pt]{\small $i$} to (0,1.5);
  \end{tikzpicture}}\endxy
  % \mspace{10mu}
  (\ldots,\underset{i}{1},\underset{i+1}{1},\ldots)
\end{equation}
Extending a label following rule \eqref{eq:foam_diag_label} may fail, as it can lead to coordinates outside of $\{0,1,2\}$. However, if the extension does work, we call the labelling \emph{legal}. An identity diagram equipped with a legal label is called a \emph{labelled identity (foam) diagram}. By the above discussion, labelled identity diagrams are in one-to-one correspondence with webs.

Given these diagrammatics for webs, it is not difficult to extend these diagrammatics to foams. Indeed, local behaviours in \cref{defn:graded_foam} are determined by the seam, together with its orientation and colour.

\begin{definition}
  \label{defn:generic_foam_diagram}
  A (generic) \emph{foam diagram}, or simply \emph{diagram}, is a diagram obtained by vertical and horizontal juxtapositions of identity diagrams and generators given below:
  \begin{gather*}
    \allowdisplaybreaks
      \begin{tabular}{*{5}{c@{\hskip 3ex}}c}
        % \small dot &\small \stackunder{\text{rightward}}{\text{cup}} & \small \stackunder{\text{leftward}}{\text{cap}} & \small \stackunder{\text{rightward}}{\text{cap}} & \small \stackunder{\text{leftward}}{\text{cup}}
        %
        \small dot &\small rightward cup & \small leftward cap & \small leftward cup & \small rightward cap
        \\*[1ex]
        %%%% dot %%%%
        \xy(0,0)*{\begin{tikzpicture}
          \node[fdot1] at (0,0) {};
          \node at (.2,-.2) {\scriptsize $i$};
        \end{tikzpicture}}\endxy
        &
        %%%% rightward cup %%%%
        \xy (0,2)*{ \begin{tikzpicture} [scale=1]
          \draw[diag1,->] (0,1) node[left=-3pt]{\scriptsize $i$} to [out=270,in=180] (.25,.5) to [out=0,in=270] (.5,1);
          % \node at (-.1,0.5) {\scriptsize $\lambda$};
        \end{tikzpicture} }\endxy
        &
        %%%% leftward cap %%%%
        \xy(0,0)*{\begin{tikzpicture}[scale=1]
          \draw[diag1,<-] (0,0) to [out=90,in=180] (.25,.5)
            to [out=0,in=90] (.5,0) node[right=-3pt]{\scriptsize $i$};
          % \node at (0.6,.5) {\scriptsize $\lambda$};
        \end{tikzpicture}}\endxy
        &
        %%%% leftward cup %%%%
        \xy(0,2)*{\begin{tikzpicture}[scale=1]
          \draw[diag1,<-] (0,1) to [out=270,in=180] (.25,.5)
            to [out=0,in=270] (.5,1) node[right=-3pt]{\scriptsize $i$};
        \end{tikzpicture}}\endxy
        &
        %%%% rightward cap %%%%
        \xy (0,0)*{\begin{tikzpicture} [scale=1]
          \draw[diag1,->] (0,0) node[left=-3pt]{\scriptsize $i$} to[out=90,in=180] (.25,.5) to[out=0,in=90] (.5,0);
          % \node at (0.6,.5) {\scriptsize $\lambda$};
        \end{tikzpicture} }\endxy
        \\*[1ex]
        $\degd$ & $\degrcu$ & $\deglca$ & $\deglcu$ & $\degrca$
      \end{tabular}
      \\[2ex]
      \begin{tabular}{*{3}{c@{\hskip 3ex}}c}
        \small downward crossing & \small leftward crossing & \small upward crossing & \small rightward crossing
        \\*[1ex]
        %%%% downward crossing %%%%
        \xy (0,0)*{\begin{tikzpicture}[scale=.4]
          \draw [diag1,->] (1,1) to (-1,-1) node[left=-3pt]{\scriptsize $i$};
          \draw [diag2,->] (-1,1) to (1,-1) node[right=-3pt]{\scriptsize $j$};
          % \node at (1.4,0) {\scriptsize $\lambda$};
        \end{tikzpicture}}\endxy
        &
        %%%% leftward crossing %%%%
        \xy (0,0)*{\begin{tikzpicture}[scale=.4]
          \draw [diag2,->] (1,1) to (-1,-1) node[left=-3pt]{\scriptsize $j$};
          \draw [diag1,<-] (-1,1) to (1,-1) node[right=-3pt]{\scriptsize $i$};
          % \node at (1.4,0) {\scriptsize $\lambda$};
        \end{tikzpicture}}\endxy
        &
        %%%% upward crossing %%%%
        \xy (0,0)*{\begin{tikzpicture}[scale=.4]
          \draw [diag1,<-] (1,1) to (-1,-1) node[left=-3pt]{\scriptsize $i$};
          \draw [diag2,<-] (-1,1) to (1,-1) node[right=-3pt]{\scriptsize $j$};
          % \node at (1.4,0) {\scriptsize $\lambda$};
        \end{tikzpicture}}\endxy
        &
        %%%% rightward crossing %%%%
        \xy (0,0)*{\begin{tikzpicture}[scale=.4]
          \draw [diag2,<-] (1,1) to (-1,-1) node[left=-3pt]{\scriptsize $j$};
          \draw [diag1,->] (-1,1) to (1,-1) node[right=-3pt]{\scriptsize $i$};
          % \node at (1.4,0) {\scriptsize $\lambda$};
        \end{tikzpicture}}\endxy
        \\*[1ex]
        $(0,0)$ & $(0,0)$ & $(0,0)$ & $(0,0)$
      \end{tabular}
  \end{gather*}
  where a crossing exists only if $\abs{i-j}>1$.
  Each generator is equipped with a $\bZ^2$-degree, which extends additively to generic diagrams.
  We assume that in a generic diagram, generators are in general position with respect to the vertical direction.
\end{definition}

In a foam diagram, we refer to strands coloured $i$ as \emph{$i$-strands}, and dots coloured $i$ as \emph{$i$-dots}.
As for an identity diagram, we can label the regions of a generic diagram with elements of $\uLam_d$. The label is said to be \emph{legal} if it satisfies condition \eqref{eq:foam_diag_label} as before, and if for each $i$-dot contained in a region labelled with $\lambda$, we have $\lambda_i=1$. This latter condition corresponds to the fact that dots can only sit on 1-facets. A \emph{labelled (foam) diagram} is a foam diagram equipped with a legal label. By the above discussion, labelled diagrams (regarded up to planar isotopies preserving the vertical positions of generators) are in one-to-one correspondence with foams (regarded up to Morse-preserving isotopies).

This provides a string diagrammatics for the graded-2\nbd-cate\-gory $\gfoam_d$:

\begin{figure}
  \def\scl{.7}
  \newcommand{\ins}{3pt}% space between diagram and text
  \newcommand{\outs}{12pt}% space between two (set of) relation(s)
  \begin{gather*}
    %%%%% BRAID-LIKE MOVES %%%%%
    \begingroup
      %%%%% Reidemeister 2 moves %%%%%
      \tikz[scale=.8*\scl,baseline={([yshift=.8ex]current bounding box.center)}]{
        \draw[diag1]  +(0,-.75) node[below] {\textcolor{black}{\scriptsize $i$}}
          .. controls (0,-.375) and (1,-.375) .. (1,0)
          .. controls (1,.375) and (0, .375) .. (0,.75);
        \draw[diag2]  +(1,-.75) node[below] {\textcolor{black}{\scriptsize $j$}}
          .. controls (1,-.375) and (0,-.375) .. (0,0)
          .. controls (0,.375) and (1, .375) .. (1,.75);
        % \node at (1.3,0) {\scriptsize $\lambda$};
      }
      \;=\;
      \tikz[scale=.8*\scl,baseline={([yshift=.8ex]current bounding box.center)}]{
        \draw[diag1] (0,-.75) node[below] {\textcolor{black}{\scriptsize $i$}} to (0,.75);
        \draw[diag2] (1,-.75) node[below] {\textcolor{black}{\scriptsize $j$}} to (1,.75);
        % \node at (1.3,0) {\scriptsize $\lambda$};
      }
      \mspace{80mu}
      %
      %%%%% Reidemeister 3 moves %%%%%
      \tikz[scale=.7*\scl,baseline={([yshift=.9ex]current bounding box.center)}]{
        \draw[diag1]  +(0,0)node[below] {\textcolor{black}{\scriptsize $i$}}
          .. controls (0,0.5) and (2, 1) ..  +(2,2);
        \draw[diag3]  +(2,0)node[below] {\textcolor{black}{\scriptsize $k$}}
          .. controls (2,1) and (0, 1.5) ..  +(0,2);
        \draw[diag2]  (1,0)node[below] {\textcolor{black}{\scriptsize $j$}}
          .. controls (1,0.5) and (0, 0.5) ..  (0,1)
          .. controls (0,1.5) and (1, 1.5) ..  (1,2);
        % \node at (2.1,1) {\scriptsize $\lambda$};
      }
      \;=\;
      \tikz[scale=.7*\scl,baseline={([yshift=.9ex]current bounding box.center)}]{
        \draw[diag1]  +(0,0)node[below] {\textcolor{black}{\scriptsize $i$}}
          .. controls (0,1) and (2, 1.5) ..  +(2,2);
        \draw[diag3]  +(2,0)node[below] {\textcolor{black}{\scriptsize $k$}}
          .. controls (2,.5) and (0, 1) ..  +(0,2);
        \draw[diag2]  (1,0)node[below]{\textcolor{black} {\scriptsize $j$}}
          .. controls (1,0.5) and (2, 0.5) ..  (2,1)
          .. controls (2,1.5) and (1, 1.5) ..  (1,2);
        % \node at (2.3,1) {\scriptsize $\lambda$};
      }
    \endgroup
    \\[\ins]
    \text{\footnotesize braid-like relations}
    \\[\outs]
    %%%%% PITCHFORKS %%%%%
    \begingroup
      \begin{tikzpicture}[baseline = 0, scale=1.2*\scl]
      	\draw[diag1](-.5,.4) to (0,-.3);
      	\draw[diag2] (0.3,-0.3)
      		to[out=90, in=0] (0,0.2)
      		to[out = -180, in = 40] (-0.5,-0.3);
      	\node at (-0.5,-.5) {\scriptsize $j$};
      	\node at (0,-.5) {\scriptsize $i$};
      	% \node at (0.4,0) {\scriptsize $\lambda$};
      \end{tikzpicture}
      \;=\;
      \begin{tikzpicture}[baseline = 0,scale=1.2*\scl]
      	\draw[diag1](.6,.4) to (.1,-.3);
      	\draw[diag2] (0.6,-0.3)
      		to[out=140, in=0] (0.1,0.2)
      		to[out = -180, in = 90] (-0.2,-0.3);
      	\node at (-0.2,-.5) {\scriptsize $j$};
      	\node at (0.1,-.5) {\scriptsize $i$};
      	% \node at (0.7,0) {\scriptsize $\lambda$};
      \end{tikzpicture}
      \mspace{80mu}
      \begin{tikzpicture}[baseline = 5, scale=1.2*\scl]
      	\draw[diag1](-.5,-.3) to (0,.4);
      	\draw[diag2] (0.3,0.4)
      		to[out=-90, in=0] (0,-0.1)
      		to[out = 180, in = -40] (-0.5,0.4);
      	\node at (-0.5,.6) {\scriptsize $j$};
      	\node at (-0.05,.6) {\scriptsize $i$};
      	% \node at (0.5,0.1) {\scriptsize $\lambda$};
      \end{tikzpicture}
      \;=\;
      \begin{tikzpicture}[baseline = 5, scale=1.2*\scl]
      	\draw[diag1](.6,-.3) to (.1,.4);
      	\draw[diag2] (0.6,0.4)
      		to[out=-140, in=0] (0.1,-0.1)
      		to[out = 180, in = -90] (-0.2,0.4);
      	\node at (-0.25,.6)  {\scriptsize $j$};
      	\node at (0.15,.6) {\scriptsize $i$};
      	% \node at (0.7,0.1) {\scriptsize $\lambda$};
      \end{tikzpicture}
    \endgroup
    \\[\ins]
    \text{\footnotesize pitchfork relations}
    \\[\outs]
    %%%%% ZIGZAGS %%%%%
    \begingroup
      \xy (0,1)*{\tikz[scale=1*\scl]{
        \draw[diag1,<-] (0,-.5) to (0,0) to[out=90,in=180] (.25,.5)
          to[out=0,in=90] (.5,0) to[out=-90,in=180] (.75,-.5)
          to[out=0,in=-90] (1,0) to (1,.5) node[left=-3pt]{\scriptsize $i$};
        % \node at (1.2,0) {\tiny $\lambda$};
      }}\endxy
      =
      \xy (0,1)*{\begin{tikzpicture}[scale=1*\scl]
      	\draw [diag1,<-] (0,0) to (0,1) node[left=-3pt]{\scriptsize $i$};
      	% \node at (.3,.5) {\tiny $\lambda$};
      \end{tikzpicture}}\endxy
      \mspace{30mu}
      \xy (0,-1)*{\tikz[scale=1*\scl]{
        \draw[diag1,<-] (0,.5) to (0,0) to[out=-90,in=180] (.25,-.5)
          to[out=0,in=-90] (.5,0) to[out=90,in=180] (.75,.5)
          to[out=0,in=90] (1,0) to (1,-.5) node[left=-3pt]{\scriptsize $i$};
        % \node at (1.2,0) {\tiny $\lambda$};
      }}\endxy
      =X
      \xy (0,-1)*{\begin{tikzpicture}[scale=1*\scl]
        \draw [diag1,->] (0,0) node[left=-3pt]{\scriptsize $i$} to (0,1);
        % \node at (.3,.5) {\tiny $\lambda$};
      \end{tikzpicture}}\endxy
      \mspace{30mu}
      \xy (0,-1)*{\tikz[scale=1*\scl]{
        \draw[diag1,->] (0,-.5) node[left=-3pt]{\scriptsize $i$} to (0,0) to[out=90,in=180] (.25,.5)
          to[out=0,in=90] (.5,0) to[out=-90,in=180] (.75,-.5)
          to[out=0,in=-90] (1,0) to (1,.5);
        % \node at (1.2,0) {\tiny $\lambda$};
      }}\endxy
      =
      \norma{Z}{Z}{Z^2}\;
      \xy (0,-1)*{\begin{tikzpicture}[scale=1*\scl]
      	\draw [diag1,->] (0,0) node[left=-3pt]{\scriptsize $i$} to (0,1);
      	% \node at (.3,.5) {\tiny $\lambda$};
      \end{tikzpicture}}\endxy
      \mspace{30mu}
      \xy (0,1)*{\tikz[scale=1*\scl]{
        \draw[diag1,->] (0,.5) node[left=-3pt]{\scriptsize $i$} to (0,0) to[out=-90,in=180] (.25,-.5)
          to[out=0,in=-90] (.5,0) to[out=90,in=180] (.75,.5)
          to[out=0,in=90] (1,0) to (1,-.5);
        % \node at (1.2,0) {\tiny $\lambda$};
      }}\endxy
      =
      \norma{YZ}{YZ}{YZ^2}\;
      \xy (0,1)*{\begin{tikzpicture}[scale=1*\scl]
        \draw [diag1,<-] (0,0) to (0,1) node[left=-3pt]{\scriptsize $i$};
        % \node at (.3,.5) {\tiny $\lambda$};
      \end{tikzpicture}}\endxy
    \endgroup
    \\[\ins]
    \text{\footnotesize zigzag relations (or adjunction relations)}
    \\[\outs]
    %%%%% DOT-TYPE RELATIONS %%%%%
    \begin{IEEEeqnarraybox}{cCcCc}
      %%%%% DOT ANNIHILATION %%%%%
      \left(\tikz[scale=.8,baseline={([yshift=0ex]current bounding box.center)}]{
      \node[fdot1] at (0,0) {};
      \node at (.3,-.2) {\scriptsize $i$};
      }\right)^2
      =
      0
      &\mspace{50mu}&
      %%%%% DOT SWAP %%%%%
      \tikz[scale=.8*\scl,baseline={([yshift=.8ex]current bounding box.center)}]{
      \node[fdot1] at (-1,0) {};
      \node at (-.7,-.2) {\scriptsize $i$};
      \draw[diag1,->] (0,-.75) node[below] {\textcolor{black}{\scriptsize $i$}} to (0,.75);
      }
      =\norma{}{-\;\;}{}
      \tikz[scale=.8*\scl,baseline={([yshift=.8ex]current bounding box.center)}]{
      \node[fdot2] at (-1.2,0) {};
      \node at (-.7,-.5) {\scriptsize $i+1$};
      \draw[diag1,->] (0,-.75) node[below] {\textcolor{black}{\scriptsize $i$}} to (0,.75);
      }
      &\mspace{50mu}&
      %%%%% DOT SLIDE %%%%%
      \tikz[scale=.8*\scl,baseline={([yshift=.8ex]current bounding box.center)}]{
      \node[fdot2] at (-.7,0) {};
      \node at (-.4,-.4) {\scriptsize $j$};
      \draw[diag1] (0,-.75) node[below] {\textcolor{black}{\scriptsize $i$}} to (0,.75);
      }
      \;=\;
      \tikz[scale=.8*\scl,baseline={([yshift=.8ex]current bounding box.center)}]{
      \node[fdot2] at (.7,0) {};
      \node at (1,-.4) {\scriptsize $j$};
      \draw[diag1] (0,-.75) node[below] {\textcolor{black}{\scriptsize $i$}} to (0,.75);
      }
      \mspace{20mu}
      \text{if }j\neq i,i+1
      \\[\ins]
      \text{\footnotesize dot annihilation}
      &&
      \text{\footnotesize dot migration}
      &&
      \text{\footnotesize dot slide}
    \end{IEEEeqnarraybox}
    \\[\outs]
    %%%%% BUBBLES %%%%%
    \begin{IEEEeqnarraybox}{cCc}
      %%%%% COUNTER-CLOCKWISE BUBBLES %%%%%
      {}\xy(0,0)*{\begin{tikzpicture}[scale=\scl]
        \draw[diag1,
        decoration={markings, mark=at position 0 with {\arrow{>}}},
        postaction={decorate}]
        (0,0) circle (.5cm);
        \node at (.8,-.3) {\scriptsize $i$};
        \node[fdot1] at (0,0) {};
        \node at (.2,-.2) {\scriptsize $i$};
      \end{tikzpicture}}\endxy
      =
      1
      \mspace{40mu}
      \xy(0,0)*{\begin{tikzpicture}[scale=\scl]
        \draw[diag1,
        decoration={markings, mark=at position 0 with {\arrow{>}}},
        postaction={decorate}]
        (0,0) circle (.5cm);
        \node at (.8,-.3) {\scriptsize $i$};
      \end{tikzpicture}}\endxy
      =
      0
      & \mspace{70mu} &
      %%%%% CLOCKWISE BUBBLES %%%%%
      \begingroup
        \xy(0,0)*{\begin{tikzpicture}[scale=\scl]
          \draw[diag1,
          decoration={markings, mark=at position 0 with {\arrow{<}}},
          postaction={decorate}]
          (0,0) circle (.5cm);
          \node at (.7,-.3) {\scriptsize $i$};
        \end{tikzpicture}}\endxy
        =
        \norma{}{}{Z}
        \mspace{10mu}
        \xy(0,-1.2)*{\begin{tikzpicture}[scale=\scl]
          \node[fdot1] at (0,0) {};
          \node at (.2,-.2) {\scriptsize $i$};
        \end{tikzpicture}}\endxy
        % \mspace{10mu}
        \norma{+XY}{-XY}{+XYZ}
        \mspace{10mu}
        \xy(0,-1.2)*{\begin{tikzpicture}[scale=\scl]
          \node[fdot2] at (0,0) {};
          \node at (.5,-.3) {\scriptsize $i+1$};
        \end{tikzpicture}}\endxy
      \endgroup
      \\[\ins]
      \text{\footnotesize evaluation of counter-clockwise bubbles}
      &&
      \text{\footnotesize evaluation of clockwise bubbles}
    \end{IEEEeqnarraybox}
    \\[\outs]
    %%%%% SHADING RELATED RELATION %%%%%
    \begin{IEEEeqnarraybox}{cCc}
      %%%%% NECK-CUTTING RELATION %%%%%
      \begingroup
        \xy(0,-2)*{\begin{tikzpicture}[scale=.8*\scl]
          \draw[diag1,<-] (0,0) node[below] {\scriptsize $i$} to (0,1.8);
          \draw[diag1,->] (1,0) to (1,1.8);
        \end{tikzpicture}}\endxy
        \;=\;
        \xy(0,-2)*{\begin{tikzpicture}[scale=.8*\scl]
          \draw[diag1,<-] (0,0) node[below] {\scriptsize $i$} to[out=90,in=180] (.5,.75) to[out=0,in=90] (1,0);
          \begin{scope}[shift={(0,.8)}]
            \draw[diag1,->] (0,1) to[out=-90,in=180] (.5,.25) to[out=0,in=-90] (1,1);
          \end{scope}
          \begin{scope}[shift={(.5,1.7)}]
            \node[fdot1] at (0,0) {};
            % \node at (.2,-.2) {\scriptsize $i$};
          \end{scope}
        \end{tikzpicture}}\endxy
        \;+\;
        \xy(0,-2)*{\begin{tikzpicture}[scale=.8*\scl]
          \draw[diag1,<-] (0,0) node[below] {\scriptsize $i$} to[out=90,in=180] (.5,.75) to[out=0,in=90] (1,0);
          \begin{scope}[shift={(0,.8)}]
            \draw[diag1,->] (0,1) to[out=-90,in=180] (.5,.25) to[out=0,in=-90] (1,1);
          \end{scope}
          \begin{scope}[shift={(.5,.1)}]
            \node[fdot1] at (0,0) {};
            % \node at (.2,-.2) {\scriptsize $i$};
          \end{scope}
        \end{tikzpicture}}\endxy
      \endgroup
      & \mspace{120mu} &
      %%%%% ISOTOPY RELATION %%%%%
      \begingroup
        \xy(0,-2)*{\begin{tikzpicture}[scale=.7*\scl]
          \draw[diag2,->] (0,0) node[below] {\scriptsize $i+1$} to (0,3);
          \draw[diag2,<-] (1,0) to (1,3);
          \draw[diag1,<-] (-.5,0) to (-.5,3);
          \draw[diag1,->] (1.5,0) node[below] {\scriptsize $i$} to (1.5,3);
        \end{tikzpicture}}\endxy
        \;=\;
        \norma{}{-\;}{Z^{-1}\;}
        \xy(0,-2)*{\begin{tikzpicture}[scale=.7*\scl]
          \draw[diag2,->] (0,0) node[below] {\scriptsize $i+1$} to[out=90,in=180] (.5,.75) to[out=0,in=90] (1,0);
          \draw[diag1,<-] (-.5,0) to[out=90,in=180] (.5,1.25) to[out=0,in=90] (1.5,0) node[below] {\scriptsize $i$};
          \begin{scope}[shift={(0,2)}]
            \draw[diag2,<-] (0,1) to[out=-90,in=180] (.5,.25) to[out=0,in=-90] (1,1);
            \draw[diag1,->] (-.5,1) to[out=-90,in=180] (.5,-.25) to[out=0,in=-90] (1.5,1);
          \end{scope}
        \end{tikzpicture}}\endxy
      \endgroup
      \\[\ins]
      \text{\footnotesize neck-cutting relation}
      &&
      \text{\footnotesize squeezing relation}
    \end{IEEEeqnarraybox}
  \end{gather*}

  \caption{Relations in $\diag_d^\Lambda$. We omitted the objects labelling the regions of each diagram: this avoids clutter and emphasizes that relations are independent on the ambient object. If no orientation is given, the relation holds for all orientations. In the case of the braid-like and pitchfork relations, colours should be so that the crossings exist.}
  \label{fig:rel_diagfoam}
\end{figure}

\begin{definition}
  \label{defn:string-2-foam}
    ${\diag}_d^\Lambda$ is the $(\bZ^2,\bilfoam)$-graded-2-category whose objects are elements of the set $\uLam_d$, 1\nbd-morphisms are labelled identity foam diagrams, and 2\nbd-morphisms are labelled foam diagrams, regarded up to the following relations:
  \begin{enumerate}[(i)]
    \item If $\varphi\colon D_1\to D_2$ is a planar isotopy such that $\varphi$ preserves the relative vertical positions of the generators, then $D_1=D_2$ in $\diag_d^\Lambda$.
    \item If $\varphi\colon D_1\to D_2$ is a planar isotopy as above \emph{except} that it interchanges the vertical positions of two generators $p$ and $q$, with $p$ vertically above $q$, then $D_1=\bilfoam(\deg p,\deg q)D_2$ in $\diag_d^\Lambda$.
    \item All the local relations in \cref{fig:rel_diagfoam} above, viewed with a legal label.
  \end{enumerate}
\end{definition}

\begin{proposition}
  $\gfoam_d$ and $\diag_d^\Lambda$ are isomorphic as $(\bZ^2,\bilfoam)$-graded-2-categories.
\end{proposition}

\begin{proof}
  The $\bZ^2$\nbd-grading is preserved by correspondence between foams and diagrams. It remains to check that foam relations in \cref{defn:cat_foam} correspond to diagrammatic relations in \cref{defn:string-2-foam}.
  
  Relations (i) in \cref{defn:cat_foam} correspond to relations (i) in \cref{defn:string-2-foam}, together with graded interchange laws that involve at least one crossing, and braid-like relations, pitchfork relations and dot slide.
  
  Relations (ii) in \cref{defn:cat_foam} correspond to graded interchange laws that only involve cups, caps and dots.

  Relations in \cref{fig:rel_foam} correspond to relations in \cref{fig:rel_diagfoam}, except braid-like relations, pitchfork relations and dot slide. For instance, the horizontal neck-cutting corresponds to the evaluation of clockwise bubbles.
\end{proof}

Note that none of the relations in \cref{defn:string-2-foam} depends on the objects, that is, the label of the foam diagrams. In the sequel, we mostly leave foam diagrams unlabelled, as the label is either irrelevant to the discussion or understood from the context.
In particular, this applies when computing secondary relations from the defining relations in $\diag_d^\Lambda$, as in the following lemma:

\begin{lemma}
\label{lem:more_rels_diagfoam}
  The following local relations hold in $\diag_d^\Lambda$ for any choice of legal label (omitted below):
  \begin{gather*}
    %%% VNC type B %%%
    {}\xy(0,-2)*{\begin{tikzpicture}[scale=.5]
      \draw[diag1,->] (0,0) node[below] {\scriptsize $i$} to[out=90,in=180] (.5,.75) to[out=0,in=90] (1,0);
      \begin{scope}[shift={(0,.8)}]
        \draw[diag1,<-] (0,1) to[out=-90,in=180] (.5,.25) to[out=0,in=-90] (1,1);
      \end{scope}
    \end{tikzpicture}}\endxy
    =
    \norma{XZ}{XZ}{XZ^2}\;
    \xy(0,-2)*{\begin{tikzpicture}[scale=.5]
      \draw[diag1,->] (0,0) node[below] {\scriptsize $i$} to (0,1.8);
      \draw[diag1,<-] (1,0) node[below] {\scriptsize $i$} to (1,1.8);
      \begin{scope}[shift={(1.7,.9)}]
        \node[fdot1] at (0,0) {};
        \node at (.3,-.2) {\scriptsize $i$};
      \end{scope}
    \end{tikzpicture}}\endxy
    \;+\norma{YZ}{YZ}{YZ^2}\;
    \xy(0,-2)*{\begin{tikzpicture}[scale=.5]
      \draw[diag1,->] (0,0) node[below] {\scriptsize $i$} to (0,1.8);
      \draw[diag1,<-] (1,0) node[below] {\scriptsize $i$} to (1,1.8);
      \begin{scope}[shift={(-.7,.9)}]
        \node[fdot1] at (0,0) {};
        \node at (.3,-.2) {\scriptsize $i$};
      \end{scope}
    \end{tikzpicture}}\endxy
    \mspace{40mu}
    %%% inverse of iso2 %%%
    {}\xy(0,-2)*{\begin{tikzpicture}[scale=.4]
      \draw[diag1=1.5,->] (0,0) node[below] {\scriptsize $i$} to[out=90,in=180] (.5,.75) to[out=0,in=90] (1,0);
      \draw[diag2=1.5,<-] (-.5,0) to[out=90,in=180] (.5,1.25) to[out=0,in=90] (1.5,0) node[below] {\scriptsize $i+1$};
      \begin{scope}[shift={(0,2)}]
        \draw[diag1=1.5,<-] (0,1) to[out=-90,in=180] (.5,.25) to[out=0,in=-90] (1,1);
        \draw[diag2=1.5,->] (-.5,1) to[out=-90,in=180] (.5,-.25) to[out=0,in=-90] (1.5,1);
      \end{scope}
    \end{tikzpicture}}\endxy
    =
    \norma{XY}{-XY}{XYZ}
    \mspace{5mu}
    \xy(0,-2)*{\begin{tikzpicture}[scale=.4]
      \draw[diag1=1.5,->] (0,0) node[below] {\scriptsize $i$} to (0,3);
      \draw[diag1=1.5,<-] (1,0) to (1,3);
      \draw[diag2=1.5,<-] (-.5,0) to (-.5,3);
      \draw[diag2=1.5,->] (1.5,0) node[below] {\scriptsize $i+1$} to (1.5,3);
    \end{tikzpicture}}\endxy
    \\[1ex]
    %%% ISO TYPE B %%%
    {}\xy(0,1)*{\begin{tikzpicture}[scale=.45]
      \draw[diag2=1.5,
        decoration={markings, mark=at position 0 with {\arrow{<}}},
        postaction={decorate}]
        (0,0) circle (.5cm);
      \node at (0,-.8) {\tiny $i+1$};
      \draw[diag1=1.5,
        decoration={markings, mark=at position 0 with {\arrow{>}}},
        postaction={decorate}]
        (0,0) circle (1.2cm);
      \node at (1.7,-.3) {\scriptsize $i$};
    \end{tikzpicture}}\endxy
    \;=\;
    \norma{1}{-1}{Z}
    \mspace{50mu}
    {}\xy(0,1)*{\begin{tikzpicture}[scale=.45]
      \draw[diag1=1.5,
        decoration={markings, mark=at position 0 with {\arrow{<}}},
        postaction={decorate}]
        (0,0) circle (.5cm);
      \node at (0,-.8) {\scriptsize $i$};
      \draw[diag2=1.5,
        decoration={markings, mark=at position 0 with {\arrow{>}}},
        postaction={decorate}]
        (0,0) circle (1.2cm);
      \node at (2.3,-.3) {\scriptsize $i+1$};
    \end{tikzpicture}}\endxy
    \;=\;
    \norma{XY}{-XY}{XYZ}
  \end{gather*}
\end{lemma}

%%%%%%%%%%%%%%%%%%%%%%%%%%%%%%%
%%%          basis          %%%
%%%%%%%%%%%%%%%%%%%%%%%%%%%%%%%
\subsection{A basis for graded \texorpdfstring{$\glt$}{gl2}-foams}
\label{subsec:basis_foam}

Each 2Hom-space of $\gfoam_d$ admits a basis with the following one-line description:
\begin{quote}
  \textsc{slogan.} \emph{A hom-basis in $\gfoam_d$ is given by picking a foam whose underlying surface is a union of disks, and considering all the ways to put dots on the disks.}
\end{quote}
This kind of family of foams is formalized below as a \emph{reduced family}, and foams appearing in such families as \emph{reduced foams}.
In this section, we restrict ourselves to the description of reduced families.
The following theorem is shown in \cite{Schelstraete_RewritingModuloDiagrammatic_2026} (see also \cite{Schelstraete_OddKhovanovHomology_2024}) using higher linear rewriting theory:

\begin{theorem}[Basis theorem for graded $\glt$-foams]
  \label{thm:basis_foam}
  Let $W,W'\colon\mu\to\lambda$ be two webs with the same source and target.
  $W$~and $W'$ admit a reduced family (\cref{defn:foam_reduced_family}) and every reduced family constitutes a basis of $\Hom_{\gfoam_d}(W,W')$.
\end{theorem}

Let $W,W'\colon\mu\to\lambda$ be two webs with the same source and target.
We denote $\undcomp(W)\sqcup_\partial \undcomp(W')$ the closed 1-manifold obtained by glueing $\undcomp(W)$ and $\undcomp(W')$ along their common boundary points.
Note that if $F\colon W\to W'$ is a foam, then $\partial (\undcomp(F))$ is homeomorphic to $\undcomp(W)\sqcup_\partial \undcomp(W')$.
A \emph{reduced foam} is a foam $F$ such that $\undcomp(F)$ is a (dotted) union of disks with at most one dot on each disk.
In that case, there is a bijection between $\pi_0(\undcomp(F))$, the connected components of $\undcomp(F)$, and $\pi_0(\undcomp(W)\sqcup_\partial \undcomp(W'))$, the connected components of $\undcomp(W)\sqcup_\partial \undcomp(W')$, mapping a disk to its boundary.
For $\delta\subset \pi_0(\undcomp(W)\sqcup_\partial \undcomp(W'))$, we say that $\undcomp(F)$ (or abusing terminology, $F$) is \emph{$\delta$-dotted} if the following holds:
\begin{quote}
  \emph{$\delta$-dotted:} there is a dot (resp.\ no dot) sitting on a disk $d$ in $\undcomp(F)$ if and only if its boundary $\partial d\in \pi_0(\undcomp(W)\sqcup_\partial \undcomp(W'))$ is in $\delta$ (resp.\ is not in $\delta$).  
\end{quote}
We say that $F$ is \emph{undotted} whenever it is $\emptyset$-dotted.

\begin{definition}
  \label{defn:foam_reduced_family}
  Let $W,W'\colon\mu\to\lambda$ be two webs with the same source and target.
  A \emph{reduced family} is a family of foams $F_\delta\colon W\to W'$, indexed by subsets $\delta\subset\pi_0(\undcomp(W)\sqcup_\partial \undcomp(W'))$, such that each $F_\delta$ is a $\delta$-dotted reduced foam.
\end{definition}

Thanks to \cref{lem:foam_all_isotopy_are_rels}, reduced families are essentially unique, in the sense that if $\{F_{\delta}\}_\delta$ and $\{F_{\delta}'\}_\delta$ are two reduced families for the same hom-space, then $F_\delta\projrel F_\delta'$ for all $\delta$.
In particular:

\begin{lemma}
  \label{lem:reduced_family_basis_implies_every_reduced_family_basis}
  Let $W,W'\colon\mu\to\lambda$ be two webs with the same source and target. If a reduced family is a basis of $\Hom_{\gfoam_d}(W,W')$, then every reduced family is a basis of $\Hom_{\gfoam_d}(W,W')$.
  \hfill\qed
\end{lemma}

As shown in \cite[Lemma~2.5]{BHP+_Functoriality$mathfraksl_2$Tangle_2023}, counting the number of closed components in $\undcomp(W)\sqcup_\partial \undcomp(W')$ defines a non-degenerate pairing on $\glt$-webs:
\begin{equation*}
  \langle W,W'\rangle\coloneqq (q+q^{-1})^{\abs{\pi_0(\undcomp(W)\sqcup_\partial \undcomp(W'))}},
\end{equation*}
where $\abs{\pi_0(\undcomp(W)\sqcup_\partial \undcomp(W'))}$ denotes the number of connected components in $\undcomp(W)\sqcup_\partial \undcomp(W')$.
This pairing coincides with the web evaluation formula given in \cite{Robert_NewWayEvaluate_2015}; see also \cite[p.~1315]{BHP+_Functoriality$mathfraksl_2$Tangle_2023}.

\begin{corollary}[Non-degeneracy of graded $\glt$-foams]
\label{cor:non-degeneracy_foam}
  Let $W,W'\colon\mu\to\lambda$ be two webs with the same source and target.
  Denote by $\#1(\lambda,\mu)$ the number of $1$'s in $\lambda$ plus the number of $1$'s in $\mu$.
  Then the space $\Hom_{\gfoam_d}(W,W')$ is a free $\ringfoam$\nbd-module and:
  \begin{equation*}
    \gdim_q\Hom_{\gfoam_d}(W,W')=q^{\#1(\lambda,\mu)/2}\langle W,W'\rangle,
  \end{equation*}
  where $\gdim_q(-)$ denotes the graded rank with respect to the $q$-grading.
  \hfill\qed
\end{corollary}

\begin{proof}
  Let $F\colon W\to W'$ be a $\delta$-dotted reduced foam.
  Following \cref{rem:foam_grading}, we have that
  $\qdeg F = 2\abs{\delta}-\abs{\pi_0(\undcomp(W)\sqcup_\partial \undcomp(W'))} + \#1(\lambda,\mu)/2$.
  The result follows from \cref{thm:basis_foam}.
\end{proof}
%%%%%%%%%%%%%%%%%%%%%%%%%%%%%%%%%%%%%%%%%%%%
%%%          decategorification          %%%
%%%%%%%%%%%%%%%%%%%%%%%%%%%%%%%%%%%%%%%%%%%%
\subsection{The categorification theorem}
\label{sec:categorification_foam}

Recall the definition of the Grothendieck ring of a graded-2\nbd-cate\-gory and related notions from \cref{subsec:more_cat_facts}.
We decategorify $\gfoam_d$ with respect to the quantum grading (see \cref{rem:foam_grading}), that is:
\begin{gather*}
  K_0(\gfoam_d)\vert_q = K_0(\underline{(\gfoam_d)}^{\oplus}_q).
\end{gather*}
As explained in \cref{subsec:more_cat_facts}, we call it the \emph{quantum Grothendieck ring} of $\gfoam_d$. It has the structure of a $\bZ[q,q^{-1}]$-linear category.

\begin{theorem}[Categorification theorem]
\label{thm:foam_categorification}
  The graded-2\nbd-cate\-gory of graded $\glt$-foams graded-categorifies the category of $\glt$-webs:
  \begin{equation*}
    K_0(\gfoam_d)\vert_q\cong \web_d,
  \end{equation*}
  where the isomorphism is an isomorphism of $\bZ[q,q^{-1}]$-linear categories.
\end{theorem}

This generalizes the analogous statement in the non-graded case, see e.g.\ \cite[Theorem~2.11]{BHP+_Functoriality$mathfraksl_2$Tangle_2023}.
The proof follows the same line of thought, borrowing the general strategy from \cite{KL_DiagrammaticApproachCategorification_2009}.

\begin{proof}[Proof of \cref{thm:foam_categorification}]
  Let $\gamma\colon\web_d\to K_0(\gfoam_d)\vert_q$ mapping a web $W$ to its image $[W]$ in $K_0(\gfoam_d)\vert_q$.
  The neck-cutting relation, the squeezing relation and the first braid-like relation in \cref{fig:rel_diagfoam} show that the web relations lift to isomorphisms in $\underline{(\gfoam_d)}^{\oplus}_q$, so that $\gamma$ is well-defined.
  By \cref{cor:non-degeneracy_foam}, $\gamma$ preserves a non-degenerate pairing, so it must be injective.
\end{proof}

\section{Covering Khovanov homology}
\label{sec:topo}
In this section, we define an invariant of oriented tangles and show that it coincides with odd Khovanov homology when restricted to links.
More generally, we define an invariant that coincides with \emph{covering Khovanov homology}, an invariant of links defined by Putyra \cite{Putyra_2categoryChronologicalCobordisms_2014}.
Both constructions are defined over the ring $\ringfoam$ given with elements $X$, $Y$ and $Z$ such that $X^2=Y^2=1$ (see \cref{defn:ring_R_bil})\footnote{More precisely, Putyra's original description is in the graded case.}.
The even case (choosing $X=Y=Z=1$) recovers (even) Khovanov homology, while the odd case (choosing $X=Z=1$ and $Y=-1$) recovers odd Khovanov homology\footnote{Equivalently, one can work in the graded case ($X$, $Y$ and $Z$ are formal parameters) and change the base ring at the level of chain complexes.}.

To distinguish the two constructions, we call Putyra's construction \emph{covering $\mathfrak{sl}_2$-Khovanov homology} and denote it $\slKh(L)$ for an oriented link $L$, and we call our construction \emph{covering $\mathfrak{gl}_2$-Khovanov homology} and denote it $\glKh(T)$ for an oriented tangle $T$. The latter coincides with \emph{not even Khovanov homology} \cite{Vaz_NotEvenKhovanov_2020} in the odd case (choosing $X=Z=1$ and $Y=-1$).
See also \cref{rem:NP_arc_algebra} for connections with the work of Naisse and Putyra \cite{NP_OddKhovanovHomology_2020}.

To state our claim precisely, we introduce the following completions of $\uLam_d$ (see \eqref{eq:defn_uLam}), $\web_d$ (\cref{defn:webs}) and $\gfoam_d$ (\cref{defn:cat_foam}):
\begin{definition}
  \label{defn:colimit_of_uLam_web_foam}
  The set $\uLam$, the $\bZ[q,q^{-1}]$-linear category $\web$ and the $(\bZ^2,\bilfoam)$-graded-2\nbd-cate\-gory $\gfoam$ are respectively defined as:
  \begin{align*}
    \uLam&\coloneqq
    \colim(\ldots\hookrightarrow\uLam_d\hookrightarrow\uLam_{d+2}\hookrightarrow\ldots),
    \\
    \web&\coloneqq
    \colim(\ldots\hookrightarrow\web_d\hookrightarrow\web_{d+2}\hookrightarrow\ldots),
    \\
    \gfoam&\coloneqq
    \colim(\ldots\hookrightarrow\gfoam_d\hookrightarrow\gfoam_{d+2}\hookrightarrow\ldots),
  \end{align*}
  where the embeddings denote the addition of a double point on the right, a double line on top, and a double facet at the back.
\end{definition}

The fact that the above indeed are embeddings follows from \cref{lem:web_spatial_isotopy} and \cref{thm:basis_foam}.

As an example, in the category $\web$ the following identity webs are identified:
\begin{gather*}
  \xy(0,0)*{\begin{tikzpicture}[scale=.7]
    \draw[web1] (0,0) to (2,0);
    \draw[web1] (0,.5) to (2,.5);
    \draw[web1] (0,1) to (2,1);
  \end{tikzpicture}}\endxy
  \;=\;
  \xy(0,0)*{\begin{tikzpicture}[scale=.7]
    \draw[web1] (0,0) to (2,0);
    \draw[web1] (0,.5) to (2,.5);
    \draw[web1] (0,1) to (2,1);
    \draw[web2] (0,1.5) to (2,1.5);
    \draw[web2] (0,2) to (2,2);
  \end{tikzpicture}}\endxy
  \qquad
  \text{in }\web.
\end{gather*}
Informally, working in $\uLam$, $\web$ and $\gfoam$ means that one can always ``add a double point on the right, a double line on top, and a double facet at the back''.

For $\cC$ a $(G,\bil)$-graded-2\nbd-cate\-gory, we denote by $\Ch(\cC)$ the $\ringfoam$\nbd-linear category of chain complexes in $\cC$ and chain morphisms. We say that a chain morphism is \emph{degree-preserving} if each of its components is degree-preserving. Note that in general, chain morphisms need not be degree-preserving.

Recall from \cref{subsec:more_cat_facts} that $\underline{(\gfoam)}^{\oplus}_q$ denotes the additive $q$\nbd-shif\-ted closure of $\gfoam$: it allows formal shifts of webs in the quantum grading, restricts graded foams to those preserving the quantum grading, and allows formal direct sums. As before, we view the quantum grading as independent from the $\bZ^2$-grading; hence $\underline{(\gfoam)}^{\oplus}_q$ is still a $(\bZ^2,\bilfoam)$-graded-2\nbd-cate\-gory (for $\bilfoam$ as in \cref{defn:ring_R_bil}), and shifts in quantum degree do not affect the $\bZ^2$-degree.

For each sliced oriented tangle diagram $D_T$ representing an oriented tangle $T$, we define in \cref{subsec:topo_defn_invariant} a chain complex
\begin{gather*}
  \glKom(D_T)\in \Ch(\underline{(\gfoam_d)}^{\oplus}_q).
\end{gather*}

Then:

\begin{theorem}
\label{thm:invariance}
  Let $D_T$ be a sliced oriented tangle diagram presenting an oriented tangle $T$, and denote $N_+$ and $N_-$ the number of positive and negative crossings, respectively.
  The homotopy type of $q^{-2N_++N_-}t^{N_+}\glKom(D_T)$, denoted $\glKh(T)$, is an invariant of the oriented tangle $T$.
\end{theorem}

The proof of \cref{thm:invariance} is given in \cref{subsec:proof_invariance}.
This construction is a graded analogue of the construction given in \cite{LQR_KhovanovHomologySkew_2015}. Crucially, it uses a graded horizontal composition (or ``object-adapted'' tensor product) of chain complexes (\cref{sec:complexes}), for which we give a minimal introduction in \cref{subsec:hypercubic_chain_complexes}.

In \cref{sec:odd_kh_homology}, we review the definition of covering $\slt$-Khovanov homology. For each oriented link diagram $D_L$, it associates a complex $\slKom(D_L)$ in $\Ch(\GModfoam)$, the category of chain complexes in $\bZ$-graded $\ringfoam$\nbd-modules. The homotopy type of $q^{-2N_++N_-}t^{N_+}\slKom(D_L)$ is an invariant of the oriented link $L$.

Finally, \cref{sec:equiv_covering_homologies} shows the equivalence between the two constructions, when restricted to links.
To state it, denote by $\emptyset\in\uLam$ the empty weight and by $\emptyset\coloneqq \id_\emptyset$ its identity, the empty web. Recall that in $\uLam$ (resp.\ in $\web$), the empty weight (resp.\ the empty web) is the same as an arbitrary vertical juxtaposition of double points (resp.\ double lines).
Denote by $\gfoam(\emptyset,\emptyset)$ the $\bZ$-graded $\ringfoam$\nbd-linear category obtained by restricting $\gfoam$ to the object $\emptyset$ (the $\bZ$\nbd-grading being the quantum grading), and let
\begin{gather*}
  \cA_\glt\colon\gfoam(\emptyset,\emptyset)\to\GModfoam
\end{gather*} 
be the representable functor $\cA_\glt\coloneqq\Hom_{\gfoam(\emptyset,\emptyset)}(\emptyset,-)$. It canonically extends to a functor
\begin{gather*}
  \cA_\glt\colon
  \Ch(\underline{(\gfoam)}^{\oplus}_q)(\emptyset,\emptyset)\to\Ch(\GModfoam).
\end{gather*} 
Finally, we need the following Technical Condition:

\begin{definition}
  \label{defn:technical_condition_R}
  We say that the ring $\ringfoam$ from \cref{defn:ring_R_bil} \emph{verifies the Technical Condition} if every monomial\footnote{More precisely, $a$ is an element of the abelian group generated by $X$, $Y$ and $Z$, viewed as elements of the abelian group $(\ringfoam^\times,\cdot)$.} $a$ in $X,Y,Z$ such that
  $(1-a)(1+XY)=0$
  is equal to either $1$ or $XY$.
\end{definition}

The Technical Condition is verified in all of the three canonical choices, either even, odd or graded.

We can now state the main result of this section:

\begin{theorem}
\label{thm:equivalence_with_odd}
  Assume that $\ringfoam$ verifies the Technical Condition (\cref{defn:technical_condition_R}).
  Let $D_L$ be a sliced link diagram presenting an oriented link $L$. We have the following degree-preserving isomorphism of chain complexes of $\ringfoam$\nbd-modules:
  \[\cA_\glt(\glKom(D_L))\cong\slKom(D_L).\]
\end{theorem}

The Technical Condition is only used in \cref{prop:non-ladybug-cases}. This theorem is the content of \cref{bigthm:equivalence_with_odd} in the odd case (choosing $X=Z=1$ and $Y=-1$).

%%%%%%%%%%%%%%%%%%%%%%%%%%%%%%%%%%%
%%%          invariant          %%%
%%%%%%%%%%%%%%%%%%%%%%%%%%%%%%%%%%%
\subsection{A covering \texorpdfstring{$\mathfrak{gl}_2$}{gl2}-Khovanov homology for oriented tangles}
\label{sec:defn_invariant}

\subsubsection{Composition of hypercubic chain complexes}
\label{subsec:hypercubic_chain_complexes}

We describe the horizontal composition of two \emph{hypercubic complexes} in a given graded-2\nbd-cate\-gory. Hypercubic complexes are special cases of homogeneous polycomplexes, introduced in full generality in \cref{defn:Hcomplex}. The horizontal composition that we describe is the specialization of the definitions \cref{defn:sigma_tensor_product} and \cref{defn:standard_tensor_product}\footnote{More precisely, their straightforward generalizations where graded-monoidal categories are replaced by graded-2-categories.}.
This is the minimal description necessary for the construction of covering $\glt$-Khovanov homology.

\begin{notation}
  \label{nota:basic_hypercubic_complex}
  Fix $N\in\bN$.
  We view $\{0,1\}^N$ as a hypercubic lattice and denote $(\ve_i)_{i\in\{1,\ldots,N\}}$ the canonical basis of $\bZ^N$. For each ${\vr\in\{0,1\}^N}$ and each $i\in\{1,\ldots,N\}$, we write $\vr\to\vr+\ve_i$ the corresponding edge in the hypercube $\{0,1\}^N$.
\end{notation}

Fix $\cC$ a $(G,\bil)$-graded-2\nbd-cate\-gory. Whenever we write a composition in $\cC$, it is tacitly assumed that the 1\nbd-morphisms or 2\nbd-morphisms involved are composable.

\begin{definition}
  A \emph{hypercubic complex} in $\cC$ and of dimension $N$ is a pair $\bA=(A,\alpha)$ consisting of the following data:
  \begin{enumerate}[(i)]
    \item for each vertex $\vr\in\{0,1\}^N$, a 1\nbd-morphism $A^\vr$ in $\cC$,
    \item for each edge $\vr\to\vr+\ve_i$, a homogeneous 2\nbd-morphism $\alpha^\vr_i\colon A^\vr\to A^{\vr+\ve_i}$ in $\cC$, such that each square anti-commutes:
    \begin{gather*}
      \alpha^{\vr+\ve_{i_1}}_{i_2}\starop_1\alpha^\vr_{i_1}=-\alpha^{\vr+\ve_{i_2}}_{i_1}\starop_1\alpha^\vr_{i_2}
    \end{gather*}
    for all suitable $\vr\in\{0,1\}^N$ and $i_1,i_2\in\{1,\ldots,N\}$.
    
    \item The grading is constant in a given direction, in the sense that for any $i\in\{1,\ldots,N\}$, either $\alpha^\vr_{i} = \alpha^\vs_{i}=0$ for all $\vr,\vs\in\{0,1\}^N$, or else $\deg\alpha^\vr_{i} = \deg\alpha^\vs_{i}$ for all $\vr,\vs\in\{0,1\}^N$.
    \end{enumerate}
\end{definition}

Given such a hypercubic complex $\bA=(A,\alpha)$, we define the following element of $G$:
\begin{gather*}
  \abs{\alpha}(\vr) \coloneqq \sum_{i\colon\vr_i=1} \deg_G(\alpha^\mathbf{0}_i),
\end{gather*}
where $\mathbf{0}\coloneqq (0,\ldots,0)\in\{0,1\}^N$ and we abused notation setting $\deg_G(0)=0$. Alternatively, the element $\abs{\alpha}(\vr)$ is the sum of the $G$-degrees along a path from $\mathbf{0}$ to $\vr$.

\begin{definition}
  \label{defn:tensor_product_hypercubic}
  Let $\bA=(A,\alpha)$ and $\bB=(B,\beta)$ be two hypercubic complexes of dimensions $N$ and $M$ respectively.
  The \emph{horizontal composition $\bA\starop_0\bB$} of $\bA$ and $\bB$ is the hypercubic chain complex of dimension $N+M$ defined by the following data:
  \begin{enumerate}[(i)]
    \item on each vertex $(\vr,\vs)\in\{0,1\}^{N+M}$, the 1\nbd-morphism $A^\vr\starop_0 B^\vs$;
    \item on each edge $(\vr,\vs)\to(\vr,\vs)+\ve_k$, the homogeneous 2\nbd-morphism
    \begin{IEEEeqnarray*}{c}
      % \label{eq:sign_tensor_product_hypercubic}
      (\alpha\starop_0\beta)^{(\vr,\vs)}_k\coloneqq
      % \mspace{400mu}
      % \\*[1ex]
      % \mspace{20mu}
      \begin{cases}
        \alpha^\vr_i\starop_0\id_{B^\vs} & \text{if }k=i\in\{1,\ldots,N\},\\
        (-1)^{\abs{\vr}}\bil\left(\abs{\alpha}(\vr),\beta^\vs_j\right)
        \id_{A^\vr}\starop_0\beta^\vs_j & \text{if }k=j\in\{N+1,\ldots,N+M\},
      \end{cases}
    \end{IEEEeqnarray*}
  \end{enumerate}
\end{definition}

The sign appearing above is the \emph{graded Koszul rule}. By \cref{thm:invariance_homotopy_classes}, this horizontal composition is coherent with homotopies (see also \cref{bigthm:intro_tensor_product}).

Note that a length-two chain complex whose differential is homogeneous is exactly a hypercubic complex of dimension one. In particular, if $\bA_1,\ldots,\bA_N$ is a family of horizontally composable length-two chain complexes with homogeneous differentials, \cref{defn:tensor_product_hypercubic} defines their $N$-fold horizontal composition.

\medbreak 

\begin{figure}[p]
  \centering
  \begin{gather*}
    \Circled{1}\mspace{40mu}
    \xy(0,0)*{\begin{tikzpicture}[scale=.5,rotate=90]
      % LHS
      \draw[web1] (1,2) to (1,1) to[out=-90,in=0] (.5,.5) to[out=180,in=-90] (0,1) to (0,4) to[out=90,in=180] ++(.5,.5) to[out=0,in=90] (1,4);
      % RHS
      \draw[web1] (2,2) to[out=-90,in=180] ++(.5,-.5) to[out=0,in=-90] (3,2) to (3,5) to[out=90,in=0] (2.5,5.5) to[out=180,in=90] (2,5) to (2,4);
      % bot=crossing
      \draw[web1] (1,2) to[out=90,in=-90] (2,3);
      \draw[web1,overdraw] (2,2) to[out=90,in=-90] (1,3);
      % top=crossing
      \begin{scope}[shift={(0,1)}]
        \draw[web1] (1,2) to[out=90,in=-90] (2,3);
        \draw[web1,overdraw] (2,2) to[out=90,in=-90] (1,3);
      \end{scope}
    \end{tikzpicture}}\endxy
    %%%%%%%%%%%%%%%%%%%%
    \mspace{120mu}
    %%%%%%%%%%%%%%%%%%%%
    \Circled{2}\mspace{40mu}
    \xy(0,0)*{\begin{tikzpicture}[scale=.5,rotate=90]
      % doubles
      \draw[web2] (.5,0) to (.5,.5);
      \draw[web2] (2.5,0) to (2.5,1.5);
      \draw[web2] (.5,4.5) to[out=90,in=-90] (2.5,6) to (2.5,7);
      \draw[web2] (1.5,6.5) to (1.5,7);
      % LHS
      \draw[web1] (1,2) to (1,1) to[out=-90,in=0] (.5,.5) to[out=180,in=-90] (0,1) to (0,4) to[out=90,in=180] ++(.5,.5) to[out=0,in=90] (1,4);
      % RHS
      \draw[web1] (2,2) to[out=-90,in=180] ++(.5,-.5) to[out=0,in=-90] (3,2) to (3,5);
      \draw[web1,overdraw] (3,5) to[out=90,in=-90] (2,6);
      \draw[web1] (2,6) to[out=90,in=0] ++(-.5,.5) to[out=180,in=90] (1,6);
      \draw[web1,overdraw] (1,6) to[out=-90,in=90] (2,4.5);
      \draw[web1] (2,4.5) to (2,4);
      % bot=crossing
      \draw[web1] (1,2) to[out=90,in=-90] (2,3);
      \draw[web1,overdraw] (2,2) to[out=90,in=-90] (1,3);
      % top=crossing
      \begin{scope}[shift={(0,1)}]
        \draw[web1] (1,2) to[out=90,in=-90] (2,3);
        \draw[web1,overdraw] (2,2) to[out=90,in=-90] (1,3);
      \end{scope}
    \end{tikzpicture}}\endxy
    %%%%%%%%%%%%%%%%%%%%
    \\[10ex]
    %%%%%%%%%%%%%%%%%%%%
    \xy(0,0)*{\begin{tikzpicture}
      \node at (-5,4) {$\Circled{3}$};
      \node (L) at (-4,0) {
        $q^2t^{-2}$\;\xy(0,0)*{\begin{tikzpicture}[scale=.4,rotate=90]
          % doubles
          \draw[web2] (.5,0) to (.5,.5);
          \draw[web2] (2.5,0) to (2.5,1.5);
          \draw[web2] (.5,4.5) to[out=90,in=-90] (2.5,6) to (2.5,7);
          \draw[web2] (1.5,6.5) to (1.5,7);
          % LHS
          \draw[web1] (1,2) to (1,1) to[out=-90,in=0] (.5,.5) to[out=180,in=-90] (0,1) to (0,4) to[out=90,in=180] ++(.5,.5) to[out=0,in=90] (1,4);
          % RHS
          \draw[web1] (2,2) to[out=-90,in=180] ++(.5,-.5) to[out=0,in=-90] (3,2) to (3,5);
          \draw[web1,overdraw] (3,5) to[out=90,in=-90] (2,6);
          \draw[web1] (2,6) to[out=90,in=0] ++(-.5,.5) to[out=180,in=90] (1,6);
          \draw[web1,overdraw] (1,6) to[out=-90,in=90] (2,4.5);
          \draw[web1] (2,4.5) to (2,4);
          % bot = 0-resolution
          \draw[web2] (1.5,2.3) to (1.5,3-.3);
          \draw[web1] (1,2) to[out=90,in=90] (2,2);
          \draw[web1] (1,3) to[out=-90,in=-90] (2,3);
          % top = 0-resolution
          \begin{scope}[shift={(0,1)}]
            \draw[web2] (1.5,2.3) to (1.5,3-.3);
            \draw[web1] (1,2) to[out=90,in=90] (2,2);
            \draw[web1] (1,3) to[out=-90,in=-90] (2,3);
          \end{scope}
        \end{tikzpicture}}\endxy
      };
      \node (R) at (4,0) {
        \begin{tikzpicture}[scale=.4,rotate=90]
          % doubles
          \draw[web2] (.5,0) to (.5,.5);
          \draw[web2] (2.5,0) to (2.5,1.5);
          \draw[web2] (.5,4.5) to[out=90,in=-90] (2.5,6) to (2.5,7);
          \draw[web2] (1.5,6.5) to (1.5,7);
          % LHS
          \draw[web1] (1,2) to (1,1) to[out=-90,in=0] (.5,.5) to[out=180,in=-90] (0,1) to (0,4) to[out=90,in=180] ++(.5,.5) to[out=0,in=90] (1,4);
          % RHS
          \draw[web1] (2,2) to[out=-90,in=180] ++(.5,-.5) to[out=0,in=-90] (3,2) to (3,5);
          \draw[web1,overdraw] (3,5) to[out=90,in=-90] (2,6);
          \draw[web1] (2,6) to[out=90,in=0] ++(-.5,.5) to[out=180,in=90] (1,6);
          \draw[web1,overdraw] (1,6) to[out=-90,in=90] (2,4.5);
          \draw[web1] (2,4.5) to (2,4);
          % bot = 1-resolution
          \draw[web1] (1,2) to (1,3);
          \draw[web1] (2,2) to (2,3);
          % top = 0-resolution
          \begin{scope}[shift={(0,1)}]
            \draw[web1] (1,2) to (1,3);
            \draw[web1] (2,2) to (2,3);
          \end{scope}
          %
          % \draw[snakecd] (-.5,0) to (-.5,7);
        \end{tikzpicture}
      };
      \node (T) at (0,3) {
        $qt^{-1}$\;\xy(0,0)*{\begin{tikzpicture}[scale=.4,rotate=90]
          % doubles
          \draw[web2] (.5,0) to (.5,.5);
          \draw[web2] (2.5,0) to (2.5,1.5);
          \draw[web2] (.5,4.5) to[out=90,in=-90] (2.5,6) to (2.5,7);
          \draw[web2] (1.5,6.5) to (1.5,7);
          % LHS
          \draw[web1] (1,2) to (1,1) to[out=-90,in=0] (.5,.5) to[out=180,in=-90] (0,1) to (0,4) to[out=90,in=180] ++(.5,.5) to[out=0,in=90] (1,4);
          % RHS
          \draw[web1] (2,2) to[out=-90,in=180] ++(.5,-.5) to[out=0,in=-90] (3,2) to (3,5);
          \draw[web1,overdraw] (3,5) to[out=90,in=-90] (2,6);
          \draw[web1] (2,6) to[out=90,in=0] ++(-.5,.5) to[out=180,in=90] (1,6);
          \draw[web1,overdraw] (1,6) to[out=-90,in=90] (2,4.5);
          \draw[web1] (2,4.5) to (2,4);
          % bot = 1-resolution
          \draw[web1] (1,2) to (1,3);
          \draw[web1] (2,2) to (2,3);
          % top = 0-resolution
          \begin{scope}[shift={(0,1)}]
            \draw[web2] (1.5,2.3) to (1.5,3-.3);
            \draw[web1] (1,2) to[out=90,in=90] (2,2);
            \draw[web1] (1,3) to[out=-90,in=-90] (2,3);
          \end{scope}
        \end{tikzpicture}}\endxy
      };
      \node (B) at (0,-3) {
        $qt^{-1}$\;\xy(0,0)*{\begin{tikzpicture}[scale=.4,rotate=90]
          % doubles
          \draw[web2] (.5,0) to (.5,.5);
          \draw[web2] (2.5,0) to (2.5,1.5);
          \draw[web2] (.5,4.5) to[out=90,in=-90] (2.5,6) to (2.5,7);
          \draw[web2] (1.5,6.5) to (1.5,7);
          % LHS
          \draw[web1] (1,2) to (1,1) to[out=-90,in=0] (.5,.5) to[out=180,in=-90] (0,1) to (0,4) to[out=90,in=180] ++(.5,.5) to[out=0,in=90] (1,4);
          % RHS
          \draw[web1] (2,2) to[out=-90,in=180] ++(.5,-.5) to[out=0,in=-90] (3,2) to (3,5);
          \draw[web1,overdraw] (3,5) to[out=90,in=-90] (2,6);
          \draw[web1] (2,6) to[out=90,in=0] ++(-.5,.5) to[out=180,in=90] (1,6);
          \draw[web1,overdraw] (1,6) to[out=-90,in=90] (2,4.5);
          \draw[web1] (2,4.5) to (2,4);
          % bot = 0-resolution
          \draw[web2] (1.5,2.3) to (1.5,3-.3);
          \draw[web1] (1,2) to[out=90,in=90] (2,2);
          \draw[web1] (1,3) to[out=-90,in=-90] (2,3);
          % top = 1-resolution
          \begin{scope}[shift={(0,1)}]
            \draw[web1] (1,2) to (1,3);
            \draw[web1] (2,2) to (2,3);
          \end{scope}
        \end{tikzpicture}}\endxy
      };
      \draw[->] (L) to[bend left=20] (T);
      \draw[->] (L) to[bend right=20] (B);
      \draw[->] (T) to[bend left=20] (R);
      \draw[->] (B) to[bend right=20] (R);
      \node[anchor=south east] at (-2.5,2.2) {
        \begin{tikzpicture}[xscale=.8,scale=.6]
          \draw[diag1=1.5,<-] (0,0) to (0,1) node[above]{\small 1};
          \draw[diag3=1.5,<-] (.5,0) to (.5,1) node[above]{\small 3};
          \draw[diag2=1.5,->] (1,0) to (1,1) node[above]{\small 2};
          \draw[diag2=1.5,<-] (1.5,0) to (1.5,1);
          \draw[diag1=1.5,->] (2,0) to (2,1);
          \draw[diag1=1.5,<-] (2.5,0) to (2.5,1);
          \draw[diag2=1.5,->] (3,0) to (3,1);
          \draw[diag2=1.5,<-] (3.5,0) to (3.5,1);
          \draw[diag2=1.5,->] (4,0) to[out=90,in=180] (4.5,.5) to[out=0,in=90] (5,0);
          \draw[diag3=1.5,->] (5.5,0) to (5.5,1);
          \draw[diag1=1.5,->] (6,0) to (6,1);
        \end{tikzpicture}
      };
      \node[anchor=north east] at (-2.5,-2.2) {
        \begin{tikzpicture}[xscale=.8,scale=.6]
          \draw[diag1=1.5,<-] (0,0) to (0,1);
          \draw[diag3=1.5,<-] (.5,0) to (.5,1);
          \draw[diag2=1.5,->] (1,0) to (1,1);
          \draw[diag2=1.5,<-] (1.5,0) to (1.5,1);
          \draw[diag1=1.5,->] (2,0) to (2,1);
          \draw[diag1=1.5,<-] (2.5,0) to (2.5,1);
          \draw[diag2=1.5,->] (3,0) to[out=90,in=180] (3.5,.5) to[out=0,in=90] (4,0);
          \draw[diag2=1.5,->] (4.5,0) to (4.5,1);
          \draw[diag2=1.5,<-] (5,0) to (5,1);
          \draw[diag3=1.5,->] (5.5,0) to (5.5,1);
          \draw[diag1=1.5,->] (6,0) to (6,1);
        \end{tikzpicture}
      };
      \node[anchor=south west] at (2.5,2.2) {
        \begin{tikzpicture}[xscale=.8,scale=.6]
          \draw[diag1=1.5,<-] (0,0) to (0,1);
          \draw[diag3=1.5,<-] (.5,0) to (.5,1);
          \draw[diag2=1.5,->] (1,0) to (1,1);
          \draw[diag2=1.5,<-] (1.5,0) to (1.5,1);
          \draw[diag1=1.5,->] (2,0) to (2,1);
          \draw[diag1=1.5,<-] (2.5,0) to (2.5,1);
          \draw[diag2=1.5,->] (3,0) to[out=90,in=180] (3.5,.5) to[out=0,in=90] (4,0);
          \draw[diag3=1.5,->] (4.5,0) to (4.5,1);
          \draw[diag1=1.5,->] (5,0) to (5,1);
        \end{tikzpicture}
      };
      \node[anchor=north west] at (2.1,-2.4) {$-Y$
        \xy(0,0)*{\begin{tikzpicture}[xscale=.8,scale=.6]
          \draw[diag1=1.5,<-] (0,0) to (0,1);
          \draw[diag3=1.5,<-] (.5,0) to (.5,1);
          \draw[diag2=1.5,->] (1,0) to (1,1);
          \draw[diag2=1.5,<-] (1.5,0) to (1.5,1);
          \draw[diag1=1.5,->] (2,0) to (2,1);
          \draw[diag1=1.5,<-] (2.5,0) to (2.5,1);
          \draw[diag2=1.5,->] (3,0) to[out=90,in=180] (3.5,.5) to[out=0,in=90] (4,0);
          \draw[diag3=1.5,->] (4.5,0) to (4.5,1);
          \draw[diag1=1.5,->] (5,0) to (5,1);
        \end{tikzpicture}}\endxy
      };
    \end{tikzpicture}}\endxy
  \end{gather*}

  \caption{Defining procedure for $\glKom$, in the case of a sliced tangle diagram presenting the Hopf link. The variables $q$ and $t$ respectively denote shift in quantum and homological degrees.
  Both differentials have $\bZ^2$-degree $(0,1)$: one checks that the graded Koszul rule in \cref{defn:tensor_product_hypercubic} adds the scalar $-Y$, as pictured in the figure.}
  \label{fig:example_invariant} 
\end{figure}

\subsubsection{Definition of covering \texorpdfstring{$\glt$}{gl2}-Khovanov homology}
\label{subsec:topo_defn_invariant}

We now define a chain complex $\glKom(D)\in\Ch(\underline{(\gfoam)}^{\oplus}_q)$ for every sliced tangle diagram $D$.
In fact, $\glKom(D)$ is independent of the orientation of $D$, and we shall not mention it again in this section.
The reader can follow the procedure on the example given in \cref{fig:example_invariant}, with $D$ pictured at step \Circled{1}.

We shall need the following definitions of ``mixed crossings'',  between a single line and a double line:
\begin{equation}
\label{eq:defn_web_rb}
  {}\xy(0,0)*{\begin{tikzpicture}[scale=.6,xscale=1.2,transform shape,rotate=90]
    \pic at (0,0) {webcr_rb};
  \end{tikzpicture}}\endxy
  \mspace{10mu}\coloneqq\mspace{10mu}
  \xy(0,0)*{\begin{tikzpicture}[scale=.6,xscale=.8,yscale=.6,transform shape,rotate=90]
    \pic at (0,0) {web-};
    \pic at (1,1) {web+};
    \draw[web1] (1.5,-.5) to (1.5,.5);
    \draw[web1] (-.5,.5) to (-.5,1.5);
  \end{tikzpicture}}\endxy
  \qquad\an\qquad
  \xy(0,0)*{\begin{tikzpicture}[scale=.6,xscale=1.2,transform shape,rotate=90]
    \pic at (0,0) {webcr_br};
  \end{tikzpicture}}\endxy
  \mspace{10mu}\coloneqq\mspace{10mu}
  \xy(0,0)*{\begin{tikzpicture}[scale=.6,xscale=-.8,yscale=.6,transform shape,rotate=90]
    \pic at (0,0) {web-};
    \pic at (1,1) {web+};
    \draw[web1] (1.5,-.5) to (1.5,.5);
    \draw[web1] (-.5,.5) to (-.5,1.5);
  \end{tikzpicture}}\endxy\;.
\end{equation}
These mixed crossings satisfy the following relations in $\web$:
\begin{gather}
  \label{eq:rel_web_rb}
  {}\xy(0,0)*{\begin{tikzpicture}[scale=.5,transform shape,rotate=90]
    \pic at (0,0) {webcr_rb};
    \pic at (0,1) {webcr_br};
  \end{tikzpicture}}\endxy
  \mspace{10mu}=\mspace{10mu}
  \xy(0,0)*{\begin{tikzpicture}[scale=.5,transform shape,rotate=90]
    \draw[web1] (.5,-.5) to (.5,1.5);
    \draw[web2] (-.5,-.5) to (-.5,1.5);
  \end{tikzpicture}}\endxy
  \qquad\an\qquad
  \xy(0,0)*{\begin{tikzpicture}[scale=.3,transform shape,rotate=90]
    \pic at (1.5,0.5) {web-};
    \pic at (0,1) {webcr_rb};
    \pic at (1,2) {webcr_rb};
    \draw[web2] (0,0) to (0,1);
    \draw[web1] (2,1) to (2,2);
    \draw[web1] (0,2) to (0,3);
  \end{tikzpicture}}\endxy
  \mspace{10mu}=\mspace{10mu}
  \xy(0,0)*{\begin{tikzpicture}[scale=.3,transform shape,rotate=90]
    \pic at (0,0) {web-};
    \draw[web2] (1.5,-1.5) to (1.5,1.5);
    \draw[web2] (0,-1.5) to (0,-.5);
    \draw[web1] (-.5,.5) to (-.5,1.5);
    \draw[web1] (.5,.5) to (.5,1.5);
  \end{tikzpicture}}\endxy\;.
\end{gather}
as well as all the relations obtained from the above by reflecting vertically and horizontally.
\medbreak
The procedure starts by telling how to assign a web to an elementary flat tangle diagram, that is, to a cap and a cup (see \cite[section~4A2]{LQR_KhovanovHomologySkew_2015}). There is more than one web $W$ such that $\undcomp(W)$ is a cup (or a cap), but we fix a choice by fixing the endpoints.
To enforce the given endpoints, we use the mixed crossings \eqref{eq:defn_web_rb}. 

Say that $\lambda\in\uLam$ is \emph{antidominant} if it is antidominant as a $\glt$-weight, that is, if it is (non-strictly) increasing. To any set of $n$ points on a line corresponds a unique antidominant weight $\lambda\in\uLam_d$ for $n\leq d$ and $n=d\mod 2$. In turn, those antidominant weights define a unique element in $\uLam$.
Given any elementary flat tangle diagram, we pick a web representative whose endpoints are antidominant by ``adding a double line to the cup or cap and sliding it to the top''. For instance:
\begin{equation*}
  \xy (0,0)*{\begin{tikzpicture}[scale=.5,rotate=90,transform shape]
    \draw[line width=1.5pt] (0+.5,0) to[out=90,in=180] (.5+.5,.5+.5) to[out=0,in=90] (1+.5,0);
    \draw[line width=1.5pt] (2,0) to (2,3);
    \draw[line width=1.5pt] (3,0) to (3,3);
  \end{tikzpicture}}\endxy
  \mspace{10mu}\mapsto\mspace{10mu}
  \xy (0,0)*{\begin{tikzpicture}[scale=.5,transform shape,rotate=90]
    \pic at (0,1) {web+};
    \pic at (0,1.5) {webcr_rb};
    \pic at (1,2.5) {webcr_rb};
    \draw[web1] (1,.5) to (1,1.5);
    \draw[web1] (2,.5) to (2,2.5);
    \draw[web1] (0,2.5) to (0,3.5);
  \end{tikzpicture}}\endxy
\end{equation*}
Note that fixing a choice for the endpoints ensures that if two elementary flat tangle diagrams are composable, then so are the corresponding webs in $\web$. We may extend this assignment to non-flat tangle diagrams by formally adjoining crossings to our web diagrammatics (see step \Circled{2} in \cref{fig:example_invariant}).

The procedure extends to chain complexes in $\Ch(\underline{(\gfoam)}^{\oplus}_q)$. For cups and caps, it is the chain complex concentrated in homological degree 0 corresponding to the previously assigned web. For crossings, this is given by the Khovanov--Blanchet bracket, generalized to the graded case:
\begin{IEEEeqnarray*}{rCrl}
  {}\xy(0,0)*{\begin{tikzpicture}[scale=-1*.7,xscale=1.2,rotate=90,transform shape]
  \draw[black,line width=1.5pt] (1,0) to[out=90,in=-90] (0,1);
  \draw[black,line width=1.5pt,overdraw=8pt] (0,0) to[out=90,in=-90] (1,1);
  \end{tikzpicture}}\endxy
  \qquad&\mapsto&\qquad
  % \left[\mspace{10mu}
  qt^{-1}\;
  \xy(0,0)*{\begin{tikzpicture}[scale=.7,rotate=90,transform shape]
    \draw[web1] (0,0) to (0,2);
    \draw[web1] (1,0) to (1,2);
    % \draw[snakecd] (-.5,0) to (-.5,2);
  \end{tikzpicture}}\endxy
  \quad&\xrightarrow{
    \xy(0,0)*{\begin{tikzpicture}[scale=.7,transform shape]
      \pic[diag1] at (0,0) {lcup};
    \end{tikzpicture}}\endxy
  }\quad
  \xy(0,0)*{\begin{tikzpicture}[scale=1*.7,transform shape,rotate=90]
    \pic at (0,0) {web+};
    \pic at (0,1) {web-};
  \end{tikzpicture}}\endxy
  % \mspace{10mu}\right]
  \\[2ex]
  {}\xy(0,0)*{\begin{tikzpicture}[scale=-1*.7,xscale=1.2,rotate=90,transform shape]
  \draw[black,line width=1.5pt] (0,0) to[out=90,in=-90] (1,1);
  \draw[black,line width=1.5pt,overdraw=8pt] (1,0) to[out=90,in=-90] (0,1);
  \end{tikzpicture}}\endxy
  \qquad&\mapsto&\qquad
  % \left[\mspace{10mu}
  qt^{-1}\;
  \xy(0,0)*{\begin{tikzpicture}[scale=1*.7,transform shape,rotate=90]
    \pic at (0,0) {web+};
    \pic at (0,1) {web-};
  \end{tikzpicture}}\endxy
  \quad&\xrightarrow{
    \xy(0,0)*{\begin{tikzpicture}[scale=.7,transform shape]
      \pic[diag1] at (0,0) {rcap};
    \end{tikzpicture}}\endxy
  }\quad
  \;\xy(0,0)*{\begin{tikzpicture}[scale=1*.7,rotate=90,transform shape]
    \draw[web1] (0,0) to (0,2);
    \draw[web1] (1,0) to (1,2);
  \end{tikzpicture}}\endxy
  % \mspace{10mu}\right]
\end{IEEEeqnarray*}
The variables $q$ and $t$ respectively denote shift in quantum and homological degrees.
Note that the differentials are degree-preserving with respect to the quantum degree. However, they are not degree-preserving with respect to the $\bZ^2$-degree: the former has $\bZ^2$-degree $(1,0)$, while the latter has $\bZ^2$-degree $(0,1)$. If one restricts to the odd case $(X=Z=1$ and $Y=-1$), the former has even parity and the latter has odd parity.

Finally, let $D$ be a sliced tangle diagram. Then $\glKom(D)$ is defined as the horizontal composition (see \cref{defn:tensor_product_hypercubic}) of the chain complexes assigned to each slice of $D$.
This ends the definition of $\glKom(D)$ (see step \Circled{3} in \cref{fig:example_invariant}).
\hfill$\diamond$

\subsubsection{Proof of invariance}
\label{subsec:proof_invariance}

In this subsection, we prove \cref{thm:invariance}.

Since the horizontal composition of chain complexes is coherent with homotopies (\cref{thm:invariance_homotopy_classes}), the proof can be done locally. It suffices to check invariance under the Reide\-mei\-ster--Turaev moves for sliced oriented tangle diagrams (see e.g.\ \cite{Ohtsuki_QuantumInvariantsStudy_2002}).
We split the proof in two lemmas, the first one (\cref{lem:covering_khovanov_invariance_planar_isotopies}) dealing with planar isotopies and the second one (\cref{lem:covering_khovanov_invariance_reidemeister_moves}) dealing with Reidemeister moves.

\begin{lemma}
  \label{lem:covering_khovanov_invariance_planar_isotopies}
  Let $D_1$ and $D_2$ be two sliced tangle diagrams.
  If $D_1$ and $D_2$ are planar isotopic, then $\glKom(D_1)$ and $\glKom(D_2)$ are isomorphic.
\end{lemma}

\begin{proof}
  It suffices to check invariance under elementary planar isotopies for sliced tangle diagrams, as given by \cref{fig:Turaev_moves_planar_isotopies}.
  \begin{figure}[h]
    \def\scl{.8}
    \centering
    \begin{gather*}
      % homotopy of diagrams 2
      % \label{eq:turaev_homotopy2}
      {}\xy (0,0)*{\begin{tikzpicture}[thin,scale=.7*\scl,rotate=90]
        \draw[very thick] (.7,0) to (.7,3);
        \draw[very thick] (1.3,0) to (1.3,3);
        \node[draw,fill=white,rectangle,rounded corners,line width=1,inner sep=7pt,align=center] at (1,.8) {$T$};
        \node[,rotate=90] at (2,1.5) {\ldots};
        \begin{scope}[shift={(2,0)}]
          \draw[very thick] (.7,0) to (.7,3);
          \draw[very thick] (1.3,0) to (1.3,3);
          \node[draw,fill=white,rectangle,rounded corners,line width=1,inner sep=7pt,align=center] at (1,2.2) {$T'$};
        \end{scope}
      \end{tikzpicture}}\endxy
      \mspace{7mu}\longleftrightarrow\mspace{7mu}
      \xy (0,0)*{\begin{tikzpicture}[thin,scale=.7*\scl,rotate=90]
        \draw[very thick] (.7,0) to (.7,3);
        \draw[very thick] (1.3,0) to (1.3,3);
        \node[draw,fill=white,rectangle,rounded corners,line width=1,inner sep=7pt,align=center] at (1,2.2) {$T$};
        \node[,rotate=90] at (2,1.5) {\ldots};
        \begin{scope}[shift={(2,0)}]
          \draw[very thick] (.7,0) to (.7,3);
          \draw[very thick] (1.3,0) to (1.3,3);
          \node[draw,fill=white,rectangle,rounded corners,line width=1,inner sep=7pt,align=center] at (1,.8) {$T'$};
        \end{scope}
    \end{tikzpicture}}\endxy
    \\[2ex]
    % zigzags
    % \label{eq:turaev_zigzag}
    {}\xy (0,0)*{\begin{tikzpicture}[very thick,xscale=0.6,yscale=.4,scale=\scl,rotate=90]
      \draw[] (0,0) to[out=90,in=180] (0.5,1.5) to[out=0,in=120] (1,1) to[out=-60,in=180] (1.5,0.5) to[out=0,in=-90] (2,2);
    \end{tikzpicture}}\endxy
    \mspace{7mu}\longleftrightarrow\mspace{7mu}
    \xy (0,0)*{\begin{tikzpicture}[very thick,scale=0.6*\scl,rotate=90]
      \draw[] (0,0) to (0,2);
    \end{tikzpicture}}\endxy
    \mspace{7mu}\longleftrightarrow\mspace{7mu}
    \xy (0,0)*{\begin{tikzpicture}[very thick,xscale=0.6,yscale=-.4,scale=\scl,rotate=90]
      \draw[] (0,0) to[out=90,in=180] (0.5,1.5) to[out=0,in=120] (1,1) to[out=-60,in=180] (1.5,0.5) to[out=0,in=-90] (2,2);
    \end{tikzpicture}}\endxy
    \\[2ex]
    {}\xy (0,0)*{\begin{tikzpicture}[scale=.5*\scl,very thick,rotate=90]
      \draw (1,0) to[out=45,in=-45] (1,2);
      \draw[overdraw] (0,0) to[out=70,in=180] (0.5,1.2) to[out=0,in=120] (2,0);
    \end{tikzpicture}}\endxy
    \mspace{7mu}\longleftrightarrow\mspace{7mu}
    \xy (0,0)*{\begin{tikzpicture}[scale=.5*\scl,yscale=-1,very thick,rotate=90]
      \draw (1,0) to[out=45,in=-45] (1,2);
      \draw[overdraw] (0,0) to[out=70,in=180] (0.5,1.2) to[out=0,in=120] (2,0);
    \end{tikzpicture}}\endxy
    \mspace{80mu}
    \xy (0,0)*{\begin{tikzpicture}[scale=.5*\scl,xscale=-1,very thick,rotate=90]
      \draw (1,0) to[out=45,in=-45] (1,2);
      \draw[overdraw] (0,0) to[out=70,in=180] (0.5,1.2) to[out=0,in=120] (2,0);
    \end{tikzpicture}}\endxy
    \mspace{7mu}\longleftrightarrow\mspace{7mu}
    \xy (0,0)*{\begin{tikzpicture}[scale=-.5*\scl,very thick,rotate=90]
      \draw (1,0) to[out=45,in=-45] (1,2);
      \draw[overdraw] (0,0) to[out=70,in=180] (0.5,1.2) to[out=0,in=120] (2,0);
    \end{tikzpicture}}\endxy
    \end{gather*}
    \caption{Elementary planar moves for sliced tangle diagrams. Here $T$ (resp.\ $T'$) denotes a crossing, a cup or a cap.}
    \label{fig:Turaev_moves_planar_isotopies}
  \end{figure}

  If the elementary planar isotopy does not involve any crossing, the complex is concentrated in a single homological degree. Finding an isomorphism of complexes reduces to finding an isomorphism of webs. Thanks to the categorification theorem (\cref{thm:foam_categorification}) and \cref{lem:web_spatial_isotopy}, two webs $W_1$ and $W_2$ are isomorphic in $\underline{(\gfoam)}^{\oplus}_q$ precisely if $\undcomp(W_1)$ and $\undcomp(W_2)$ are isotopic.

  If the elementary planar isotopy contains exactly one crossing, the assigned complexes are length-two complexes with isotopic webs as vertices. Fix a pair of isomorphisms between these webs. Thanks to \cref{lem:foam_all_isotopy_are_rels} (if two foams have the same underlying surface, they are equal up to invertible scalar), this pair defines a chain morphism up to invertible scalar. Renormalizing provides a genuine chain isomorphism.

  Finally, the only elementary planar isotopy with at least two crossings is the planar isotopy interchanging two crossings.
  The associated chain complexes have the form $F\starop_1 F'$ and $F'\starop_1 F$ respectively, for some pair of length-two complexes $F$ and $F'$.
  A chain isomorphism $F\starop_1 F'\cong F'\starop_1 F$ is given by suitable pointwise compositions of foam crossings, exchanging $F$ with $F'$.
  This commutes thanks to braid-like relations.
\end{proof}

To show invariance under Reidemeister moves, we follow Bar-Natan's strategy for the even case \cite{Bar-Natan_KhovanovsHomologyTangles_2005,Bar-Natan_FastKhovanovHomology_2007}, using delooping and Gaussian elimination.
\emph{Delooping} denotes using circle evaluation in $\web$ to remove a circle,
while \emph{Gaussian elimination} is the following general homological fact \cite[Lemma~3.2]{Bar-Natan_FastKhovanovHomology_2007}:

\begin{samepage}
  \begin{lemma}
  \label{lem:iso_in_square_complex}
    In any additive category, if $\alpha$ is an isomorphism in the complex $C_\bullet$ given by
    \[\begin{tikzcd}[row sep=3,cramped]
    & X \ar[dr,"\gamma"] \ar[dd,"\oplus"{description},phantom,dashed] & \\
    W \ar[ur,"\alpha"] \ar[dr,"\beta"'] &  & Z \\
    & Y \ar[ur,"\delta"'] &
    \end{tikzcd}\;,\]
    then there exists a homotopy equivalence of complexes $C_\bullet\to C'_\bullet$ where $C'_\bullet$ is the complex $\begin{tikzcd}[cramped]
    Y \ar[r,"\delta"] & Z
    \end{tikzcd}$. Moreover, this homotopy equivalence is a strong deformation retract.
    Similarly, if $\gamma$ is an isomorphism there exists a strong deformation retract from $C_\bullet$ into $\begin{tikzcd}[cramped]
    W \ar[r,"\beta"] & Y
    \end{tikzcd}$.
    \hfill\qed
  \end{lemma}
\end{samepage}

One need not know the definition of a strong deformation retract, except for its appearance in \cref{lem:cone_and_strongdef} below.

\begin{figure}[ht]
  \def\scl{.9}
  \centering
  \begin{gather*}
    % R1
    % \label{eq:turaev_R1}
    {}\xy (0,1.3)* {\begin{tikzpicture}[very thick,scale=.7*\scl,rotate=-90]
      \draw (.8,1) to[out=-90,in=0] (.5,.6) to[out=180,in=-90] (0,2);
      \draw[<-,overdraw] (0,0) to[out=90,in=180] (.5,1.4) to[out=0,in=90] (.8,1);
    \end{tikzpicture}}\endxy
    \mspace{7mu}\longleftrightarrow\mspace{7mu}
    \xy (0,0)* {\begin{tikzpicture}[very thick,scale=.7*\scl,rotate=-90]
      \draw[<-] (0,0) to (0,2);
    \end{tikzpicture}}\endxy
    \mspace{7mu}\longleftrightarrow\mspace{7mu}
    \xy (0,0)* {\begin{tikzpicture}[very thick,scale=.7*\scl,rotate=-90]
      \draw[<-] (0,0) to[out=90,in=180] (.5,1.4) to[out=0,in=90] (.8,1);
      \draw[overdraw] (.8,1) to[out=-90,in=0] (.5,.6) to[out=180,in=-90] (0,2);
    \end{tikzpicture}}\endxy
    \\[2ex]
    % R2
    % \label{eq:turaev_R2}
    {}\xy (0,0)* {\begin{tikzpicture}[very thick,yscale=.5,xscale=.7,scale=\scl,rotate=-90]
      \draw[<-] (0,0) to[out=45,in=-90] (1,1) to[out=90,in=-45] (0,2);
      \draw[<-,overdraw] (1,0) to[out=135,in=-90] (0,1) to[out=90,in=-135] (1,2);
    \end{tikzpicture}}\endxy
    \mspace{7mu}\longleftrightarrow\mspace{7mu}
    \xy (0,0)* {\begin{tikzpicture}[very thick,yscale=.5,xscale=.7,scale=\scl,rotate=-90]
        \draw[<-] (0,0) to (0,2);
        \draw[<-] (1,0) to (1,2);
    \end{tikzpicture}}\endxy
    \mspace{80mu}
    % R2
    % \label{eq:turaev_R2}
    {}\xy (0,0)* {\begin{tikzpicture}[very thick,yscale=.5,xscale=.7,scale=\scl,rotate=-90]
      \draw[->] (0,0) to[out=45,in=-90] (1,1) to[out=90,in=-45] (0,2);
      \draw[<-,overdraw] (1,0) to[out=135,in=-90] (0,1) to[out=90,in=-135] (1,2);
    \end{tikzpicture}}\endxy
    \mspace{7mu}\longleftrightarrow\mspace{7mu}
    \xy (0,0)* {\begin{tikzpicture}[very thick,yscale=.5,xscale=.7,scale=\scl,rotate=-90]
        \draw[->] (0,0) to (0,2);
        \draw[<-] (1,0) to (1,2);
    \end{tikzpicture}}\endxy
    \\[2ex]
    % R3
    % \label{eq:turaev_R3}
    \xy (0,0)* {\begin{tikzpicture}[very thick,scale=.7*\scl,rotate=-90]
        \draw[<-] (1.6,0) to[out=120,in=-60] (0,2);
        \draw[<-,overdraw] (.8,0) to[out=120,in=-90] (0,1) to[out=90,in=-120] (.8,2);
        \draw[<-,overdraw] (0,0) to[out=60,in=-120] (1.6,2);
    \end{tikzpicture}}\endxy
    \mspace{7mu}\longleftrightarrow\mspace{7mu}
    \xy (0,0)* {\begin{tikzpicture}[very thick,scale=-.7*\scl,rotate=-90]
        \draw[->] (1.6,0) to[out=120,in=-60] (0,2);
        \draw[->,overdraw] (.8,0) to[out=120,in=-90] (0,1) to[out=90,in=-120] (.8,2);
        \draw[->,overdraw] (0,0) to[out=60,in=-120] (1.6,2);
    \end{tikzpicture}}\endxy
  \end{gather*}
  \caption{Reidemeister moves for sliced oriented tangle diagrams.}
  \label{fig:reidemeister_moves}
\end{figure}

\begin{lemma}
  \label{lem:covering_khovanov_invariance_reidemeister_moves}
  Let $D_1$ and $D_2$ be two oriented sliced tangle diagrams.
  Denote $N_+^i$ (resp.\ $N_-^i$) the number of positive (resp.\ negative) crossings in $D_i$.
  If $D_1$ and $D_2$ are related by a Reidemeister move (see \cref{fig:reidemeister_moves}), then
  \[
    q^{-2N_+^1+N_-^1}t^{N_+^1}\glKom(D_1)
    \an
    q^{-2N_+^2+N_-^2}t^{N_+^2}\glKom(D_2)
  \]
  are homotopic.
\end{lemma}

\begin{proof}
  Consider first the Reidemeister I move. The proof of invariance is essentially contained in \cref{fig:invariance_R1}. On the left, the complex associated to the left-hand side of the move. On the right, the same complex after delooping in homological degree zero, and simplifying in homological degree one. By \cref{lem:more_rels_diagfoam}, the bottom arrow is an isomorphism, so Gaussian elimination gives a homotopy equivalence with the left-top web.
  It is shifted in quantum and homological degrees, but this is fixed by the renormalization.
  Invariance under the other Reidemeister I move is proved similarly.

  \begin{figure}[h]
    \begin{gather*}
      \begingroup
        % {}\left(
          qt^{-1}\;\xy(0,0)*{\begin{tikzpicture}[scale=.5,rotate=90]
          \draw[web1] (0,-.5) to (0,1.5);
          \pic[transform shape] at (1,0) {web-};
          \pic[transform shape] at (1,1) {web+};
          % \draw[snakecd] (-.5,-.5) to (-.5,1.5);
        \end{tikzpicture}}\endxy
        \;\xrightarrow{
        \xy(0,0)*{\begin{tikzpicture}[scale=.7,xscale=.4,yscale=.7]
          \pic[transform shape,diag1=1.2] at (0,0) {rcup};
          \node[above] at (0,1) {\scriptsize $1$};
          \draw[diag2=1.2,<-] (-1,0) to (-1,1)
          node[above] {\scriptsize $2$};
          \draw[diag2=1.2,->] (2,0) to (2,1);
        \end{tikzpicture}}\endxy
        }\;
        \xy(0,0)*{\begin{tikzpicture}[scale=.4,yscale=-1,transform shape,rotate=90]
          \pic at (0,0) {web-};
          \pic at (1,1) {web+};
          \pic at (1,2) {web-};
          \pic at (0,3) {web+};
          \draw[web1] (1.5,-.5) to (1.5,.5);
          \draw[web1] (-.5,.5) to (-.5,2.5);
          \draw[web1] (1.5,2.5) to (1.5,3.5);
        \end{tikzpicture}}\endxy
        % \right)
      \endgroup
      \\[3ex]
      \mspace{10mu}\simeq\mspace{10mu}
      \\[3ex]
      \begingroup
        % \left(
          \xy(0,0)*{\begin{tikzpicture}
          \node[anchor=east] (A1) at (0,1) {
            $q^2t^{-1}$\;
            \xy(0,0)*{\begin{tikzpicture}[scale=.5,rotate=90]
              \draw[web1] (0,0) to (0,2);
              \draw[web2] (1,0) to (1,2);
              % \draw[snakecd] (-.5,0) to (-.5,2);
            \end{tikzpicture}}\endxy
          };
          \node[anchor=east] (A2) at (0,-1) {
            $t^{-1}$\;
            \xy(0,0)*{\begin{tikzpicture}[scale=.5,rotate=90]
              \draw[web1] (0,0) to (0,2);
              \draw[web2] (1,0) to (1,2);
              % \draw[snakecd] (-.5,0) to (-.5,2);
            \end{tikzpicture}}\endxy
          };
          \node (B) at (2.5,0) {
            \xy(0,0)*{\begin{tikzpicture}[scale=.5,rotate=90]
              \draw[web1] (0,0) to (0,2);
              \draw[web2] (1,0) to (1,2);
            \end{tikzpicture}}\endxy
          };
          \draw[->] (A1) to node[above]{
            \xy(0,1)*{\begin{tikzpicture}[scale=.4]
              \draw[diag1=1.2,
                decoration={markings, mark=at position 0 with {\arrow{<}}},
                postaction={decorate}]
                (0,0) circle (.5cm);
              \draw[diag2=1.2,
                decoration={markings, mark=at position 0 with {\arrow{>}}},
                postaction={decorate}]
                (0,0) circle (1.2cm);
              \node[fdot2=1.5pt,anchor=center] at (0,-.8){};
            \end{tikzpicture}}\endxy
          } (B);
          \draw[->] (A2) to node[below]{
            \xy(0,1)*{\begin{tikzpicture}[scale=.4]
              \draw[diag1=1.2,
                decoration={markings, mark=at position 0 with {\arrow{<}}},
                postaction={decorate}]
                (0,0) circle (.5cm);
              \draw[diag2=1.2,
                decoration={markings, mark=at position 0 with {\arrow{>}}},
                postaction={decorate}]
                (0,0) circle (1.2cm);
            \end{tikzpicture}}\endxy
          } (B);
        \end{tikzpicture}}\endxy
        % \right)
      \endgroup
    \end{gather*}
    \caption{Proof of invariance under Reidemeister I}
    \label{fig:invariance_R1}
  \end{figure}

  Reidemeister II is proved similarly. The complex associated to the non-trivial side of this move is pictured in \cref{fig:invariance_R2}.
  Delooping can be applied to the bottom web, and Gaussian elimination (together with zigzag relations) shows that this complex is homotopy equivalent to the top web. Renormalization concludes.

  \begin{figure}[h]
    \begin{equation*}
      \begin{tikzpicture}
        % LADDERS
        \node (A) at (-.5+.5,0) {
          $q^2t^{-2}$\;
          \xy (0,0)*{\begin{tikzpicture}[scale=.5,transform shape,rotate=90]
            \pic at (0,0) {web+};
            \pic at (0,1) {web-};
          \end{tikzpicture}}\endxy\;
        };
        \node (B) at (3.5,-1) {
          $qt^{-1}$\;
          \xy (0,-1)*{\begin{tikzpicture}[scale=.5,transform shape,rotate=90]
            \pic at (0,0) {web+};
            \pic at (0,1) {web-};
            \pic at (0,2) {web+};
            \pic at (0,3) {web-};
            % \draw[snakecd] (-1,-.5) to (-1,3.5);
          \end{tikzpicture}}\endxy\;
        };
        \node (C) at (5,1) {
          $qt^{-1}$\;
          \xy (0,-.5)*{\begin{tikzpicture}[scale=.5,transform shape,rotate=90]
            \draw[web1] (0,0) to (0,2);
            \draw[web1] (1,0) to (1,2);
            % \draw[snakecd] (-.5,0) to (-.5,2);
          \end{tikzpicture}}\endxy\;
        };
        \node (D) at (9.5-.5,0.1) {
          \xy (0,0)*{\begin{tikzpicture}[scale=.5,transform shape,rotate=90]
          \pic at (0,0) {web+};
          \pic at (0,1) {web-};
          \end{tikzpicture}}\endxy\;
        };
        % STRINGS
        \draw[->] (A) to node[above]{
        \begin{tikzpicture}[scale=.5,transform shape]
          \pic[diag1,line width=1.2] at (0,0) {rcap};
        \end{tikzpicture}
        } (C);
        \draw[->] (A) to node[below=2pt,pos=0]{
        \begin{tikzpicture}[scale=.5,transform shape]
          \pic[diag1,line width=1.2] at (0,0) {lcup};
          \draw[diag1,line width=1.2,->] (1.5,0) to (1.5,1);
          \draw[diag1,line width=1.2,<-] (2,0) to (2,1);
        \end{tikzpicture}
        } (B);
        \draw[->] (B) to node[below=1pt,pos=.6]{
        $-Z^{-1}\;$\xy(0,1)*{\begin{tikzpicture}[scale=.5,transform shape]
        \pic[diag1,line width=1.2] at (2.5,0) {rcap};
        \draw[diag1,line width=1.2,->] (1.5,0) to (1.5,1);
        \draw[diag1,line width=1.2,<-] (2,0) to (2,1);
        \end{tikzpicture}}\endxy
        } (D);
        \draw[->] (C) to node[above]{
        \begin{tikzpicture}[scale=.5,transform shape]
          \pic[diag1,line width=1.2] at (0,0) {lcup};
        \end{tikzpicture}
        } (D);
      \end{tikzpicture}
    \end{equation*}
    \caption{Proof of invariance under Reidemeister II}
    \label{fig:invariance_R2}
  \end{figure}

  Following Bar-Natan \cite{Bar-Natan_KhovanovsHomologyTangles_2005}, we use cones to simplify the proof of invariance under Reidemeister III. Denote $\Gamma(\psi)$ the cone associated to a morphism of complex $\psi$.
  Any hypercube $T$ of dimension $n$ can be seen as a cone. Indeed, choose a direction $k$ in the hypercube. Ignoring the $k$-edges, $T$ breaks in two hypercubes $T_1$ and $T_2$ of dimension $n-1$. Then, switch the signs of all differentials in $T_1$. The hypercube $T_1$ is still a complex, but the $k$-faces of $T$ now commute instead of anti-commuting: the $k$-edges form a morphism $\psi\colon T_1\to T_2$. It is then easy to see that $\Gamma(\psi) = T$.
  We shall use the following lemma:

  \begin{lemma}[{\cite[Lemma 4.5]{Bar-Natan_KhovanovsHomologyTangles_2005}}]
    \label{lem:cone_and_strongdef}
    The cone construction is invariant, up to homotopy, under compositions with strong deformation retracts. That is, if $\psi$ and $F$ are composable morphisms of complexes and $F$ is a deformation retract, then $\Gamma(F\psi)\simeq\Gamma(\psi)$.
    \hfill\qed
  \end{lemma}

  \begin{figure}[h]
    \begin{equation*}
      \begin{tikzpicture}
        % NODES
        \node at (.4+.4,.2) {\tiny 3};
        \node at (.4,-.8+.2) {\tiny 2};
        \node at (-.2,-1+.4) {\tiny 1};
        \begingroup
        \node (000) at (0,0) {\begin{tikzpicture}[thick,rotate=90,transform shape]
        \pic at (0,0) {III_web};
        \pic at (0,.4) {III_web};
        \pic at (0,.8) {III_web};
        \end{tikzpicture}};
        \node (001) at (3,0) {\begin{tikzpicture}[thick,rotate=90,transform shape]
        \pic at (0,0) {III_web};
        \pic at (0,.4) {III_web};
        \pic at (0,.8) {EFI_web};
        \end{tikzpicture}};
        \node (010) at (1,-1) {\begin{tikzpicture}[thick,rotate=90,transform shape]
        \pic at (0,0) {III_web};
        \pic at (0,.4) {IEF_web};
        \pic at (0,.8) {III_web};
        \end{tikzpicture}};
        \node (100) at (0,-3) {\begin{tikzpicture}[thick,rotate=90,transform shape]
        \pic at (0,0) {EFI_web};
        \pic at (0,.4) {III_web};
        \pic at (0,.8) {III_web};
        \end{tikzpicture}};
        \node (011) at (4,-1) {\begin{tikzpicture}[thick,rotate=90,transform shape]
        \pic at (0,0) {III_web};
        \pic at (0,.4) {IEF_web};
        \pic at (0,.8) {EFI_web};
        \end{tikzpicture}};
        \node (101) at (3,-3) {\begin{tikzpicture}[thick,rotate=90,transform shape]
        \pic at (0,0) {EFI_web};
        \pic at (0,.4) {III_web};
        \pic at (0,.8) {EFI_web};
        \end{tikzpicture}};
        \node (110) at (1,-4) {\begin{tikzpicture}[thick,rotate=90,transform shape]
        \pic at (0,0) {EFI_web};
        \pic at (0,.4) {IEF_web};
        \pic at (0,.8) {III_web};
        \end{tikzpicture}};
        \node (111) at (4,-4) {\begin{tikzpicture}[thick,rotate=90,transform shape]
        \pic at (0,0) {EFI_web};
        \pic at (0,.4) {IEF_web};
        \pic at (0,.8) {EFI_web};
        \end{tikzpicture}};
        \node (singleton) at (2,-7) {\begin{tikzpicture}[thick,rotate=90,transform shape]
        \pic at (0,0) {EFI_web};
        \pic at (0,.4) {IEF_web};
        \end{tikzpicture}};
        \endgroup
        % TOP
        \draw[->] (000) to (001);
        \draw[->] (000) to (010);
        \draw[->] (000) to (100);
        \draw[->] (001) to (011);
        \draw[->] (001) to node[left,pos=.6]{\tiny $-X$} (101);
        \draw[->] (100) to (110);
        \draw[->] (100) to (101);
        \node[fill=white] at (3,-1) {};
        \node[fill=white] at (1,-3) {};
        \draw[->] (010) to node[above,pos=.4]{\tiny $-X$} (011);
        \draw[->] (010) to node[left,pos=.4]{\tiny $-X$} (110);
        \draw[->] (101) to (111);
        % BOTTOM
        \draw[->,bend left=10] (101) to (singleton);
        \node[fill=white] at (2.9,-4) {};
        \draw[->] (011) to (111);
        \draw[->] (110) to node[above,pos=.4]{\tiny $-X$} (111);
        \node at (3.5,-5.7) {\xy (0,1)*{\begin{tikzpicture}[xscale=.3,yscale=.4]
        \draw[diag1,thick,<-] (-3.1,0) to (-3.1,-1);
        \draw[diag1,thick,<-] (-2.4,-1) to[out=90,in=180] (-1.9,0) to[out=0,in=90] (-1.4,-1);
        \draw[diag1,thick,->] (-.7,0) to (-.7,-1);
        \draw[diag2,thick,<-] (-4.8,0) to[out=-90,in=180] (-4.3,-1) to[out=0,in=-90] (-3.8,0);
        \end{tikzpicture}}\endxy
        };
        \draw[->] (110) to node[left] {$X\cdot \Id$} (singleton);
        % MORPHISMS
        \draw[<-,very thick] (4.5,-2.5) to node [right] {$\psi_1$} (4.5,-1.5);
        \draw[<-,very thick] (4.5,-2.5-3.5) to node [right] {$F_1$} (4.5,-1.5-3.5);
      \end{tikzpicture}
      \mspace{10mu}
      \begin{tikzpicture}
        % NODES
        \node at (.4+.4,.2) {\tiny 2};
        \node at (.4,-.8+.2) {\tiny 1};
        \node at (-.2,-1+.4) {\tiny 3};
        \begingroup
        \node (000) at (0,0) {\begin{tikzpicture}[thick,rotate=90,transform shape]
        \pic at (0,0) {III_web};
        \pic at (0,.4) {III_web};
        \pic at (0,.8) {III_web};
        \end{tikzpicture}};
        \node (001) at (3,0) {\begin{tikzpicture}[thick,rotate=90,transform shape]
        \pic at (0,0) {III_web};
        \pic at (0,.4) {EFI_web};
        \pic at (0,.8) {III_web};
        \end{tikzpicture}};
        \node (010) at (1,-1) {\begin{tikzpicture}[thick,rotate=90,transform shape]
        \pic at (0,0) {IEF_web};
        \pic at (0,.4) {III_web};
        \pic at (0,.8) {III_web};
        \end{tikzpicture}};
        \node (100) at (0,-3) {\begin{tikzpicture}[thick,rotate=90,transform shape]
        \pic at (0,0) {III_web};
        \pic at (0,.4) {III_web};
        \pic at (0,.8) {IEF_web};
        \end{tikzpicture}};
        \node (011) at (4,-1) {\begin{tikzpicture}[thick,rotate=90,transform shape]
        \pic at (0,0) {IEF_web};
        \pic at (0,.4) {EFI_web};
        \pic at (0,.8) {III_web};
        \end{tikzpicture}};
        \node (101) at (3,-3) {\begin{tikzpicture}[thick,rotate=90,transform shape]
        \pic at (0,0) {III_web};
        \pic at (0,.4) {EFI_web};
        \pic at (0,.8) {IEF_web};
        \end{tikzpicture}};
        \node (110) at (1,-4) {\begin{tikzpicture}[thick,rotate=90,transform shape]
        \pic at (0,0) {IEF_web};
        \pic at (0,.4) {III_web};
        \pic at (0,.8) {IEF_web};
        \end{tikzpicture}};
        \node (111) at (4,-4) {\begin{tikzpicture}[thick,rotate=90,transform shape]
        \pic at (0,0) {IEF_web};
        \pic at (0,.4) {EFI_web};
        \pic at (0,.8) {IEF_web};
        \end{tikzpicture}};
        \node (singleton) at (2,-7) {\begin{tikzpicture}[thick,rotate=90,transform shape]
        \pic at (0,0) {EFI_web};
        \pic at (0,.4) {IEF_web};
        \end{tikzpicture}};
        \endgroup
        % TOP
        \draw[->] (000) to (001);
        \draw[->] (000) to (010);
        \draw[->] (000) to (100);
        \draw[->] (001) to (011);
        \draw[->] (001) to node[left,pos=.6]{\tiny $-X$} (101);
        \draw[->] (100) to (110);
        \draw[->] (100) to (101);
        \node[fill=white] at (3,-1) {};
        \node[fill=white] at (1,-3) {};
        \draw[->] (010) to node[above,pos=.4]{\tiny $-X$} (011);
        \draw[->] (010) to node[left,pos=.4]{\tiny $-X$} (110);
        \draw[->] (101) to (111);
        % BOTTOM
        \draw[->,bend left=10] (101) to (singleton);
        \node[fill=white] at (2.9,-4) {};
        \draw[->] (011) to (111);
        \draw[->] (110) to node[above,pos=.4]{\tiny $-X$} (111);
        \node at (3,-5.7) {$\Id$};
        \draw[->] (110) to node[left] {
        \begin{tikzpicture}[xscale=.3,yscale=.4]
          \draw[diag2,thick,<-] (-3.1,0) to (-3.1,-1);
          \draw[diag2,thick,<-] (-2.4,-1) to[out=90,in=180] (-1.9,0) to[out=0,in=90] (-1.4,-1);
          \draw[diag2,thick,->] (-.7,0) to (-.7,-1);
          \draw[diag1,thick,<-] (0,0) to[out=-90,in=180] (.5,-1) to[out=0,in=-90] (1,0);
        \end{tikzpicture}
        } (singleton);
        % MORPHISMS
        \draw[<-,very thick] (4.5,-2.5) to node [right] {$\psi_2$} (4.5,-1.5);
        \draw[<-,very thick] (4.5,-2.5-3.5) to node [right] {$F_2$} (4.5,-1.5-3.5);
      \end{tikzpicture}
    \end{equation*}
    \caption{Proof of invariance under Reidemeister III}
    \label{fig:invariance_R3}
  \end{figure}

  The proof of invariance is then essentially contained in \cref{fig:invariance_R3}, where we discarded quantum and homological shifts for clarity.
  On the top, one sees the hypercubes associated to each side of the Reidemeister III move, viewed as cones over morphisms respectively denoted $\psi_1$ and $\psi_2$. The ordering 1, 2 and 3 of the three directions corresponds to the ordering of the crossings, reading from right to left. Note that the top faces are identical.
  The bottom of the picture shows how the bottom faces are simplified using delooping and Gaussian elimination, giving homotopy equivalences $F_1$ and $F_2$.

  Thanks to \cref{lem:iso_in_square_complex}, $F_1$ and $F_2$ are strong deformation retracts, so by the above lemma, $\Gamma(F_i\psi_i)\simeq\Gamma(\psi_i)$ for $i=1,2$.
  To conclude, it only remains to compare $F_1\psi_1$ and $F_2\psi_2$. Delving into the proof of \cref{lem:iso_in_square_complex} (or computing by hand), we get explicit $F_1$ and $F_2$; the result is shown in \cref{fig:invariance_R3}. We get that $F_1\psi_1=F_2\psi_2$, and hence $\Gamma(F_1\psi_1)\simeq\Gamma(F_2\psi_2)$.
\end{proof}

\begin{proof}[Proof of \cref{thm:invariance}]
  The renormalization with $q^{-2N_++N_-}$ does not affect the result of \cref{lem:covering_khovanov_invariance_planar_isotopies}: if $D_1$ and $D_2$ are two planar isotopic oriented sliced tangle diagrams, then $q^{-2N_+^1+N_-^1}\glKom(D_1)$ and $q^{-2N_+^2+N_-^2}\glKom(D_2)$ are isomorphic, and in particular homotopic (here we used the notations of \cref{lem:covering_khovanov_invariance_reidemeister_moves}). We conclude with invariance under Reidemeister moves (\cref{lem:covering_khovanov_invariance_reidemeister_moves}).
\end{proof}

%%%%%%%%%%%%%%%%%%%%%%%%%%%%%%%%%%%%%%%%%%%%%%%
%%%          odd khovanov homology          %%%
%%%%%%%%%%%%%%%%%%%%%%%%%%%%%%%%%%%%%%%%%%%%%%%
\subsection{Review of covering \texorpdfstring{$\slt$}{sl2}-Khovanov homology for links}
\label{sec:odd_kh_homology}

We review the construction of covering $\slt$-Khovanov homology as defined by Putyra~\cite{Putyra_2categoryChronologicalCobordisms_2014}.
His construction uses a 2\nbd-cate\-gory of ``chronological cobordisms'', but for our purpose, we give here a ``low-tech'' definition of covering $\slt$-Khovanov homology, directly generalizing the original definition of odd Khovanov homology of Ozsváth, Rasmussen and Szabó \cite{ORS_OddKhovanovHomology_2013}.

\medbreak

We first give some preliminary definitions.
Recall the ring $\ringfoam$ from \cref{defn:ring_R_bil} with distinguished elements $X,Y,Z$ such that $X^2=Y^2=1$.
For $n\in\bN$, let $\wedge_{\ringfoam}(a_1,\ldots,a_n)$ be the $\ringfoam$\nbd-algebra generated by variables $a_1,\ldots,a_n$ and subject to the following relations:
\begin{IEEEeqnarray*}{rCll}
  a_ia_j&=&XYa_ja_i \qquad& \text{for }1\leq i,j\leq n,\\
  a_i^2&=&0 & \text{for }1\leq i\leq n.
\end{IEEEeqnarray*}
Denote by $\wedge_{\ringfoam}^r(a_1,\ldots,a_n)$ the $\ringfoam$\nbd-submodule of $\wedge_{\ringfoam}(a_1,\ldots,a_n)$ generated by words of length $r$ in the letters $a_1,\ldots,a_n$. We endow $\wedge_{\ringfoam}(a_1,\ldots,a_n)$ with a $\bZ$\nbd-grading, the \emph{$q$-grading}, setting $\qdeg p = 2r-n$ whenever $p\in\wedge_{\ringfoam}^r(a_1,\ldots,a_n)$.
Define also the following linear maps:
\begin{align*}
  m_{a_1,a_2;a}\colon \wedge_{\ringfoam}(a_1,a_2,x_1,\ldots,x_n) \to& \wedge_{\ringfoam}(a,x_1,\ldots,x_n)\\
  p\mapsto& p\vert_{a_1,a_2\mapsto a}\\[1ex]
  \Delta_{a;a_1,a_2}\colon \wedge_{\ringfoam}(a,x_1,\ldots,x_n) 
  \to& \wedge_{\ringfoam}(a_1,a_2,x_1,\ldots,x_n)\\
  p\mapsto& (a_1+XY a_2)p\vert_{a\mapsto a_1}\\
   &= (a_1+XY a_2)p\vert_{a\mapsto a_2}
\end{align*}
Here $a_1,a_2\mapsto a$ means that one should replace every instance of $a_1$ and $a_2$ by $a$ in $p$, and similarly for $a\mapsto a_1$ and $a\mapsto a_2$. Note that in the latter case, replacing $a$ by either $a_1$ or $a_2$ gives the same result since $(a_1+XY a_2)a_1=(a_1+XY a_2)a_2$.

With respect to the $q$-grading, these maps are graded maps with $q$-degree
\begin{gather*}
  \qdeg (m_{a_1,a_2;a})=\qdeg (\Delta_{a;a_1,a_2})=1.
\end{gather*}
Note that one recovers the Frobenius algebra $\bZ[a_1,\ldots,a_n]/{(a_1^2=\ldots=a_n^2=0)}$ with its product and coproduct in the even case (choosing $X=Y=Z=1$), and the exterior algebra in variables $a_1,\ldots,a_n$ in the odd case ($X=Z=1$ and $Y=-1$).

Recall \cref{nota:basic_hypercubic_complex}. For $\vr\in\{0,1\}^N$ and $k,l\in\{1,\ldots,N\}$ where $k<l$, the square
\begin{equation}
  \label{eq:square_in_hypercube}
  \begin{tikzcd}
    \vr 
    \drar["\circlearrowright",phantom]
    \rar \dar & \vr+\ve_k \dar
    \\
    \vr+\ve_l \rar & \vr+\ve_k+\ve_l
  \end{tikzcd}
\end{equation}
is given an orientation as depicted, and we denote it by $\square^{\vr}_{k,l}$.

\medbreak

Let then $D$ be a link diagram with $N$ crossings. The hypercubic complex $\slKom(D)$ is constructed through the following steps:
\begin{enumerate}[(i)]
  \item \emph{Hypercube of resolutions:} fix an arbitrary order on the crossings of $D$. Each crossing can be \emph{resolved} into two possible planar diagrams, respectively the \emph{0-resolution} (on the left) or the \emph{1-resolution} (on the right):
  \begin{equation*}
    \def\scl{2}
    \xy(0,0)*{
    \begin{tikzpicture}[scale=.5*\scl,rotate=90]
      \draw[black,very thick] (1,0) to (.65,.35);
      \draw[black,very thick] (.35,.65) to (0,1);
      \draw[black,very thick] (0,0) to (1,1);
    \end{tikzpicture}
    }\endxy
    \mspace{15mu}\mapsto\mspace{15mu}
    \;\xy(0,0)*{
    \begin{tikzpicture}[scale=.5*\scl,rotate=90]
      \draw[black,very thick] (1,0) to[out=135,in=-135] (1,1);      
      \draw[black,very thick] (0,0) to[out=45,in=-45] (0,1);
      \draw[very thick,draw=ocre,-stealth] (0.5,.5) to (.8,.5);
      \draw[very thick,draw=metallic_blue,stealth-] (0.2,.5) to (.5,.5);
    \end{tikzpicture}
    }\endxy
    \mspace{15mu}
    \xrightarrow{
    \xy(0,0)*{\begin{tikzpicture}[scale=-.3,xscale=1.2]
      \def\sh{.3}
      \def\hrf{2}% height of the red foam
      \def\shtwo{.08}
      \draw[draw_foam1] (\sh,3-\sh) to (-\sh,3-\sh) to (-\sh,-\sh)
        .. controls (1,-\sh) and (1.5,-\sh) .. (1.5,0)
        .. controls (1.5,\hrf) and (3,\hrf) .. (3,0)
        .. controls (3,-\sh) and (3.5,-\sh) .. (4.5-\sh,-\sh) to (4.5-\sh,\sh);
      \draw[draw_foam1,dashed] (+\sh,3-\sh) to (4.5-\sh,3-\sh) to (4.5-\sh,+\sh);
      \fill[fill_foam1=\shop+.15] (4.5-\sh,-\sh) to (4.5-\sh,3-\sh) to (-\sh,3-\sh) to (-\sh,-\sh)
        .. controls (1,-\sh) and (1.5,-\sh) .. (1.5,0)
        .. controls (1.5,\hrf) and (3,\hrf) .. (3,0)
        .. controls (3,-\sh) and (3.5,-\sh) .. (4.5-\sh,-\sh);
      % blue below
      \draw[draw_foam1] (4.5-\sh,+\sh) to (4.5+\sh,+\sh) to (4.5+\sh,3+\sh) to (+\sh,3+\sh) to (+\sh,3-\sh);
      \draw[draw_foam1] (+\sh,3-\sh) to (+\sh,+\sh)
        .. controls (1,+\sh) and (1.5,+\sh) .. (1.5,0)
        .. controls (1.5,\hrf) and (3,\hrf) .. (3,0)
        .. controls (3,+\sh) and (3.5,+\sh) .. (4.5+\sh,+\sh);
      \fill[fill_foam1=\shop+.25] (4.5+\sh,+\sh) to (4.5+\sh,3+\sh) to (+\sh,3+\sh) to (+\sh,+\sh)
        .. controls (1,+\sh) and (1.5,+\sh) .. (1.5,0)
        .. controls (1.5,\hrf) and (3,\hrf) .. (3,0)
        .. controls (3,+\sh) and (3.5,+\sh) .. (4.5+\sh,+\sh);
    \end{tikzpicture}}\endxy
    }\mspace{15mu}
    \xy(0,0)*{
      \begin{tikzpicture}[scale=.5*\scl,rotate=-90+90]
        \draw[black,very thick] (1,0) to[out=135,in=-135] (1,1);      
        \draw[black,very thick] (0,0) to[out=45,in=-45] (0,1);
        \draw[very thick,draw=ocre,-stealth] (0.5,.5) to (.8,.5);
        \draw[very thick,draw=metallic_blue,stealth-] (0.2,.5) to (.5,.5);
      \end{tikzpicture}
    }\endxy
  \end{equation*}
  A \emph{resolution} of $D$ is a choice of resolutions for each crossing. The resolutions of $D$ can be pictured as sitting on the vertices of a hypercube $\{0,1\}^N$, where for $\vr\in\{0,1\}^N$ the binary $\vr_i$ encodes the chosen resolution for the $i$-th crossing.
  Each edge $\vr\to\vr+\ve_i$ of the hypercube connects two resolutions that only differ at the $i$-crossing. This edge is decorated with a saddle cobordism, which can either be a merge or a split depending on the global context.
  Finally, for each crossing one must choose an orientation on the arcs of the two resolutions: the red or the blue orientation. We call it the \emph{arc orientation}. Equivalently, an arc orientation is a choice of arc orientation for the 0-resolution, which induces an arc orientation for the 1-resolution by rotating a quarter of a turn clockwise. 
  % We denote $H_\slt(D)$ the hypercube $\{0,1\}^N$ decorated as above.

  \item \emph{Algebraization:} we turn the hypercube of resolutions into a hypercube in the category of $\bZ$-graded $\ringfoam$\nbd-modules. To each vertex $\vr\in\{0,1\}^N$ we associate the $\ringfoam$\nbd-module
  \begin{gather*}
    {V_\vr\coloneqq q^{N-\abs{\vr}}t^{-(N-\abs{\vr})}\;\wedge_{\ringfoam}(a_1,\ldots,a_n)},    
  \end{gather*}
  where $q$ and $t$ denote shifts in quantum and homological degree respectively, and $n$ is the number of connected components in the corresponding resolution. One should think of each variable as attached to one connected component. In addition, each edge is replaced by an $\ringfoam$\nbd-linear map between relevant $\ringfoam$\nbd-modules:
  \begin{equation*}
    \xy(0,0)*{\begin{tikzpicture}[scale=-.7]
      \draw (0,0) to[out=-90,in=-90] (1,0);
      \draw (0,0) to[out=90,in=90] (1,0);
      \draw (0,0) to[out=90,in=-90] (-1,2);
      \draw (1,0) to[out=90,in=-90] (2,2);
      \draw (0,2) to[out=-90,in=180] (.5,1) to[out=0,in=-90] (1,2);
      \draw[dashed] (-1,2) to[out=-90,in=-90] (0,2);
      \draw (0,2) to[out=90,in=90] (-1,2);
      \draw[thick,dashed] (1,2) to[out=-90,in=-90] (2,2);
      \draw (2,2) to[out=90,in=90] (1,2);
      \draw[line width=.5mm,stealth-] (0,2) to (1,2);
      \node at (.5,0) {\small $a$};
      \node at (-.5,2) {\small $a_2$};
      \node at (1.5,2) {\small $a_1$};
    \end{tikzpicture}}\endxy
    \mapsto m_{a_1,a_2;a}
    \qquad\an\qquad
    \xy(0,0)*{\begin{tikzpicture}[scale=.7]
      \begin{scope}
        \draw (0,0) to[out=-90,in=-90] (1,0);
        \draw[dashed] (0,0) to[out=90,in=90] (1,0);
        \draw (0,0) to[out=90,in=-90] (-1,2);
        \draw (1,0) to[out=90,in=-90] (2,2);
        \draw (0,2) to[out=-90,in=180] (.5,1) to[out=0,in=-90] (1,2);
        \draw (-1,2) to[out=-90,in=-90] (0,2) to[out=90,in=90] (-1,2);
        \draw (1,2) to[out=-90,in=-90] (2,2) to[out=90,in=90] (1,2);
      \end{scope}
      \draw[line width=.5mm,-stealth] (0,2) to (1,2);
      \node at (.5,0) {\small $a$};
      \node at (-.5,2) {\small $a_1$};
      \node at (1.5,2) {\small $a_2$};
    \end{tikzpicture}}\endxy
    \mapsto \Delta_{a;a_1,a_2}
  \end{equation*}
  Note the importance of the extra arrows, which give a preferred choice of ordering between the two circles corresponding to the variables $a_1$ and $a_2$.
  We denote $H_\slt(D)$ the resulting hypercube.

  \item \emph{Commutativity:} as defined, squares in the algebraized hypercube do not necessarily commute. In fact, if we consider a generic square
  \begin{equation*}
    \begin{tikzcd}
      \vr 
      \drar["\circlearrowright",phantom]
      \rar["F_{*0}"] \dar["F_{0*}"'] & \vr+\ve_k \dar["F_{1*}"]
      \\
      \vr+\ve_l \rar["F_{*1}"'] & \vr+\ve_k+\ve_l
    \end{tikzcd},
  \end{equation*}
  we have:
  \[F_{1*}\circ F_{*0}=\psi_\slt(\square^{\vr}_{k,l})F_{*1}\circ F_{0*}.\]
  Here $\psi_\slt$ is the $\ringfoam^\times$-valued 2-cochain on the hypercube defined by \cref{tab:sign_assigments}.
  As shown in \cite{ORS_OddKhovanovHomology_2013,Putyra_2categoryChronologicalCobordisms_2014}, $\psi_\slt$ is a cocycle.
  An \emph{$\slt$-scalar assignment} is a choice of a 1-cochain $\epsilon_\slt$ such that $\partial \epsilon_\slt=\psi_\slt$. Such a choice always exists: by contractibility of the hypercube, a 2-cocycle is always a 2-coboundary.
  Given a choice of $\slt$-scalar assignment $\epsilon_\slt$, we multiply each edge $e$ of the hypercube by $\epsilon_\slt(e)$: this makes each square commute.
  We denote $H_\slt(D,\epsilon_\slt)$ the resulting hypercube.

  \begin{table}[t]
    \centering
    \def\scl{.3}
    \def\wcl{.5}

    \begin{tabular}{c@{\hskip 7ex}c}
      \begin{tabular}[t]{c}
        {\LARGE $1$}
        \\[3ex]
        \begin{tikzpicture}[scale=\scl]
          \draw (0,0) circle (1cm);
          \draw (3,0) circle (1cm);
          \draw[line width= \wcl mm] (1,0) to (2,0);
          \begin{scope}[shift={(6,0)}]
            \draw (0,0) circle (1cm);
            \draw (3,0) circle (1cm);
            \draw[line width= \wcl mm] (1,0) to (2,0);
          \end{scope}
        \end{tikzpicture}
        \\[2ex]
        \begin{tikzpicture}[scale=\scl]
          \draw (0,0) circle (1cm);
          \draw (3,0) circle (1cm);
          \draw (6,0) circle (1cm);
          \draw[line width= \wcl mm] (1,0) to (2,0);
          \draw[line width= \wcl mm] (4,0) to (5,0);
        \end{tikzpicture}
        \\[2ex]
        \begin{tikzpicture}[scale=\scl]
          \draw (0,0) circle (1cm);
          \draw (3.5,0) circle (1cm);
          \draw[line width= \wcl mm,-stealth] (0.94,0.342) to (2.5+1-0.94,0.342);
          \draw[line width= \wcl mm,stealth-] (0.94,-0.342) to (2.5+1-0.94,-0.342);
        \end{tikzpicture}
        \\[2ex]
        \begin{tikzpicture}[scale=\scl]
          \draw (0,0) circle (1cm);
          \draw (3,0) circle (1cm);
          \draw[line width= \wcl mm] (1,0) to (2,0);
          \draw (6.5,0) ellipse (1.3cm and 1cm);
          \draw[line width= \wcl mm] (6.5,-1) to (6.5,1);
        \end{tikzpicture}
        \\[2ex]
        \begin{tikzpicture}[scale=\scl]
          \draw (0,0) circle (1cm);
          \draw[line width= \wcl mm] (1,0) to (2,0);
          \draw (3.3,0) ellipse (1.3cm and 1cm);
          \draw[line width= \wcl mm] (3.3,-1) to (3.3,1);
        \end{tikzpicture}
      \end{tabular}
      &
      \begin{tabular}[t]{c}
        {\LARGE $XY$}
        \\[3ex]
        \begin{tikzpicture}[scale=\scl]
          \draw (0,0) ellipse (1.3cm and 1cm);
          \draw[line width= \wcl mm] (0,-1) to (0,1);
          \draw (3.5,0) ellipse (1.3cm and 1cm);
          \draw[line width= \wcl mm] (3.5,-1) to (3.5,1);
        \end{tikzpicture}
        \\[2ex]
        \begin{tikzpicture}[scale=\scl]
          \draw (0,0) ellipse (1.7cm and 1cm);
          \draw[line width= \wcl mm] (-.7,-.9) to (-.7,.9);
          \draw[line width= \wcl mm] (.7,-.9) to (.7,.9);
        \end{tikzpicture}
        \\[2ex]
        \begin{tikzpicture}[scale=\scl]
          \draw (0,0) circle (1cm);
          \draw (3.5,0) circle (1cm);
          \draw[line width= \wcl mm,-stealth] (0.94,0.342) to (2.5+1-0.94,0.342);
          \draw[line width= \wcl mm,-stealth] (0.94,-0.342) to (2.5+1-0.94,-0.342);
        \end{tikzpicture}
        \\[7ex]
        \begin{tabular}[t]{c@{\hskip 8ex}c}
          {\LARGE $1$}
          &
          {\LARGE $XY$}
          \\[3ex]
          \begin{tikzpicture}[scale=\scl]
            \draw (0,0) circle (1cm);
            \draw[line width= \wcl mm,-stealth] (0,-1) to (0,1);
            \draw[line width= \wcl mm,-stealth] (1,0) .. controls (2,0) and (2,1.5) .. (0,1.5) .. controls (-2,1.5) and (-2,0) .. (-1,0);
          \end{tikzpicture}
          &
          \begin{tikzpicture}[scale=\scl]
            \draw (0,0) circle (1cm);
            \draw[line width= \wcl mm,-stealth] (0,-1) to (0,1);
            \draw[line width= \wcl mm,stealth-] (1,0) .. controls (2,0) and (2,1.5) .. (0,1.5) .. controls (-2,1.5) and (-2,0) .. (-1,0);
          \end{tikzpicture}
        \end{tabular}
      \end{tabular}
    \end{tabular}
    
    \caption{Definition of $\psi_\slt$ for covering $\slt$-Khovanov homology. Each square is uniquely represented by the (relevant local piece of) resolution at the initial point. If no orientation on the arrows is given, then the value of $\psi_\slt$ is independent of the choice of orientations.
    The last two cases are called the \emph{ladybugs}.%
    % Note that we do not need to specify the orientation of the squares, as the image of $\psi_\slt$ is isomorphic to $\bZ/2\bZ$.
    }
    \label{tab:sign_assigments}
  \end{table}
  
  \item \emph{Koszul rule:} we apply the Koszul rule (multiplying each edge $\vr\to\vr+\ve_i$ by $(-1)^{\#\{r_j=1\mid j<i\}}$) to turn every commutative square into an anti-commutative square.
  We denote $\slKom(D)$ the resulting hypercubic complex.
\end{enumerate}
% This ends the definition of $\slKom(D)$.\hfill$\diamond$

\begin{theorem}[{\cite{ORS_OddKhovanovHomology_2013,Putyra_2categoryChronologicalCobordisms_2014}}]
\label{thm:invariance_ORS}
  Let $D$ be an oriented link diagram with respectively $N_+$ and $N_-$ positive and negative crossings, presenting an oriented link $L$.
  The isomorphism class of $\slKom(D)$ is independent of the choice of ordering on crossings, the choice of arc orientations, and the choice of $\slt$-scalar assignment. Moreover, the homotopy type of $q^{-2N_++N_-}t^{N_+}\slKom(D)$, denoted $\slKh(L)$, is an invariant of $L$.
\end{theorem}

\begin{remark}
  \label{rem:type_XY_sign_assignment}
  The 2-cocycle $\psi_\slt$ is not the only choice that makes the construction above work. Indeed, for the last two cases of \cref{tab:sign_assigments}, called the \emph{ladybugs},\footnote{This terminology is borrowed from \cite{LS_KhovanovStableHomotopy_2014}.} we both have
  \begin{equation*}
    F_{1*}\circ F_{*0}=\psi_\slt(\square^{\vr}_{k,l})F_{*1}\circ F_{0*}
    \an
    F_{1*}\circ F_{*0}=XY\psi_\slt(\square^{\vr}_{k,l})F_{*1}\circ F_{0*}.
  \end{equation*}
  Define $\overline{\psi}_\slt$ to be the 2-cochain defined as $\psi_\slt$ except for the ladybugs, where instead we set
  \begin{equation*}
      \overline{\psi}_\slt\left(
      \xy(0,0)*{\begin{tikzpicture}[scale=.4]
        \draw (0,0) circle (1cm);
        \draw[line width= .6 mm,-stealth] (0,-1) to (0,1);
        \draw[line width= .6 mm,-stealth] (1,0) .. controls (2,0) and (2,1.5) .. (0,1.5) .. controls (-2,1.5) and (-2,0) .. (-1,0);
      \end{tikzpicture}}\endxy
      \right)=XY
      \qquad\an\qquad
      \overline{\psi}_\slt\left(
      \xy(0,0)*{\begin{tikzpicture}[scale=.4]
        \draw (0,0) circle (1cm);
        \draw[line width= .6 mm,-stealth] (0,-1) to (0,1);
        \draw[line width= .6 mm,stealth-] (1,0) .. controls (2,0) and (2,1.5) .. (0,1.5) .. controls (-2,1.5) and (-2,0) .. (-1,0);
      \end{tikzpicture}}\endxy
      \right)=1.
  \end{equation*}
  The results of \cref{thm:invariance_ORS} still hold in this case. Let us use the notations $\slKh(L,\psi_\slt)$ and $\slKh(L,\overline{\psi}_\slt)$ to distinguish the two constructions. In \cite{ORS_OddKhovanovHomology_2013}, they are respectively called ``type X'' and ``type Y'' (choosing $X=Z=1$ and ${Y=-1}$; no analogy between the scalar $X$ and ``type X'' intended). It is shown in Putyra \cite{Putyra_2categoryChronologicalCobordisms_2014} that $\slKh(L,\psi_\slt)$ and $\slKh(L,\overline{\psi}_\slt)$ are in fact isomorphic. 
\end{remark}

%%%%%%%%%%%%%%%%%%%%%%%%%%%%%%%%%%%%%
%%%          equivalence          %%%
%%%%%%%%%%%%%%%%%%%%%%%%%%%%%%%%%%%%%
\subsection{Covering \texorpdfstring{$\slt$}{sl2}- and \texorpdfstring{$\glt$}{gl2}-Khovanov homologies are isomorphic}
\label{sec:equiv_covering_homologies}

This section is devoted to the proof of \cref{thm:equivalence_with_odd}.
Here is a quick summary:
\begin{enumerate}[(i)]
  \item In \cref{subsec:proof_equiv_hypercube}, we restate the definition of covering $\glt$-Khovanov homology using a \emph{$\glt$-hypercube of resolutions}. The rest of the proof consists in comparing this $\glt$-hypercube with the $\slt$-hypercube defined above.
  \item To compare the hypercubes, we need to compare the $\ringfoam$\nbd-modules at each vertex. This requires a choice of basis for each $\Hom_\gfoam(\emptyset,W)$, called \emph{cup foams}, that we describe in \cref{subsec:proof_equiv_cup_foam}.
  \item In \cref{subsec:proof_equiv_up_to_scalar}, we use the above basis to define a family of isomorphisms on the level of vertices. This defines a proper morphism of hypercubes only \emph{up to invertible scalar}; we call it a \emph{projective morphism}. We state a certain 2-cocycle condition such that, if satisfied, the aforementioned family of isomorphisms can be rescaled into a genuine isomorphism of hypercubes.
  \item The proof of \cref{thm:equivalence_with_odd} then reduces to the analysis of this 2-cocycle condition. \Cref{subsec:proof_equiv_local_analysis} shows that this can be done locally, looking only at the cases pictured in \cref{tab:sign_assigments}. In most cases, general considerations show that the 2-cocycle condition is necessarily verified.
  \item However, these general considerations do not work for the ladybugs (see \cref{rem:type_XY_sign_assignment}). To deal with these two cases, we require finer results on the independence on all choices involved in the above family of isomorphisms. This is done in \cref{subsec:proof_equiv_choices}, which concludes the proof.
\end{enumerate}

% ==============================
\subsubsection{Cup foams}
\label{subsec:proof_equiv_cup_foam}

Call a web $W$ \emph{closed} if $\undcomp(W)$ is a closed 1-manifold. Recall the basis described in \cref{subsec:basis_foam}. In the special case where the domain is the empty web and the codomain is a closed web $W$ (recall that in $\underline{(\gfoam)}^{\oplus}_q$, the empty web is equal to a juxtaposition of double lines), an undotted reduced foam has the following form:
\begin{gather*}
  \undcomp(F) \;=\;
  \xy(0,0)*{\begin{tikzpicture}[scale=1*.8]
    \draw (1,0) to (-.5,1) to (2.5,1) to (1,0);
    \node at (1,.5) {cups};
    \node at (1,-.3) {$\emptyset$};
    \node at (1,1.3) {$\undcomp(W)$};
  \end{tikzpicture}}\endxy
\end{gather*}
We call such an $F$ an \emph{undotted cup foam on $W$}. If we wish to allow $F$ to (possibly) carry dots, we simply say that $F$ is a \emph{cup foam on $W$}.
We write $\pi_0(\undcomp(W))$ for the set of closed components of $\undcomp(W)$.
In this context, \cref{thm:basis_foam} is restated as follows:

\begin{proposition}
  \label{prop:cup_foam_generic_basis}
  Let $W$ be a closed web. Let $B$ be a set containing precisely one cup foam
  \[\beta_\delta\colon\emptyset\to W\]
  for each subset $\delta\subset \pi_0(\undcomp(W))$, so that for each $c\in \pi_0(\undcomp(W))$, the corresponding disk in $\undcomp(\beta_\delta)$ is dotted if and only if $c\in\delta$. Then $B$ is a basis for the $\ringfoam$\nbd-module $\Hom_\gfoam(\emptyset,W)$.
  \hfill\qed
\end{proposition}

Fix an undotted cup foam $\beta^W$ for $W$ and pick a total order on $\pi_0(\undcomp(W))$. For each subset ${\delta\subset \undcomp(W)}$, denote by $\id_W^\delta\colon W\to W$ the foam identical to $\id_W$ except for an additional dot on the connected component $c$ for each $c\in\delta$, ordering the dots increasingly with respect to the total order on $\pi_0(\undcomp(W))$, reading from bottom to top on $\id_W^\delta$. This defines $\id_W^\delta$ uniquely. We denote $\beta^W_\delta\coloneqq\id^\delta_W\circ\beta^W$.
Schematically:
\begin{gather*}
  \undcomp(\beta^W_\delta) \;=\;
  \xy(0,0)*{\begin{tikzpicture}[scale=1*.8]
    \draw (1,0) to (-.5,1) to (2.5,1) to (1,0);
    \draw (-.5,1) rectangle (2.5,2);
    \node at (1,1.5) {dots on $\delta$};
    \node at (1,.5) {cups};
    \node at (1,-.3) {$\emptyset$};
    \node at (1,1.3+1) {$\undcomp(W)$};
  \end{tikzpicture}}\endxy
\end{gather*}
It follows from \cref{prop:cup_foam_generic_basis} that:

\begin{corollary}
  \label{cor:basis_cup_foam}
  For every choice of undotted cup foam $\beta^W$ and total order on $\pi_0(\undcomp(W))$, the family $\{\beta^W_\delta\}_{\delta\subset \undcomp(W)}$ defines a basis for the $\ringfoam$\nbd-module $\Hom_\gfoam(\emptyset,W)$.
  \hfill\qed
\end{corollary}

Let $V$ and $W$ be two closed webs. If $F\colon V\to W$ is a zip or an unzip, then $\undcomp(F)$ is either a merge or a split. Below we sometimes speak of the closed components in the domain and codomain of $\undcomp(F)$ to refer only to the closed components making up the boundary of the closed component in $\undcomp(F)$ containing the saddle. The distinction should be clear by the context.

Recall the symbol $\projrel$ to mean ``equal up to multiplying by an invertible scalar'', defined in the paragraph before \cref{lem:foam_all_isotopy_are_rels}.

\begin{proposition}
  \label{prop:properties_merge_split}
  Let $V$ and $W$ be closed webs and $\beta^V$ (resp.\ $\beta^{W}$) be a choice of undotted cup foam for $V$ (resp.\ $W$). Let $F\colon V\to W$ be a zip or an unzip. Then:
  \begin{enumerate}[(i)]
    \item If $\undcomp(F)$ is a merge, then (as depicted in \eqref{eq:properties_merge})
    \[
      F\circ\beta^V\projrel \beta^{W}.
    \]
    \item If $\undcomp(F)$ is a split, then (as depicted in \eqref{eq:properties_split})
    \[
      F\circ\beta^V\projrel \beta^{W}_{i_1}+XY\beta^{W}_{i_2},
    \]
    where $i_1$ and $i_2$ are the connected components of $W$ corresponding to the codomain of $\undcomp(F)$ (note that the ordering of $i_1$ and $i_2$ is irrelevant in that statement).
  \end{enumerate}
\end{proposition}

Here is the schematic for \cref{prop:properties_merge_split}:
\begin{gather}
  \label{eq:properties_merge}
  {}\xy(0,0)*{\begin{tikzpicture}[scale=.9*.8]
    \draw (1,0) to (-.5,1) to (2.5,1) to (1,0);
    \node at (1,.5) {cups};
    \node at (1,-.3) {$\emptyset$};
    \draw (-.5,1) rectangle (2.5,2);
    \node at (1,1.5) {merge};
    \node at (1,1.3+1) {$W$};
    \node[fill=white] at (1,1) {\small $V$};
  \end{tikzpicture}}\endxy
  \;\projrel\;
  \xy(0,0)*{\begin{tikzpicture}[scale=.9*.8]
    \draw (1,0) to (-.5,1) to (2.5,1) to (1,0);
    \node at (1,.5) {cups};
    \node at (1,-.3) {$\emptyset$};
    \node at (1,1.3) {$W$};
  \end{tikzpicture}}\endxy
  \\[1ex]
  \label{eq:properties_split}
  {}\xy(0,0)*{\begin{tikzpicture}[scale=.9*.8]
    \draw (1,0) to (-.5,1) to (2.5,1) to (1,0);
    \node at (1,.5) {cups};
    \node at (1,-.3) {$\emptyset$};
    \draw (-.5,1) rectangle (2.5,2);
    \node at (1,1.5) {split};
    \node at (1,1.3+1) {$W$};
    \node[fill=white] at (1,1) {\small $V$};
  \end{tikzpicture}}\endxy
  \;\projrel\;
  \xy(0,0)*{\begin{tikzpicture}[scale=.9*.8]
    \draw (1,0) to (-.5,1) to (2.5,1) to (1,0);
    \node at (1,.5) {cups};
    \node at (1,-.3) {$\emptyset$};
    \draw (-.5,1) rectangle (2.5,2);
    \node at (1,1.5) {dot on $i_1$};
    \node at (1,1.3+1) {$W$};
    \node[fill=white] at (1,1) {\small $V$};
  \end{tikzpicture}}\endxy
  \;+XY\;
  \xy(0,0)*{\begin{tikzpicture}[scale=.9*.8]
    \draw (1,0) to (-.5,1) to (2.5,1) to (1,0);
    \node at (1,.5) {cups};
    \node at (1,-.3) {$\emptyset$};
    \draw (-.5,1) rectangle (2.5,2);
    \node at (1,1.5) {dot on $i_2$};
    \node at (1,1.3+1) {$W$};
    \node[fill=white] at (1,1) {\small $V$};
  \end{tikzpicture}}\endxy
\end{gather}

\begin{proof}[Proof of \cref{prop:properties_merge_split}]
  Fix a total order on $\pi_0(\undcomp(V))$ and $\pi_0(\undcomp(W))$.
  Using the neck-cutting successively on the identity of $W$, one can decompose it as
  \[\id_W = \sum_{\delta\subset \pi_0(\undcomp(W))}\beta^W_\delta\circ\beta_{\delta^c}^c,\]
  where $\delta^c\coloneqq \pi_0(\undcomp(W))\setminus \delta$, and each $\beta^{c}_\delta\colon W\to \emptyset$ is a \emph{cap foam} in the sense that $\undcomp(\beta^{c}_\delta)$ is a union of disks, each disk being dotted depending on $\delta$ as in \cref{prop:cup_foam_generic_basis}.
  This allows us to write:
  \begin{equation*}
    F\circ\beta^V=\sum_{\delta\subset \pi_0(\undcomp(W))}\beta^W_\delta\circ\left(\beta_{\delta^c}^c\circ F\circ\beta^V\right),
  \end{equation*}
  where each $\beta_{\delta^c}^c\circ F\circ\beta^V$ is a \emph{closed foam}, that is, a foam with domain and codomain the empty web.
  Then, we apply the following result on the evaluation of closed foams, which is a consequence of \cref{thm:basis_foam}:

  \begin{lemma}
    \label{lem:evaluation_closed_foams}
    Let $U\colon\emptyset\to\emptyset$ be a closed foam in $\gfoam$. Then:
    \begin{center}
      \hfill
      $\xy(0,4)*{U\projrel
      \begin{cases}
        \id_\emptyset & \text{if each closed component of $\undcomp(U)$ is a sphere with a single dot},\\
        0 & \text{otherwise.}
      \end{cases}}\endxy$
      \hfill\qed
    \end{center}
  \end{lemma}

  By the above lemma, there exist invertible scalars $\tau$, $\tau_1$ and $\tau_2$ such that:
  \begin{equation*}
    F\circ\beta^V=\tau \beta^{W}
    \quad\text{or}\quad
    F\circ\beta^V=\tau_1 \beta^{W}_{i_1}+\tau_2\beta^{W}_{i_2},
  \end{equation*}
  depending on whether $\undcomp(F)$ is a merge or a split. It remains to show that $\tau_1/\tau_2=XY$ in the latter case. For that, we use \cref{lem:moving_dots}. Assume $i_1<i_2$ for the purpose of the computation, so that $\id_{i_1,i_2}=\id_{i_2}\circ\id_{i_1}$:
  \begin{IEEEeqnarray*}{rCl}
    \tau_1 \beta^{W}_{i_1,i_2}
    &=& \id_{i_2}^{W}\circ (\tau_1\beta^{W}_{i_1}+\tau_2\beta^{W}_{i_2})
    = \id_{i_2}^{W} \circ F\circ \beta^V \\
    &\overset{\ref{lem:moving_dots}}{=}&  \id_{i_1}^{W}\circ F\circ \beta^V
    =  \id_{i_1}^{W}\circ (\tau_1\beta^{W}_{i_1}+\tau_2\beta^{W}_{i_2})
    = XY\tau_2 \beta^{W}_{i_1,i_2},
  \end{IEEEeqnarray*}
  where in the last equality we used that two dots interchange at the cost of the scalar $XY$.
  Since $\beta^{W}_{i_1,i_2}\neq 0$ belongs to a free family by \cref{prop:cup_foam_generic_basis}, we must have $(\tau_1-XY\tau_2)=0$, which concludes.
\end{proof}

% ==============================
\subsubsection{The $\glt$-hypercube of resolutions}
\label{subsec:proof_equiv_hypercube}

We reformulate the definition of covering $\glt$-Khovanov homology to emphasize the similarities with covering $\slt$-Khovanov homology.
Recall
\[\cA_\glt\coloneqq\Hom_{\gfoam(\emptyset,\emptyset)}(\emptyset,-)\]
defined in the introduction of this section.

Let $D$ be a sliced oriented link diagram with $N$ crossings. As described in \cref{subsec:topo_defn_invariant}, we can associate with $D$ a knotted web $W_D$.
Then, starting with $W_D$, the complex $\cA_\glt(\glKom(D))$ can be defined as follows:
\begin{enumerate}[(i)]
  \item \emph{Hypercube of resolutions:} each crossing can be resolved into a \emph{web 0\nbd-resolution} or a \emph{web 1\nbd-resolution}:
  \begin{IEEEeqnarray*}{c}
    \def\scl{1.2}
    \xy(0,0)*{
    \begin{tikzpicture}[scale=.5*\scl,transform shape]
      \draw[black,very thick] (1,0) to (0,1);
      \draw[black,very thick] (.35,.35) to (0,0);
      \draw[black,very thick] (.65,.65) to (1,1);
    \end{tikzpicture}
    }\endxy
    \mspace{15mu}\mapsto\mspace{15mu}
    qt^{-1}\;\xy(0,0)*{
    \begin{tikzpicture}[scale=.3*\scl,transform shape,rotate=90]
      \draw[web1] (0,0) to (0,2);
      \draw[web1] (1,0) to (1,2);
    \end{tikzpicture}
    }\endxy
    \mspace{5mu}
    \xrightarrow{
    % ZIP
    \xy(0,0)*{\begin{tikzpicture}[scale=-.2,xscale=1.2]
      \def\sh{.3}
      \def\hrf{2}% height of the red foam
      \def\shtwo{.08}
      % red
      \path[draw_foam2,fill_foam2] (3,0-\shtwo) to (1.5,0-\shtwo)
        .. controls (1.5,\hrf) and (3,\hrf) .. (3,0-\shtwo);
      \path[draw_foam2,fill_foam2] (3,0+\shtwo) to (1.5,0+\shtwo)
        .. controls (1.5,\hrf) and (3,\hrf) .. (3,0+\shtwo);
      % \draw[draw_foam2] (1.5,0) to (3,0);
      % blue above
      \draw[draw_foam1] (\sh,3-\sh) to (-\sh,3-\sh) to (-\sh,-\sh)
        .. controls (1,-\sh) and (1.5,-\sh) .. (1.5,0)
        .. controls (1.5,\hrf) and (3,\hrf) .. (3,0)
        .. controls (3,-\sh) and (3.5,-\sh) .. (4.5-\sh,-\sh) to (4.5-\sh,\sh);
      \draw[draw_foam1,dashed] (+\sh,3-\sh) to (4.5-\sh,3-\sh) to (4.5-\sh,+\sh);
      \fill[fill_foam1=\shop+.15] (4.5-\sh,-\sh) to (4.5-\sh,3-\sh) to (-\sh,3-\sh) to (-\sh,-\sh)
        .. controls (1,-\sh) and (1.5,-\sh) .. (1.5,0)
        .. controls (1.5,\hrf) and (3,\hrf) .. (3,0)
        .. controls (3,-\sh) and (3.5,-\sh) .. (4.5-\sh,-\sh);
      % blue below
      \draw[draw_foam1] (4.5-\sh,+\sh) to (4.5+\sh,+\sh) to (4.5+\sh,3+\sh) to (+\sh,3+\sh) to (+\sh,3-\sh);
      \draw[draw_foam1] (+\sh,3-\sh) to (+\sh,+\sh)
        .. controls (1,+\sh) and (1.5,+\sh) .. (1.5,0)
        .. controls (1.5,\hrf) and (3,\hrf) .. (3,0)
        .. controls (3,+\sh) and (3.5,+\sh) .. (4.5+\sh,+\sh);
      \fill[fill_foam1=\shop+.25] (4.5+\sh,+\sh) to (4.5+\sh,3+\sh) to (+\sh,3+\sh) to (+\sh,+\sh)
        .. controls (1,+\sh) and (1.5,+\sh) .. (1.5,0)
        .. controls (1.5,\hrf) and (3,\hrf) .. (3,0)
        .. controls (3,+\sh) and (3.5,+\sh) .. (4.5+\sh,+\sh);
    \end{tikzpicture}}\endxy
    }\mspace{5mu}
    \xy(0,0)*{
    \begin{tikzpicture}[scale=.3*\scl,,transform shape,rotate=90]
      \pic at (0,0) {web+};
      \pic at (0,1) {web-};
    \end{tikzpicture}
    }\endxy
    \quad\an\quad
    \def\scl{1.2}
    \xy(0,0)*{
    \begin{tikzpicture}[scale=.5*\scl,transform shape,rotate=90]
      \draw[black,very thick] (0,0) to (.35,.35);
      \draw[black,very thick] (.65,.65) to (1,1);
      \draw[black,very thick] (1,0) to (0,1);
    \end{tikzpicture}
    }\endxy
    \mspace{15mu}\mapsto\mspace{15mu}
    qt^{-1}\;\xy(0,0)*{
    \begin{tikzpicture}[scale=.3*\scl,,transform shape,rotate=90]
      \pic at (0,0) {web+};
      \pic at (0,1) {web-};
    \end{tikzpicture}
    }\endxy
    \mspace{5mu}
    \xrightarrow{
    % UNZIP
    \xy(0,0)*{\begin{tikzpicture}[scale=.2,xscale=1.2]
      \def\sh{.3}
      \def\shtwo{.08}
      \def\hrf{2}% height of the red foam
      % red
      \path[draw_foam2,fill_foam2] (3,0-\shtwo) to (1.5,0-\shtwo)
      .. controls (1.5,\hrf) and (3,\hrf) .. (3,0-\shtwo);
      \path[draw_foam2,fill_foam2] (3,0+\shtwo) to (1.5,0+\shtwo)
      .. controls (1.5,\hrf) and (3,\hrf) .. (3,0+\shtwo);
      % \draw[draw_foam2] (1.5,0) to (3,0);
      % blue above
      \draw[draw_foam1] (\sh,3-\sh) to (-\sh,3-\sh) to (-\sh,-\sh)
        .. controls (1,-\sh) and (1.5,-\sh) .. (1.5,0);
      \draw[draw_foam1] (1.5,0) .. controls (1.5,\hrf) and (3,\hrf) .. (3,0);
      \draw[draw_foam1] (3,0)  .. controls (3,-\sh) and (3.5,-\sh) .. (4.5-\sh,-\sh) to (4.5-\sh,\sh);
      \draw[draw_foam1] (+\sh,3-\sh) to (4.5-\sh,3-\sh) to (4.5-\sh,+\sh);
      \fill[fill_foam1=\shop+.25] (4.5-\sh,-\sh) to (4.5-\sh,3-\sh) to (-\sh,3-\sh) to (-\sh,-\sh)
        .. controls (1,-\sh) and (1.5,-\sh) .. (1.5,0)
        .. controls (1.5,\hrf) and (3,\hrf) .. (3,0)
        .. controls (3,-\sh) and (3.5,-\sh) .. (4.5-\sh,-\sh);
      % blue below
      \draw[draw_foam1] (4.5-\sh,+\sh) to (4.5+\sh,+\sh) to (4.5+\sh,3+\sh) to (+\sh,3+\sh) to (+\sh,3-\sh);
      \draw[draw_foam1,dashed] (+\sh,3-\sh) to (+\sh,+\sh)
        .. controls (1,+\sh) and (1.5,+\sh) .. (1.5,0)
        .. controls (1.5,\hrf) and (3,\hrf) .. (3,0)
        .. controls (3,+\sh) and (3.5,+\sh) .. (4.5+\sh,+\sh);
      \fill[fill_foam1=\shop+.15] (4.5+\sh,+\sh) to (4.5+\sh,3+\sh) to (+\sh,3+\sh) to (+\sh,+\sh)
        .. controls (1,+\sh) and (1.5,+\sh) .. (1.5,0)
        .. controls (1.5,\hrf) and (3,\hrf) .. (3,0)
        .. controls (3,+\sh) and (3.5,+\sh) .. (4.5+\sh,+\sh);
    \end{tikzpicture}}\endxy
    }\mspace{5mu}
    \xy(0,0)*{
    \begin{tikzpicture}[scale=.3*\scl,,transform shape,rotate=90]
      \draw[web1] (0,0) to (0,2);
      \draw[web1] (1,0) to (1,2);
     \end{tikzpicture}
    }\endxy
  \end{IEEEeqnarray*}
  A \emph{web resolution} is a choice of web resolutions for each crossing. Fixing an ordering on the crossings, they can be pictured as sitting on the vertices of the hypercube $\{0,1\}^N$, whose edges are decorated with a zip or an unzip, depending on whether the associated crossing is positive or negative.

  \item \emph{Algebraization:} we apply the functor $\cA_\glt$ to the hypercube.
  Denote $H_\glt(D)$ the decorated hypercube so obtained.
  
  \item \emph{Commutativity:} a \emph{$\glt$-scalar assignment}\footnote{The notion of scalar assignment here is slightly different from \cref{sec:complexes}, as the latter already includes the Koszul rule.} is an $\ringfoam^\times$-valued 1-cochain $\epsilon_\glt$ on the hypercube $\{0,1\}^N$, such that $\partial\epsilon_\glt = \psi_\glt$ where $\psi_\glt$ is a 2\nbd-cocycle defined as
  \begin{equation*}
    \psi_\glt(\square^\vr_{k,l})\coloneqq \bilfoam(\deg F_{0*},\deg F_{*0})^{-1} = \bilfoam(\deg F_{*1},\deg F_{1*}).
  \end{equation*}
  We multiply each edge $e$ by $\epsilon_\glt(e)$. This makes each square commute. This defines a hypercube $H_\glt(D;\epsilon_\glt)$.

  \item \emph{Koszul rule:} we apply the Koszul rule to turn every commutative square into an anti-commu\-tative square.
\end{enumerate}
As a consequence of \cref{lem:tensor_product_uniqueness}, the isomorphism class of the complex obtained does not depend on the choice of $\glt$-scalar assignment.
It is easily checked that the construction above coincides with the definition of covering $\cA_\glt(\glKom(D))$ given in \cref{subsec:hypercubic_chain_complexes}. In particular, it is shown in \cref{defn:standard_tensor_product} that the graded Koszul rule is a $\glt$-scalar assignment.\footnote{With the caveat of the previous footnote.}

% ==============================
\subsubsection{A projective isomorphism of hypercubes}
\label{subsec:proof_equiv_up_to_scalar}

Let then $D$ be a sliced link diagram with $N$ crossings, with a fixed choice of ordering on crossings and arc orientations.
To show \cref{thm:equivalence_with_odd}, it suffices to exhibit an isomorphism between the two hypercubes $H_\slt(D;\epsilon_\slt)$ and $H_\glt(D;\epsilon_\glt)$, where an \emph{isomorphism of hypercubes} is a family of isomorphisms at each vertex, such that all squares involved commute.
We abuse notation and denote $H_\slt(D)$ (resp.\ $H_\glt(D)$) both the hypercube of resolutions (step (i)) and its algebraization (step (ii)): the relevant hypercube should be clear by the context.
We also use $\slt$ and $\glt$ as subscripts to distinguish features of the two constructions. For instance, we write $\vr_\glt$ to denote the decoration on the vertex $\vr\in\{0,1\}^N$ of the hypercube $H_\glt(D)$ (that is, a web $W$, or the graded $\ringfoam$\nbd-module $\Hom_\gfoam(\emptyset,W)$ depending on the context).

Looking at step (i) in the construction of the hypercubes, it is clear that $\undcomp(H_\glt(D))=H_\slt(D)$: that is, $\undcomp(\vr_\glt)=\vr_\slt$ for each $\vr\in\{0,1\}^N$, and similarly for edges.
For each vertex $W=\vr_\glt$ of $H_\glt(D)$, fix a choice of undotted cup foam $\beta^W$ and a choice of total ordering on the set of connected components $\pi_0(\undcomp(W))$.
By \cref{cor:basis_cup_foam}, $\vr_\glt$ has basis given by the set $\{\beta^W_\delta\}_{\delta\subset \pi_0(\undcomp(W))}$. On the other hand, $\vr_\slt$ has basis given by the set $\{a_\delta\}_{\delta\subset \pi_0(\undcomp(W))}$, where $a_\delta = a_{i_k}\ldots a_{i_1}$ with $\delta=\{i_1,\ldots,i_k\}$ and $i_1<\ldots<i_k$.
Hence, the map
\[\iota_W\colon \beta^W_\delta\mapsto \bilfoam(\abs{\delta}\degd,\beta^W) a_\delta\]
defines an isomorphism of graded $\ringfoam$\nbd-modules.\footnote{One can think of $a_\delta$ as corresponding to $\beta^W\circ\id_\delta^W$: this makes no sense when it comes to composition, but this is coherent with the interchange relation. Thinking this way makes some of the identities below clearer.} Note that $\iota_W$ depends on the choice of $\beta^W$, but not on the choice of ordering on $\pi_0(\undcomp(W))$.
Moreover, by \cref{prop:properties_merge_split} the family of those isomorphisms at each vertex defines a \emph{projective} isomorphism of hypercubes, in the sense that for each edge ${e\colon V\to W}$ in $H_\glt(D)$, the square
\begin{equation}
  \label{eq:square_e}
  \square_e\coloneqq\quad
  \begin{tikzcd}
    V \rar["e"] \arrow[d,"\iota_V"'] 
    \drar["\circlearrowright",phantom]
    & W \arrow[d,"\iota_W"] & H_\glt(D)
    \\
    \undcomp(V) \rar["\undcomp(e)"'] & \undcomp(W) & H_\slt(D)
  \end{tikzcd}
\end{equation}
commutes up to invertible scalar (the symbol $\circlearrowright$ denotes a choice of orientation).
This scalar is computed in the following lemma:

\begin{lemma}
  \label{lem:iota_commute_up_to_scalar}
  Denote by $\tau\in \ringfoam^\times$ the invertible scalar given by \cref{prop:properties_merge_split}. Then either
  \begin{equation*}
    e\circ \beta^V=\tau\beta^{W}
    \quad\text{or}\quad
    e\circ \beta^V=\tau\left(\beta^{W}_{i_1}+XY\beta^{W}_{i_2}\right),
  \end{equation*}
  depending on whether $\undcomp(e)$ is a merge or a split.
  We distinguish $i_1$ from $i_2$ using the arc orientation: it goes from $i_1$ to $i_2$. Define
  \begin{equation*}
    \psi_\beta(e)\coloneqq
    \begin{cases}
      \tau & \text{if $\undcomp(e)$ is a merge},\\
      \tau\bilfoam(\degd,\beta^W) & \text{if $\undcomp(e)$ is a split}.
   \end{cases}
  \end{equation*}
  Then $\iota_W\circ e = \psi_\beta(e) (\undcomp(e)\circ\iota_V)$.
\end{lemma}

Note that $\psi_\beta(e)$ depends in general on the choices of arc orientation on $e$ and undotted cup foams $\beta^V$ and $\beta^W$, but not on the choice of ordering on closed components. We postpone the proof of \cref{lem:iota_commute_up_to_scalar} until after this discussion. 

Ultimately, we are interested in comparing the hypercubes $H_\slt(D;\epsilon_\slt)$ and $H_\glt(D;\epsilon_\glt)$. The above suggests the strategy of finding a 0-cochain $\varphi$ on $\{0,1\}^N$ such that the square
\begin{equation*}
  \begin{tikzcd}
    \undcomp(V) \rar["\epsilon_\glt(e)e"] \arrow[d,"\varphi(V)\iota_V"']
    \drar["\circlearrowright",phantom]
    & W \arrow[d,"\varphi(W)\iota_W"] & H_\glt(D;\epsilon_\glt)
    \\
    \undcomp(V) \rar["\epsilon_\slt(e)\undcomp(e)"'] & \undcomp(W) & H_\slt(D;\epsilon_\slt)
  \end{tikzcd}
\end{equation*}
commutes. That would define an isomorphism of hypercubes, and prove \cref{thm:equivalence_with_odd}.

We can rephrase the problem as follows. Denote by $H_\iota(D)$ the $(N+1)$-dimensional hypercube decorated as $H_\glt(D)$ on $\{0,1\}^N\times\{0\}$, as $H_\slt(D)$ on $\{0,1\}^N\times\{1\}$, and decorated with $\iota$ on edges $\vr\times\{0\}\to\vr\times\{1\}$.
Define also the following 2-cochain $\psi$ on $H_\iota(D)$:
\begin{equation*}
  \psi \coloneqq \begin{cases}
    \psi_\glt & \text{on }\{0,1\}^N\times\{0\},\\
    \psi_\slt & \text{on }\{0,1\}^N\times\{1\},\\
    \psi_\beta(e)^{-1} & \text{on }\square_e.
  \end{cases}
\end{equation*}
where we recall $\square_e$ from \eqref{eq:square_e}.

Assume that $\psi$ is a 2-cocycle. Then by contractibility of the hypercube, it is a coboundary, and there exists some 1-cochain $\epsilon$ such that $\partial\epsilon=\psi$. Denote by $H_\iota(D;\epsilon)$ the hypercube $H_\iota(D)$ obtained by multiplying each edge by its value on $\epsilon$.
By definition, $\epsilon_\glt\coloneqq\epsilon\vert_{\{0,1\}^N\times\{0\}}$ (resp.\ $\epsilon_\slt\coloneqq\epsilon\vert_{\{0,1\}^N\times\{1\}}$) is a $\glt$-scalar assignment (resp.\ a $\slt$-scalar assignment).
In other words, the hypercube $\{0,1\}^N\times\{0\}$ (resp.\ $\{0,1\}^N\times\{1\}$) in $H_\iota(D;\epsilon)$ coincides with $H_\glt(D;\epsilon_\glt)$ (resp.\ $H_\slt(D;\epsilon_\slt)$).
Moreover, by \cref{lem:iota_commute_up_to_scalar} all squares of the kind $\square_e$ commute, so that the $(N+1)$-direction in the hypercube $H_\iota(D;\epsilon)$ defines an isomorphism of hypercubes between $H_\glt(D;\epsilon_\glt)$ and $H_\slt(D;\epsilon_\slt)$.

\begin{figure}
  \begin{equation*}
    \cube_S\coloneqq\quad
    \xy(0,0)*{\begin{tikzpicture}
      % NODES
      \node (000) at (0,0) {\small $W_{00}$};
      \node (001) at (3,0) {\small $W_{01}$};
      \node (010) at (1,-1) {\small $W_{10}$};
      \node (011) at (4,-1) {\small $W_{11}$};
      \node (100) at (0,-3) {\small $\undcomp(W_{00})$};
      \node (101) at (3,-3) {\small $\undcomp(W_{01})$};
      \node (110) at (1,-4) {\small $\undcomp(W_{10})$};
      \node (111) at (4,-4) {\small $\undcomp(W_{11})$};
      \draw[->] (001) to node[left,pos=.6]{\tiny $\iota_{W_{01}}$} (101);
      \node[fill=white] at (3,-1) {};
      \draw[->] (000) to node[above]{\tiny $e_{0*}$} (001);
      \draw[->] (000) to node[right]{\tiny $e_{*0}$} (010);
      \draw[->] (001) to node[right]{\tiny $e_{*1}$} (011);
      \draw[->] (010) to node[above]{\tiny $e_{1*}$} (011);
      \draw[->] (100) to node[below,pos=.55]{\tiny $\undcomp(e_{0*})$} (101);
      \draw[->] (100) to node[left]{\tiny $\undcomp(e_{*0})$} (110);
      \draw[->] (101) to node[left]{\tiny $\undcomp(e_{*1})$} (111);
      \draw[->] (110) to node[below]{\tiny $\undcomp(e_{1*})$} (111);
      \draw[->] (000) to node[left]{\tiny $\iota_{W_{00}}$} (100);
      \node[fill=white] at (1,-3) {};
      \draw[->] (010) to node[right,pos=.4]{\tiny $\iota_{W_{10}}$} (110);
      \draw[->] (011) to node[right]{\tiny $\iota_{W_{11}}$} (111);
      \node at (6,-.5) {$S$};
      \node at (6,-3.5) {$\undcomp(S)$};
    \end{tikzpicture}}\endxy
  \end{equation*}
  \caption{The 3-dimensional cube $\cube_S$, where $S$ is a square of $H_\glt(D)$, pictured on the top.}
  \label{fig:proof_equiv_3D_cube}
\end{figure}

When is $\psi$ a 2-cocycle? Note that it suffices that $\psi$ is a 2-cocycle on every 3-dimensional cube of $H_\iota(D)$. As $\psi_\glt$ and $\psi_\slt$ are already 2-cocycles, it is only necessary that $\psi$ is a 2-cocycle on every 3-dimensional cube $\cube_S\coloneqq S\times\{0,1\}$ for $S$ a square of $H_\glt(D)$ (see \cref{fig:proof_equiv_3D_cube}). We give an orientation on $\cube_S$ such that it agrees with the orientation of $e_{0*}$. 

To sum up, we have shown that:

\begin{proposition}
  \label{prop:cycle_relation_for_odd_equiv}
  Let $D$ be a sliced link diagram with $N$ crossings. Assume given an ordering on crossings, a choice of arc orientations, and a choice of undotted cup foams on the webs decorating $H_\glt(D)$.
  If for every square $S$ of $H_\glt(D)$, the identity
  \[\partial\psi(\cube_S)=1\]
  holds,
  then there exist scalar assignments $\epsilon_\glt$ and $\epsilon_\slt$ such that $H_\glt(D;\epsilon_\glt)$ and $H_\slt(D;\epsilon_\slt)$ are isomorphic.\hfill$\square$
\end{proposition}

We end this subsection with the proof of \cref{lem:iota_commute_up_to_scalar}.

\begin{proof}[Proof of \cref{lem:iota_commute_up_to_scalar}]
  We have two cases: either $\undcomp(e)$ is a merge, or it is a split. In both cases, we compare where the basis element $\beta^W_\delta$ is mapped through the two paths defining the square. The first case is depicted below, where the two end-results are separated by a dashed line:
  \begin{equation*}
    \begin{tikzpicture}[xscale=5.5,yscale=1.5]
      \node (A) at (0,0) {$\beta^V_\delta$};
      \node (B) at (1,0) {$\tau\bilfoam(\deg e,\abs{\delta}\degd)\;\beta^W_\delta$};
      \node (C) at (0,-2) {$\bilfoam(\abs{\delta}\degd,\deg\beta^V)\;a^V_\delta$};
      \node (D) at (1,-1) {$\tau\bilfoam(\deg e,\abs{\delta}\degd)\bilfoam(\abs{\delta}\degd,\deg\beta^W)\;a^W_\delta$};
      \node (E) at (1,-2) {$\bilfoam(\abs{\delta}\degd,\deg\beta^V)\;a^W_\delta$};
      \draw[->] (A) to (B);
      \draw[->] (A) to (C);
      \draw[->] (B) to (D);
      \draw[->] (C) to (E);
      \draw[dashed] (.5,-1.5) to (1.5,-1.5);
    \end{tikzpicture}
  \end{equation*}
  We use symmetry of $\bilfoam$ and the identity $\deg \beta^{W}=\deg \beta^{V}+\deg e$ to conclude.
  Similarly, the computation for the second case gives:
  \begin{equation*}
    \begin{tikzpicture}[xscale=6,yscale=1.5]
      \node (A) at (0,0) {$\beta^V_\delta$};
      \node (B) at (1,0) {$\tau\bilfoam(\deg e,\abs{\delta}\degd)\;\id_W^\delta\circ\left(\beta^{W}_{i_1}+XY\beta^{W}_{i_2}\right)$};
      \node (C) at (0.25,-2) {$\bilfoam(\abs{\delta}\degd,\deg\beta^V)\;a^V_\delta$};
      \node (Cp) at (0,-1.8){};
      \node (D) at (1,-1) {$\tau\bilfoam(\deg e,\abs{\delta}\degd)\bilfoam((\abs{\delta}+1)\degd,\deg\beta^W)\;a^W_\delta \left(a^{W}_{i_1}+XYa^{W}_{i_2}\right)$};
      \node (E) at (1.35,-2) {$\bilfoam(\abs{\delta}\degd,\deg\beta^V)\;\left(a^{W}_{i_1}+XYa^{W}_{i_2}\right) a^W_\delta$};
      \draw[->] (A) to (B);
      \draw[->] (A) to (Cp);
      \draw[->] (B) to (D);
      \draw[->] (C) to (E);
      \draw[dashed] (.5,-1.5) to (1.5,-1.5);
    \end{tikzpicture}
  \end{equation*}
  and the identity $\deg\beta^W+\degd=\deg\beta^V+\deg e$ concludes.
\end{proof}

% ==============================
\subsubsection{Local analysis}
\label{subsec:proof_equiv_local_analysis}

For a sliced link diagram $D$ together with choices as in \cref{prop:cycle_relation_for_odd_equiv}, we need to verify
\begin{gather*}
  \partial\psi(\cube_S)=1
\end{gather*}
for every square $S$ in $H_\glt(D)$.

The following is clear:

\begin{lemma}
  The value of $\partial\psi(\cube_S)$ does not depend on the choice of ordering on crossings.
\end{lemma}

We implicitly assume such a choice in the sequel.

Then, note that the value of $\partial\psi(\cube_S)$ only depends on choices relevant to $\cube_S$; namely, the choice of arc orientations for crossings associated to $S$, and the choice of undotted cup foams for webs associated to $S$.
In other words, whether $S$ is viewed as a square in $H_\glt(D)$ or as a square in $H_\glt(D')$ for some other diagram $D'$ does not affect the value of $\partial\psi(\cube_S)$, provided the same choices relevant to $\cube_S$ are made.
This leads to the following lemma:

\begin{lemma}
  Assume that for every sliced link diagram $S$ with exactly two crossings, and for all choices of arc orientations and undotted cup foams, the identity $\partial\psi(\cube_S)=1$ holds. Then \cref{thm:equivalence_with_odd} holds.\hfill$\square$
\end{lemma}

Recall the pictures of \cref{tab:sign_assigments}: they describe each possible isotopy class, together with the data of arc orientations, associated with such a sliced link diagram $S$.
They are obtained by recording only the 0-resolution for both crossings together with their arc orientation.
More precisely, pictures in \cref{tab:sign_assigments} only picture the non-trivial local part, obtained by removing the simple closed loops that do not contain the boundary of an arc.

Recall the ladybug local arc presentation from \cref{rem:type_XY_sign_assignment}: this was the only case where the value of the $\slt$-2-cocycle $\psi_\slt$ could be set differently. For the other cases, a generic argument is sufficient:

\begin{proposition}
  \label{prop:non-ladybug-cases}
  Assume the Technical Condition for $\ringfoam$ (\cref{defn:technical_condition_R}).
  Let $S$ be a sliced link diagram with exactly two crossings, together with choices of arc orientations and undotted cup foams. Then:
  \begin{enumerate}[(i)]
    \item if the local arc presentation of $S$ is not a ladybug, then $\partial\psi(\cube_S)=1$,
    \item if the local arc presentation of $S$ is a ladybug, then $\partial\psi(\cube_S)=1$ or $\partial\psi(\cube_S)=XY$.
  \end{enumerate}
\end{proposition}

\begin{proof}
  Recall the notations of \cref{fig:proof_equiv_3D_cube}. As $\psi$ is the 2-cochain controlling the commutativity in $\cube_S$, by definition we have that $p=\partial\psi(\cube_S)p$ for ${p=\undcomp(e_{*1})\circ \undcomp(e_{0*})}$ (see \cref{lem:iota_commute_up_to_scalar}). In particular:
  \begin{equation}
    \label{eq:proof_not_ladybug}
  (1-\partial\psi(\cube_S))p(1)=0. 
  \end{equation}
  The element $p(1)$ admits a unique decomposition into basis elements: $p(1) = \sum_{i=1}^n\lambda_i a_{\delta_i}$ for some scalars $\lambda_i\in \ringfoam$ and some subsets $\delta_i\subset \undcomp(W_{11})$. The above relation and the uniqueness of the decomposition imply that $(1-\partial\psi(\cube_S))\lambda_i=0$ for all $i=1,\ldots,n$. Hence, if any of the $\lambda_i$'s is invertible, we automatically get that $\partial\psi(\cube_S)=1$.

  The only case where we get non-invertible coefficients is when $\undcomp(e_{0*})$ is a split and $\undcomp(e_{*1})$ is a merge, that is, in the ladybug cases. In these cases, we have $p(1)=(1+XY)a_j$ for some $j\in \undcomp(W_{11})$. The Technical Condition forces either $\partial\psi(\cube_S)=1$ or $\partial\psi(\cube_S)=XY$.
\end{proof}

If $\partial\psi(\cube_S)=1$ holds for all choices with a ladybug local arc presentation, then \cref{thm:equivalence_with_odd} holds.
If on the contrary $\partial\psi(\cube_S)=XY$ holds for all choices with a ladybug local arc presentation, then \cref{thm:equivalence_with_odd} still holds, as we can apply the same reasoning using instead the $\slt$\nbd-2-cocycle $\ov{\psi}_\slt$ defined in \cref{rem:type_XY_sign_assignment}.
This amounts to comparing our $\glt$\nbd-construction with the $\slt$\nbd-construction of type Y. In either case, we construct an isomorphism between the $\glt$\nbd-hypercube and the $\slt$\nbd-hypercubes of type X or type Y, the latter two being isomorphic.

In other words, what matters is that, whatever the value of $\partial\psi(\cube_S)$, it remains the same for all the choices involved. This is given by the following proposition:

\begin{proposition}
  \label{prop:independence_on_choices}
  Let $S$ be a sliced link diagram with exactly two crossings, together with choices of arc orientations and undotted cup foams. Then the value of $\partial\psi(\cube_S)$ only depends on the local arc presentation of $S$.
\end{proposition}

The proof of \cref{prop:independence_on_choices} is given in \cref{subsec:proof_equiv_choices}.

\begin{remark}
  A direct computation shows that in fact, we have $\partial\psi(\cube_S)=XY$ in the case of a ladybug local arc presentation.
  It is interesting to note that, if we change the defining zigzag relations as follows:
  \begin{equation*}
        %%%%% ZIGZAGS %%%%%
        \begingroup
        \xy (0,1)*{\tikz[scale=.7]{
          \draw[diag1,<-] (0,-.5) to (0,0) to[out=90,in=180] (.25,.5)
            to[out=0,in=90] (.5,0) to[out=-90,in=180] (.75,-.5)
            to[out=0,in=-90] (1,0) to (1,.5) node[left=-3pt]{\scriptsize $i$};
          % \node at (1.2,0) {\tiny $\lambda$};
        }}\endxy
        \,=\,
        \xy (0,1)*{\begin{tikzpicture}[scale=.7]
          \draw [diag1,<-] (0,0) to (0,1) node[left=-3pt]{\scriptsize $i$};
          % \node at (.3,.5) {\tiny $\lambda$};
        \end{tikzpicture}}\endxy
        \mspace{30mu}
        \xy (0,-1)*{\tikz[scale=.7]{
          \draw[diag1,<-] (0,.5) to (0,0) to[out=-90,in=180] (.25,-.5)
            to[out=0,in=-90] (.5,0) to[out=90,in=180] (.75,.5)
            to[out=0,in=90] (1,0) to (1,-.5) node[left=-3pt]{\scriptsize $i$};
          % \node at (1.2,0) {\tiny $\lambda$};
        }}\endxy
        \,=\,X\;
        \xy (0,-1)*{\begin{tikzpicture}[scale=.7]
          \draw [diag1,->] (0,0) node[left=-3pt]{\scriptsize $i$} to (0,1);
          % \node at (.3,.5) {\tiny $\lambda$};
        \end{tikzpicture}}\endxy
        \mspace{30mu}
        \xy (0,-1)*{\tikz[scale=.7]{
          \draw[diag1,->] (0,-.5) node[left=-3pt]{\scriptsize $i$} to (0,0) to[out=90,in=180] (.25,.5)
            to[out=0,in=90] (.5,0) to[out=-90,in=180] (.75,-.5)
            to[out=0,in=-90] (1,0) to (1,.5);
          % \node at (1.2,0) {\tiny $\lambda$};
        }}\endxy
        \,=\,
        XYZ^{2}\;
        \xy (0,-1)*{\begin{tikzpicture}[scale=.7]
          \draw [diag1,->] (0,0) node[left=-3pt]{\scriptsize $i$} to (0,1);
          % \node at (.3,.5) {\tiny $\lambda$};
        \end{tikzpicture}}\endxy
        \mspace{30mu}
        \xy (0,1)*{\tikz[scale=.7]{
          \draw[diag1,->] (0,.5) node[left=-3pt]{\scriptsize $i$} to (0,0) to[out=-90,in=180] (.25,-.5)
            to[out=0,in=-90] (.5,0) to[out=90,in=180] (.75,.5)
            to[out=0,in=90] (1,0) to (1,-.5);
          % \node at (1.2,0) {\tiny $\lambda$};
        }}\endxy
        \,=\,
        XZ^{2}\;
        \xy (0,1)*{\begin{tikzpicture}[scale=.7]
          \draw [diag1,<-] (0,0) to (0,1) node[left=-3pt]{\scriptsize $i$};
          % \node at (.3,.5) {\tiny $\lambda$};
        \end{tikzpicture}}\endxy
      \endgroup
  \end{equation*}
  leaving the rest of the definition identical, we get $\partial\psi(\cube_S)=1$ in all cases. Following \cite[Section~4.5]{Schelstraete_RewritingModuloDiagrammatic_2026} (see also \cite[Section~7.4]{Schelstraete_OddKhovanovHomology_2024}), this variant satisfies the same basis theorem as $\gfoam_d$, and the same proof goes through. It is not clear to the authors whether this other graded-2-category of graded $\glt$-foams is isomorphic to $\gfoam_d$.
\end{remark}

% ==============================
\subsubsection{Independence on choices}
\label{subsec:proof_equiv_choices}

We conclude the proof of \cref{thm:equivalence_with_odd} by proving \cref{prop:independence_on_choices}. With the notation of \cref{fig:proof_equiv_3D_cube}, recall that $S$ is oriented as follows:
\begin{equation*}
  S=\quad
  \begin{tikzcd}
    W_{00} \rar["e_{*0}"] \arrow[d,"e_{0*}"'] 
    \drar["\circlearrowright",phantom]
    & W_{10} \arrow[d,"e_{1*}"]
    \\
    W_{01} \rar["e_{*1}"'] & W_{11}
  \end{tikzcd}
\end{equation*}
and similarly for $\undcomp(S)$. We compute that:
\begin{equation*}
  \partial\psi(\cube_S)=\psi_\glt(S)\psi_\slt(\undcomp(S))^{-1}\partial\psi_\beta(S).
\end{equation*}

\begin{lemma}
  \label{lem:independence_on_arc_cup}
  Let $S$ be a sliced link diagram with exactly two crossings, together with choices of arc orientations and undotted cup foams. Then the value of $\partial\psi(\cube_S)$ does not depend on the choices of arc orientations and undotted cup foams.
\end{lemma}

\begin{proof}
  Assume we swap the arc orientation of the $k$th crossing. Then $\partial\psi_\beta(S)$ will contribute an additional factor $XY$ for every split in direction $k$. Looking case by case at \cref{tab:sign_assigments}, one checks that this change is exactly compensated by the contribution of $\psi_\slt(\undcomp(S))$. The value of $\psi_\glt(S)$ does not change. Hence, changing the arc orientations does not change the value of $\partial\psi(\cube_S)$.

  Assume we change the choice of undotted cup foams instead. Let $\beta$ and $\wt{\beta}$ be two choices of undotted cup foams, identical everywhere except at some vertex $W$. Denote by $\tau$ the invertible scalar such that $\wt{\beta}^W=\tau\beta^W$.\footnote{The existence of such a scalar can be deduced from \cref{thm:basis_foam}.} Then:
  \begin{equation*}
    \psi_\wt{\beta}(\to W) = \psi_\beta(\to W)/\tau
    \quad\an\quad
    \psi_\wt{\beta}(W\to) = \psi_\beta(W\to)\cdot\tau,
  \end{equation*}
  where $\to W$ (resp.\ $W\to$) denotes an incoming (resp.\ outgoing) edge in $H_\glt(S)$. In particular, $\partial\psi_\beta(S) = \partial\psi_\wt{\beta}(S)$. One also checks that the values of $\psi_\glt(S)$ and $\psi_\slt(\undcomp(S))$ do not change. Hence, changing the choice of undotted cup foams does not change the value of $\partial\psi(\cube_S)$.
\end{proof}

Next we check \cref{prop:independence_on_choices} for planar isotopies. Actually, we show a bit more, namely that the value of $\partial\psi(\cube_S)$ is independent of the choice of spatial representative for $S$. A \emph{spatial sliced representative} of a link diagram $S$ is a sliced diagram that is a representative of $S$ up to spatial isotopies (see \cref{defn:spatial_isotopy}).

\begin{lemma}
  \label{lem:independence_on_spatial_iso}
  Let $S$ be a sliced link diagram with exactly two crossings, together with choices of arc orientations and undotted cup foams. Then the value of $\partial\psi(\cube_S)$ does not depend on the spatial sliced representative.
\end{lemma}

\begin{proof}
  Let $S$ and $\wt{S}$ be two sliced link diagrams in the same spatial isotopy class. Throughout the proof we use the notation $\wt{(-)}$ to distinguish features related to $\wt{S}$. By \cref{lem:independence_on_arc_cup}, we may freely choose arc orientations and undotted cup foams.
  
  Pick a choice of arc orientations on $S$. Then there are arc orientations on $\wt{S}$ naturally associated with the one on $S$: for instance, one can use an orientation on $S$ to record arc orientations, and orientation of link diagrams is preserved by spatial isotopies. With this choice, we have $\wt{\psi}_\slt=\psi_\slt$.

  Recall that if $S$ and $\wt{S}$ are planar isotopic, then by \cref{lem:covering_khovanov_invariance_planar_isotopies} there exist scalar assignments $\epsilon$ and $\wt{\epsilon}$ such that $H_\glt(S;\epsilon)$ and $H_\glt(\wt{S};\wt{\epsilon})$ are isomorphic. This extends to the spatial case, thanks to the following pair of isomorphisms in $\gfoam$:
  \begin{equation*}
    \begin{tikzpicture}
      \node[anchor=east] (A) at (0,0) {
        \begin{tikzpicture}[scale=.5,xscale=.8]
          \draw[diag2,<-] (0,0) to (0,2);
          \draw[diag2,->] (1,0) to (1,2);
        \end{tikzpicture}
      };
      \node[anchor=west] (B) at (7,0) {
        \begin{tikzpicture}[scale=.5,xscale=.8]
          \draw[diag2,<-] (0,0) to (0,2);
          \draw[diag1,->] (1,0) to (1,2);
          \draw[diag1,<-] (2,0) to (2,2);
          \draw[diag1,->] (3,0) to (3,2);
          \draw[diag1,<-] (4,0) to (4,2);
          \draw[diag2,->] (5,0) to (5,2);
        \end{tikzpicture}
      };
      \draw[->] (.5,.15) to (6.5,.15);
      \draw[<-] (.5,-.15) to (6.5,-.15);
      \node[anchor=north] at (3.5,-.1) {
        $\xy(0,0)*{\begin{tikzpicture}[scale=.25]
          \draw[diag1=1.5,<-] (-1,0) to[out=90,in=180] (0,1) to[out=0,in=90] (1,0);
          \draw[diag1=1.5,->] (-2,0) to[out=90,in=180] (0,2) to[out=0,in=90] (2,0);
          \draw[diag2=1.5,<-] (-3,0) to[out=90,in=180] (0,3) to[out=0,in=90] (3,0);
          \draw[diag2=1.5,->] (-1,4.5) to[out=-90,in=180] (0,3.5) to[out=0,in=-90] (1,4.5);
          \node[fdot2=1.5pt,anchor=north] at (0,4.5){};
        \end{tikzpicture}}\endxy
        +
        \xy(0,0)*{\begin{tikzpicture}[scale=.25]
          \draw[diag1=1.5,<-] (-1,0) to[out=90,in=180] (0,1) to[out=0,in=90] (1,0);
          \draw[diag1=1.5,->] (-2,0) to[out=90,in=180] (0,2) to[out=0,in=90] (2,0);
          \draw[diag2=1.5,<-] (-3,0) to[out=90,in=180] (0,3) to[out=0,in=90] (3,0);
          \draw[diag2=1.5,->] (-1,4.5) to[out=-90,in=180] (0,3.5) to[out=0,in=-90] (1,4.5);
          \node[fdot1=1.5pt,anchor=south] at (0,0){};
        \end{tikzpicture}}\endxy$
      };
      \node[anchor=south] at (3.5,.1) {
        $XY\left(\;\xy(0,0)*{\begin{tikzpicture}[scale=.25,rotate=180]
          \draw[diag1=1.5,<-] (-1,0) to[out=90,in=180] (0,1) to[out=0,in=90] (1,0);
          \draw[diag1=1.5,->] (-2,0) to[out=90,in=180] (0,2) to[out=0,in=90] (2,0);
          \draw[diag2=1.5,<-] (-3,0) to[out=90,in=180] (0,3) to[out=0,in=90] (3,0);
          \draw[diag2=1.5,->] (-1,4.5) to[out=-90,in=180] (0,3.5) to[out=0,in=-90] (1,4.5);
          \node[fdot2=1.5pt,anchor=south] at (0,4.5){};
        \end{tikzpicture}}\endxy
        +
        \xy(0,0)*{\begin{tikzpicture}[scale=.25,rotate=180]
          \draw[diag1=1.5,<-] (-1,0) to[out=90,in=180] (0,1) to[out=0,in=90] (1,0);
          \draw[diag1=1.5,->] (-2,0) to[out=90,in=180] (0,2) to[out=0,in=90] (2,0);
          \draw[diag2=1.5,<-] (-3,0) to[out=90,in=180] (0,3) to[out=0,in=90] (3,0);
          \draw[diag2=1.5,->] (-1,4.5) to[out=-90,in=180] (0,3.5) to[out=0,in=-90] (1,4.5);
          \node[fdot1=1.5pt,anchor=north] at (0,0){};
        \end{tikzpicture}}\endxy\;\right)$
      };
    \end{tikzpicture}
  \end{equation*}
  Denote $\varphi\colon H_\glt(S;\epsilon)\to H_\glt(\wt{S};\wt{\epsilon})$ such an isomorphism. If $\beta$ is a choice of undotted cup foams for $S$, then $\varphi\circ \beta$ is a choice of undotted cup foams for $\wt{S}$. Then, for each edge $e\colon\vr\to\vs$ in $\{0,1\}^2$, denoting $F_e$ and $\wt{F}_e$ the corresponding foams in $H_\glt(S)$ and $H_\glt(\wt{S})$ respectively, we have that:
  \begin{align*}
    \wt{\epsilon}(e) \left(\wt{F}_e\circ\wt{\beta}^\vr \right)
    &= \epsilon(e)\left(\varphi_\vs\circ F_e\circ \beta^\vr\right)\\
    &= \epsilon(e)\psi_\beta(e)
    \begin{cases}
      \beta^\vs & \text{ if $F_e$, $\wt{F}_e$ are merges},\\
      (\beta_{i_1}^\vs+XY\beta_{i_2}^\vs) & \text{ if $F_e$, $\wt{F}_e$ are splits}.
    \end{cases}
  \end{align*}
  Hence, we have $\psi_\wt{\beta}\wt{\epsilon}=\psi_\beta\epsilon$. That implies $\partial(\psi_\wt{\beta})\wt{\psi}_\glt=\partial(\psi_\beta)\psi_\glt$. This concludes the proof.
\end{proof}

Finally, we need to check that $\partial\psi(\cube_S)$ only depends on the local arc presentation. Up to spatial isotopies, we can slide away all closed simple loops. The next lemma concludes the proof of \cref{prop:independence_on_choices}.

\begin{lemma}
  Let $D_0$ and $D_1$ be two sliced link diagrams, such that $D_1$ has two crossings and $D_0$ none. Then $\partial\psi(\cube_{D_1})=\partial\psi(\cube_{D_0\star_0 D_1})$.
\end{lemma}

\begin{proof}
  If $W_0$ is the web corresponding to $D_0$ and $\beta^{W_0}$ is a choice of undotted cup foam for $W_0$, then $(\id_{W_0}\star_0\beta^W)\star_1\beta^{W_0}$ is an undotted cup foam for every undotted cup foam $\beta^W$. This defines a choice of undotted cup foams on $D_0\star_0 D_1$, given one on $D_1$. It is also clear how to define arc orientations on $D_0\star_0 D_1$ given such a choice on $D_1$.
  With these choices, $\psi_\glt$, $\psi_\slt$ and $\psi_\beta$ remain identical, and \cref{lem:independence_on_arc_cup} concludes.
\end{proof}

\section{A graded-categorification of the \texorpdfstring{$q$}{q}-Schur algebra of level 2}
\label{sec:schur}
This section introduces a diagrammatic graded-2\nbd-cate\-gory that categorifies the $q$-Schur algebra of level 2. Then, we define a graded foamation 2-functor that relates this construction to graded $\glt$-foams.
This can be seen as a partial super analogue of \cite{LQR_KhovanovHomologySkew_2015} in the $\glt$ case. See also \cite{Mackaay_Sl3FoamsKhovanovLaudaCategorification_2009} for earlier work.

\medbreak

Diagrammatic categorification of quantum groups was independently introduced by Khovanov--Lauda \cite{KL_DiagrammaticApproachCategorification_2009} and Rouquier \cite{Rouquier_2KacMoodyAlgebras_2008}. A super analogue of this construction was given by Brundan and Ellis \cite{BE_SuperKacMoody_2017}, building on earlier work of Kang, Kashiwara and Tsuchioka \cite{KKT_QuiverHeckeSuperalgebras_2016}. For odd $\mathfrak{sl}_2$, this was already studied in \cite{EL_OddCategorification$U_qmathfraksl_2$_2016,EKL_OddNilHeckeAlgebra_2014}.
See also \cite{KKO_SupercategorificationQuantumKac_2014,KKO_SupercategorificationQuantumKacMoody_2013,HW_CategorificationQuantumKacMoody_2015} for related work.

In the even case, a categorification of the $q$-Schur algebra appeared in \cite{MSV_DiagrammaticCategorification$q$Schur_2013}. In \cite{Vaz_NotEvenKhovanov_2020}, the second author defined a supercategorification of the negative half of the $q$-Schur algebra of level 2. A graded version of this construction was given in \cite{NP_OddKhovanovHomology_2020}. In the same paper, Naisse and Putyra also defined a ``1-map'' from this graded version to their construction. This 1-map has similarities with our graded foamation 2-functor: we expect the two to coincide once an equivalence between \cite{NP_OddKhovanovHomology_2020} and our category of graded $\glt$-foams is found.

Our presentation of the categorification of the $q$-Schur algebra of level 2 is similar to the presentation of super Kac--Moody 2-algebras in \cite{BE_SuperKacMoody_2017}. However, and contrary to the even case, it is not obtained as a quotient of their construction. See also \cref{rem:intro_relation_to_supercategorified_QG}.

\medbreak

\Cref{subsec:qschur} reviews the definition of $\qschur_{n,d}$, the idempotented $q$-Schur algebra of level 2. Its graded-categorification, that we call the \emph{graded 2-Schur algebra} $\catschur_{n,d}$, is introduced in \cref{subsec:graded_2_schur}. \Cref{sec:foamation_functor} then defines the \emph{graded foamation 2-functor} from $\catschur_{n,d}$ into $\gfoam_d$, our graded-2\nbd-cate\-gory of $\glt$-foams defined in \cref{subsec:graded_foams}. Finally, we show in \cref{subsec:schur_categorification} that $\catschur_{n,d}$ categorifies $\qschur_{n,d}$.

%%%%%%%%%%%%%%%%%%%%%%%%%%%%%%%%%
%%%          q-Schur          %%%
%%%%%%%%%%%%%%%%%%%%%%%%%%%%%%%%%
\subsection{The \texorpdfstring{$q$}{q}-Schur algebra of level 2}
\label{subsec:qschur}

\newcommand{\schursuper}[2]{#1}

\begin{definition}
  The \emph{(idempotented) $q$-Schur algebra of level 2}, or simply \emph{Schur algebra}, is the $\bZ[q,q^{-1}]$-linear category $\qschur_{n,d}$ such that:
  \begin{itemize}
    \item Objects are weights in the set
    \begin{equation*}
      \Lnd\coloneqq\{\lambda\in\{0,1,2\}^n\mid \lambda_1+\ldots+\lambda_n=d\}.
    \end{equation*}
    \item Morphisms are $\bZ[q,q^{-1}]$-linear combinations of iterated compositions of identity morphisms $1_\lambda\colon\lambda\to\lambda$ and generating morphisms
    \begin{equation*}
      e_i1_\lambda\colon \lambda\to\lambda+\alpha_i
      \quad\an\quad
      f_i1_\lambda\colon \lambda\to\lambda-\alpha_i\qquad i=1,\ldots,n-1,
    \end{equation*}
    where $\alpha_i\coloneqq(0,\ldots,1,-1,\ldots,0)\in\bZ^n$ with 1 being on the $i$-th coordinate. Morphisms are subject to the \emph{Schur quotient}
    \begin{gather*}
      1_\lambda = 0 \qquad\text{ if }\lambda\not\in\Lnd
    \end{gather*}
    and to the following relations:
    \begin{gather}
      \label{eq:quantum_algebra_rel}
      \begin{cases}
        (e_if_j-f_je_i)1_\lambda = 
        \schursuper{
        \delta_{ij}[\overline{\lambda}_i]_q1_\lambda 
        }{
        \delta_{ij}\pi^{\lambda_{i+1}}[\overline{\lambda}_i]_{q,\pi}1_\lambda 
        }
        &\\
        (e_ie_j-e_je_i)1_\lambda = 0 & \text{ for }\abs{i-j}>1\\
        (f_if_j-f_jf_i)1_\lambda = 0 & \text{ for }\abs{i-j}>1
      \end{cases}
    \end{gather}
    where $\overline{\lambda}_i\coloneqq\lambda_i-\lambda_{i+1}$, the symbol $\delta_{ij}$ is the Kronecker delta and
    \[
    \schursuper{
      [m]_q\coloneqq \frac{q^m-q^{-m}}{q-q^{-1}} =
      \begin{cases}
        q^{m-1}+q^{m-3}+\ldots+q^{3-m}+q^{1-m}&\text{if }m\geq 0,\\
        -(q^{m-1}+q^{m-3}+\ldots+q^{3-m}+q^{1-m})&\text{if }m< 0,
      \end{cases}
    }{
      [m]_{\pi,q}\coloneqq \frac{q^m-(\pi q)^{-m}}{q-(\pi q)^{-1}}
    }
    \]
    is the $m$th quantum integer.
  \end{itemize}
\end{definition}

Recall that a $\bZ[q,q^{-1}]$-linear category is the same as a $\bZ[q,q^{-1}]$-algebra with a distinguished set of idempotents, so that $\qschur_{n,d}$ is indeed an algebra.
The ``level 2'' stands for the fact that the value of coordinates is at most two.
The Schur quotient implies that a morphism that factors through a weight not in $\Lambda_{n,d}$ is set to zero.
In the sequel, it is understood that an expression involving a weight that does not belong to $\Lnd$ is set to zero.

\begin{remark}
  \label{rem:schur_relation_with_RP}
  The $q$-Schur algebra of level 2 is an integral form for the fundamental representations of $U_q(\glt)$, analogous to the role of the Temperley--Lieb algebra for $U_q(\slt)$.
  Let $\bigwedge^k\coloneqq\bigwedge^k(\bC(q)^2)$ for $k=0,1,2$ denote the fundamental representations of $U_q(\glt)$, with $\bC(q)^2$ the standard representation. Let $\mathrm{Fund}_{n,d}(U_q(\glt))$ denote the full subcategory of $U_q(\glt)$-representations consisting of $n$-fold tensor products $\bigwedge^{k_1}\otimes\ldots\otimes\bigwedge^{k_n}$ with $k_1+\ldots+ k_n=d$.
  Then:
  \begin{equation*}
    \qschur_{n,d}\otimes\bC(q)\cong\mathrm{Fund}_{n,d}(U_q(\glt)),
  \end{equation*}
  where the isomorphism is an isomorphism of $\bC(q)$-linear categories (see \cite{CKM_WebsQuantumSkew_2014}).
\end{remark}

\subsubsection{Relationship with \texorpdfstring{$\glt$}{gl2}-webs}
\label{subsec:schur_ladder_to_web}

Each generating 1-morphism of $\qschur_{n,d}$ can be represented as a \emph{ladder diagram}:
\begin{equation*}
  e_i1_{(\lambda_i,\lambda_{i+1})} \;\mapsto\;
  \xy (0,0)*{\begin{tikzpicture}[scale=.8]
    \draw[line width=1.5pt,directed=.25,directed=.75] (2,0) to (0,0);
    \draw[line width=1.5pt,directed=.35,directed=.85] (2,1) to (0,1);
    \draw[line width=1.5pt,directed] (1.1,0) to (.9,1);
    \node[right] at (2.2,0) {$\lambda_i$};
    \node[right] at (2.2,1) {$\lambda_{i+1}$};
    \node[left] at (-.2,0) {$\lambda_i+1$};
    \node[left] at (-.2,1) {$\lambda_{i+1}-1$};
  \end{tikzpicture}}\endxy
  \qquad
  f_i1_{(\lambda_i,\lambda_{i+1})} \;\mapsto\;
  \xy (0,0)*{\begin{tikzpicture}[scale=.8]
    \draw[line width=1.5pt,directed=.25,directed=.75] (2,1) to (0,1);
    \draw[line width=1.5pt,directed=.35,directed=.85] (2,0) to (0,0);
    \draw[line width=1.5pt,rdirected] (.9,0) to (1.1,1);
    \node[right] at (2.2,0) {$\lambda_i$};
    \node[right] at (2.2,1) {$\lambda_{i+1}$};
    \node[left] at (-.2,0) {$\lambda_i-1$};
    \node[left] at (-.2,1) {$\lambda_{i+1}+1$};
  \end{tikzpicture}}\endxy
\end{equation*}
Representing coordinates 0, 1 and 2 respectively with a dotted line $\xy (0,0)*{\begin{tikzpicture}[thick,rotate=90]
  \draw[web0] (0,0) to (0,.5);
\end{tikzpicture}}\endxy$ , 
a single line $\xy (0,0)*{\begin{tikzpicture}[thick,rotate=90]
  \draw[web1] (0,0) to (0,.5);
\end{tikzpicture}}\endxy$
and a double line $\xy (0,0)*{\begin{tikzpicture}[,rotate=90]
  \draw[web2] (0,0) to (0,.5);
\end{tikzpicture}}\endxy$, ladder diagrams take the following local form:
\begin{gather*}
  e1_{(1,1)} \;\mapsto\;
  \xy (0,0)*{\begin{tikzpicture}[thick,xscale=.5,yscale=.5,rotate=90]
  \draw[web0] (1,1) to (1,2);
  \draw[web2] (0,1) to (0,2);
  \draw[web1] (0,0) to (0,1) to (1,1) to (1,0);
  \end{tikzpicture}}\endxy
  \mspace{30mu}
  e1_{(0,1)} \;\mapsto\;
  \xy (0,0)*{\begin{tikzpicture}[thick,xscale=.5,yscale=.5,rotate=90]
  \draw[web0] (0,0) to (0,1);\draw[web0] (1,1) to (1,2);
  \draw[web1] (1,0) to (1,1) to (0,1) to (0,2);
  \end{tikzpicture}}\endxy
  \mspace{30mu}
  e1_{(1,2)} \;\mapsto\;
  \xy (0,0)*{\begin{tikzpicture}[thick,xscale=.5,yscale=.5,rotate=90]
  \draw[web2] (1,0) to (1,1);\draw[web2] (0,1) to (0,2);
  \draw[web1] (0,0) to (0,1) to (1,1) to (1,2);
  \end{tikzpicture}}\endxy
  \mspace{30mu}
  e1_{(0,2)} \;\mapsto\;
  \xy (0,0)*{\begin{tikzpicture}[thick,xscale=.5,yscale=.5,rotate=90]
  \draw[web0] (0,0) to (0,1);
  \draw[web2] (1,0) to (1,1);
  \draw[web1] (0,2) to (0,1) to (1,1) to (1,2);
  \end{tikzpicture}}\endxy
  \\[1ex]
  f1_{(1,1)} \;\mapsto\;
  \xy (0,0)*{\begin{tikzpicture}[thick,xscale=-.5,yscale=.5,rotate=90]
  \draw[web0] (1,1) to (1,2);
  \draw[web2] (0,1) to (0,2);
  \draw[web1] (0,0) to (0,1) to (1,1) to (1,0);
  \end{tikzpicture}}\endxy
  \mspace{30mu}
  f1_{(1,0)} \;\mapsto\;
  \xy (0,0)*{\begin{tikzpicture}[thick,xscale=-.5,yscale=.5,rotate=90]
  \draw[web0] (0,0) to (0,1);\draw[web0] (1,1) to (1,2);
  \draw[web1] (1,0) to (1,1) to (0,1) to (0,2);
  \end{tikzpicture}}\endxy
  \mspace{30mu}
  f1_{(2,1)} \;\mapsto\;
  \xy (0,0)*{\begin{tikzpicture}[thick,xscale=-.5,yscale=.5,rotate=90]
  \draw[web2] (1,0) to (1,1);\draw[web2] (0,1) to (0,2);
  \draw[web1] (0,0) to (0,1) to (1,1) to (1,2);
  \end{tikzpicture}}\endxy
  \mspace{30mu}
  f1_{(2,0)} \;\mapsto\;
  \xy (0,0)*{\begin{tikzpicture}[thick,xscale=-.5,yscale=.5,rotate=90]
  \draw[web0] (0,0) to (0,1);
  \draw[web2] (1,0) to (1,1);
  \draw[web1] (0,2) to (0,1) to (1,1) to (1,2);
  \end{tikzpicture}}\endxy
\end{gather*}
The ladder diagrammatics can be understood as a ``rigidification'' of web diagrammatics (see \cref{subsec:webs}). Forgetting zero entries defines a mapping $\lambda\mapsto\ul{\lambda}$ from $\Lambda_{n,d}$ to $\uLam_d$, and forgetting dotted lines and smoothing corners defines a functor $F_{n,d}$ from $\qschur_{n,d}$ to $\web_d$.

The following is a special case of the combination of Theorem~4.4.1 and Theorem~3.3.1 from \cite{CKM_WebsQuantumSkew_2014}:

\begin{lemma}
  \label{lem:faithful_schur_web}
  The functor $F_{n,d}\colon\qschur_{n,d}\to\web_d$ is faithful.
  \hfill\qed
\end{lemma}

%%%%%%%%%%%%%%%%%%%%%%%%%%%%%%%%%%%%%
%%%          cat q-schur          %%%
%%%%%%%%%%%%%%%%%%%%%%%%%%%%%%%%%%%%%
\subsection{The graded 2-Schur algebra}
\label{subsec:graded_2_schur}

Recall the definitions of $\ringfoam$ and $\bilfoam\colon\bZ^2\times\bZ^2\to \ringfoam^\times$ in \cref{defn:ring_R_bil}. We use the following notation:
\begin{equation*}
  p_{ij}\coloneqq
  ((\alpha_j)_{i+1},-(\alpha_j)_i)  
  =\begin{cases}
    (0,1) & \text{if }j=i-1,\\
    (-1,-1) & \text{if }j=i,\\
    (1,0) & \text{if }j=i+1,\\
    (0,0) & \text{otherwise.}
  \end{cases}
\end{equation*}

\begin{definition}
  \label{defn:two_schur}
  The \emph{graded 2-Schur algebra} $\catschur_{n,d}$ is the $\bZ$-graded $(\bZ^2,\bilfoam)$-graded-2\nbd-cate\-gory such that:
  \begin{itemize}
    \item Objects are the elements of $\Lnd$.
    \item 1\nbd-morphisms are compositions of the generating 1\nbd-morphisms
    \begin{equation*}
      1_\lambda\colon\lambda\to\lambda,\quad\sF_i1_\lambda\colon \lambda\to\lambda-\alpha_i
      \quad\an\quad
      \sE_i1_\lambda\colon \lambda\to\lambda+\alpha_i,
    \end{equation*}
    whenever both $\lambda$ and $\lambda-\alpha_i$ (resp.\ $\lambda$ and $\lambda+\alpha_i$) are elements of $\Lnd$. Using string diagrammatics, the identity $1_\lambda$ is not pictured, and the non-trivial 1-generators are pictured as follows:
    \begin{equation*}
      \xy (0,0)*{\begin{tikzpicture}[scale=.25]
        \draw [schur1, ->] (-2,0) -- (-2,-4);
        \node at (-.5,-1.5) {\tiny $\lambda$};
        \node at (-4.5, -1.5) {\tiny $\lambda-\alpha_i$};
        \node at (-2.5,0) {\tiny $i$};
      \end{tikzpicture}}\endxy
      = \id_{\sF_i1_\lambda}
      \mspace{80mu}
      \xy (0,0)*{\begin{tikzpicture}[scale=.25]
        \draw [schur1, ->] (-2,-4) to (-2,0);
          \node at (-.5,-1.5) {\tiny $\lambda$};
          \node at (-4.5, -1.5) {\tiny $\lambda+\alpha_i$};
        \node at (-2.5,-4) {\tiny $i$};
      \end{tikzpicture}}\endxy
      = \id_{\sE_i1_\lambda}
    \end{equation*}
    Note that we read from bottom to top and from right to left.
    \item 2\nbd-morphisms are $\ringfoam$\nbd-linear combinations of string diagrams generated by formal vertical and horizontal compositions of the following generating 2\nbd-morphisms:
    \begin{gather*}
      \begin{array}{c@{\hskip 10ex}c}
        %%%% leftward cup %%%%
        \xy (0,0)*{\begin{tikzpicture}[scale=1]
          \draw[schur1,<-] (0,1) to [out=270,in=180] (.25,.5)
            to [out=0,in=270] (.5,1) node[right=-3pt]{\tiny $i$};
          \node at (.6,0.5) {\tiny $\lambda$};
        \end{tikzpicture}}\endxy
        \colon 1_\lambda\to \sE_i\sF_i1_\lambda
        &
        %%%% leftward cap %%%%
        \xy(0,0)*{\begin{tikzpicture}[scale=1]
          \draw[schur1,<-] (0,0) to [out=90,in=180] (.25,.5)
            to [out=0,in=90] (.5,0) node[right=-3pt]{\tiny $i$};
          \node at (0.6,.5) {\tiny $\lambda$};
        \end{tikzpicture}}\endxy
        \colon \sE_i\sF_i1_\lambda\to 1_\lambda
        \\[2ex]
        (\lambda_{i+1},-\lambda_i)+(0,1) & -(\lambda_{i+1},-\lambda_i)+(1,0)
      \end{array}
      \\[3ex]
      \begin{array}{c@{\hskip 10ex}c}
        %%%% downward dot %%%%
        \xy (0,1.5)*{\begin{tikzpicture}[scale=.8]
          \draw [schur1,<-] (0,0) to node[sdot1]{} (0,1) node[left=-3pt]{\tiny $i$};
          \node at (.3,.5) {\tiny $\lambda$};
        \end{tikzpicture}}\endxy
        \colon \sF_i1_\lambda\to \sF_i1_\lambda
        &
        %%%% downward crossing %%%%
        \xy (0,0)*{\begin{tikzpicture}[scale=.4]
          \draw [schur1,->] (1,1) to (-1,-1) node[left=-3pt]{\tiny $i$};
          \draw [schur2,->] (-1,1) to (1,-1) node[right=-3pt]{\tiny $j$};
          \node at (1.4,0) {\tiny $\lambda$};
        \end{tikzpicture}}\endxy
        \colon \sF_i\sF_j1_\lambda\to \sF_j\sF_i1_\lambda
        \\[2ex]
        \degd & p_{ij}
      \end{array}
    \end{gather*}
    where the $\bZ^2$-degree $\deg_{\bZ^2}$ is given below each generator. Such string diagrams are called \emph{Schur diagrams}.
  \end{itemize}
  2\nbd-morphisms are further subject to axioms described below. The \emph{quantum grading} on $\catschur_{n,d}$ is the $\bZ$\nbd-grading defined by $\qdeg(D)\coloneqq q(\deg_{\bZ^2}(D))$ (where $q(a,b)=a+b$; see \cref{rem:foam_grading}).
\end{definition}

%   \caption{Primary generators of $\catschur_{n,d}$ and their $\bZ^2$-grading.}
%   \label{tab:gen_catschur}
% \end{table}

One only needs to label one region in a given diagram, as it determines the label of all the other regions. If this forces a region to be labelled by a weight $\lambda$ that does not belong to $\Lnd$, we set this diagram to zero. In that case, we say that the diagram is zero ``due to the Schur quotient''.

We assume that the generators in a string diagram are always in generic position, in the sense that the vertical projection defines a separative Morse function. Generating 2\nbd-morphisms are subject to the following local relations:
\begin{enumerate}[(1)]
  %%%% superinterchange law %%%%%
  \item As in any graded-2\nbd-cate\-gory, we have the graded interchange law:
  \begin{equation*}\label{eq:chronology}
    \begin{tikzpicture}[very thick, scale=.4, baseline={([yshift=.7ex]current bounding box.center)}]
        \draw (-1,-.5) -- (1,-.5);\draw (-1, .5) -- (1, .5);
        \draw (-1,-.5) -- (-1,.5);\draw (1,-.5) -- (1,.5);
        \node at (0,0) {\small $f$};
        \draw (-.75,-2.5)node[below]{\tiny $i_1$} -- (-.75,-.5);
        \node at (0.05, -2.3) {$\dotsm$};
        \draw (.75,-2.5)node[below]{\tiny $i_k$} -- (.75,-.5);
        \draw (-.75,1) -- (-.75,.5);
        \node at (0.05, .85) {$\dotsm$};
        \draw (.75,1) -- (.75,.5);
    \end{tikzpicture}
    \mspace{15mu}
    \begin{tikzpicture}[very thick,scale=.4,baseline={([yshift=.8ex]current bounding box.center)}]
        \draw (1.5,-.5) -- (3.5,-.5);\draw (1.5, .5) -- (3.5, .5);
        \draw (1.5,-.5) -- (1.5,.5);\draw (3.5,-.5) -- (3.5,.5);
        \node at (2.5,0) {\small $g$};
        \draw (1.75,-1)node[below]{\tiny $i_1$} -- (1.75,-.5);
        \node at (2.55, -.85) {$\dotsm$};
        \draw (3.25,-1)node[below]{\tiny $i_k$} -- (3.25,-.5);
        \draw (1.75,2.5) -- (1.75,.5);
        \node at (2.55, 2.3) {$\dotsm$};
        \draw (3.25,2.5) -- (3.25,.5);
    \end{tikzpicture}
    = \bilfoam(\deg_{\bZ^2} f,\deg_{\bZ^2}g)\;
    \begin{tikzpicture}[very thick,scale=.4,baseline={([yshift=.8ex]current bounding box.center)}]
        \draw (1.5,-.5) -- (3.5,-.5);\draw (1.5, .5) -- (3.5, .5);
        \draw (1.5,-.5) -- (1.5,.5);\draw (3.5,-.5) -- (3.5,.5);
        \node at (2.5,0) {\small $f$};
        \draw (1.75,-1)node[below]{\tiny $i_1$} -- (1.75,-.5);
        \node at (2.55, -.85) {$\dotsm$};
        \draw (3.25,-1)node[below]{\tiny $i_k$} -- (3.25,-.5);
        \draw (1.75,2.5) -- (1.75,.5);
        \node at (2.55, 2.3) {$\dotsm$};
        \draw (3.25,2.5) -- (3.25,.5);
    \end{tikzpicture}
    \mspace{15mu}
    \begin{tikzpicture}[very thick,scale=.4,baseline={([yshift=.8ex]current bounding box.center)}]
        \draw (-1,-.5) -- (1,-.5);\draw (-1, .5) -- (1, .5);
        \draw (-1,-.5) -- (-1,.5);\draw (1,-.5) -- (1,.5);
        \node at (0,0) {\small $g$};
        \draw (-.75,-2.5) node[below]{\tiny $i_1$} -- (-.75,-.5);
        \node at (0.05, -2.3) {$\dotsm$};
        \draw (.75,-2.5) node[below]{\tiny $i_k$} -- (.75,-.5);
        \draw (-.75,1) -- (-.75,.5);
        \node at (0.05, .85) {$\dotsm$};
        \draw (.75,1) -- (.75,.5);
    \end{tikzpicture}
  \end{equation*}

  %%%%% level 2 Schur quotient %%%%%
  % \item The (level 2) Schur quotient. The identity 2\nbd-morphism $1_\lambda$ of the 1\nbd-morphism $\lambda$ is zero if $\lambda\notin\Lambda_{n,d}$.
  % This means that we set to zero all diagrams containing a region with a label not in $\Lambda_{n,d}$.

  %%%%% dot annihilation %%%%%
  \item Two dots annihilate:
  \begin{equation}
    \label{eq:dotnil}
    %%%%% dot annihilation %%%%%
      \tikz[scale=1.3,baseline={([yshift=.8ex]current bounding box.center)}]{
        \draw[schur1,<-] (0,0) node[below] {\textcolor{black}{\tiny $i$}} to
          node [pos=.3,sdot1]{}
          node [pos=.7,sdot1]{} (0,1);
        \node at (.3,.5) {\tiny $\lambda$};
      } = 0
  \end{equation}

  %%%%% KLR relations %%%%%
  \item Graded KLR algebra relations for downward crossings:
  \begin{gather}
    \label{eq:klrR2}
    %%%%% Reidemeister 2 moves %%%%%
    \xy(0,-1)*{\begin{tikzpicture}[xscale=.7,yscale=.8]
      \draw[schur1,<-]  +(0,-.75) node[below] {\textcolor{black}{\tiny $i$}}
        .. controls (0,-.375) and (1,-.375) .. (1,0)
        .. controls (1,.375) and (0, .375) .. (0,.75);
      \draw[schur2,<-]  +(1,-.75) node[below] {\textcolor{black}{\tiny $j$}}
        .. controls (1,-.375) and (0,-.375) .. (0,0)
        .. controls (0,.375) and (1, .375) .. (1,.75);
      \node at (1.3,0) {\tiny $\lambda$};
    \end{tikzpicture}}\endxy
    =
    \left\{\begin{array}{cl}
      %%%%% i = j %%%%%
      0 & \text{ if }i = j
      %%%%% distant %%%%%
      \\[2ex]
      \xy(0,0)*{\begin{tikzpicture}[scale=.8,yscale=0.75,xscale=.7]
        \draw[schur1,<-] (0,-.75) node[below] {\textcolor{black}{\tiny $i$}} to (0,.75);
        \draw[schur2,<-] (1,-.75) node[below] {\textcolor{black}{\tiny $j$}} to (1,.75);
        \node at (1.3,0) {\tiny $\lambda$};
      \end{tikzpicture}}\endxy
      & \text{ if }\vert i-j \vert > 1
      %%%%% j=i-1 %%%%%
      \\[2ex]
      -XYZ\;
      \xy(0,-1)*{\begin{tikzpicture}[scale=.8,yscale=0.75,xscale=.7]
        \draw[schur1,<-] (0,-.75) node[below] {\textcolor{black}{\tiny $i$}}
          to node [sdot1]{} (0,.75);
        \draw[schur2,<-] (1,-.75) node[below] {\textcolor{black}{\tiny $i+1$}} to (1,.75);
        \node at (1.3,0) {\tiny $\lambda$};
      \end{tikzpicture}}\endxy
      +\ XYZ\;
      \xy(0,-1)*{\begin{tikzpicture}[scale=.8,yscale=0.75,xscale=.7]
        \draw[schur1,<-]  (0,-.75) node[below] {\textcolor{black}{\tiny $i$}} to (0,.75);
        \draw[schur2,<-]  (1,-.75) node[below] {\textcolor{black}{\tiny $i+1$}}
          to node[sdot2]{} (1,.75);
        \node at (1.3,0) {\tiny $\lambda$};
      \end{tikzpicture}}\endxy
      & \text{ if }j=i+1
      %%%%% j=i+1 %%%%%
      \\[2ex]
      YZ^2\;
      \xy(0,-1)*{\begin{tikzpicture}[scale=.8,yscale=0.75,xscale=.7]
        \draw[schur1,<-] (0,-.75) node[below] {\textcolor{black}{\tiny $i$}}
          to node [sdot1]{} (0,.75);
        \draw[schur2,<-] (1,-.75) node[below] {\textcolor{black}{\tiny $i-1$}} to (1,.75);
        \node at (1.3,0) {\tiny $\lambda$};
      \end{tikzpicture}}\endxy
      -YZ^2\;
      \xy(0,-1)*{\begin{tikzpicture}[scale=.8,yscale=0.75,xscale=.7]
        \draw[schur1,<-]  (0,-.75) node[below] {\textcolor{black}{\tiny $i$}} to (0,.75);
        \draw[schur2,<-]  (1,-.75) node[below] {\textcolor{black}{\tiny $i-1$}}
          to node[sdot2]{} (1,.75);
        \node at (1.3,0) {\tiny $\lambda$};
      \end{tikzpicture}}\endxy
      & \text{ if }j=i-1
    \end{array}\right.
    \\[4ex]\label{eq:dotslides_first}
    %%%%% dot slides %%%%%
    \begin{tikzpicture}[scale=.7,baseline={([yshift=.9ex]current bounding box.center)}]
      \draw[schur1,<-] (0,-.5) node[below] {\textcolor{black}{\tiny $i$}}
        .. controls (0,0) and (1,0) .. (1,.5);
      \draw[schur2,<-] (1,-.5) node[below] {\textcolor{black}{\tiny $j$}}
        .. controls (1,0) and (0,0) .. (0,.5) node [near end,sdot2]{};
      \node at (1.2,0) {\tiny $\lambda$};
    \end{tikzpicture}
    \mspace{0mu} =  \bilfoam(\degd,p_{ij})
    \begin{tikzpicture}[scale=.7,baseline={([yshift=.9ex]current bounding box.center)}]
      \draw[schur1,<-] (0,-.5) node[below] {\textcolor{black}{\tiny $i$}}
        .. controls (0,0) and (1,0) .. (1,.5);
      \draw[schur2,<-] (1,-.5) node[below] {\textcolor{black}{\tiny $j$}}
        .. controls (1,0) and (0,0) .. (0,.5) node [near start,sdot2]{};
      \node at (1.2,0) {\tiny $\lambda$};
    \end{tikzpicture}
    \mspace{50mu}
    \text{ for $i\neq j$}
    %%%%%
    \\[1ex]\label{eq:dotslides_second}
    \begin{tikzpicture}[scale=.7,baseline={([yshift=.9ex]current bounding box.center)}]
      \draw[schur1,<-] (0,-.5) node[below] {\textcolor{black}{\tiny $i$}}
        .. controls (0,0) and (1,0) .. (1,.5) node [near start,sdot1]{};
      \draw[schur2,<-] (1,-.5) node[below] {\textcolor{black}{\tiny $j$}}
        .. controls (1,0) and (0,0) .. (0,.5);
      \node at (1.2,0) {\tiny $\lambda$};
    \end{tikzpicture}
    \mspace{0mu} = \bilfoam(p_{ij},\degd)
    \begin{tikzpicture}[scale=.7,baseline={([yshift=.9ex]current bounding box.center)}]
      \draw[schur1,<-] (0,-.5) node[below] {\textcolor{black}{\tiny $i$}}
        .. controls (0,0) and (1,0) .. (1,.5) node [near end,sdot1]{};
      \draw[schur2,<-] (1,-.5) node[below] {\textcolor{black}{\tiny $j$}}
        .. controls (1,0) and (0,0) .. (0,.5);
      \node at (1.2,0) {\tiny $\lambda$};
    \end{tikzpicture}
    \mspace{50mu}
    \text{ for $i\neq j$}
    \\[4ex]\label{eq:dotslide-nilH}
    \begin{tikzpicture}[scale=.7,baseline={([yshift=.9ex]current bounding box.center)}]
      \draw[schur1,<-] (0,-.5) node[below] {\textcolor{black}{\tiny $i$}}
        .. controls (0,0) and (1,0) .. (1,.5);
      \draw[schur1,<-] (1,-.5) node[below] {\textcolor{black}{\tiny $i$}}
        .. controls (1,0) and (0,0) .. (0,.5) node [near end,sdot1]{};
      \node at (1.3,0) {\tiny $\lambda$};
    \end{tikzpicture}
    -XY
    \begin{tikzpicture}[scale=.7,baseline={([yshift=.9ex]current bounding box.center)}]
      \draw[schur1,<-] (0,-.5) node[below] {\textcolor{black}{\tiny $i$}}
        .. controls (0,0) and (1,0) .. (1,.5);
      \draw[schur1,<-] (1,-.5) node[below] {\textcolor{black}{\tiny $i$}}
        .. controls (1,0) and (0,0) .. (0,.5) node [near start,sdot1]{};
      \node at (1.3,0) {\tiny $\lambda$};
    \end{tikzpicture}
    =
    \begin{tikzpicture}[scale=.7,baseline={([yshift=.9ex]current bounding box.center)}]
      \draw[schur1,<-] (0,-.5) node[below] {\textcolor{black}{\tiny $i$}} to (0,.5);
      \draw[schur1,<-] (1,-.5) node[below] {\textcolor{black}{\tiny $i$}} to (1,.5);
      \node at (1.3,0) {\tiny $\lambda$};
    \end{tikzpicture}
    =
    \begin{tikzpicture}[scale=.7,baseline={([yshift=.9ex]current bounding box.center)}]
      \draw[schur1,<-] (0,-.5) node[below] {\textcolor{black}{\tiny $i$}}
        .. controls (0,0) and (1,0) .. (1,.5) node [near start,sdot1]{};
      \draw[schur1,<-] (1,-.5) node[below] {\textcolor{black}{\tiny $i$}}
        .. controls (1,0) and (0,0) .. (0,.5);
      \node at (1.3,0) {\tiny $\lambda$};
    \end{tikzpicture}
    -XY
    \begin{tikzpicture}[scale=.7,baseline={([yshift=.9ex]current bounding box.center)}]
      \draw[schur1,<-] (0,-.5) node[below] {\textcolor{black}{\tiny $i$}}
        .. controls (0,0) and (1,0) .. (1,.5) node [near end,sdot1]{};
      \draw[schur1,<-] (1,-.5) node[below] {\textcolor{black}{\tiny $i$}}
        .. controls (1,0) and (0,0) .. (0,.5);
      \node at (1.3,0) {\tiny $\lambda$};
    \end{tikzpicture}
    \\[4ex]\label{eq:klrR3}
    %%%%% Reidemeister 3 moves %%%%%
    \begin{array}{c}
      \xy(0,-1)*{\begin{tikzpicture}[scale=.6]
        \draw[schur1,<-]  +(0,0)node[below] {\textcolor{black}{\tiny $i$}}
          .. controls (0,0.5) and (2, 1) ..  +(2,2);
        \draw[schur2,<-]  +(2,0)node[below] {\textcolor{black}{\tiny $k$}}
          .. controls (2,1) and (0, 1.5) ..  +(0,2);
        \draw[schur3,<-]  (1,0)node[below] {\textcolor{black}{\tiny $j$}}
          .. controls (1,0.5) and (0, 0.5) ..  (0,1)
          .. controls (0,1.5) and (1, 1.5) ..  (1,2);
        \node at (2.1,1) {\tiny $\lambda$};
      \end{tikzpicture}}\endxy
      \mspace{5mu}=
      \bilfoam(p_{jk},p_{ij})\bilfoam(p_{ik},p_{ij})\bilfoam(p_{jk},p_{ik})
      \xy(0,-1)*{\begin{tikzpicture}[scale=.6]
        \draw[schur1,<-]  +(0,0)node[below] {\textcolor{black}{\tiny $i$}}
          .. controls (0,1) and (2, 1.5) ..  +(2,2);
        \draw[schur2,<-]  +(2,0)node[below] {\textcolor{black}{\tiny $k$}}
          .. controls (2,.5) and (0, 1) ..  +(0,2);
        \draw[schur3,<-]  (1,0)node[below]{\textcolor{black} {\tiny $j$}}
          .. controls (1,0.5) and (2, 0.5) ..  (2,1)
          .. controls (2,1.5) and (1, 1.5) ..  (1,2);
        \node at (2.3,1) {\tiny $\lambda$};
      \end{tikzpicture}}\endxy
      \\[-2ex]
      \text{ unless }i=k\text{ and }\vert i-j\vert = 1
    \end{array}
    %%%%%
    \\[4ex]\label{eq:R3serre_first_case}
    -YZ^{-2}
    \xy(0,-1)*{\begin{tikzpicture}[scale=.6]
      \draw[schur1,<-]  +(0,0)node[below] {\textcolor{black}{\tiny $i$}}
        .. controls (0,0.5) and (2, 1) ..  +(2,2);
      \draw[schur1,<-]  +(2,0)node[below] {\textcolor{black}{\tiny $i$}}
        .. controls (2,1) and (0, 1.5) ..  +(0,2);
      \draw[schur2,<-]  (1,0)node[below] {\textcolor{black}{\tiny $i+1$}}
        .. controls (1,0.5) and (0, 0.5) ..  (0,1)
        .. controls (0,1.5) and (1, 1.5) ..  (1,2);
      \node at (2.1,1) {\tiny $\lambda$};
    \end{tikzpicture}}\endxy
    +
    Z^{-1}
    \xy(0,-1)*{\begin{tikzpicture}[scale=.6]
      \draw[schur1,<-]  +(0,0)node[below] {\textcolor{black}{\tiny $i$}}
        .. controls (0,1) and (2, 1.5) ..  +(2,2);
      \draw[schur1,<-]  +(2,0)node[below] {\textcolor{black}{\tiny $i$}}
        .. controls (2,.5) and (0, 1) ..  +(0,2);
      \draw[schur2,<-]  (1,0)node[below] {\textcolor{black}{\tiny $i+1$}}
        .. controls (1,0.5) and (2, 0.5) ..  (2,1)
        .. controls (2,1.5) and (1, 1.5) ..  (1,2);
      \node at (2.3,1) {\tiny $\lambda$};
    \end{tikzpicture}}\endxy
    =
    \xy(0,-1)*{\begin{tikzpicture}[scale=.6]
      \draw[schur1,<-] (0,0) node[below] {\textcolor{black}{\tiny $i$}} to (0,2);
      \draw[schur2,<-] (1,0) node[below] {\textcolor{black}{\tiny $i+1$}} to (1,2);
      \draw[schur1,<-] (2,0) node[below] {\textcolor{black}{\tiny $i$}} to (2,2);
      \node at (2.3,1) {\tiny $\lambda$};
    \end{tikzpicture}}\endxy
    %%%%%
    \\\label{eq:R3serre_second_case}
    XYZ^{-1}
    \xy(0,-1)*{\begin{tikzpicture}[scale=.6]
      \draw[schur1,<-]  +(0,0)node[below] {\textcolor{black}{\tiny $i$}}
        .. controls (0,0.5) and (2, 1) ..  +(2,2);
      \draw[schur1,<-]  +(2,0)node[below] {\textcolor{black}{\tiny $i$}}
        .. controls (2,1) and (0, 1.5) ..  +(0,2);
      \draw[schur2,<-]  (1,0)node[below] {\textcolor{black}{\tiny $i-1$}}
        .. controls (1,0.5) and (0, 0.5) ..  (0,1)
        .. controls (0,1.5) and (1, 1.5) ..  (1,2);
      \node at (2.1,1) {\tiny $\lambda$};
    \end{tikzpicture}}\endxy
    -
    XZ^{-2}
    \xy(0,-1)*{\begin{tikzpicture}[scale=.6]
      \draw[schur1,<-]  +(0,0)node[below] {\textcolor{black}{\tiny $i$}}
        .. controls (0,1) and (2, 1.5) ..  +(2,2);
      \draw[schur1,<-]  +(2,0)node[below] {\textcolor{black}{\tiny $i$}}
        .. controls (2,.5) and (0, 1) ..  +(0,2);
      \draw[schur2,<-]  (1,0)node[below] {\textcolor{black}{\tiny $i-1$}}
        .. controls (1,0.5) and (2, 0.5) ..  (2,1)
        .. controls (2,1.5) and (1, 1.5) ..  (1,2);
      \node at (2.3,1) {\tiny $\lambda$};
    \end{tikzpicture}}\endxy
    =
    \xy(0,-1)*{\begin{tikzpicture}[scale=.6]
      \draw[schur1,<-] (0,0) node[below] {\textcolor{black}{\tiny $i$}} to (0,2);
      \draw[schur2,<-] (1,0) node[below] {\textcolor{black}{\tiny $i-1$}} to (1,2);
      \draw[schur1,<-] (2,0) node[below] {\textcolor{black}{\tiny $i$}} to (2,2);
      \node at (2.3,1) {\tiny $\lambda$};
    \end{tikzpicture}}\endxy
  \end{gather}

  %%%%% adjunction relations %%%%%
  \item Graded adjunction relations:
  \begin{equation}
    \label{eq:adjrels}
    \xy (0,0)*{\begin{tikzpicture}[scale=1*.8]
    	\draw [schur1,<-] (0,0) to (0,1) node[left=-3pt]{\tiny $i$};
    	\node at (.3,.5) {\tiny $\lambda$};
    \end{tikzpicture}}\endxy
    =
    \xy (0,0)*{\tikz[scale=1*.8]{
      \draw[schur1,<-] (0,-.5) to (0,0) to[out=90,in=180] (.25,.5)
        to[out=0,in=90] (.5,0) to[out=-90,in=180] (.75,-.5)
        to[out=0,in=-90] (1,0) to (1,.5) node[left=-3pt]{\tiny $i$};
      \node at (1.2,0) {\tiny $\lambda$};
    }}\endxy
    \mspace{80mu}
    \xy (0,-2)*{\begin{tikzpicture}[scale=1*.8]
    	\draw [schur1,->] (0,0) node[left=-3pt]{\tiny $i$} to (0,1);
    	\node at (.3,.5) {\tiny $\lambda$};
    \end{tikzpicture}}\endxy
    =
    X^{1+\lambda_{i+1}}Y^{\lambda_i}
    \xy (0,-2)*{\tikz[scale=1*.8]{
      \draw[schur1,<-] (0,.5) to (0,0) to[out=-90,in=180] (.25,-.5)
        to[out=0,in=-90] (.5,0) to[out=90,in=180] (.75,.5)
        to[out=0,in=90] (1,0) to (1,-.5) node[left=-3pt]{\tiny $i$};
      \node at (1.2,0) {\tiny $\lambda$};
    }}\endxy
  \end{equation}
\end{enumerate}
Finally, we require invertibility axioms. To define them, we introduce the following shorthands:
\begin{equation*}
  %%%% multiple dots %%%%
  \xy (0,0)*{\begin{tikzpicture}[scale=1*.8]
    \draw [schur1,<-] (0,0) to node[sdot1]{}
    node[above left=-2pt]{\tiny $n$} (0,1) node[left=-3pt]{\tiny $i$};
    \node at (.3,.5) {\tiny $\lambda$};
  \end{tikzpicture}}\endxy
  \coloneqq
  \left(
  \xy (0,0)*{\begin{tikzpicture}[scale=1*.8]
    \draw [schur1,<-] (0,0) to node[sdot1]{} (0,1) node[left=-3pt]{\tiny $i$};
    \node at (.3,.5) {\tiny $\lambda$};
  \end{tikzpicture}}\endxy\right)^{\circ n}
  \mspace{100mu}
  %%%% leftward crossing %%%%
  \xy (0,-.5)*{\begin{tikzpicture}[scale=.4*.8]
    \draw [schur1,->] (1,1) to (-1,-1) node[left=-3pt]{\tiny $i$};
    \draw [schur2,<-] (-1,1) to (1,-1) node[right=-3pt]{\tiny $j$};
    \node at (1.4,0) {\tiny $\lambda$};
  \end{tikzpicture}}\endxy
  \coloneqq
  \xy (0,0)*{\begin{tikzpicture}[scale=1*.8]
    \draw[schur2,<-] (0.3,-0.5) to (-0.3,0.5);
    \draw[schur1,<-] (-0.8,0.5) to [out=-90,in=180] (-0.5,-0.5) to [out=0,in=230] (-0.2,-0.3) to [out=50,in=180] (0.5,0.5) to [out=0,in=90] (0.8,-0.5);
    \node at (.9,-.5) {\tiny $i$};
    \node at (-0.4,.5) {\tiny $j$};
    \node at (1.05,.05) {\tiny $\lambda$};
  \end{tikzpicture}}\endxy
\end{equation*}
Recall the notation $\overline{\lambda}_i\coloneqq\lambda_i-\lambda_{i+1}$.
Then:
\begin{enumerate}[(1),resume]
  %%%%% isomorphisms %%%%%
  \item Except if they are zero due to the Schur quotient, the following 2\nbd-mor\-phisms are isomorphisms in the graded additive envelope of $\catschur_{n,d}$ (see \cref{subsec:more_cat_facts}):
  \begin{align}
    %%%% i\neq j %%%%
    \label{eq:iso-ij}
    \xy (0,-.5)*{\begin{tikzpicture}[scale=.4*.8]
      \draw [schur1,->] (1,1) to (-1,-1) node[left=-3pt]{\tiny $i$};
      \draw [schur2,<-] (-1,1) to (1,-1) node[right=-3pt]{\tiny $j$};
      \node at (1.4,0) {\tiny $\lambda$};
    \end{tikzpicture}}\endxy
    &\colon\sF_i\sE_j(\lambda)\to \sE_j\sF_i(\lambda)
    && \text{if }i\neq j
    \\
    %%%% \lambda_i>0 %%%%
    \label{eq:isoii-pos}
    \xy (0,-.5)*{\begin{tikzpicture}[scale=.4*.8]
      \draw [schur1,->] (1,1) to (-1,-1) node[left=-3pt]{\tiny $i$};
      \draw [schur1,<-] (-1,1) to (1,-1) node[right=-3pt]{\tiny $i$};
      \node at (1.4,0) {\tiny $\lambda$};
    \end{tikzpicture}}\endxy
    \oplus\bigoplus\limits_{n=0}^{\overline{\lambda}_i-1}
    \xy (0,0)*{\begin{tikzpicture}[scale=1.2*.8]
      \draw[schur1,<-] (0,1) to [out=270,in=180] (.25,.5)
        to [out=0,in=270] node[sdot1]{}
        node[below right=-2pt]{\tiny $n$} (.5,1) node[right=-3pt]{\tiny $i$};
      \node at (-.1,0.5) {\tiny $\lambda$};
    \end{tikzpicture}}\endxy
    &\colon \sF_i\sE_i(\lambda)\oplus \lambda^{\oplus[\overline{\lambda}_i]}\to\sE_i\sF_i(\lambda)
    && \text{if }\overline{\lambda}_i \geq 0
    \\
    %%%% \lambda_i<0 %%%%
    \label{eq:isoii-neg}
    \xy (0,-.5)*{\begin{tikzpicture}[scale=.4*.8]
      \draw [schur1,->] (1,1) to (-1,-1) node[left=-3pt]{\tiny $i$};
      \draw [schur1,<-] (-1,1) to (1,-1) node[right=-3pt]{\tiny $i$};
      \node at (1.4,0) {\tiny $\lambda$};
    \end{tikzpicture}}\endxy
    \oplus\bigoplus\limits_{n=0}^{-\overline{\lambda}_i-1}
    \xy(0,0)*{\begin{tikzpicture}[scale=1.2*.8]
      \draw[schur1,<-] (0,0) to[out=90,in=180] node[sdot1]{}
        node[above left=-2pt]{\tiny $n$} (.25,.5)
        to[out=0,in=90] (.5,0) node[right=-3pt]{\tiny $i$};
      \node at (0.6,.5) {\tiny $\lambda$};
    \end{tikzpicture}}\endxy
    &\colon \sF_i\sE_i(\lambda)\to \sE_i\sF_i(\lambda)\oplus \lambda^{\oplus[-\overline{\lambda}_i]}
    && \text{if }\overline{\lambda}_i \leq 0
  \end{align} 
\end{enumerate}
This ends the definition of the relations on the graded 2-Schur algebra.\hfill$\diamond$

\begin{remark}
  Let us elaborate on the invertibility axioms above.
  They are equivalent to the existence of some unnamed generators which are entries of the inverse matrices of \eqref{eq:iso-ij}, \eqref{eq:isoii-pos} and \eqref{eq:isoii-neg}, and some unnamed relations that precisely encompass the fact that those generators form inverse matrices.
  This definition follows Rouquier's approach \cite{Rouquier_2KacMoodyAlgebras_2008} to 2-Kac--Moody algebras (categorified quantum groups) and Brundan and Ellis' approach \cite{BE_SuperKacMoody_2017} to super 2\nbd-Kac--Moody algebras.
  Unravelling the definition would lead to a more explicit (but heavier) definition, similar to Khovanov and Lauda's approach \cite{KL_DiagrammaticApproachCategorification_2009} to categorified quantum groups.
\end{remark}

\begin{remark}
  Some scalars in the relations above can be understood as artefacts of interchanging the vertical positions of the generators. For instance, the scalar in relation \eqref{eq:dotslides_first} is precisely the scalar that appears when interchanging a dot and a $(i,j)$-crossing. A similar reasoning applies to relations \eqref{eq:dotslides_second} and \eqref{eq:klrR3}.
\end{remark}

In \cite[Definition~3.2]{MSV_DiagrammaticCategorification$q$Schur_2013}, a categorification of the $q$-Schur algebra was constructed, denoted $\cS_{n,d}$. Let us write $\cS_{n,d}^\bullet$ for the $\ringfoam$\nbd-linear\footnote{The linear 2\nbd-cate\-gory $\cS_{n,d}$ is defined over $\bQ$ in \cite{MSV_DiagrammaticCategorification$q$Schur_2013}, but it can be defined over $\bZ$ and hence over any unital commutative ring.} 2\nbd-cate\-gory obtained from $\cS_{n,d}$ by further imposing relation \eqref{eq:dotnil}.
Then:

\begin{proposition}
  Let $\ul{(\catschur_{n,d})}^{\oplus}_q$ be the additive shifted closure of $\catschur_{n,d}$ with respect to the quantum grading (see \cref{subsec:more_cat_facts} and \cref{defn:two_schur}). Then:
  \begin{equation*}
    \cS_{n,d}^\bullet \cong \left.\ul{(\catschur_{n,d})}^{\oplus}_q\right\vert_{X=Y=Z=1}.
  \end{equation*}
\end{proposition}

\begin{proof}[Sketch of proof]
  Defining the unnamed generators implied by the invertibility axioms and doing a relation chase exhibits the missing relations. One only needs to rescale the $(i,i)$-downward crossing by $(-1)$.
  This exactly follows the proof from \cite{Brundan_DefinitionKacMoody_2016} showing the equivalence between Rouquier's definition of Kac--Moody 2-algebras and Khovanov--Lauda categorified quantum groups. See also \cite{BE_SuperKacMoody_2017} for a related statement in the super case.
\end{proof}

%%%%%%%%%%%%%%%%%%%%%%%%%%%%%%%%%%%%%
%%%          equivalence          %%%
%%%%%%%%%%%%%%%%%%%%%%%%%%%%%%%%%%%%%
\subsection{The graded foamation 2-functor}
\label{sec:foamation_functor}

This section exhibits $\gfoam_d$ as a 2-representation of $\catschur_{n,d}$. More precisely, for each $n,d\in\bN$ there exists a $(\bZ^2,\bilfoam)$-graded-2-functor
\begin{equation*}
  \cF_{n,d}\colon \catschur_{n,d}\to\gfoam_d.
\end{equation*}
This categorifies the functor $F_{n,d}\colon\qschur_{n,d}\to\web_d$ defined in \cref{subsec:schur_ladder_to_web}.
% This categorifies the ladder diagrammatics of the Schur algebra, which we now recall.

On the level of objects, the functor $F_{n,d}$ maps a weight $\lambda\in\Lnd$ to the weight $\ul{\lambda}\in\uLam_d$, obtained by forgetting all zero entries in $\lambda$. Recall the colour of a coordinate defined in \cref{subsec:webs}. For $i\in\{1,\ldots,n\}$ such that $\lambda_i\neq 0$, we denote $\underline{i}_\lambda$ the colour of the coordinate of the ``image'' of $i$ in $\ul{\lambda}$.
For instance:
\begin{gather*}
  \underline{(1,0,1,2,0,1)}=(1,1,2,1)
  \quad\an\quad \underline{1}_\lambda=1,\;\underline{3}_\lambda=2,\;\underline{4}_\lambda=3\an\underline{6}_\lambda=5.
\end{gather*}
In the string diagrammatics of foams, the functor $F_{n,d}$ is given by
\begin{gather}
  \label{eq:defn_schur_to_foam_downward_strand}
  \begin{array}{*{2}{c@{\hskip 3ex}}*{6}{c@{\hskip 2ex}}c}
    \xy (0,4)*{\begin{tikzpicture}[scale=1.3]
      \draw[schur1,<-] (0,0) to (0,1) node[above right=-4pt]{\scriptsize $i$};
      \node at (.2,0.5) {\tiny $\lambda$};
    \end{tikzpicture}}\endxy
    &
    \mapsto
    &
    \xy (0,4)*{\begin{tikzpicture}[scale=1]
      \node at (.3,0.5) {$\underline{\lambda}$};
    \end{tikzpicture}}\endxy
    &
    ,
    &
    \xy (0,4)*{\begin{tikzpicture}[scale=1]
      \draw[diag1,->] (0,0) to (0,1) node[above right=-4pt]{\scriptsize $\ul{i}_\lambda$};
      \node at (.3,0.5) {\tiny $\underline{\lambda}$};
    \end{tikzpicture}}\endxy
    &
    ,
    &
    \xy (0,4)*{\begin{tikzpicture}[scale=1]
      \draw[diag1,<-] (0,0) to (0,1) node[above right=-4pt]{\scriptsize $\ul{i}_\lambda$};
      \node at (.3,0.5) {\tiny $\underline{\lambda}$};
    \end{tikzpicture}}\endxy
    &
    ,
    &
    \xy (0,4)*{\begin{tikzpicture}[scale=1]
      \draw[diag2,<-] (-.7,0) to (-.7,1) node[above left=-4pt]{\scriptsize $\ul{i}_\lambda{+}1$};
      \draw[diag1,->] (0,0) to (0,1) node[above right=-4pt]{\scriptsize $\ul{i}_\lambda$};
      \node at (.3,0.5) {\tiny $\underline{\lambda}$};
    \end{tikzpicture}}\endxy
    \\[3ex]
    (\lambda_i,\lambda_{i+1}) && (1,0) && (2,0) && (1,1) && (2,1)
  \end{array}
  \\
  \label{eq:defn_schur_to_foam_upward_strand}
  \begin{array}{*{2}{c@{\hskip 5ex}}*{6}{c@{\hskip 2ex}}c}
    \xy (0,4)*{\begin{tikzpicture}[scale=1.3]
      \draw[schur1,->] (0,0) to (0,1) node[above right=-4pt]{\scriptsize $i$};
      \node at (.2,0.5) {\tiny $\lambda$};
    \end{tikzpicture}}\endxy
    &
    \mapsto
    &
    \xy (0,4)*{\begin{tikzpicture}[scale=1]
      \node at (.3,0.5) {$\underline{\lambda}$};
    \end{tikzpicture}}\endxy
    &
    ,
    &    
    \xy (0,4)*{\begin{tikzpicture}[scale=1]
      \draw[diag1,->] (0,0) to (0,1) node[above right=-4pt]{\scriptsize $\ul{i}_\lambda$};
      \node at (.3,0.5) {\tiny $\underline{\lambda}$};
    \end{tikzpicture}}\endxy
    &
    ,
    &
    \xy (0,4)*{\begin{tikzpicture}[scale=1]
      \draw[diag1,<-] (0,0) to (0,1) node[above right=-4pt]{\scriptsize $\ul{i}_\lambda$};
      \node at (.3,0.5) {\tiny $\underline{\lambda}$};
    \end{tikzpicture}}\endxy
    &
    ,
    &
    \xy (0,4)*{\begin{tikzpicture}[scale=1]
      \draw[diag1,<-] (-.7,0) to (-.7,1) node[above left=-4pt]{\scriptsize $\ul{i}_\lambda$};
      \draw[diag2,->] (0,0) to (0,1) node[above right=-4pt]{\scriptsize $\ul{i}_\lambda{+}1$};
      \node at (.3,0.5) {\tiny $\underline{\lambda}$};
    \end{tikzpicture}}\endxy
    \\[3ex]
    (\lambda_i,\lambda_{i+1}) && (0,1) && (0,2) && (1,1) && (1,2)
  \end{array}
\end{gather}
The local data of $(\lambda_i,\lambda_{i+1})$ is given below each case.

Following \cite[p.~59]{NP_OddKhovanovHomology_2020}, we shall use the scalar
\begin{gather*}
  \Gamma_\lambda(i)\coloneqq (-XY)^{\#\{\lambda_j=1\mid j\leq i\}}
\end{gather*}
to normalize the graded foamation 2-functor.

\begin{proposition}
  \label{prop:foamation_two_functor}
  There exists a $(\bZ^2,\bilfoam)$-graded-2-functor
  \begin{equation*}
    \cF_{n,d}\colon \catschur_{n,d}\to\gfoam_d
  \end{equation*}
  defined on generating 2\nbd-morphisms as in \cref{fig:defn_foamation_functor}. We call $\cF_{n,d}$ the \emph{graded foamation 2-functor}.
\end{proposition}

\begin{figure}
  \def\SclSch{.6}
  \def\SclFoam{.7}
  \def\hsin{0ex}
  \def\hsout{2ex}
  \begin{gather*}
    \begin{array}{*{2}{c@{\hskip 1ex}}*{6}{c@{\hskip 1ex}}c}
      %%%%% DOWNWARD DOT %%%%%
      \xy (0,4)*{\begin{tikzpicture}[scale=1.3*\SclSch]
        \draw[schur1,<-] (0,0) node[below=-1pt]{\scriptsize $i$} to node[sdot1]{} (0,1);
        \node at (.3,.5) {\scriptsize $\lambda$};
      \end{tikzpicture}}\endxy
      \quad
      &
      \overset{\raise.7em\hbox{\scalebox{1.1}{$\Gamma_\lambda(i)$}}}{\mapsto}
      &
      \quad
      \xy (0,4)*{\begin{tikzpicture}[scale=1*\SclFoam]
        \node[fdot1] at (0,0) {};
        \node at (.5,0) {\scriptsize $\ul{\lambda}$};
        \node at (-.3,.3) {\scriptsize $\ul{i}_\lambda$};
      \end{tikzpicture}}\endxy
      &&
      \xy (0,4)*{\begin{tikzpicture}[scale=1*\SclFoam]
        \draw[diag1,->] (0,0) to (0,1) node[above]{\scriptsize $\ul{i}_\lambda$};
        \node[fdot2] at (-.6,0.5) {};
        \node at (.5,.5) {\scriptsize $\ul{\lambda}$};
      \end{tikzpicture}}\endxy
      &&
      \xy (0,4)*{\begin{tikzpicture}[scale=1*\SclFoam]
        \draw[diag1,<-] (0,0) to (0,1) node[above]{\scriptsize $\ul{i}_\lambda$};
        \node[fdot1] at (.6,0.5) {};
        \node at (1.3,.5) {\scriptsize $\ul{\lambda}$};
      \end{tikzpicture}}\endxy
      &&
      \xy (0,4)*{\begin{tikzpicture}[scale=1*\SclFoam]
        \draw[diag2,<-] (-.7,0) to (-.7,1) node[above,shift={(-.3,0)}]{\scriptsize $\ul{i}_\lambda{+}1$};
        \draw[diag1,->] (0,0) to (0,1) node[above]{\scriptsize $\ul{i}_\lambda$};
        \node[fdot2] at (-.35,0.5) {};
        \node at (.5,.5) {\scriptsize $\ul{\lambda}$};
      \end{tikzpicture}}\endxy
      \\[\hsin]
      (\lambda_i,\lambda_{i+1}) && (1,0) && (2,0) && (1,1) && (2,1)
      \\[\hsout]
      \xy (0,4)*{\begin{tikzpicture}[scale=1*\SclSch]
        \draw[schur1,<-] (-.5,0) to[out=-90,in=180] (0,-.6) to[out=0,in=-90] (.5,0)node[above=-1pt]{\scriptsize $i$};
        \node at (.7,-.4) {\scriptsize $\lambda$};
      \end{tikzpicture}}\endxy
      &\mapsto&
      \xy (0,4)*{\begin{tikzpicture}[scale=1*\SclFoam]
        \node at (0,0) {$\underline{\lambda}$};
      \end{tikzpicture}}\endxy
      &&
      Z^{-2}\mspace{-15mu}\xy (0,4)*{\begin{tikzpicture}[scale=.8*\SclFoam]
        \draw[diag1,->] (-.5,0) node[above]{\scriptsize $\ul{i}_\lambda$} to[out=-90,in=180] (0,-.6) to[out=0,in=-90] (.5,0);
        \node at (.7,-.7) {\scriptsize $\ul{\lambda}$};
      \end{tikzpicture}}\endxy
      &&
      \xy (0,4)*{\begin{tikzpicture}[scale=.8*\SclFoam]
        \draw[diag1,<-] (-.5,0) to[out=-90,in=180] (0,-.6) to[out=0,in=-90] (.5,0)  node[above]{\scriptsize $\ul{i}_\lambda$};
        \node at (.7,-.7) {\scriptsize $\ul{\lambda}$};
      \end{tikzpicture}}\endxy
      &&
      Z^{-1}\mspace{-15mu}\xy (0,4)*{\begin{tikzpicture}[scale=.8*\SclFoam]
        \draw[diag2,<-] (-.5,0) to[out=-90,in=180] (0,-.6) to[out=0,in=-90] (.5,0)node[above]{\scriptsize $\ul{i}_\lambda{+}1$};
        \draw[diag1,->] (-.5*2,0) node[above]{\scriptsize $\ul{i}_\lambda$} to[out=-90,in=180] (0,-.6*2) to[out=0,in=-90] (.5*2,0);
        \node at (1.2,-1.2) {\scriptsize $\ul{\lambda}$};
      \end{tikzpicture}}\endxy
      \\[\hsin]
      (\lambda_i,\lambda_{i+1}) && (1,0) && (2,0) && (1,1) && (2,1)
      \\[\hsout]
      %%%%% LEFTWARD CAP %%%%%
      \xy (0,4)*{\begin{tikzpicture}[scale=1*\SclSch]
        \draw[schur1,<-] (-.5,0) to[out=90,in=180] (0,.6) to[out=0,in=90] (.5,0)node[below=-1pt]{\scriptsize $i$};
        \node at (.7,.4) {\scriptsize $\lambda$};
      \end{tikzpicture}}\endxy
      &\mapsto&
      \xy (0,4)*{\begin{tikzpicture}[scale=1*\SclFoam]
        \node at (0,0) {$\underline{\lambda}$};
      \end{tikzpicture}}\endxy
      &&
      \xy (0,4)*{\begin{tikzpicture}[scale=.8*\SclFoam]
        \draw[diag1,<-] (-.5,0) to[out=90,in=180] (0,.6) to[out=0,in=90] (.5,0) node[below]{\scriptsize $\ul{i}_\lambda$};
        \node at (.7,.7) {\scriptsize $\ul{\lambda}$};
      \end{tikzpicture}}\endxy
      &&
      \xy (0,4)*{\begin{tikzpicture}[scale=.8*\SclFoam]
        \draw[diag1,->] (-.5,0) node[below]{\scriptsize $\ul{i}_\lambda$} to[out=90,in=180] (0,.6) to[out=0,in=90] (.5,0);
        \node at (.7,.7) {\scriptsize $\ul{\lambda}$};
      \end{tikzpicture}}\endxy
      &&
      \xy (0,4)*{\begin{tikzpicture}[scale=.8*\SclFoam]
        \draw[diag1,->] (-.5,0) node[below]{\scriptsize $\ul{i}_\lambda$} to[out=90,in=180] (0,.6) to[out=0,in=90] (.5,0);
        \draw[diag2,<-] (-.5*2,0) to[out=90,in=180] (0,.6*2) to[out=0,in=90] (.5*2,0) node[below,shift={(.3,0)}]{\scriptsize $\ul{i}_\lambda{+}1$};
        \node at (1.2,1.2) {\scriptsize $\ul{\lambda}$};
      \end{tikzpicture}}\endxy
      \\[\hsin]
      (\lambda_i,\lambda_{i+1}) && (0,1) && (0,2) && (1,1) && (1,2)
      \\[\hsout]
      %%%%% ADJACENT CROSSING %%%%%
      \xy (0,3)*{\begin{tikzpicture}[scale=1.1*\SclSch]
        \draw[schur1,<-] (0,0) node[below=-1pt]{\scriptsize $i$} to (1,1);
        \draw[schur2,<-] (1,0) node[below=-1pt]{\scriptsize $i-1$} to (0,1);
        \node at (1,.5) {\scriptsize $\lambda$};
      \end{tikzpicture}}\endxy
      &\mapsto&
      \xy (0,1)*{\begin{tikzpicture}[scale=.8*\SclFoam]
        \draw[diag2,->,shift={(.5,0)}] (-.5,0) to[out=90,in=180] (0,.6) to[out=0,in=90] (.5,0) node[below]{\scriptsize $\ul{i}_\lambda{-}1$};
        % \pic at (0,0) {rcap=diag2}node[below right=-1pt,shift={(-.4,0)}]{\scriptsize $\ul{i}_\lambda{-}1$};
        \node at (1.3,.7) {\scriptsize $\underline{\lambda}$};
      \end{tikzpicture}}\endxy
      &&
      \xy (0,3)*{\begin{tikzpicture}[scale=.8*\SclFoam]
        \draw[diag2,->,shift={(.5,0)}] (-.5,0) to[out=90,in=180] (0,.6) to[out=0,in=90] (.5,0) node[below]{\scriptsize $\ul{i}_\lambda{-}1$};
        \draw[diag1,->] (1.5,0) to (1.5,1) node[above=-1pt]{\scriptsize $\ul{i}_\lambda{-}2$};
        \node at (2,.5) {\scriptsize $\underline{\lambda}$};
      \end{tikzpicture}}\endxy
      &&
      \xy (0,3)*{\begin{tikzpicture}[scale=.8*\SclFoam]
        \draw[diag2,->,shift={(.5,0)}] (-.5,0) to[out=90,in=180] (0,.6) to[out=0,in=90] (.5,0) node[below]{\scriptsize $\ul{i}_\lambda{-}1$};
        \draw[diag3,<-] (-.5,0) to (-.5,1) node[above=-1pt]{\scriptsize $\ul{i}_\lambda$};
        \node at (1.3,.7) {\scriptsize $\underline{\lambda}$};
      \end{tikzpicture}}\endxy
      &&
      \xy (0,5)*{\begin{tikzpicture}[scale=.8*\SclFoam]
        \draw[diag2,->,shift={(.5,0)}] (-.5,0) to[out=90,in=180] (0,.6) to[out=0,in=90] (.5,0) node[below]{\scriptsize $\ul{i}_\lambda{-}1$};
        \draw[diag3,<-] (-.5,0) to[out=90,in=-90] (1,2) node[above=-1pt]{\scriptsize $\ul{i}_\lambda$};
        \draw[diag1,->] (1.5,0) to[out=90,in=-90] (0,2) node[above=-1pt]{\scriptsize $\ul{i}_\lambda{-}2$};
        \node at (1.8,1) {\scriptsize $\underline{\lambda}$};
      \end{tikzpicture}}\endxy
      \\
      (\lambda_{i-1},\lambda_{i},\lambda_{i+1}) && (1,1,0) && (2,1,0) && (1,1,1) && (2,1,1)
      \\[\hsout]
      \xy (0,1)*{\begin{tikzpicture}[scale=1.1*\SclSch]
        \draw[schur1,<-] (0,0) node[below=-1pt]{\scriptsize $i-1$} to (1,1);
        \draw[schur2,<-] (1,0) node[below=-1pt]{\scriptsize $i$} to (0,1);
        \node at (1,.5) {\scriptsize $\lambda$};
      \end{tikzpicture}}\endxy
      &
      \overset{\raise.7em\hbox{\scalebox{1.1}{$\Gamma_\lambda(i)$}}}{\mapsto}
      &
      \xy (0,5)*{\begin{tikzpicture}[scale=.8*\SclFoam]
        \draw[diag2,<-,shift={(-.5,1)}] (-.5,0) node[above]{\scriptsize $\ul{i}_\lambda{-}1$} to[out=-90,in=180] (0,-.6) to[out=0,in=-90] (.5,0);
        \node at (.4,.3) {\scriptsize $\underline{\lambda}$};
      \end{tikzpicture}}\endxy
      &&
      \xy (0,3)*{\begin{tikzpicture}[scale=.8*\SclFoam]
        \draw[diag2,<-,shift={(.5,1)}] (-.5,0) node[above]{\scriptsize $\ul{i}_\lambda{-}1$} to[out=-90,in=180] (0,-.6) to[out=0,in=-90] (.5,0);
        \draw[diag1,->] (1.5,0) node[below=-1pt]{\scriptsize $\ul{i}_\lambda{-}2$} to (1.5,1);
        \node at (2,.4) {\scriptsize $\underline{\lambda}$};
      \end{tikzpicture}}\endxy
      &&
      \xy (0,3)*{\begin{tikzpicture}[scale=.8*\SclFoam]
        \draw[diag2,<-,shift={(.5,1)}] (-.5,0) node[above]{\scriptsize $\ul{i}_\lambda{-}1$} to[out=-90,in=180] (0,-.6) to[out=0,in=-90] (.5,0);
        \draw[diag3,<-] (-.5,0) node[below=-1pt]{\scriptsize $\ul{i}_\lambda$} to (-.5,1);
        \node at (1,.2) {\scriptsize $\underline{\lambda}$};
      \end{tikzpicture}}\endxy
      &&
      \xy (0,1)*{\begin{tikzpicture}[scale=.8*\SclFoam]
        \draw[diag2,<-,shift={(.5,1)}] (-.5,0) node[above]{\scriptsize $\ul{i}_\lambda{-}1$} to[out=-90,in=180] (0,-.6) to[out=0,in=-90] (.5,0);
        \draw[diag3,<-] (1,-1) node[below=-1pt]{\scriptsize $\ul{i}_\lambda$} to[out=90,in=-90] (-.5,1);
        \draw[diag1,->] (0,-1) node[below=-1pt]{\scriptsize $\ul{i}_\lambda{-}2$} to[out=90,in=-90] (1.5,1);
        \node at (1.8,0) {\scriptsize $\underline{\lambda}$};
      \end{tikzpicture}}\endxy
      \\
      (\lambda_{i-1},\lambda_{i},\lambda_{i+1}) && (1,1,0) && (2,1,0) && (1,1,1) && (2,1,1)
    \end{array}
    \\[\hsout]
    %%%%% CROSSING ii %%%%%
    \xy (0,0)*{\begin{tikzpicture}[scale=1.1*\SclSch]
      \draw[schur1,<-] (0,0) to (1,1);
      \draw[schur1,<-] (1,0) node[below=-1pt]{\scriptsize $i$} to (0,1);
      \node at (1,.5) {\scriptsize $\lambda$};
    \end{tikzpicture}}\endxy
    % \overset{\raise.7em\hbox{\scalebox{1.1}{$(-XY)\Gamma_\lambda(i)$}}}{\mapsto}
    \mapsto(-XY)\Gamma_\lambda(i)
    \xy (0,-2)*{\begin{tikzpicture}[scale=.8*\SclFoam]
      \pic[transform shape] at (0,0) {lcap=diag1}node[below=-1pt]{\scriptsize $\ul{i}_\lambda$};
      \pic[transform shape] at (0,.8) {rcup=diag1};
      \node at (1.3,.8) {\scriptsize $\underline{\lambda}$};
    \end{tikzpicture}}\endxy
    \mspace{30mu}
    %%%%% DISTANT CROSSING %%%%%
    \xy (0,-1)*{\begin{tikzpicture}[scale=1.1*\SclSch]
      \draw[schur1,<-] (0,0) node[below=-1pt]{\scriptsize $i$} to (1,1);
      \draw[schur2,<-] (1,0) node[below=-1pt]{\scriptsize $j$} to (0,1);
      \node at (1,.5) {\scriptsize $\lambda$};
    \end{tikzpicture}}\endxy
    \mspace{0mu}\mapsto\mspace{0mu}
    \xy (0,0)*{\begin{tikzpicture}[scale=1*\SclFoam]
      \draw[diag1,dashed,line width=3pt] (0,0) to (1,1);
      \draw[diag2,dashed,line width=3pt] (1,0) to (0,1);
      \node at (1.3,.5) {\scriptsize $\underline{\lambda}$};
    \end{tikzpicture}}\endxy
    \mspace{10mu}\text{ if $\abs{i-j}>1$}
  \end{gather*}
  
  \caption{
    Definition of $\cF_{n,d}$ on generating 2\nbd-morphisms. For dots, rightward cups and caps, and adjacent crossings, it depends on the local value of $\lambda$, which is given below each case.
    A symbol $\Gamma_\lambda(i)$ above a ``mapsto'' arrow means that the codomain should be multiplied by $\Gamma_\lambda(i)$.
    For distant crossings (last picture), the picture means that one should replace each strand of the Schur crossing with the corresponding strand or pair of strands as prescribed by \eqref{eq:defn_schur_to_foam_downward_strand}. This defines a foam diagram consisting of one, two or four crossings.
    %As foam crossings have trivial grading, one need not worry about their respective vertical positions (this happens in the case with four crossings, i.e.\ when $(\lambda_i,\lambda_i+1)=(\lambda_j,\lambda_j+1)=(2,1)$).
  }
  \label{fig:defn_foamation_functor}
\end{figure}

\begin{proof}
  One checks that $\cF_{n,d}$ preserves the $\bZ^2$-grading.
  We need to check that the images through $\cF_{n,d}$ of the defining relations of the graded 2-Schur algebra are relations in $\gfoam_d$. This is analogous to the proof of \cite[Proposition~3.3]{LQR_KhovanovHomologySkew_2015}. The main additional work is to check that the scalars match. For readability, we leave implicit the label of regions for foam diagrams.
  
  The fact that $\cF_{n,d}$ respects relation~\eqref{eq:dotnil} follows from dot annihilation in foams. For relation~\eqref{eq:klrR2}, the case $i=j$ follows from the evaluation of an undotted counterclockwise bubble and the case $\abs{i-j}>1$ follows from the Reidemeister II braid-like relation.
  Consider then the case $j=i-1$. For $\lambda_i=2$, both sides of~\eqref{eq:klrR2} are zero due to the Schur quotient. If $(\lambda_{i-1},\lambda_{i},\lambda_{i+1}) = (1,1,0)$, we have:
  \begin{align*}
    \cF_{n,d}\left(
    \tikz[scale=.8*.8,baseline={([yshift=.8ex]current bounding box.center)}]{
      \draw[schur1,<-]  +(0,-.75) node[below] {\textcolor{black}{\tiny $i$}}
        .. controls (0,-.375) and (1,-.375) .. (1,0)
        .. controls (1,.375) and (0, .375) .. (0,.75);
      \draw[schur2,<-]  +(1,-.75) node[below] {\textcolor{black}{\tiny $i-1$}}
        .. controls (1,-.375) and (0,-.375) .. (0,0)
        .. controls (0,.375) and (1, .375) .. (1,.75);
      \node at (1.3,0) {\tiny $\lambda$};
    }\right)
    \mspace{10mu}&=\mspace{10mu}
    \Gamma_\lambda(i)
    \xy (0,-2)*{\begin{tikzpicture}[scale=.7*.8]
      \pic[transform shape] at (0,0) {rcap=diag2}node[below=-1pt]{\scriptsize $\ul{i}_\lambda{-}1$};
      \pic[transform shape] at (0,.8) {lcup=diag2};
      % \node at (1.3,.8) {\scriptsize $\underline{\lambda}$};
    \end{tikzpicture}}\endxy
    \\
    &=
    \Gamma_\lambda(i)
    \left(
    \;YZ^2\;
    \xy(0,-2)*{\begin{tikzpicture}[scale=.7*.8,transform shape]
      \draw[diag2,->] (0,0) node[below] {\scriptsize $\ul{i}_\lambda{-}1$} to (0,1.8);
      \draw[diag2,<-] (1,0) to (1,1.8);
      \begin{scope}[shift={(-.7,.9)}]
        \node[fdot2=3pt] at (0,0) {};
        %\node at (.2,-.2) {\scriptsize $i$};
      \end{scope}
    \end{tikzpicture}}\endxy
    +XZ^2\;
    \xy(0,-2)*{\begin{tikzpicture}[scale=.7*.8,transform shape]
      \draw[diag2,->] (0,0) node[below] {\scriptsize $\ul{i}_\lambda{-}1$} to (0,1.8);
      \draw[diag2,<-] (1,0) to (1,1.8);
      \begin{scope}[shift={(1.7,.9)}]
        \node[fdot2=3pt] at (0,0) {};
        %\node at (.2,-.2) {\scriptsize $i$};
      \end{scope}
    \end{tikzpicture}}\endxy
    \right)
    \\
    &=
    \cF_{n,d}\left(YZ^2\;
    \tikz[scale=.8*.8,yscale=0.75,baseline={([yshift=.8ex]current bounding box.center)}]{
      \draw[schur1,<-] (0,-.75) node[below] {\textcolor{black}{\tiny $i$}}
        to node [sdot1]{} (0,.75);
      \draw[schur2,<-] (1,-.75) node[below] {\textcolor{black}{\tiny $i-1$}} to (1,.75);
      \node at (1.3,0) {\tiny $\lambda$};
    }
    -YZ^2\;
     \tikz[scale=.8*.8,yscale=0.75,baseline={([yshift=.8ex]current bounding box.center)}]{
      \draw[schur1,<-]  (0,-.75) node[below] {\textcolor{black}{\tiny $i$}} to (0,.75);
      \draw[schur2,<-]  (1,-.75) node[below] {\textcolor{black}{\tiny $i-1$}}
        to node[sdot2]{} (1,.75);
      \node at (1.3,0) {\tiny $\lambda$};
    }\right)
  \end{align*}
  where we used \cref{lem:more_rels_diagfoam} and $\Gamma_{\lambda-\alpha_{i-1}}(i)=\Gamma_{\lambda}(i)=(-XY)\Gamma_\lambda(i-1)$.
  Other cases for which $\lambda_i=1$ are computed similarly. If $\lambda_i=0$, the left-hand side of~\eqref{eq:klrR2} is zero due to the Schur quotient. If $(\lambda_{i-1},\lambda_{i},\lambda_{i+1}) = (2,0,1)$, the image of the right-hand side is
  \begin{gather*}
    \Gamma_{\lambda-\alpha_{i-1}}(i)\;YZ^2\;
    \xy (0,2)*{\begin{tikzpicture}[scale=1]
      \draw[diag2,<-] (-.7,0) to (-.7,1) node[above,shift={(-.3,0)}]{\scriptsize $\ul{i}_\lambda{+}1$};
      \draw[diag1,->] (0,0) to (0,1) node[above]{\scriptsize $\ul{i}_\lambda$};
      \node[fdot2] at (-.35,0.5) {};
      % \node at (.5,.5) {\scriptsize $\ul{\lambda}$};
    \end{tikzpicture}}\endxy
    -\Gamma_\lambda(i-1)YZ^2\;
    \xy (0,2)*{\begin{tikzpicture}[scale=1]
      \draw[diag2,<-] (-.7,0) to (-.7,1) node[above,shift={(-.3,0)}]{\scriptsize $\ul{i}_\lambda{+}1$};
      \draw[diag1,->] (0,0) to (0,1) node[above]{\scriptsize $\ul{i}_\lambda$};
      \node[fdot2] at (-.35,0.5) {};
      % \node at (.5,.5) {\scriptsize $\ul{\lambda}$};
    \end{tikzpicture}}\endxy
    =0,
  \end{gather*}
  which follows from $\Gamma_{\lambda-\alpha_{i-1}}(i)=\Gamma_\lambda(i-1)$.
  Other cases for which $\lambda_i=0$ are computed similarly.
  When $j=i+1$, both sides are zero due to the Schur quotient when $\lambda_j=0$, and only the left-hand side when $\lambda_j=2$. For the representative case $(\lambda_{j-1},\lambda_{j},\lambda_{j+1}) = (1,1,0)$, one gets:
  \begin{align*}
    \cF_{n,d}\left(
    \xy(0,0)*{\begin{tikzpicture}[scale=.8*.8,baseline={([yshift=.8ex]current bounding box.center)}]
      \draw[schur1,<-]  +(0,-.75) node[below] {\textcolor{black}{\tiny $j-1$}}
        .. controls (0,-.375) and (1,-.375) .. (1,0)
        .. controls (1,.375) and (0, .375) .. (0,.75);
      \draw[schur2,<-]  +(1,-.75) node[below] {\textcolor{black}{\tiny $j$}}
        .. controls (1,-.375) and (0,-.375) .. (0,0)
        .. controls (0,.375) and (1, .375) .. (1,.75);
      \node at (1.3,0) {\tiny $\lambda$};
    \end{tikzpicture}}\endxy\right)
    \mspace{10mu}&=\Gamma_\lambda(j)\mspace{10mu}
    \xy(0,0)*{\begin{tikzpicture}[scale=.6]
      \draw[diag2,
      decoration={markings, mark=at position 0 with {\arrow{<}}},
      postaction={decorate}]
      (0,0) circle (.5cm);
      \node at (1.5,-.3) {\scriptsize $l(\ul{j})-1$};
    \end{tikzpicture}}\endxy
    \\
    &=\Gamma_\lambda(j)
    \left(Z\mspace{10mu}
    \xy(0,-1.2)*{\begin{tikzpicture}
      \node[fdot2] at (0,0) {};
      \node at (.8,-.2) {\scriptsize $l(\ul{j})-1$};
    \end{tikzpicture}}\endxy
    \mspace{10mu}
    +\;XYZ
    \mspace{10mu}
    \xy(0,-1.2)*{\begin{tikzpicture}
      \node[fdot3] at (0,0) {};
      \node at (.5,-.2) {\scriptsize $l(\ul{j})$};
    \end{tikzpicture}}\endxy
    \right)
    \\
    &=
    \cF_{n,d}\left(
    -XYZ\xy(0,-1)*{\begin{tikzpicture}[scale=.8,yscale=0.75]
      \draw[schur1,<-] (0,-.75) node[below] {\textcolor{black}{\tiny $j-1$}}
        to node [sdot1]{} (0,.75);
      \draw[schur2,<-] (1,-.75) node[below] {\textcolor{black}{\tiny $j$}} to (1,.75);
      \node at (1.3,0) {\tiny $\lambda$};
    \end{tikzpicture}}\endxy
    +XYZ 
    \xy(0,-1)*{\begin{tikzpicture}[scale=.8,yscale=0.75]
      \draw[schur1,<-]  (0,-.75) node[below] {\textcolor{black}{\tiny $j-1$}} to (0,.75);
      \draw[schur2,<-]  (1,-.75) node[below] {\textcolor{black}{\tiny $j$}}
        to node[sdot2]{} (1,.75);
      \node at (1.3,0) {\tiny $\lambda$};
    \end{tikzpicture}}\endxy\right)
  \end{align*}
  which follows from $\Gamma_{\lambda-\alpha_{j}}(j-1)=(-XY)\Gamma_\lambda(j)$.
  For the representative case $(\lambda_{i-1},\lambda_{i},\lambda_{i+1}) = (1,2,0)$, the image of the right-hand side of \eqref{eq:klrR2} gives:
  \begin{gather*}
    -\Gamma_{\lambda-\alpha_{j}}(j-1)\;XYZ\;
    \xy (0,2)*{\begin{tikzpicture}[scale=1*.8]
      \draw[diag1,<-] (-.7,0) to (-.7,1) node[above,shift={(-.3,0)}]{\scriptsize $\ul{j}_\lambda-1$};
      \draw[diag2,->] (0,0) to (0,1) node[above]{\scriptsize $\ul{j}_\lambda$};
      \node[fdot1] at (-.35,0.5) {};
      % \node at (.5,.5) {\scriptsize $\ul{\lambda}$};
    \end{tikzpicture}}\endxy
    +\Gamma_\lambda(j)XYZ\;
    \xy (0,2)*{\begin{tikzpicture}[scale=1]
      \draw[diag1,<-] (-.7,0) to (-.7,1) node[above,shift={(-.3,0)}]{\scriptsize $\ul{j}_\lambda-1$};
      \draw[diag2,->] (0,0) to (0,1) node[above]{\scriptsize $\ul{j}_\lambda$};
      \node[fdot3] at (-.35,0.5) {};
      % \node at (.5,.5) {\scriptsize $\ul{\lambda}$};
    \end{tikzpicture}}\endxy
    =0,
  \end{gather*}
  which follows from $\Gamma_{\lambda-\alpha_{j}}(j-1)=\Gamma_\lambda(j)$ and dot migration.
  
  For relations~\eqref{eq:dotslides_first} and~\eqref{eq:dotslides_second},
  %, \eqref{eq:dotslides_second} and~\eqref{eq:dotjumpneib}
  the claim follows from the graded interchange law and dot migration.
  The relation~\eqref{eq:dotslide-nilH} follows directly from the (vertical) neck-cutting relation.

  Consider now relation~\eqref{eq:klrR3}. If all colours are pairwise equal or adjacent (i.e.\ cases $i=j=k$ and $i=j=k\pm 1$ together with permutations), then either the case is excluded by assumption, or both sides are zero by the Schur quotient. In particular, we can disregard the normalization with the $\Gamma$'s (this follows from the fact that $\Gamma_{\lambda-\alpha_j}(i)=\Gamma_\lambda(i)$ whenever $i\neq j$).
  If a colour is distant from the two others (i.e.\ cases $\abs{i-j}>1$ and $\abs{i-k}>1$ together with permutations), then the image under $\cF_{n,d}$ of relation~\eqref{eq:klrR3} is a composition of graded interchanges, consisting in moving vertically the image of the only crossing with a (possibly) non-trivial $\bZ^2$-grading.
  The last six cases also consist of graded interchanges, each interchanging two saddles (zip or unzip) sharing a common 1-facet. We picture below the domain and codomain of the six cases:
  \begin{gather*}
    \begin{tikzcd}[column sep=10em,ampersand replacement=\&]
      \xy(0,0)*{\begin{tikzpicture}[xscale=.7,yscale=.5,rotate=90,scale=.7]
        \draw[web2] (1,.5) to ++(0,2);
        \draw[web1] (1,0) to ++(0,.5);
        \draw[web1] (1,2.5) to ++(0,.5);
        \draw[web1] (2,0) to ++(0,1.5);
        \draw[web1,web0] (2,1.5) to ++(0,1);
        \draw[web1] (2,2.5) to ++(0,.5);
        \draw[web1] (0,.5) to ++(1,0);
        \draw[web1] (1,2.5) to ++(1,0);
        \draw[web1] (2,1.5) to ++(1,0);
        \draw[dashed] (0,0) to ++(0,3);
        \draw[dashed] (3,0) to ++(0,3);
      \end{tikzpicture}}\endxy
      \arrow[r,"{(i,j,k)=(l+1,l+2,l)}",shift left=2]
      \&
      \xy(0,0)*{\begin{tikzpicture}[xscale=.7,yscale=.5,rotate=90,scale=.7]
        \draw[web1,web0] (1,.5) to ++(0,2);
        \draw[web1] (1,0) to ++(0,.5);
        \draw[web1] (1,2.5) to ++(0,.5);
        \draw[web1] (2,0) to ++(0,.5);
        \draw[web2] (2,.5) to ++(0,1);
        \draw[web1] (2,1.5) to ++(0,1.5);
        \draw[web1] (0,2.5) to ++(1,0);
        \draw[web1] (1,0.5) to ++(1,0);
        \draw[web1] (2,1.5) to ++(1,0);
        \draw[dashed] (0,0) to ++(0,3);
        \draw[dashed] (3,0) to ++(0,3);
      \end{tikzpicture}}\endxy
      \arrow[l,"{(i,j,k)=(l,l+2,l+1)}",shift left=2]
    \end{tikzcd}
    \\[1ex]
    \begin{tikzcd}[column sep=10em,ampersand replacement=\&]
      \xy(0,0)*{\begin{tikzpicture}[xscale=.7,yscale=.5,rotate=90,scale=.7]
        \draw[web1] (1,0) to ++(0,.5);
        \draw[web1,web0] (1,.5) to ++(0,1);
        \draw[web1] (1,1.5) to ++(0,1.5);
        \draw[web1] (2,0) to ++(0,.5);
        \draw[web2] (2,.5) to ++(0,2);
        \draw[web1] (2,2.5) to ++(0,.5);
        \draw[web1] (0,1.5) to ++(1,0);
        \draw[web1] (1,0.5) to ++(1,0);
        \draw[web1] (2,2.5) to ++(1,0);
        \draw[dashed] (0,0) to ++(0,3);
        \draw[dashed] (3,0) to ++(0,3);
      \end{tikzpicture}}\endxy
      \arrow[r,"{(i,j,k)=(l+2,l,l+1)}",shift left=2]
      \&
      \xy(0,0)*{\begin{tikzpicture}[xscale=.7,yscale=.5,rotate=90,scale=.7]
        \draw[web1] (1,0) to ++(0,1.5);
        \draw[web2] (1,1.5) to ++(0,1);
        \draw[web1] (1,2.5) to ++(0,.5);
        \draw[web1] (2,0) to ++(0,.5);
        \draw[web1,web0] (2,.5) to ++(0,2);
        \draw[web1] (2,2.5) to ++(0,.5);
        \draw[web1] (0,1.5) to ++(1,0);
        \draw[web1] (1,2.5) to ++(1,0);
        \draw[web1] (2,0.5) to ++(1,0);
        \draw[dashed] (0,0) to ++(0,3);
        \draw[dashed] (3,0) to ++(0,3);
      \end{tikzpicture}}\endxy
      \arrow[l,"{(i,j,k)=(l+1,l,l+2)}",shift left=2]
    \end{tikzcd}
    \\[1ex]
    \begin{tikzcd}[column sep=10em,ampersand replacement=\&]
      \xy(0,0)*{\begin{tikzpicture}[xscale=.7,yscale=.5,rotate=90,scale=.7]
        \draw[web1] (1,0) to ++(0,.5);
        \draw[web2] (1,.5) to ++(0,1);
        \draw[web1] (1,1.5) to ++(0,1.5);
        \draw[web1] (2,0) to ++(0,1.5);
        \draw[web2] (2,1.5) to ++(0,1);
        \draw[web1] (2,2.5) to ++(0,.5);
        \draw[web1] (0,0.5) to ++(1,0);
        \draw[web1] (1,1.5) to ++(1,0);
        \draw[web1] (2,2.5) to ++(1,0);
        \draw[dashed] (0,0) to ++(0,3);
        \draw[dashed] (3,0) to ++(0,3);
      \end{tikzpicture}}\endxy
      \arrow[r,"{(i,j,k)=(l+2,l+1,l)}",shift left=2]
      \&
      \xy(0,0)*{\begin{tikzpicture}[xscale=.7,yscale=.5,rotate=90,scale=.7]
        \draw[web1] (1,0) to ++(0,1.5);
        \draw[web1,web0] (1,1.5) to ++(0,1);
        \draw[web1] (1,2.5) to ++(0,.5);
        \draw[web1] (2,0) to ++(0,.5);
        \draw[web1,web0] (2,.5) to ++(0,1);
        \draw[web1] (2,1.5) to ++(0,1.5);
        \draw[web1] (0,2.5) to ++(1,0);
        \draw[web1] (1,1.5) to ++(1,0);
        \draw[web1] (2,0.5) to ++(1,0);
        \draw[dashed] (0,0) to ++(0,3);
        \draw[dashed] (3,0) to ++(0,3);
      \end{tikzpicture}}\endxy
      \arrow[l,"{(i,j,k)=(l,l+1,l+2)}",shift left=2]
    \end{tikzcd}
  \end{gather*}
  This concludes the proof for the relation~\eqref{eq:klrR3}.

  A case-by-case analysis of the possible values for $(\lambda_i,\lambda_{i+1},\lambda_{i+2})$ shows that except for four possibilities, the three diagrams in \eqref{eq:R3serre_first_case} are all set to zero by the Schur quotient. In the remaining cases, exactly one summand on the left-hand side is set to zero by the Schur quotient.
  If $(\lambda_i,\lambda_{i+1},\lambda_{i+2})=(2,0,0)$, then
  \begin{align*}
    \cF_{n,d}\left(-YZ^{-2}
    \xy(0,0)*{\begin{tikzpicture}[scale=.7*.7]
      \draw[schur1,<-]  +(0,0)node[below] {\textcolor{black}{\tiny $i$}}
        .. controls (0,0.5) and (2, 1) ..  +(2,2);
      \draw[schur1,<-]  +(2,0)node[below] {\textcolor{black}{\tiny $i$}}
        .. controls (2,1) and (0, 1.5) ..  +(0,2);
      \draw[schur2,<-]  (1,0)node[below] {\textcolor{black}{\tiny $i+1$}}
        .. controls (1,0.5) and (0, 0.5) ..  (0,1)
        .. controls (0,1.5) and (1, 1.5) ..  (1,2);
      \node at (2.1,1) {\tiny $\lambda$};
    \end{tikzpicture}}\endxy\right)
    & = XZ^{-2}\;\;
    \begin{tikzpicture}[scale=.35*.8,baseline={([yshift=0ex]current bounding box.center)}]
      \pic[transform shape] at (0,0) {cup=diag1};
      \pic[transform shape] at (1,1) {cap=diag1};
      \pic[transform shape] at (0,3) {cap=diag1};
      \pic[transform shape] at (1,2) {cup=diag1};
      \draw[diag1] (0,1) to (0,3);
      \draw[diag1] (2,0) node[below]{\scriptsize $\ul{i}_\lambda$} to (2,1);
      \draw[diag1,->] (2,3) to (2,4);
    \end{tikzpicture}
    =
    \begin{tikzpicture}[scale=.35*.8,baseline={([yshift=0ex]current bounding box.center)}]
      \draw[diag1,->] (0,0) node[below]{\scriptsize $\ul{i}_\lambda$} to (0,4);
    \end{tikzpicture}
    =
    \cF_{n,d}\left(\;\tikz[scale=.7*.7,xscale=.8,,baseline={([yshift=0ex]current bounding box.center)}]{
      \draw[schur1,<-] (0,0) node[below] {\textcolor{black}{\tiny $i$}} to (0,2);
      \draw[schur2,<-] (1,0) node[below] {\textcolor{black}{\tiny $i+1$}} to (1,2);
      \draw[schur1,<-] (2,0) node[below] {\textcolor{black}{\tiny $i$}} to (2,2);
      \node at (2.3,1) {\tiny $\lambda$};
    }\right),
  \end{align*}
  using $\Gamma_{\lambda-\alpha_i}(i+1)=\Gamma_\lambda(i)$.
  The case $(\lambda_i,\lambda_{i+1},\lambda_{i+2})=(2,0,1)$ is similar. If $(\lambda_i,\lambda_{i+1},\lambda_{i+2})=(2,1,0)$, then
  \begin{align*}
    \cF_{n,d}\left(Z^{-1}
    \xy(0,0)*{\begin{tikzpicture}[scale=.7*.7]
      \draw[schur1,<-]  +(0,0)node[below] {\textcolor{black}{\tiny $i$}}
        .. controls (0,1) and (2, 1.5) ..  +(2,2);
      \draw[schur1,<-]  +(2,0)node[below] {\textcolor{black}{\tiny $i$}}
        .. controls (2,.5) and (0, 1) ..  +(0,2);
      \draw[schur2,<-]  (1,0)node[below] {\textcolor{black}{\tiny $i+1$}}
        .. controls (1,0.5) and (2, 0.5) ..  (2,1)
        .. controls (2,1.5) and (1, 1.5) ..  (1,2);
      \node at (2.3,1) {\tiny $\lambda$};
    \end{tikzpicture}}\endxy\right)
    & =Z^{-1}\mspace{-10mu}
    \xy(0,-2)*{\begin{tikzpicture}[scale=.6*.8]
      \draw[diag2,->] (0,0) node[below] {\scriptsize $\ul{i}_\lambda{+}1$} to[out=90,in=180] (.5,.75) to[out=0,in=90] (1,0);
      \draw[diag1,<-] (-.5,0) to[out=90,in=180] (.5,1.25) to[out=0,in=90] (1.5,0) node[below] {\scriptsize $\ul{i}_\lambda$};
      \begin{scope}[shift={(0,2)}]
        \draw[diag2,<-] (0,1) to[out=-90,in=180] (.5,.25) to[out=0,in=-90] (1,1);
        \draw[diag1,->] (-.5,1) to[out=-90,in=180] (.5,-.25) to[out=0,in=-90] (1.5,1);
      \end{scope}
    \end{tikzpicture}}\endxy
    =
    \xy(0,-2)*{\begin{tikzpicture}[scale=.6*.8]
      \draw[diag2,->] (0,0) node[below] {\scriptsize $\ul{i}_\lambda{+}1$} to (0,3);
      \draw[diag2,<-] (1,0) to (1,3);
      \draw[diag1,<-] (-.5,0) to (-.5,3);
      \draw[diag1,->] (1.5,0) node[below] {\scriptsize $\ul{i}_\lambda$} to (1.5,3);
    \end{tikzpicture}}\endxy
    =
    \cF_{n,d}\left(\;\tikz[scale=.7*.7,xscale=.8,baseline={([yshift=0ex]current bounding box.center)}]{
      \draw[schur1,<-] (0,0) node[below] {\textcolor{black}{\tiny $i$}} to (0,2);
      \draw[schur2,<-] (1,0) node[below] {\textcolor{black}{\tiny $i+1$}} to (1,2);
      \draw[schur1,<-] (2,0) node[below] {\textcolor{black}{\tiny $i$}} to (2,2);
      \node at (2.3,1) {\tiny $\lambda$};
    }\right),
  \end{align*}
  using $\Gamma_{\lambda-\alpha_{i+1}}(i)=(-XY)\Gamma_\lambda(i+1)$.
  The case $(\lambda_i,\lambda_{i+1},\lambda_{i+2})=(2,1,1)$ is similar.

  A similar analysis can be done for the relation \eqref{eq:R3serre_second_case}, giving four non-trivial cases, with only two computations needed. We depict respectively the cases $(\lambda_{i-1},\lambda_i,\lambda_{i+1})=(1,2,0)$ and $(\lambda_{i-1},\lambda_i,\lambda_{i+1})=(1,1,0)$:
  \begin{align*}
    \cF_{n,d}\left(
    XYZ^{-1}
    \xy(0,0)*{\begin{tikzpicture}[scale=.7*.7]
      \draw[schur1,<-]  +(0,0)node[below] {\textcolor{black}{\tiny $i$}}
        .. controls (0,0.5) and (2, 1) ..  +(2,2);
      \draw[schur1,<-]  +(2,0)node[below] {\textcolor{black}{\tiny $i$}}
        .. controls (2,1) and (0, 1.5) ..  +(0,2);
      \draw[schur2,<-]  (1,0)node[below] {\textcolor{black}{\tiny $i-1$}}
        .. controls (1,0.5) and (0, 0.5) ..  (0,1)
        .. controls (0,1.5) and (1, 1.5) ..  (1,2);
      \node at (2.1,1) {\tiny $\lambda$};
    \end{tikzpicture}}\endxy\right)
    & = XYZ^{-1}\mspace{-10mu}
    \xy(0,-2)*{\begin{tikzpicture}[scale=.5*.8]
      \draw[diag1=1.6,->] (0,0) node[below] {\scriptsize $\ul{i}_\lambda{-}1$} to[out=90,in=180] (.5,.75) to[out=0,in=90] (1,0);
      \draw[diag2=1.6,<-] (-.5,0) to[out=90,in=180] (.5,1.25) to[out=0,in=90] (1.5,0) node[below] {\scriptsize $\ul{i}_\lambda$};
      \begin{scope}[shift={(0,2)}]
        \draw[diag1=1.6,<-] (0,1) to[out=-90,in=180] (.5,.25) to[out=0,in=-90] (1,1);
        \draw[diag2=1.6,->] (-.5,1) to[out=-90,in=180] (.5,-.25) to[out=0,in=-90] (1.5,1);
      \end{scope}
    \end{tikzpicture}}\endxy
    \mspace{-7mu}=\mspace{-7mu}
    \xy(0,-2)*{\begin{tikzpicture}[scale=.5*.8]
      \draw[diag1=1.6,->] (0,0) node[below] {\scriptsize $\ul{i}_\lambda{-}1$} to (0,3);
      \draw[diag1=1.6,<-] (1,0) to (1,3);
      \draw[diag2=1.6,<-] (-.5,0) to (-.5,3);
      \draw[diag2=1.6,->] (1.5,0) node[below] {\scriptsize $\ul{i}_\lambda$} to (1.5,3);
    \end{tikzpicture}}\endxy
    \mspace{-7mu}=
    \cF_{n,d}\left(\;
    \xy(0,0)*{\begin{tikzpicture}[scale=.7*.7,xscale=.8,]
      \draw[schur1,<-] (0,0) node[below] {\textcolor{black}{\tiny $i$}} to (0,2);
      \draw[schur2,<-] (1,0) node[below] {\textcolor{black}{\tiny $i-1$}} to (1,2);
      \draw[schur1,<-] (2,0) node[below] {\textcolor{black}{\tiny $i$}} to (2,2);
      \node at (2.3,1) {\tiny $\lambda$};
    \end{tikzpicture}}\endxy\right)
    \\
    \cF_{n,d}\left(-XZ^{-2}
    \xy(0,0)*{\begin{tikzpicture}[scale=.7*.7]
      \draw[schur1,<-]  +(0,0)node[below] {\textcolor{black}{\tiny $i$}}
        .. controls (0,1) and (2, 1.5) ..  +(2,2);
      \draw[schur1,<-]  +(2,0)node[below] {\textcolor{black}{\tiny $i$}}
        .. controls (2,.5) and (0, 1) ..  +(0,2);
      \draw[schur2,<-]  (1,0)node[below] {\textcolor{black}{\tiny $i-1$}}
        .. controls (1,0.5) and (2, 0.5) ..  (2,1)
        .. controls (2,1.5) and (1, 1.5) ..  (1,2);
      \node at (2.3,1) {\tiny $\lambda$};
    \end{tikzpicture}}\endxy\right)
    & = YZ^{-2}\mspace{-10mu}
    \xy(0,-1)*{\begin{tikzpicture}[scale=-.35*.8]
      \pic[transform shape] at (0,0) {cup=diag1};
      \pic[transform shape] at (1,1) {cap=diag1};
      \pic[transform shape] at (0,3) {cap=diag1};
      \pic[transform shape] at (1,2) {cup=diag1};
      \draw[diag1] (0,1) to (0,3);
      \draw[diag1] (2,0) to (2,1);
      \draw[diag1,->] (2,3) to (2,4) node[below]{\scriptsize $\ul{i}_\lambda$};
    \end{tikzpicture}}\endxy
    =
    \xy(0,-1)*{\begin{tikzpicture}[scale=.35*.8]
      \draw[diag1,<-] (0,0) node[below]{\scriptsize $\ul{i}_\lambda$} to (0,4);
    \end{tikzpicture}}\endxy
    =
    \cF_{n,d}\left(\;
    \xy(0,0)*{\begin{tikzpicture}[scale=.7*.7,xscale=.8,]
      \draw[schur1,<-] (0,0) node[below] {\textcolor{black}{\tiny $i$}} to (0,2);
      \draw[schur2,<-] (1,0) node[below] {\textcolor{black}{\tiny $i-1$}} to (1,2);
      \draw[schur1,<-] (2,0) node[below] {\textcolor{black}{\tiny $i$}} to (2,2);
      \node at (2.3,1) {\tiny $\lambda$};
    \end{tikzpicture}}\endxy\right)
  \end{align*}
  using \cref{lem:more_rels_diagfoam}, and $\Gamma_{\lambda-\alpha_{i}}(i)=(-XY)\Gamma_\lambda(i)$ and $\Gamma_{\lambda-\alpha_{i-1}}(i)=\Gamma_\lambda(i)$, respectively.
  For relation~\eqref{eq:adjrels}, the claim follows from the zigzag relations for foams:
  \begin{equation*}
    \begin{array}{*{2}{c@{\hskip 1ex}}*{6}{c@{\hskip 1ex}}c}
      \xy (0,3)*{\tikz[scale=1*.7]{
        \draw[schur1,<-] (0,-.5) to (0,0) to[out=90,in=180] (.25,.5)
          to[out=0,in=90] (.5,0) to[out=-90,in=180] (.75,-.5)
          to[out=0,in=-90] (1,0) to (1,.5) node[left=-3pt]{\tiny $i$};
        \node at (1.2,0) {\tiny $\lambda$};
      }}\endxy
      &\mapsto&
      \xy (0,1)*{\tikz[scale=1]{
        \node at (0,0) {$\ul{\lambda}$};
      }}\endxy
      &&
      Z^{-2}\;
      \xy (0,4)*{\tikz[scale=.5*.7]{
        \pic[transform shape] at (1,0) {cup=diag1};
        \pic[transform shape] at (0,1) {cap=diag1};
        \draw[diag1] (0,0) to (0,1);
        \draw[diag1,->] (2,1) to (2,2) node[above]{\scriptsize $\ul{i}_\lambda$};
      }}\endxy
      &&
      \xy (0,4)*{\tikz[scale=.5*.7]{
        \pic[transform shape] at (1,0) {cup=diag1};
        \pic[transform shape] at (0,1) {cap=diag1};
        \draw[diag1,<-] (0,0) to (0,1);
        \draw[diag1] (2,1) to (2,2) node[above]{\scriptsize $\ul{i}_\lambda$};
      }}\endxy
      &&
      Z^{-1}\mspace{-5mu}
      \xy (0,4)*{\tikz[scale=.3*.7]{
        \draw[diag1,->] (1,0) to ++(0,2) to[out=90,in=180] ++(.5,.5) to[out=0,in=90] ++(.5,-.5) to[out=-90,in=180] ++(1.5,-1.5) to[out=0,in=-90] ++(1.5,1.5) to ++(0,2) node[above=-1pt,shift={(.3,0)}]{\scriptsize $\ul{i}_\lambda$};
        \draw[diag2,<-] (0,0) to ++(0,2) to[out=90,in=180] ++(1.5,1.5) to[out=0,in=90] ++(1.5,-1.5) to[out=-90,in=180] ++(.5,-.5) to[out=0,in=-90] ++(.5,.5) to ++(0,2) node[above=-1pt,shift={(-.3,0)}]{\scriptsize $\ul{i}_\lambda{+}1$};
      }}\endxy
      \\[0ex]
      && (1,0) && (2,0) && (1,1) && (2,1)
      \\[0ex]
      \xy (0,0)*{\tikz[scale=1*.7]{
        \draw[schur1,<-] (0,.5) to (0,0) to[out=-90,in=180] (.25,-.5)
          to[out=0,in=-90] (.5,0) to[out=90,in=180] (.75,.5)
          to[out=0,in=90] (1,0) to (1,-.5) node[left=-3pt]{\tiny $i$};
        \node at (1.2,0) {\tiny $\lambda$};
      }}\endxy
      &\mapsto&
      \xy (0,1)*{\tikz[scale=1*.7]{
        \node at (0,0) {$\ul{\lambda}$};
      }}\endxy
      &&
      \xy (0,3)*{\tikz[scale=.5*.7,yscale=-1]{
        \pic[transform shape] at (1,0) {cup=diag1};
        \pic[transform shape] at (0,1) {cap=diag1};
        \draw[diag1,<-] (0,0) node[above]{\scriptsize $\ul{i}_\lambda$} to (0,1);
        \draw[diag1] (2,1) to (2,2);
      }}\endxy
      &&
      Z^{-2}\mspace{-5mu}
      \xy (0,3)*{\tikz[scale=.5*.7,yscale=-1]{
        \pic[transform shape] at (1,0) {cup=diag1};
        \pic[transform shape] at (0,1) {cap=diag1};
        \draw[diag1] (0,0) node[above]{\scriptsize $\ul{i}_\lambda$} to (0,1);
        \draw[diag1,->] (2,1) to (2,2);
      }}\endxy
      &&
      Z^{-1}\mspace{-30mu}
      \xy (0,3)*{\tikz[scale=.3*.7,yscale=-1]{
        \draw[diag1,<-] (1,0) node[above=-1pt,shift={(.3,0)}]{\scriptsize $\ul{i}_\lambda$} to ++(0,2) to[out=90,in=180] ++(.5,.5) to[out=0,in=90] ++(.5,-.5) to[out=-90,in=180] ++(1.5,-1.5) to[out=0,in=-90] ++(1.5,1.5) to ++(0,2);
        \draw[diag2,->] (0,0) node[above=-1pt,shift={(-.3,0)}]{\scriptsize $\ul{i}_\lambda{+}1$} to ++(0,2) to[out=90,in=180] ++(1.5,1.5) to[out=0,in=90] ++(1.5,-1.5) to[out=-90,in=180] ++(.5,-.5) to[out=0,in=-90] ++(.5,.5) to ++(0,2);
      }}\endxy
      \\[-1ex]
      && (0,1) && (0,2) && (1,1) && (1,2)
    \end{array}
  \end{equation*}
  % The scalars are $1$, $Z^2$, $1$ and $Z$ in the first case, and $1$, $X$, $YZ^2$ and $XYZ$ in the second case, before normalizing with $Z$.

  It remains to check that $\cF_{n,d}$ preserves the three inverse axioms~\eqref{eq:iso-ij}, \eqref{eq:isoii-pos} and \eqref{eq:isoii-neg}. Note that this is independent of the coefficients in the defining relations of $\gfoam_d$, as those inverse axioms only impose existence.
  For relation \eqref{eq:iso-ij}, the claim follows from graded isotopies and zigzag relations. For $\ov{\lambda}_i=0$, one checks that the image of the leftward monochromatic crossing is the identity, up to multiplication by an invertible scalar. Otherwise, it is zero by the Schur quotient. For $\ov{\lambda}_i=\pm 1$, the image of the leftward cup (resp.\ cap) is invertible due to the squeezing relation. Finally, for $\ov{\lambda}_i=\pm 2$ this is precisely the (vertical) neck-cutting relation.
\end{proof}

In \cite{Vaz_NotEvenKhovanov_2020}, the second author defined a ``thick calculus'' for the negative half of the super version of the graded 2-Schur $\catschur_{n,d}$. This amounts to working in a subcategory of the Karoubi envelope. One can similarly define a thick calculus for the full graded 2-Schur, defining a graded-2\nbd-cate\-gory $\check{\catschur}_{n,d}$.
Then, following the line of the proof of the analogous result in the non-graded case as given by Queffelec and Rose \cite[Theorem~3.9]{QR_$mathfraksl_n$Foam2category_2016}, one can show that the graded foamation 2-functor factors through $\check{\catschur}_{n,d}$.
The inclusion
\begin{equation*}
  \check{\catschur}_{n,d}\hookrightarrow\check{\catschur}_{n+1,d}
\end{equation*}
is defined on objects as $(\lambda_1,\ldots,\lambda_n)\mapsto (\lambda_1,\ldots,\lambda_n,0)$ and similarly for 1\nbd-morphisms and 2\nbd-mor\-phisms.
Finally, following the line of proof of \cite[Proposition~3.22]{QR_$mathfraksl_n$Foam2category_2016}\footnote{Or rather, the analogue of this proposition when one imposes \eqref{eq:dotnil} and dot annihilation respectively.} in the $\glt$-case, one can show the following proposition. We leave the details to the reader.

\begin{proposition}
  We have the following equivalence of $(\bZ^2,\bilfoam)$-graded-2-categories:
\begin{equation*}
  \colim_{n\in\bN}\left(\ldots\hookrightarrow\check{\catschur}_{n,d}\hookrightarrow\check{\catschur}_{n+1,d}\hookrightarrow\ldots\right) \cong \gfoam_d.
\end{equation*}
\end{proposition}

%%%%%%%%%%%%%%%%%%%%%%%%%%%%%%%%%%%%%%%%%%
%%%          categorification          %%%
%%%%%%%%%%%%%%%%%%%%%%%%%%%%%%%%%%%%%%%%%%
\subsection{The categorification theorem}
\label{subsec:schur_categorification}

Recall the notion of quantum Grothendieck ring from \cref{subsec:more_cat_facts}.
As claimed, the graded 2-Schur algebra categorifies the Schur algebra of level 2:

\begin{theorem}
  \label{thm:schur_categorification}
  The graded 2-Schur algebra categorifies the $q$-Schur algebra of level 2:
  \begin{equation*}
    K_0(\catschur_{n,d})\vert_q\cong \qschur_{n,d},
  \end{equation*}
  where the isomorphism is an isomorphism of $\bZ[q,q^{-1}]$-linear categories.
\end{theorem}

Note that \cref{thm:schur_categorification} relies on \cref{thm:foam_categorification} and hence on \cref{thm:basis_foam}, whose proof is given in \cite{Schelstraete_RewritingModuloDiagrammatic_2026}.

\begin{proof}[Proof of \cref{thm:schur_categorification}]
  Relations \eqref{eq:quantum_algebra_rel} become 2-isomorphisms in $\catschur_{n,d}$. The first relation is categorified by the invertibility axioms \eqref{eq:iso-ij}, \eqref{eq:isoii-pos} and \eqref{eq:isoii-neg}. The second relation is categorified by \eqref{eq:klrR2} in the case $\abs{i-j}>1$. Finally, the third relation is categorified by the analogue of \eqref{eq:klrR2}, case $\abs{i-j}>1$, for upward strands.
  We define the \emph{upward crossing} as:
  \begin{equation*}
    %%%% upward crossing %%%%
    \xy (0,-.5)*{\begin{tikzpicture}[scale=.3]
      \draw [schur1,<-] (1,1) to (-1,-1) node[left=-3pt]{\tiny $i$};
      \draw [schur2,<-] (-1,1) to (1,-1) node[right=-3pt]{\tiny $j$};
      \node at (1.4,0) {\tiny $\lambda$};
    \end{tikzpicture}}\endxy
    \coloneqq
    \xy (0,-.5)*{\begin{tikzpicture}[scale=.6]
      \draw[schur2,<-] (-1.3,1.5) to (-1.3,.5) to ++(0,-1) to[out=-90,in=180] ++(.8,-.6) to[out=0,in=-90] (0.3,-0.5) to[out=90,in=-90] (-0.3,0.5) to[out=90,in=180] (.5,1.1) to[out=0,in=90] (1.3,.5) to (1.3,-1.5);
      \draw[schur1,<-] (-0.8,1.5) to (-0.8,0.5) to [out=-90,in=180] (-0.5,-0.5) to [out=0,in=230] (-0.2,-0.3) to [out=50,in=180] (0.5,0.5) to [out=0,in=90] (0.8,-0.5) to (0.8,-1.5);
      \node at (.9,-1.5) {\tiny $i$};
      \node at (1.4,-1.5) {\tiny $j$};
      \node at (1.55,0) {\tiny $\lambda$};
    \end{tikzpicture}}\endxy
  \end{equation*}
  It then follows from adjunction relations \eqref{eq:adjrels} and \eqref{eq:klrR2}, case $\abs{i-j}>1$, that:
  \begin{equation*}
    \xy(0,-.8)*{\begin{tikzpicture}[scale=.6]
      \draw[schur1,->]  +(0,-.75) node[below] {\textcolor{black}{\tiny $i$}}
      .. controls (0,-.375) and (1,-.375) .. (1,0)
      .. controls (1,.375) and (0, .375) .. (0,.75);
      \draw[schur2,->]  +(1,-.75) node[below] {\textcolor{black}{\tiny $j$}}
        .. controls (1,-.375) and (0,-.375) .. (0,0)
        .. controls (0,.375) and (1, .375) .. (1,.75);
      \node at (1.3,0) {\tiny $\lambda$};
    \end{tikzpicture}}\endxy
    \mspace{10mu}=\mspace{10mu}
    \xy(0,-.8)*{\begin{tikzpicture}[scale=.6]
      \draw[schur1,->] (0,-.75) node[below] {\textcolor{black}{\tiny $i$}} to (0,.75);
      \draw[schur2,->] (1,-.75) node[below] {\textcolor{black}{\tiny $j$}} to (1,.75);
      \node at (1.3,0) {\tiny $\lambda$};
    \end{tikzpicture}}\endxy
    \qquad \text{ if }\vert i-j \vert > 1.
  \end{equation*}
  This implies that there exists a $\bZ[q,q^{-1}]$-linear functor $\qschur_{n,d}\rightarrow K_0(\catschur_{n,d})\vert_q$, full and surjective, fitting into the following commutative diagram:
  \begin{equation*}
    \begin{tikzcd}
      \qschur_{n,d} & K_0(\catschur_{n,d})\vert_q
      \\
      \web_d & K_0(\gfoam_d)\vert_q
      \arrow[from=1-1,to=1-2]
      \arrow[from=1-1,to=2-1,"F_{n,d}"']
      \arrow[from=1-2,to=2-2,"K_0(\cF_{n,d})"]
      \arrow[from=2-1,to=2-2,"\cong"]
      \arrow[from=1-1,to=2-2,phantom,"\circlearrowright"]
    \end{tikzcd}
  \end{equation*}
  The bottom arrow is an isomorphism thanks to \cref{thm:foam_categorification} and the left arrow is faithful thanks to \cref{lem:faithful_schur_web}. We conclude that the top arrow is also faithful, and hence it is an isomorphism.
\end{proof}

% The proof relies on the graded foamation 2-functor and the analogous result for $\gfoam_d$ (\cref{thm:foam_categorification}), and the following lemma:

% However, in the example on the left for instance, there is another obvious way to show that the two diagrams are equal, without further embeddings: evaluate the circles to $q+q^{-1}$. A similar strategy works in the second case, and in fact we show in the proof below that the ability the evaluate circles is enough to capture all isotopies.

\section{Chain complexes in graded-monoidal categories}
\label{sec:complexes}
This section introduces a notion of tensor product for chain complexes in a given graded-monoidal category (\cref{defn:graded_monoidal_cat}); more precisely, we restrict our study to \emph{homogeneous polycomplexes}. The special case of super-2-categories first appeared in the first author's master's thesis \cite{Schelstraete_SupercategorificationKhovanovlikeTangle_2020}.
Crucially, we then show that this tensor product leaves homotopy classes invariant. Although we work in the context of graded-monoidal categories for simplicity, all definitions and results extend in a straightforward way to graded-2-categories.

This section is not conceptually difficult, but it is technical.
We refer the reader to \cref{subsec:hypercubic_chain_complexes} for a minimalistic version, sufficient for the purpose of \cref{sec:topo}.

\begin{notation}
  \label{nota:complexes_basics}
  Fix $n\in\bN$ and $\cI\coloneqq\{i_1<\ldots<i_n\}$ a set of ordered indices, isomorphic to $\{1<\ldots<n\}$ as an ordered set.
  We view $\bZ^\cI\cong\bZ^n$ as an $n$-dimensional lattice and write elements $\vr\in\bZ^\cI$ as ordered tuples $\vr=(r_{i_1},\ldots,r_{i_n})$; moreover, we abbreviate $r_{i_k}$ as $r_k$, so that $\vr=(r_1,\ldots,r_n)$.
  We write $\abs{\vr}\coloneqq r_1+\ldots+r_n$ and denote $(\ve_i)_{i\in\cI}$ the canonical basis of $\bZ^\cI\cong\bZ^n$:
  \[\ve_i\coloneqq(\underset{1}{0},\ldots,\underset{i-1}{0},\underset{i}{1},\underset{i+1}{0},\ldots,\underset{n}{0}).\]
  Hereabove the underlying numbers denote coordinates.
  For each ${\vr\in\{0,1\}^\cI}$ and each $i\in\cI$, we write $\vr\to\vr+\ve_i$ the corresponding edge in the lattice $\bZ^\cI$.
\end{notation}

\begin{notation}
\label{nota:complexes_graded_str}
  Fix $\ring$, $G$ and $\bil$ as in \cref{subsec:graded_two_cat}, and fix a $(G,\bil)$-graded-monoidal ca\-te\-gory~$\cC$.
  To reduce clutter, we often abuse notation and write $f$ instead of $\deg f$, where $f$ is a morphism of $\cC$. The distinction should be clear by the context. We also write $*$ instead of $\bil$. For instance, for $f$, $g$ and $h$ morphisms, we may write $(f+g)* h$ for $\bil(\deg f+\deg g,\deg h)$. We sometimes use the subscripted equal signs $=_\ring$ and $=_G$ to emphasize where equality holds.
\end{notation}

\begin{quote}
  \textbf{Assumption.} To simplify this already technical section, we assume that $\bil$ is symmetric (\cref{defn:bil_symmetric}) throughout.
\end{quote}

%%%%%%%%%%%%%%%%%%%%%%%%%%%%%%%%%
%%%          summary          %%%
%%%%%%%%%%%%%%%%%%%%%%%%%%%%%%%%%
\subsection{A graded Koszul rule}

Recall that for a graph such as $\bZ^\cI$, a \emph{cycle} is an oriented loop of edges, and a 1-cochain on a graph is said to be a \emph{1-cocycle} if it is zero on all cycles. For $\bZ^\cI$, a 1-cochain is a 1-cocycle if and only if it is zero on every square.
A graph 1-cocycle is always the boundary of a graph 0-cochain.

Recall that if $V=\oplus_{g\in G} V_g$ is a $G$\nbd-graded $\ring$-module, an element $v\in V$ is said to be \emph{homogeneous} if $v\in V_g$ for some $g\in G$. If moreover $v\neq 0$, then $v$ has a well-defined degree $\deg v =g$.

\begin{definition}
\label{defn:Hcomplex}
  A \emph{homogeneous $n$-polycomplex}
  % (or simply \emph{$H$-complex})
  $\bA=(A,\alpha,\psi_\bA)$ is the data of:
  \begin{enumerate}[label=(\roman*)]
    \item a family $A\coloneqq (A^\vr)_{\vr\in\bZ^\cI}$ of objects $A^\vr\in\cC$%
    %, such that $A^{\vr}=0$ for $\abs{\vr}$ both small and big enough
    ,

    \item a family $\alpha\coloneqq(\alpha_i^{\vr})_{\vr\in\bZ^\cI,i\in\cI}$ of \emph{homogeneous} morphisms $\alpha_i^{\vr} \colon A^{\vr}\to A^{\vr+e_i}$, such that ${\alpha_i^{\vr+\ve_i}\circ\alpha_i^\vr=0}$ for all $i\in\cI$, and such that each square anti-commutes:
    \[
      \alpha_j^{\vr+\ve_i}\circ\alpha_i^{\vr}
      =-\alpha_i^{\vr+\ve_j}\circ\alpha_j^{\vr},
    \]

    \item a $G$-valued 1-cocycle $\psi_\bA$ on $\bZ^\cI$ such that $\psi_\bA(\vr\to\vr+\ve_i)= \deg{\alpha_i^{\vr}}$ whenever $\alpha_i^{\vr}\neq 0$.
  \end{enumerate}
\end{definition}

If a square in $\bA$ is non-zero in the sense that
\[
  \alpha_j^{\vr+\ve_i}\circ\alpha_i^{\vr}
  =-\alpha_i^{\vr+\ve_j}\circ\alpha_j^{\vr}\neq 0,
\]
then each of the four maps involved is non-zero, hence has a well-defined degree, and the following condition holds:
\begin{equation}
  \label{eq:polyhomogeneous_cocycle_cond}
  \deg\alpha_j^{\vr+\ve_i}+\deg\alpha_i^{\vr}
  =\deg\alpha_i^{\vr+\ve_j}+\deg\alpha_j^{\vr}.
\end{equation}
In other words, the partially-defined 1-cochain $\deg(-)$ is a 1-cocycle on non-zero squares.
Condition (iii) states the existence of a 1-cochain $\psi_\bA$ that both extends $\deg(-)$ to zero maps and the 1-cocycle condition to zero squares.
In particular, a square could be zero albeit having four non-zero edges: in which case, although the equation \eqref{eq:polyhomogeneous_cocycle_cond} makes sense, it may not hold. Condition (iii) ensures that it does.

\Cref{nota:complexes_basics} and \cref{defn:Hcomplex} introduced standard notations for a generic homogeneous polycomplex $\bA$. To continue our discussion, we give standard notations for another generic homogeneous polycomplex $\bB$.

\begin{notation}
  \label{nota:complex_B}
  Let $m\in\bN$ and $\cJ={\{j_1<\ldots<j_m\}}$ a set of ordered indices isomorphic to ${\{1<\ldots<m\}}$.
  We use generically the letter $j$ for an index $j\in\cJ$, and denote $\vs=(s_{j_1},\ldots,s_{j_m})=(s_1,\ldots,s_m)$ a vertex in the lattice $\bZ^\cJ$.
  Furthermore, we write
  \begin{gather*}
    \cI\sqcup\cJ\coloneqq\{i_1<\ldots<i_n<j_1<\ldots<j_m\}.
  \end{gather*}
  We generically use the letter $k$ for an index $k\in\cI\sqcup\cJ$ and the notation $(\vr,\vs)$ for a vertex in $\bZ^{\cI\sqcup\cJ}$.
  We let $\bB=(B,\beta,\psi_\bB)$ be a homogeneous $m$-polycomplex.
\end{notation}

Thanks to condition (ii), we associate to $\bA$ a chain complex, its \emph{total complex}, denoted as ${(\Tot(\bA),\Tot(\alpha))}$ and given by
\begin{equation*}
  \Tot(\bA)_t \coloneqq \bigoplus_{\vr\in\bZ^\cI,\abs{\vr}=t} A^{\vr}
  \quad\an\quad
  \Tot(\alpha)\vert_{A^\vr} \coloneqq \sum_{1\leq i\leq n} \alpha_i^{\vr}.
\end{equation*}
Let $\bB$ be as in \cref{nota:complex_B} and let $F\coloneqq\left(F^{\vr,\vs}\colon A^\vr\to B^\vs\right)_{\vr\in\bZ^\cI,\vs\in\bZ^\cJ}$ be a family of morphisms in $\cC$.
Similarly to the above, it gathers as a morphism $\Tot(F)\colon\Tot(\bA)\to\Tot(\bB)$.

\begin{definition}
  \label{defn:tensor_product_chain_map_chain_homotopies}
  Let $\bA=(A,\alpha,\psi_\bA)$ and $\bB=(B,\beta,\psi_\bB)$ be homogeneous polycomplexes.
  A family of morphisms in $\cC$
  \[F\coloneqq\left(F^{\vr,\vs}\colon A^\vr\to B^\vs\right)_{\vr\in\bZ^\cI,\vs\in\bZ^\cJ}\]
  is a \emph{chain morphism between $\bA$ and $\bB$} if $\Tot(F)$ is a chain morphism between $\Tot(\bA)$ and $\Tot(\bB)$.

  Let $F,G\colon\bA\to\bB$ be two chain morphisms between $\bA$ and $\bB$.
  A family of morphisms in $\cC$
  \[H\coloneqq\left(H^{\vr,\vs}\colon A^\vr\to B^\vs\right)_{\vr\in\bZ^\cI,\vs\in\bZ^\cJ}\]
  is a \emph{chain homotopy between $F$ and $G$} if $\Tot(H)$ is a chain homotopy between $\Tot(F)$ and $\Tot(G)$.
\end{definition}

In other words, while we use a restricted notion of complexes (namely, homogeneous polycomplexes), chain maps and chain homotopies between them are the usual notions of chain maps and chain homotopies.
In particular, these chain maps and chain homotopies have no homogeneity conditions.
Note also that if $F$ is a chain map (resp.\ $H$ is a chain homotopy), then $F^{\vr,\vs}=0$ whenever $\abs{\vr}\neq\abs{\vs}$ (resp.\ $H^{\vr,\vs}=0$ whenever $\abs{\vr}\neq\abs{\vs}+1$).

\medbreak

We now introduce the notion of tensor product of homogeneous polycomplexes, after setting some further notations.

For each vertex $\vr\in\bZ^\cI$, let $p$ be a path in $\bZ^\cI$ from $\mathbf{0}$ to $\vr$. Since $\psi_\bA$ is a cocycle, the value ${\abs{\alpha}(\vr)\coloneqq\psi_\bA(p)}$ does not depend on the choice of path. Most importantly, it verifies the following:
\begin{equation}
\label{eq:abs_alpha}
  \abs{\alpha}(\vr+e_i)-\abs{\alpha}(\vr)
  =_G \psi_\bA(\vr\to\vr+\ve_i)
  =_G \deg{\alpha_i^{\vr}},
\end{equation}
where the last equality holds whenever it makes sense, that is, whenever $\alpha_i^{\vr}\neq 0$.

\begin{definition}
\label{defn:sigma_tensor_product}
  Let $\epsilon$ be a $\ring^\times$-valued 1-cochain on $\bZ^{\cI\sqcup\cJ}$. With the notations above, the \emph{$\epsilon$-tensor product} of $\bA$ and $\bB$, denoted $(\bA\otimes \bB)(\epsilon)$, is the following data:
  \begin{enumerate}[label=(\roman*)]
    \item a family $A\otimes B=((A\otimes B)^{(\vr,\vs)})$ of objects $(A\otimes B)^{(\vr,\vs)} \coloneqq A^\vr\otimes B^\vs$,

    \item a family $\alpha\otimes\beta=((\alpha\otimes\beta)^{(\vr,\vs)}_k)$ of homogeneous morphisms
    \[
      (\alpha\otimes\beta)^{(\vr,\vs)}_k \coloneqq
      \begin{cases}
        \epsilon(e)(\alpha^{\vr}_i\otimes\id_{B^{\vs}})
        & k=i\in\cI,
        \\
        \epsilon(e)(\id_{A^{\vr}}\otimes\beta^{\vs}_{j})
        \quad& k=j\in\cJ,
      \end{cases}
    \]
    where $e$ denotes the edge $(\vr,\vs)\to(\vr,\vs)+\ve_k$,

    \item a 1-cocycle $\psi_\bA\otimes\psi_\bB$ on $\bZ^{\cI\sqcup\cJ}$ given on $e=(\vr,\vs)\to(\vr,\vs)+\ve_k$ by
    \[
      (\psi_\bA\otimes\psi_\bB)(e)
      \coloneqq
      \begin{cases}
        \psi_\bA(\vr\to\vr+\ve_i)
        & k=i\in\cI,
        \\
        \psi_\bB(\vs\to\vs+\ve_j)
        \quad& k=j\in\cJ.
      \end{cases}
    \]
  \end{enumerate}
\end{definition}

For $(\bA\otimes\bB)(\epsilon)$ to define a homogeneous $(n+m)$-polycomplex, $\epsilon$ needs to be such that squares anti-commute. This is encapsulated in the following lemma, where $\square^{(\vr,\vs)}_{k,l}$ with $k<l$ denotes the following oriented square in the lattice $\bZ^{\cI\sqcup\cJ}$:
\begin{equation*}
  \begin{tikzcd}
    (\vr,\vs) 
    \drar["\circlearrowright",phantom]
    \rar \dar & (\vr,\vs)+\ve_k \dar
    \\
    (\vr,\vs)+\ve_l \rar & (\vr,\vs)+\ve_k+\ve_l
  \end{tikzcd}
\end{equation*}

\begin{lemma}
  \label{lem:compatible_implies_complex}
  Say that a 1-cochain $\epsilon$ on $\bZ^{\cI\sqcup\cJ}$ is \emph{compatible} if
  \begin{IEEEeqnarray*}{c}
    \partial\epsilon(\square^{(\vr,\vs)}_{k,l})=
    % \\[1ex]
    % \mspace{20mu}
    \begin{cases}
      0 & \text{ if }k,l\in\cI\text{ or }k,l\in\cJ,
      \\
      -\psi_\bA(\vr\to\vr+\ve_i)*\psi_\bB(\vs\to\vs+\ve_j)
      & \text{ if }k=i\in\cI\text{ and }l=j\in\cJ.
    \end{cases}
  \end{IEEEeqnarray*}
  If $\epsilon$ is compatible, then $(\bA\otimes\bB)(\epsilon)$ defines a homogeneous $(n+m)$-polycomplex.
\end{lemma}

\begin{proof}
  The case $k,l\in\cI\text{ or }k,l\in\cJ$ is clear. In the case $k=i\in\cI$ and ${l=j\in\cJ}$, the square $\square^{(\vr,\vs)}_{i,j}$, decorated by the data of $\bA\otimes\bB$, has the following form:
  \begin{equation*}
    \begin{tikzcd}
      A^\vr\otimes B^\vs 
      \drar["\circlearrowright",phantom]
      \rar["\alpha^{\vr}_i\otimes\id_{B^{\vs}}"] 
      \dar["\id_{A^{\vr}}\otimes\beta^{\vs}_{j}"'] 
      & A^{\vr+\ve_i}\otimes B^\vs
      \dar["\id_{A^{\vr+\ve_i}}\otimes\beta^{\vs}_{j}"]
      \\
      A^\vr\otimes B^{\vs+\ve_j} 
      \rar["\alpha^{\vr}_i\otimes\id_{B^{\vs+\ve_j}}"',shift right=1] 
      & A^{\vr+\ve_i}\otimes B^{\vs+\ve_j} 
    \end{tikzcd}
  \end{equation*}
  The morphism corresponding to the path
  \[(\vr,\vs)\to(\vr,\vs)+\ve_i\to(\vr,\vs)+\ve_i+\ve_j\]
  is the morphism  ${(\beta^{\vs}_{j}*\alpha^{\vr}_i)\alpha^{\vr}_i\otimes\beta^{\vs}_{j}}$, while the morphism corresponding to the other path is $\alpha^{\vr}_i\otimes\beta^{\vs}_{j}$.
  If the square is zero, it automatically anti-commutes. Otherwise, it is sufficient to have
  \[\partial\epsilon(\square^{(\vr,\vs)}_{i,j})
  =_\ring -(\beta^{\vs}_{j}*\alpha^{\vr}_i)^{-1}
  =_\ring -\psi_\bA(\vr\to\vr+\ve_i)*\psi_\bB(\vs\to\vs+\ve_j),\]
  using symmetry of $\bil$ in the second equality.
\end{proof}

\begin{remark}
\label{rem:compatible_cochain}
  Note that the compatibility condition is a sufficient condition to ensure that all squares anti-commute, but a priori not a necessary one. Indeed, there might be some liberty on squares that are zero. In particular, we could have $\alpha^\vr_i\otimes\beta^\vs_j=0$ even if $\alpha^\vr_i\neq 0$ and $\beta^\vs_j\neq 0$.
  This happens for super $\glt$-foams $\cC=\sfoam_d$ (here $\otimes$ should be understood as a horizontal composition): for instance, one can pick $\alpha^\vr_i$ and $\beta^\vs_j$ elementary saddles (zip or unzip) such that $\undcomp(\id_{s(\alpha^\vr_i)}\otimes\beta^\vs_j)$ splits two circles and $\undcomp(\alpha^\vr_i\otimes\id_{t(\beta^\vs_j)})$ merges the two circles back into one circle.
\end{remark}

\Cref{defn:sigma_tensor_product} defines a tensor product of homogeneous polycomplexes, \emph{provided} that a compatible 1-cochain exists. \Cref{defn:standard_tensor_product} and \cref{lem:standard_tensor_product_is_compatible} below ensure existence.

\begin{definition}[graded Koszul rule]
\label{defn:standard_tensor_product}
  The \emph{standard 1-cochain $\epsilon_{\bA\otimes \bB}$} is the 1-cochain on $\bZ^{\cI\sqcup\cJ}$ given on $e=(\vr,\vs)\to(\vr,\vs)+\ve_k$ by
  \[
    \epsilon_{\bA\otimes \bB}(e)
    =
    \begin{cases}
      1
      & k=i\in\cI,
      \\
      (-1)^{\abs{\vr}}\abs{\alpha}(\vr)*\psi_\bB(\vs\to\vs+\ve_j)
      \quad& k=j\in\cJ.
    \end{cases}
  \]
\end{definition}

Note that if $\psi_\bA=\psi_\bB=0$, then $\abs{\alpha}(\vr)*\psi_\bB(\vs\to\vs+\ve_j) = 1$ and we recover the usual Koszul rule.

\begin{lemma}
  \label{lem:standard_tensor_product_is_compatible}
  The standard 1-cochain $\epsilon_{\bA\otimes \bB}$ is compatible in the sense of \cref{lem:compatible_implies_complex}.
\end{lemma}

\begin{proof}
  Consider the square $\square^{(\vr,\vs)}_{k,l}$ as above. If $k,l\in\cI$, the compatibility condition follows directly, and if $k,l\in\cJ$, it follows from the fact that $\psi_\bB$ is a cocycle.
  Finally, if $k=i\in\cI$ and $l=j\in\cJ$, we get (using property \eqref{eq:abs_alpha}):
  \begin{align*}
    \partial\epsilon(\square^{(\vr,\vs)}_{i,j})
    &=
    \Big((-1)^{\abs{\vr}+1}\abs{\alpha}(\vr+\ve_i)*\psi_\bB(\vs\to\vs+\ve_j)\Big)
    \\
    &\mspace{50mu}\Big((-1)^{\abs{\vr}}\abs{\alpha}(\vr)*\psi_\bB(\vs\to\vs+\ve_j)^{-1}\Big)
    \\
    &=-\big(\abs{\alpha}(\vr+\ve_i)-\abs{\alpha}(\vr)\big)*\psi_\bB(\vs\to\vs+\ve_j)
    \\
    &=-\psi_\bA(\vr\to\vr+\ve_i)*\psi_\bB(\vs\to\vs+\ve_j).\qedhere
  \end{align*}
\end{proof}

We next check that this choice of compatible 1-cochain is essentially unique, \emph{among compatible cochains} (see \cref{rem:compatible_cochain}).%
\footnote{The proof given below follows closely the proof of Lemma 2.2 in \cite{ORS_OddKhovanovHomology_2013}.}

\begin{lemma}
\label{lem:tensor_product_uniqueness}
  Let $\epsilon$ and $\epsilon'$ be two compatible 1-cochains. Then $(\bA\otimes \bB)(\epsilon)$ and $(\bA\otimes \bB)(\epsilon')$ are isomorphic as chain complexes.
\end{lemma}

\begin{proof}
  If $\epsilon$ and $\epsilon'$ are both compatible, then $\partial(\epsilon(\epsilon')^{-1})$ is zero on all squares, and hence is a \mbox{1-cocycle}. Let $\varphi$ be a 0-cochain such that $\partial\varphi=\epsilon(\epsilon')^{-1}$.
  Consider the map
  \[\Phi\colon (\bA\otimes \bB)(\epsilon)\to (\bA\otimes \bB)(\epsilon')
  \]
  corresponding to multiplication by $\varphi(\vr,\vs)$ when restricted to $(\vr,\vs)$. It is straightforward to check that this map defines an isomorphism for the associated total complexes.
\end{proof}

From now on, we simply call $(\bA\otimes \bB)(\epsilon_{\bA\otimes \bB})$ the \emph{tensor product} of $\bA$ and $\bB$ and denote it
$$
  \bA\otimes \bB=(A\otimes B,\alpha\otimes\beta,\psi_\bA\otimes\psi_\bB).
$$

\begin{example}
\label{ex:nfold_tensor_product}
  A chain complex $\bC=(C_\bullet,\partial_\bullet)$ whose differentials $\partial_t$ are all homogeneous is a homogeneous 1\nbd-polycomplex. An $n$-fold tensor product of homogeneous 1\nbd-polycomplexes is a homogeneous $n$-polycomplex. This is the specific construction used in \cref{subsec:hypercubic_chain_complexes} to define covering $\glt$-Khovanov homology.
\end{example}

We now come to the main result of this section:

\begin{theorem}
\label{thm:invariance_homotopy_classes}
  Let $\bA_1$, $\bA_2$, $\bB_1$ and $\bB_2$ be homogeneous polycomplexes. Then:
  \[
    \bA_1\simeq \bB_1\quad\an\quad \bA_2\simeq \bB_2
    \qquad\text{implies}\qquad
    \bA_1\otimes \bA_2\simeq \bB_1\otimes \bB_2,
  \]
  where $\simeq$ denotes homotopy equivalence (see \cref{defn:tensor_product_chain_map_chain_homotopies}).
\end{theorem}

The idea of the proof of \cref{thm:invariance_homotopy_classes} is straightforward: define induced morphisms and induced homotopies on the tensor product.
Precisely, \cref{thm:invariance_homotopy_classes} holds if the following holds:
\begin{enumerate}[label=(\arabic*)]
  \item Given homogeneous polycomplexes $\bA_k,\bB_k$ and chain maps $F_k\colon \bA_k\to \bB_k$ ($k=1,2$), there exists a chain map $F_1\otimes F_2\colon \bA_1\otimes \bA_2\to \bB_1\otimes \bB_2$.
  This definition is such that $F_1\otimes F_2=\Id_{\bA_1\otimes \bA_2}$ if $F_1=\Id_{\bA_1}$ and $F_2=\Id_{\bA_2}$.

  \item Given homogeneous polycomplexes $\bA_k,\bB_k$, chain maps $F_k,G_k\colon \bA_k\to \bB_k$ and homotopies $H_k\colon F_k\to G_k$ ($k=1,2$), there exists a homotopy $H_1\otimes H_2\colon F_1\otimes F_2\to G_1\otimes G_2$.
\end{enumerate}
This is the content of \cref{sec:induced_morphism_homotopy}: part (1) is shown by \cref{prop:induced_morphism}, while part (2) is shown by \cref{prop:induced_homotopy}.

\medbreak

\subsubsection{Further directions}
\label{subsec:complexes_further_directions}

We did not investigate the categorical properties of the tensor product: for our purpose, the mere existence of a homotopy equivalence in \cref{thm:invariance_homotopy_classes} is sufficient.
It remains unclear whether the construction of chain maps can be made functorial, or functorial up to homotopy.

Other possible questions include:
\begin{itemize}
  \item Can we extend the definition of the tensor product to higher structures? That is, can we define induced $n$-fold homotopies on the tensor product?

  \item Are the definitions of induced morphisms and induced homotopies unique up to higher structures, similarly to \cref{lem:tensor_product_uniqueness}?

  \item What is the most general case where one can define a (sensible) tensor product on chain complexes in a graded-monoidal category? In particular, how far can we weaken condition (iii) in \cref{defn:Hcomplex}?
\end{itemize}

\subsection{Induced morphisms and homotopies on the tensor product}
\label{sec:induced_morphism_homotopy}

%%%%%%%%%%%%%%%%%%%%%%%%%%%%%%%%%%%%%
%%%          conventions          %%%
%%%%%%%%%%%%%%%%%%%%%%%%%%%%%%%%%%%%%
\subsubsection{Preliminary remarks}

Recall the conventions of \cref{nota:complexes_basics} and \cref{nota:complex_B}. 
We denote the graded Koszul rule defined in \cref{defn:standard_tensor_product} as:
\begin{gather*}
  \epsilon_{\bA\otimes \bB}^{\vr,\vs,k}
  \coloneqq 
  % \epsilon_{\bA\otimes \bB}\big(e=(\vr,\vs)\to(\vr,\vs)+\ve_k\big)
  % =
    \begin{cases}
      1
      & k=i\in\cI,
      \\
      (-1)^{\abs{\vr}}\abs{\alpha}(\vr)*\psi_\bB(\vs\to\vs+\ve_j)
      \quad& k=j\in\cJ.
    \end{cases}
\end{gather*}

\begin{remark}
  \label{rem:complexes_expression_with_zero}
  Note that in \cref{defn:sigma_tensor_product}, for $k=j\in\cJ$ the scalar $\epsilon_{\bA\otimes \bB}^{\vr,\vs,j}$ always appears in front of an expression involving $\beta_j^\vs$.
  Thus, either $\beta_j^\vs=0$ and the value of $\epsilon_{\bA\otimes \bB}^{\vr,\vs,j}$ is irrelevant, or $\beta_j^\vs\neq 0$ and therefore $\psi_\bB(\vs\to\vs+\ve_j)=\deg\beta_j^\vs$.
  In general, if an expression contains a homogeneous element $f$, we can safely write $\deg f$ when defining or computing scalars appearing before the given expression; if $f$ turns out to be zero, this will be irrelevant anyway.
  See the next remark for a similar statement when $f$ is inhomogeneous.
\end{remark}

\begin{remark}
  \label{rem:complexes_inhomogeneous_formula}
  Write $[v]_a$ for the degree $a$ homogeneous component of a vector $v$ in a $G$-graded vector space; in particular, $v=\sum_{a\in G}[v]_a$.
  Recall that in graded-monoidal categories, a formula involving degrees of inhomogeneous morphisms should be understood by extending it additively.
  For instance, the graded interchange law
  \[(f\otimes 1)\circ (1\otimes g)=(\deg f*\deg g)(1\otimes g)\circ(f\otimes 1)\]
  for inhomogeneous $f$ and $g$ should really be understood as
  \[(f\otimes 1)\circ (1\otimes g)
  =\sum_{a,b\in G}(a*b)(1\otimes [g]_b)\circ([f]_a\otimes 1).\]
  We will encounter similar situations in what follows.
\end{remark}

\begin{notation}
  \label{nota:complexes_with_subscripts}
  We extend notations using subscripts to $\bA_k=(A_k,\alpha_k,\psi_{\bA_k})$ and $\bB_k=(B_k,\beta_k,\psi_{\bB_k})$ for $k=1,2$. From now on, we denote $\bA=\bA_1\otimes \bA_2$ and $\bB=\bB_1\otimes \bB_2$.
This includes sets of indices $\cI=\cI_1\sqcup\cI_2$ and $\cJ=\cJ_1\sqcup\cJ_2$.
We also use different shortcuts, such as $\vr=(\vr_1,\vr_2)$ and $\vs=(\vs_1,\vs_2)$, that should be clear from the context.
\end{notation}

Before entering the main proofs, we define a generic $\ring^\times$-valued 0-cochain, satisfying generic properties. In fact, \cref{defn:induced_morphism} and \cref{defn:induced_homotopy} solely depend on this generic definition, and all ``degree-wise'' computations follow from those generic properties.

\begin{definition}
\label{defn:generic_sigma}
  For $\lambda_k=\left(\lambda_k^{\vr_k,\vs_k}\right)_{(\vr_k,\vs_k)\in A_k\otimes B_k}$ a pair of $G$\nbd-valued 0\nbd-co\-chains ($k=1,2$), let $\epsilon^{\vr,\vs}_{\lambda_1,\lambda_2}$ be the following $\ring^\times$-valued 0-cochain:
  \begin{equation*}
    \epsilon^{\vr,\vs}_{\lambda_1,\lambda_2}
    \coloneqq
    \Big(\big[\lambda_1^{\vr_1,\vs_1}
    + \abs{\alpha_1}(\vr_1)
    - \abs{\beta_1}(\vs_1)\big]*\abs{\beta_2}(\vs_2)
    \Big)^{-1}
    \;
    \Big(\abs{\alpha_1}(\vr_1)*\lambda_2^{\vr_2,\vs_2}\Big).
  \end{equation*}
\end{definition}

\begin{lemma}
\label{lem:sigma_mu} 
  The generic $\ring^\times$-valued 0-cochain defined above satisfies the following:
  \begin{enumerate}[(i)]
    \item Assume $
    \nu_1\coloneqq_G\lambda_1^{\vr_1+\ve_{i_1},\vs_1}+\alpha_{i_1}^{\vr_1}
     =_G \lambda_1^{\vr_1,\vs_1-\ve_{j_1}}+\beta_{j_1}^{\vs_1-\ve_{j_1}}$ holds for all $i_1\in\cI_1$ and $j_1\in\cJ_1$; in particular, the $G$-valued 0-cochain $\nu_1$ is defined independently of the choice of $i_1$ and $j_1$. Then the following identities hold, for all $i_1\in\cI_1$ and $j_1\in\cJ_1$:
    \begin{equation*}
      \epsilon^{\vr,\vs}_{\nu_1,\lambda_2}
       =_\ring
      \epsilon^{\vr+e_{i_1},\vs}_{\lambda_1,\lambda_2}
      \left(\lambda_2^{\vr_2,\vs_2}*\alpha_{i_1}^{\vr_1}\right)
       =_\ring
      \epsilon^{\vr,\vs-e_{j_1}}_{\lambda_1,\lambda_2}.
    \end{equation*}

    \item Assume $
    \nu_2\coloneqq_G\lambda_2^{\vr_2+\ve_{i_2},\vs_2}+\alpha_{i_2}^{\vr_2}
     =_G \lambda_2^{\vr_2,\vs_2-\ve_{j_2}}+\beta_{j_2}^{\vs_2-\ve_{j_2}}$ holds for all $i_2\in\cI_2$ and $j_2\in\cJ_2$; in particular, the $G$-valued 0-cochain $\nu_2$ is defined independently of the choice of $i_2$ and $j_2$. Then the following identities hold, for all $i_2\in\cI_2$ and $j_2\in\cJ_2$:
    \begin{align*}
      \epsilon^{\vr,\vs}_{\lambda_1,\nu_2}
       &=_\ring
      (-1)^{\abs{\vr_1}}
      \epsilon_{A_1\otimes A_2}^{\vr_1,\vr_2,j_2}
      \epsilon^{\vr+\ve_{i_2},\vs}_{\lambda_1,\lambda_2}
       \\
       &=_\ring
      (-1)^{\abs{\vs_1}}
      \epsilon_{B_1\otimes B_2}^{\vs_1,\vs_2-\ve_{j_2},j_2}
      \epsilon^{\vr,\vs-\ve_{j_2}}_{\lambda_1,\lambda_2}
      \left(\beta^{\vs_2-\ve_{j_2}}_{j_2}*\lambda_1^{\vr_1,\vs_1}\right).
    \end{align*}
  \end{enumerate}
\end{lemma}

\begin{proof}
  It follows from direct computations, using relation \eqref{eq:abs_alpha}. Symmetry of $\bil$ gives the expressions as in the lemma.
  \begin{IEEEeqnarray*}{+rCl+x*}
    %%%%%%%%%%%%%%%%%%%%%%%
    \epsilon^{\vr+e_{i_1},\vs}_{\lambda_1,\lambda_2}
    &=_\ring&
    \Big(
    \big[\lambda_1^{\vr_1+\ve_{i_1},\vs_1}+\alpha_{i_1}^{\vr_1}
    + \abs{\alpha_1}(\vr_1)
    - \abs{\beta_1}(\vs_1)\big]*\abs{\beta_2}(\vs_2)
    \Big)^{-1}
    \\
    &&
    \mspace{100mu}\cdot
    \Big(\big[\abs{\alpha_1}(\vr_1)+\alpha_{i_1}^{\vr_1}\big]*\lambda_2^{\vr_2,\vs_2}\Big)
    \\
    &=_\ring&
    \epsilon^{\vr,\vs}_{\nu_1,\lambda_2}
    \left(\alpha_{i_1}^{\vr_1}* \lambda_2^{\vr_2,\vs_2}\right)
    %%%%%%%%%%%%%%%%%%%%%%%
    \\[2ex]
    \epsilon^{\vr,\vs-e_{j_1}}_{\lambda_1,\lambda_2}
    &=_\ring&
    \Big(
    \big[\lambda_1^{\vr_1,\vs_1-\ve_{j_1}}+\beta_{j_1}^{\vs_1-\ve_{j_1}}
    + \abs{\alpha_1}(\vr_1)
    - \abs{\beta_1}(\vs_1)\big]*\abs{\beta_2}(\vs_2)
    \Big)^{-1}
    \\
    &&
    \mspace{100mu}\cdot
    \Big(\abs{\alpha_1}(\vr_1)*\lambda_2^{\vr_2,\vs_2}\Big)
    \\
    &=_\ring&
    \epsilon^{\vr,\vs}_{\nu_1,\lambda_2}
    % here we use symmetry of \mu to get the expression as in the lemma; otherwise, we don't
    %%%%%%%%%%%%%%%%%%%%%%%
  \end{IEEEeqnarray*}
  \begin{IEEEeqnarray*}{+rCl+}
    \IEEEeqnarraymulticol{3}{l}{
      (-1)^{\abs{\vr_1}}
      \epsilon_{A_1\otimes A_2}^{\vr_1,\vr_2,j_2}
      \epsilon^{\vr+\ve_{i_2},\vs}_{\lambda_1,\lambda_2}
    }
    \\*
    \mspace{10mu}&=_\ring&
    \big(\abs{\alpha_1}(\vr_1)*\alpha_{i_2}^{\vr_2,\vs_2}\big)
    \Big(\big[\lambda_1^{\vr_1,\vs_1}
      + \abs{\alpha_1}(\vr_1)
      - \abs{\beta_1}(\vs_1)\big]*\abs{\beta_2}(\vs_2)
      \Big)^{-1}
    \\
    &&
    \mspace{100mu}\cdot
    \Big(\abs{\alpha_1}(\vr_1)*\lambda_2^{\vr_2+\ve_{i_2},\vs_2}\Big)
    \\
    &=_\ring&
    \Big(\big[\lambda_1^{\vr_1,\vs_1}
    + \abs{\alpha_1}(\vr_1)
    - \abs{\beta_1}(\vs_1)\big]*\abs{\beta_2}(\vs_2)
    \Big)^{-1}
    \\
    &&
    \mspace{100mu}\cdot
    \Big(\abs{\alpha_1}(\vr_1)*\big[\lambda_2^{\vr_2+\ve_{i_2},\vs_2}
    + \alpha_{i_2}^{\vr_2,\vs_2}\big]\Big)
    \\
    &=_\ring&
    \epsilon^{\vr,\vs}_{\lambda_1,\nu_2}
    %%%%%%%%%%%%%%%%%%%%%%%
    \\[2ex]
    \IEEEeqnarraymulticol{3}{l}{
      (-1)^{\abs{\vs_1}}
      \epsilon_{B_1\otimes B_2}^{\vs_1,\vs_2-\ve_{j_2},j_2}
      \epsilon^{\vr,\vs-\ve_{j_2}}_{\lambda_1,\lambda_2}
      \left(\lambda_1^{\vr_1,\vs_1}*\beta^{\vs_2-\ve_{j_2}}_{j_2}\right)^{-1}
    }
    \\*
    &=_\ring&
    \big(\abs{\beta_1}(\vs_1)*\beta_{j_2}^{\vr_2,\vs_2-\ve_{j_2}}\big)
    \\
    &&
    \mspace{100mu}\cdot
    \Big(\big[\lambda_1^{\vr_1,\vs_1}
    + \abs{\alpha_1}(\vr_1)
    - \abs{\beta_1}(\vs_1)\big]*\abs{\beta_2}(\vs_2-\ve_{j_2})
    \Big)^{-1}
    \\*
    \IEEEeqnarraymulticol{3}{r}{
      {}\cdot\Big(\abs{\alpha_1}(\vr_1)*\lambda_2^{\vr_2,\vs_2-\ve_{j_2}}\Big)
      \left(\lambda_1^{\vr_1,\vs_1}*\beta^{\vs_2-\ve_{j_2}}_{j_2}\right)^{-1}
    }
    \\
    &=_\ring&
    \Big(\big[\lambda_1^{\vr_1,\vs_1}
    + \abs{\alpha_1}(\vr_1)
    - \abs{\beta_1}(\vs_1)\big]*\abs{\beta_2}(\vs_2)
    \Big)^{-1}
    \\
    &&
    \mspace{100mu}\cdot
    \Big(\abs{\alpha_1}(\vr_1)*\big[\lambda_2^{\vr_2,\vs_2-\ve_{j_2}}+\beta^{\vs_2-\ve_{j_2}}_{j_2}\big]\Big)
    \\
    &=_\ring&
    \epsilon^{\vr,\vs}_{\lambda_1,\nu_2}.
    \hfill\qedhere
    % here again we use symmetry of \mu to get the expression as in the lemma; otherwise, we don't
  \end{IEEEeqnarray*}
\end{proof}

%%%%%%%%%%%%%%%%%%%%%%%%%%%%%%%%%%%
%%%          morphisms          %%%
%%%%%%%%%%%%%%%%%%%%%%%%%%%%%%%%%%%
\subsubsection{Induced morphism on tensor product}

Recall from \cref{defn:tensor_product_chain_map_chain_homotopies} that a chain map $F\colon\bA\to\bB$ is a family of morphisms $F^{\vr,\vs}\colon A^{\vr}\to B^{\vs}$ for all $\abs{\vr}=\abs{\vs}$ such that $\beta\circ F = F\circ\alpha$, that is:
\begin{equation*}
  \sum_{j\in\cJ} \beta_j^{\vs-\ve_j}
  \circ F^{\vr,\vs-\ve_j}
  =
  \sum_{i\in\cI} F^{\vr+\ve_i,\vs}\circ \alpha_i^{\vr}.
\end{equation*}
Recalling that $\beta_j^{\vs-\ve_j}$ and $\alpha_i^{\vr}$ are assumed to be homogeneous, the above identity is equivalent to the family of identities
\begin{equation*}
  \sum_{j\in\cJ} \beta_j^{\vs-\ve_j}
  \circ \left[F^{\vr,\vs-\ve_j}\right]_{g-\deg\beta_j^{\vs-\ve_j}}
  =
  \sum_{i\in\cI} \left[F^{\vr+\ve_i,\vs}\right]_{g-\deg\alpha_i^{\vr}}\circ \alpha_i^{\vr}
  \qquad \forall g\in G.
\end{equation*}
Here we use the notation $[-]_g$ from \cref{rem:complexes_inhomogeneous_formula} to denote the degree $g$ homogeneous component; see also \cref{rem:complexes_expression_with_zero}.

\begin{definition}
\label{defn:induced_morphism}
  Let $F_k\colon \bA_k\to \bB_k$ be chain maps for each $k=1,2$. Their \emph{induced morphism} $F=F_1\otimes F_2$ is the chain map $F\colon \bA_1\otimes \bA_2\to \bB_1\otimes \bB_2$ given by the data:
  \begin{gather*}
    F^{\vr,\vs}
    =
    \epsilon_{F_1,F_2}^{\vr,\vs}
    F_1^{\vr_1,\vs_1}\otimes F_2^{\vr_2,\vs_2},
  \end{gather*}
  where $\epsilon_{F_1,F_2}^{\vr,\vs}$ is as in \cref{defn:generic_sigma}, with the abuse of notation $F_i=\deg{F_i}$ and recalling \cref{rem:complexes_inhomogeneous_formula}.
\end{definition}

\begin{proposition}
\label{prop:induced_morphism}
  \Cref{defn:induced_morphism} gives a well-defined chain map. Moreover, if $F_1$ and $F_2$ are identities then $F$ is the identity.
\end{proposition}

\begin{proof}
  Note that if both $F_1$ and $F_2$ are identities, then
  \begin{equation*}
    \epsilon_{F_1,F_2}^{\vr,\vs}  =
    \Big(\big[0
    + \abs{\alpha_1}(\vr_1)
    + \abs{\alpha_1}(\vr_1)\big]*\abs{\alpha_2}(\vs_2)\Big)^{-1}
    \big(\abs{\alpha_1}(\vr_1)*0\big) = 1,
  \end{equation*}
  so that two identities induce the identity on the tensor product. To show the first part of the statement, we must check that $\beta\circ F = F\circ \alpha$. First, we unfold both sides of the equation:
  \begin{IEEEeqnarray*}{rClr}
    % LEFT-HAND SIDE
    \IEEEeqnarraymulticol{4}{l}{
      \sum_{j\in\cJ} \beta_j^{\vs-\ve_j}
      \circ F^{\vr,\vs-\ve_j}
    }\\*
    \quad & = &
      \sum_{j_1\in\cJ_1} \left[ \beta_{j_1}^{\vs_1-\ve_{j_1}}\otimes \Id\right]
      \circ
      \left[
        \epsilon_{F_1,F_2}^{\vr,\vs-\ve_{j_1}} 
        \;\;
        F_1^{\vr_1,\vs_1-\ve_{j_1}}\otimes F_2^{\vr_2,\vs_2}
      \right]
    &\\
    && {}+
      \sum_{j_2\in\cJ_2} \left[
      \epsilon_{B_1\otimes B_2}^{\vs_1,\vs_2-\ve_{j_2},j_2}
      \;\;
      \Id\otimes \beta_{j_2}^{\vs_2-\ve_{j_2}}\right]
    \\
    \IEEEeqnarraymulticol{3}{r}{
      \circ
      \left[
        \epsilon_{F_1,F_2}^{\vr,\vs-\ve_{j_2}} 
        \;\;
        F_1^{\vr_1,\vs_1}\otimes F_2^{\vr_2,\vs_2-\ve_{j_2}}
      \right]
    }
    &\\
    & = &
      \left[\sum_{j_1\in\cJ_1}
      \epsilon_{F_1,F_2}^{\vr,\vs-\ve_{j_1}}
      \;\;
      \beta_{j_1}^{\vs_1-\ve_{j_1}}
      \circ
      F_1^{\vr_1,\vs_1-\ve_{j_1}}\right]
      \otimes F_2^{\vr_2,\vs_2}
    & \qquad\Circled{1}
    \\
    && {}+
      F_1^{\vr_1,\vs_1}\otimes
      \Big[\sum_{j_2\in\cJ_2}
      \left(\beta_{j_2}^{\vs_2-\ve_{j_2}}*F_1^{\vr_1,\vs_1}\right)&
      \\
      \IEEEeqnarraymulticol{3}{r}{
        \epsilon_{B_1\otimes B_2}^{\vs_1,\vs_2-\ve_{j_2},j_2}
        \epsilon_{F_1,F_2}^{\vr,\vs-\ve_{j_2}}
        \;\;
        \beta_{j_2}^{\vs_2-\ve_{j_2}}
        \circ
        F_2^{\vr_2,\vs_2-\ve_{j_2}}\Big]
      }
      &
      \qquad\Circled{2}
    \\
    &&&
    \\
    % RIGHT-HAND SIDE
    \IEEEeqnarraymulticol{4}{l}{
        \sum_{i\in\cI} F^{\vr+\ve_i,\vs}\circ \alpha_i^{\vr}
    }\\*
    & = &
      \sum_{i_1\in\cI_1} \left[
        \epsilon_{F_1,F_2}^{\vr+\ve_{i_1},\vs} 
        \;\;
        F_1^{\vr_1+\ve_{i_1},\vs_1}\otimes F_2^{\vr_2,\vs_2}
      \right]
      \circ
      \left[ \alpha_{i_1}^{\vr_1}\otimes \Id\right]
    &\\*
    && {}+
      \sum_{i_2\in\cI_2} \left[
        \epsilon_{F_1,F_2}^{\vr+\ve_{i_2},\vs} 
        \;\;
        F_1^{\vr_1,\vs_1}\otimes F_2^{\vr_2+\ve_{i_2},\vs_2}
      \right]
      \circ
      \left[\epsilon_{A_1\otimes A_2}^{\vr_1,\vr_2,i_2}
      \;\;
      \Id\otimes \alpha_{i_2}^{\vr_2}\right]
    &\\
    & = &
      \left[ \sum_{i_1\in\cI_1}
      \left(F_2^{\vr_2,\vs_2}*\alpha_{i_1}^{\vr_1}\right)
      \epsilon_{F_1,F_2}^{\vr+\ve_{i_1},\vs}
      \;\;
      F_1^{\vr_1+\ve_{i_1},\vs_1}
      \circ
      \alpha_{i_1}^{\vr_1}\right]
      \otimes F_2^{\vr_2,\vs_2}
    & \Circled{1}
    \\
    && {}+
      F_1^{\vr_1,\vs_1}\otimes
      \left[\sum_{i_2\in\cI_2}
      \epsilon_{F_1,F_2}^{\vr+\ve_{i_2},\vs}
      \epsilon_{A_1\otimes A_2}^{\vr_1,\vr_2,i_2}
      \;\;
      F_2^{\vr_2+\ve_{i_2},\vs_2}
      \circ
      \alpha_{i_2}^{\vr_2}\right]
    & \Circled{2}
  \end{IEEEeqnarray*}
  We want to show that the two terms labelled $\Circled{1}$ (resp.\ $\Circled{2}$) are equal, using the chain map relation of $F_1$ (resp.\ $F_2$).
  Let us detail case $\Circled{1}$; the general strategy will be similar for case $\Circled{2}$, and for the proof of \cref{prop:induced_homotopy}.

  Restricting to homogeneous components of $F_2^{\vr_2,\vs_2}$, case $\Circled{1}$ reduces to showing
  \[
    \sum_{j_1\in\cJ_1}
      \epsilon_{F_1,F_2}^{\vr,\vs-\ve_{j_1}}
      \;\;
      \beta_{j_1}^{\vs_1-\ve_{j_1}}
      \circ
      F_1^{\vr_1,\vs_1-\ve_{j_1}}
    =
    \sum_{i_1\in\cI_1}
      \left(\left[F_2^{\vr_2,\vs_2}\right]_g*\alpha_{i_1}^{\vr_1}\right)
      \epsilon_{F_1,F_2}^{\vr+\ve_{i_1},\vs}
      \;\;
      F_1^{\vr_1+\ve_{i_1},\vs_1}
      \circ
      \alpha_{i_1}^{\vr_1}
  \]
  for each element $g\in G$ for which $\left[F_2^{\vr_2,\vs_2}\right]_g\neq 0$.
  Note that in the above equation, both $\epsilon_{F_1,F_2}^{\vr,\vs-\ve_{j_1}}$ and $\epsilon_{F_1,F_2}^{\vr+\ve_{i_1},\vs}$ depend on $F_2$ only through the value of $\deg\left[F_2^{\vr_2,\vs_2}\right]_g=g$.
  For simplicity, we abuse notation and write $F_2^{\vr_2,\vs_2}=\left[F_2^{\vr_2,\vs_2}\right]_g$; this has the same effect as assuming that $F_2^{\vr_2,\vs_2}$ is homogeneous and non-zero.
  With this convention, the equation becomes:
  \[
    \sum_{j_1\in\cJ_1}
      \epsilon_{F_1,F_2}^{\vr,\vs-\ve_{j_1}}
      \;\;
      \beta_{j_1}^{\vs_1-\ve_{j_1}}
      \circ
      F_1^{\vr_1,\vs_1-\ve_{j_1}}
    =
    \sum_{i_1\in\cI_1}
      \left(F_2^{\vr_2,\vs_2}*\alpha_{i_1}^{\vr_1}\right)
      \epsilon_{F_1,F_2}^{\vr+\ve_{i_1},\vs}
      \;\;
      F_1^{\vr_1+\ve_{i_1},\vs_1}
      \circ
      \alpha_{i_1}^{\vr_1}.
  \]
  In turn, the above can be shown by projecting onto each degree $h$ homogeneous component, for $h\in G$. Since $\beta_{j_1}^{\vs_1-\ve_{j_1}}$ (resp.\ $\alpha_{i_1}^{\vr_1}$) is assumed to be homogeneous, this enforces the degree of $F_1^{\vr_1,\vs_1-\ve_{j_1}}$ (resp.\ $F_1^{\vr_1+\ve_{i_1},\vs_1}$) and gives:
  \begin{gather*}
    \sum_{j_1\in\cJ_1}
      \epsilon_{F_1,F_2}^{\vr,\vs-\ve_{j_1}}
      \;\;
      \beta_{j_1}^{\vs_1-\ve_{j_1}}
      \circ
      \left[F_1^{\vr_1,\vs_1-\ve_{j_1}}\right]_{h-\deg\beta_{j_1}^{\vs_1-\ve_{j_1}}}
    \mspace{150mu}
    \\
    \mspace{150mu}
    =
    \sum_{i_1\in\cI_1}
      \left(F_2^{\vr_2,\vs_2}*\alpha_{i_1}^{\vr_1}\right)
      \epsilon_{F_1,F_2}^{\vr+\ve_{i_1},\vs}
      \;\;
      \left[F_1^{\vr_1+\ve_{i_1},\vs_1}\right]_{h-\deg\alpha_{i_1}^{\vr_1}}
      \circ
      \alpha_{i_1}^{\vr_1}.
  \end{gather*}
  We abuse notation and write $F_1^{\vr_1,\vs_1-\ve_{j_1}}=\left[F_1^{\vr_1,\vs_1-\ve_{j_1}}\right]_{h-\deg\beta_{j_1}^{\vs_1-\ve_{j_1}}}$. This has the same effect as assuming that $F_1^{\vr_1,\vs_1-\ve_{j_1}}$ is homogeneous with the property that if both $\beta_{j_1}^{\vs_1-\ve_{j_1}}$ and $F_1^{\vr_1,\vs_1-\ve_{j_1}}$ are non-zero, then:
  \[h=\deg\beta_{j_1}^{\vs_1-\ve_{j_1}}+\deg F_1^{\vr_1,\vs_1-\ve_{j_1}}.\]
  For our purpose, we can assume the latter holds whenever we compute with $\epsilon_{F_1,F_2}^{\vr,\vs-\ve_{j_1}}$ (see \cref{rem:complexes_expression_with_zero}).
  \Cref{lem:sigma_mu} then implies that $\epsilon_{F_1,F_2}^{\vr,\vs-\ve_{j_1}}$ is independent of $j_1$.

  A similar reasoning applies to $F_1^{\vr_1+\ve_{i_1},\vs_1}$, with the following degree condition:
  \[h=\deg F_1^{\vr_1+\ve_{i_1},\vs_1}
    +\deg\alpha_{i_1}^{r_{i_1}}.\]
  In fact, \cref{lem:sigma_mu} also shows that, whenever the relevant morphisms are non-zero, we have:
  \[\epsilon_{F_1,F_2}^{\vr,\vs-\ve_{j_1}}=
  \left(F_2^{\vr_2,\vs_2}*\alpha_{i_1}^{\vr_1}\right)
      \epsilon_{F_1,F_2}^{\vr+\ve_{i_1},\vs}.\]
  This allows us to use the chain map relation for $F_1$ and concludes.

  The proof of case $\Circled{1}$ essentially comes down to using \cref{lem:sigma_mu} with the following conditions:
  \begin{enumerate}[label=\protect\Circled{\arabic*}]
    \item $\beta_{j_1}^{s_{j_1}-1}+  F_1^{\vr_1,\vs_1-\ve_{j_1}} =_G F_1^{\vr_1+\ve_{i_1},\vs_1}
    +\alpha_{i_1}^{r_{i_1}}$ for all $i_1\in \cI_1$ and $j_1\in\cJ_1$
  \end{enumerate}
  Here we use the abuse of notation that leaves $\deg(-)$ implicit.
  Case $\Circled{2}$ is dealt with similarly, using the following conditions:
  \begin{enumerate}[label=\protect\Circled{\arabic*},resume]
    \item $\beta_{j_2}^{s_{j_2}-1}+  F_2^{\vr_2,\vs_2-\ve_{j_2}} =_G F_2^{\vr_2+\ve_{i_2},\vs_2}
    +\alpha_{i_2}^{r_{i_2}}$ for all $i_2\in \cI_2$ and $j_2\in\cJ_2$, and $\abs{\vr_1}=\abs{\vs_1}$.
  \end{enumerate}
  Again, these are exactly the assumptions needed to apply \cref{lem:sigma_mu}.\qedhere
\end{proof}

%%%%%%%%%%%%%%%%%%%%%%%%%%%%%%%%%%%%
%%%          homotopies          %%%
%%%%%%%%%%%%%%%%%%%%%%%%%%%%%%%%%%%%
\subsubsection{Induced homotopies on the tensor product}

Recall from \cref{defn:tensor_product_chain_map_chain_homotopies} that a chain homotopy $H$ between chain maps $F\colon \bA\to \bB$ and $G\colon \bA\to \bB$ is a family of morphisms
$H^{\vr,\vs}\colon A^{\vr}\to B^{\vs}$
for all $\abs{\vr}=\abs{\vs}+1$ such that
\begin{equation*}
  F^{\vr,\vs} - G^{\vr,\vs}
  =
  \sum_{i\in\cI} H^{\vr+\ve_i,\vs}
  \circ \alpha_{i}^{\vr}
  +
  \sum_{j\in\cJ} \beta_j^{\vs-\ve_j}
  \circ H^{\vr,\vs-\ve_j}.
\end{equation*}

\begin{definition}
\label{defn:induced_homotopy}
  Let $F_k\colon \bA_k\to \bB_k$ (resp.\ $G_k$) be chain maps and $H_k\colon F_k\to G_k$ homotopies for each $k=1,2$. Denote $F$ (resp.\ $G$) the chain map induced by $F_1$ and $F_2$ (resp.\ $G_1$ and $G_2$). Their \emph{induced homotopy} is the homotopy $H\colon F\to G$ given by the data:
  \begin{align*}
    H^{\vr,\vs}
    = &
      {\epsilon_{H_1,F_2}^{\vr,\vs}}
      H_1^{\vr_1,\vs_1} \otimes F_2^{\vr_2,\vs_2}
      + (-1)^{\abs{\vr_1}}\epsilon_{G_1,H_2}^{\vr,\vs}
      G_1^{\vr_1,\vs_1} \otimes H_2^{\vr_2,\vs_2}
  \end{align*}
  where the $\epsilon$'s are defined as in \cref{defn:generic_sigma}.
\end{definition}

\begin{proposition}
\label{prop:induced_homotopy}
 \Cref{defn:induced_homotopy} gives a well-defined homotopy.
\end{proposition}

\begin{proof}
  We assume the reader is familiar with the proof of \cref{prop:induced_morphism}.
  We must check that $F-G=\beta\circ H+H\circ\alpha$.
  We unfold the two terms on the right-hand side (RHS), with first the computation of $\beta\circ H$ restricted to the paths from $A^\vr$ to $B^\vs$:
  \begin{IEEEeqnarray*}{+rCl+x*}
    % H AND BETA (1)
      (\beta\circ H)\vert_{A^{\vr}}^{B^{\vs}}& = &
        \sum_{j\in\cJ} \beta_j^{\vs-\ve_j}
        \circ H^{\vr,\vs-\ve_j}
    \\ \\
    & = &
        % sum over j_1
        \sum_{j_1\in\cJ_1} \Big[\beta_{j_1}^{\vs_1-\ve_{j_1}}\otimes \Id\Big]
        \\ \\
        &&\mspace{30mu}
        ~\circ~
        \Big[
        \epsilon_{H_1,F_2}^{\vr,\vs-\ve_{j_1}} H_1^{\vr_1,\vs_1-\ve_{j_1}} \otimes F_2^{\vr_2,\vs_2}
        \\
        &&\mspace{100mu}
        +
        (-1)^{\abs{\vr_1}}\epsilon_{G_1,H_2}^{\vr,\vs-\ve_{j_1}} G_1^{\vr_1,\vs_1-\ve_{j_1}} \otimes H_2^{\vr_2,\vs_2}
        \Big]
    \\ \\
    && +
        % sum over j_2
        \sum_{j_2\in\cJ_2} \Big[
        \epsilon_{B_1\otimes B_2}^{\vs_1,\vs_2-\ve_{j_2},j_2}
        \Id \otimes \beta_{j_2}^{\vs_2-\ve_{j_2}}\Big]
    \\ \\
    &&\mspace{30mu}
    ~\circ~
        \Big[
        \epsilon_{H_1,F_2}^{\vr,\vs-\ve_{j_2}} H_1^{\vr_1,\vs_1} \otimes F_2^{\vr_2,\vs_2-\ve_{j_2}}
        \\
        &&\mspace{100mu}
        +
        (-1)^{\abs{\vr_1}}\epsilon_{G_1,H_2}^{\vr,\vs-\ve_{j_2}} G_1^{\vr_1,\vs_1} \otimes H_2^{\vr_2,\vs_2-\ve_{j_2}} \Big]
    \\ \\
    % H AND BETA (2)
    & = &
        \sum_{j_1\in\cJ_1}
        \epsilon_{H_1,F_2}^{\vr,\vs-\ve_{j_1}}
        \Big[
            \beta_{j_1}^{\vs_1-\ve_{j_1}}
            \circ
            H_1^{\vr_1,\vs_1-\ve_{j_1}}
        \Big] \otimes F_2^{\vr_2,\vs_2}
    & \Circled{1}
    \\
    && \mspace{30mu}+~
        (-1)^{\abs{\vr_1}}\epsilon_{G_1,H_2}^{\vr,\vs-\ve_{j_1}} 
        \Big[
            \beta_{j_1}^{\vs_1-\ve_{j_1}}
            \circ
            G_1^{\vr_1,\vs_1-\ve_{j_1}}
        \Big] \otimes H_2^{\vr_2,\vs_2}
    & \Circled{2}
    \\ \\
    && +
        % sum over j_2
        \sum_{j_2\in\cJ_2}
        \left(\beta_{j_2}^{\vs_2-\ve_{j_2}}*H_1^{\vr_1,\vs_1}\right)
        \epsilon_{B_1\otimes B_2}^{\vs_1,\vs_2-\ve_{j_2},j_2}
        \epsilon_{H_1,F_2}^{\vr,\vs-\ve_{j_2}}
        \\
        &&\mspace{100mu}
        ~~H_1^{\vr_1,\vs_1} \otimes \Big[
            \beta_{j_2}^{\vs_2-\ve_{j_2}}
            \circ
            F_2^{\vr_2,\vs_2-\ve_{j_2}}
        \Big]
    & \Circled{3}
    \\
    && \mspace{30mu}+~
        \left(\beta_{j_2}^{\vs_2-\ve_{j_2}}*G_1^{\vr_1,\vs_1}\right)
        (-1)^{\abs{\vr_1}}
        \epsilon_{B_1\otimes B_2}^{\vs_1,\vs_2-\ve_{j_2},j_2}
        \epsilon_{G_1,H_2}^{\vr,\vs-\ve_{j_2}}
        \\
        &&\mspace{100mu}
        ~~G_1^{\vr_1,\vs_1} \otimes \Big[
            \beta_{j_2}^{\vs_2-\ve_{j_2}}
            \circ
            H_2^{\vr_2,\vs_2-\ve_{j_2}}
        \Big]
    & \Circled{4}
  \end{IEEEeqnarray*}

  \noindent The computation of $H\circ\alpha$ restricted to the paths from $A^\vr$ to $B^\vs$ gives:
  \begin{IEEEeqnarray*}{+rCl+x*}
    % ALPHA AND H (1)
    \IEEEeqnarraymulticol{3}{l}{
      (H\circ\alpha)\vert_{A^{\vr}}^{B^{\vs}}
    }
    \\[2ex]
    & = &
        \sum_{i\in\cI} H^{\vr+\ve_i,\vs}
        \circ \alpha_{i}^{\vr}
    \nonumber
    \\ \nonumber \\
    & = &
        % sum over i_1
        \sum_{i_1\in\cI_1} \Big[
        \epsilon_{H_1,F_2}^{\vr+\ve_{i_1},\vs} H_1^{\vr_1+\ve_{i_1},\vs_1} \otimes F_2^{\vr_2,\vs_2}
        \\
        &&\mspace{100mu}
        +
        (-1)^{\abs{\vr_1}+1}
        \epsilon_{G_1,H_2}^{\vr+\ve_{i_1},\vs} G_1^{\vr_1+\ve_{i_1},\vs_1} \otimes H_2^{\vr_2,\vs_2} \Big]
    \\ \\
    &&\mspace{30mu}
        \circ
        \Big[ \alpha_{i_1}^{\vr_{1}}\otimes \Id\Big]
    \nonumber
    \\ \nonumber \\
    && +
        % sum over i_2
        \sum_{i_2\in\cI_2} \Big[
        \epsilon_{H_1,F_2}^{\vr+\ve_{i_2},\vs} H_1^{\vr_1,\vs_1} \otimes F_2^{\vr_2+\ve_{i_2},\vs_2}
        \\
        &&\mspace{100mu}
        +
        (-1)^{\abs{\vr_1}}
        \epsilon_{G_1,H_2}^{\vr+\ve_{i_2},\vs} G_1^{\vr_1,\vs_1} \otimes H_2^{\vr_2+\ve_{i_2},\vs_2} \Big]
    \\ \\
    &&\mspace{30mu}
        \circ
        \Big[
        \epsilon_{A_1\otimes A_2}^{\vr_1,\vr_2,i_2}
        \Id\otimes\alpha_{i_2}^{\vr_{2}}\Big]
    \nonumber
    \\ \nonumber \\
    % ALPHA AND H (2)
    & = &
        % sum over i_1
        \sum_{i_1\in\cI_1}
        \left(F_2^{\vr_2,\vs_2}*\alpha_{i_1}^{\vr_{1}}\right)
        \epsilon_{H_1,F_2}^{\vr+\ve_{i_1},\vs}
        \Big[
            H_1^{\vr_1+\ve_{i_1},\vs_1}
            \circ
            \alpha_{i_1}^{\vr_{1}}
        \Big] \otimes F_2^{\vr_2,\vs_2}
    & \Circled{1}
    \\
    && \mspace{40mu}+~
        \left(H_2^{\vr_2,\vs_2}*\alpha_{i_1}^{\vr_{1}}\right)
        (-1)^{\abs{\vr_1}+1}
        \epsilon_{G_1,H_2}^{\vr+\ve_{i_1},\vs}
        \Big[
            G_1^{\vr_1+\ve_{i_1},\vs_1}
            \circ
            \alpha_{i_1}^{\vr_{1}}
        \Big] \otimes H_2^{\vr_2,\vs_2}
    & \Circled{2}
    \\ \\
    && +
        % sum over i_2
        \sum_{i_2\in\cI_2}
        \epsilon_{A_1\otimes A_2}^{\vr_1,\vr_2,i_2}
        \epsilon_{H_1,F_2}^{\vr+\ve_{i_2},\vs}
        ~~H_1^{\vr_1,\vs_1} \otimes \Big[
            F_2^{\vr_2+\ve_{i_2},\vs_2}
            \circ
            \alpha_{i_2}^{\vr_{2}}
        \Big]
    & \Circled{3}
    \\
    && \mspace{50mu}+~
        (-1)^{\abs{\vr_1}}
        \epsilon_{A_1\otimes A_2}^{\vr_1,\vr_2,i_2}
        \epsilon_{G_1,H_2}^{\vr+\ve_{i_2},\vs}
        ~~G_1^{\vr_1,\vs_1} \otimes \Big[
            H_2^{\vr_2+\ve_{i_2},\vs_2}
            \circ
            \alpha_{i_2}^{\vr_{2}}
        \Big]
    & \Circled{4}
  \end{IEEEeqnarray*}

  \noindent We find four different pairs of terms, labelled $\Circled{1}$ to $\Circled{4}$, that we simplify pairwise using chain map or chain homotopy relations:
  \begin{enumerate}
    \item[\Circled{1}] We can assume that for all $i_1\in\cI_1$ and all $j_1\in\cJ_1$:
    \begin{equation*}
      \beta_{j_1}^{\vs_1-\ve_{j_1}}+H_1^{\vr_1,\vs_1-e_{j_1}}
       =_G H_1^{\vr_1+e_{i_1},\vs_1}+\alpha_{i_1}^{\vr_1}
       =_G F_1^{\vr_1,\vs_1}
       =_G G_1^{\vr_1,\vs_1}.
    \end{equation*}
    Thanks to \cref{lem:sigma_mu} (i), the sum of the two $\Circled{1}$-terms is:
    \begin{equation*}
      \left[
      \epsilon_{F_1,F_2}^{\vr_1,\vs_1} F_1^{\vr_1,\vs_1}
      -
      \epsilon_{G_1,F_2}^{\vr_1,\vs_1} G_1^{\vr_1,\vs_1}
      \right]\otimes F_2^{\vr_2,\vs_2}.
    \end{equation*}
    The computation is similar for case $\Circled{4}$.

    \item[\Circled{3}]
    We can assume $\abs{\vr_1}+1=\abs{\vs_1}$
    % and $\abs{\vr_2}=\abs{\vs_2}+1$.
    and that:
    \begin{equation*}
      \beta_{j_2}^{\vs_2-\ve_{j_2}}+F_2^{\vr_2,\vs_2-e_{j_2}}
       =_G F_2^{\vr_2+e_{i_2},\vs_2}+\alpha_{i_2}^{\vr_2}.
    \end{equation*}
    Then, \cref{lem:sigma_mu} (ii) shows that:
    \begin{equation*}
      \left(H_1^{\vr_1,\vs_1}*\beta_{j_2}^{\vs_2-\ve_{j_2}}\right)
        \epsilon_{B_1\otimes B_2}^{\vs_1,\vs_2-\ve_{j_2},j_2}
        \epsilon_{H_1,F_2}^{\vr,\vs-\ve_{j_2}}
      =
      -
      \epsilon_{A_1\otimes A_2}^{\vr_1,\vr_2,i_2}
        \epsilon_{H_1,F_2}^{\vr+\ve_{i_2},\vs}
    \end{equation*}
    independently of $j_2$ and $i_2$. We conclude that the sum of the pair $\Circled{3}$ is zero. The computation is similar for case $\Circled{2}$.
  \end{enumerate}
  We conclude:
  \begin{IEEEeqnarray*}{+rCl+x*}
    \text{RHS} &=&
    \left[
    {\epsilon_{F_1,F_2}^{\vr_1,\vs_1}} F_1^{\vr_1,\vs_1}
    -
    {\epsilon_{G_1,F_2}^{\vr_1,\vs_1}} G_1^{\vr_1,\vs_1}
    \right]
    \otimes F_2^{\vr_2,\vs_2}
    % \\*
    % &&{}
    \\
    &&
    {}+
    G_1^{\vr_1,\vs_1}\otimes
    \left[
    {\epsilon_{G_1,F_2}^{\vr_1,\vs_1}} F_2^{\vr_2,\vs_2}
    -
    {\epsilon_{G_1,G_2}^{\vr_1,\vs_1}} G_2^{\vr_2,\vs_2}
    \right]
    \\[1ex]
    &=&
    {\epsilon_{F_1,F_2}^{\vr_1,\vs_1}}
    F_1^{\vr_1,\vs_1}\otimes F_2^{\vr_2,\vs_2}
    -
    {\epsilon_{G_1,G_2}^{\vr_1,\vs_1}}
    G_1^{\vr_1,\vs_1}\otimes G_2^{\vr_2,\vs_2}
    \\[1ex]
    &=& (F-G)\vert_{A^{\vr}}^{B^{\vs}}.&\qedhere
  \end{IEEEeqnarray*}
\end{proof}

\printbibliography[heading=bibintoc]

\end{document}